\documentclass[francais,10pt,leqno,a4paper,psamsfonts]{smfart}
\RequirePackage{smfthm}%
\usepackage[francais, english]{babel}%
\newenvironment{preuvede}[1]{\begin{proof}[D\'emonstration #1]}{\end{proof}}%
\let\CAL=\mathcal%
\def\cal#1{{\CAL#1}}%
\let\goth=\mathfrak%
\let\Bbb = \mathbb%
\def\bibartp#1#2#3#4#5#6#7#8
{\bibitem{#1}
{\sc #2},
{\it #3},
{#4},
{\bf #5}, pages {#6} {\`a} {#7}, (#8)}
\def\bibart#1#2#3#4#5#6
{\bibitem{#1}
{\sc #2},
{\it #3},
{#4},
{\bf #5}, (#6)}
\def\bibliv#1#2#3#4#5
{\bibitem{#1}
{\sc #2},
{\it #3},
{#4},
(#5)}
\def\bibaart#1#2#3#4
{\bibitem{#1}
{\sc #2},
{\it #3},
{#4}}
\def\ssi{\Longleftrightarrow}% si et seulement si%
\def\fl{\rightarrow}%
\def\fle{\longrightarrow}%
\def\ifle{\hookrightarrow}%
\def\efle{\longmapsto}%
\def\iso{\stackrel{\sim}{\fle}}%
\def\is{\stackrel{\sim}{\rightarrow}}%
\def\dd#1{\frac{\displaystyle \partial \phantom{#1}}{\displaystyle \partial {#1}}}%
\def\DD#1#2{\frac{\displaystyle \partial {#1}}{\displaystyle \partial {#2}}}%
\def\efrac#1#2{\frac{\scriptstyle #1}{\scriptstyle #2}}%
\def\inte#1{\stackrel{{}_{\scriptstyle\circ}}{{#1}}}%
\def\ptsur#1{{\stackrel{{\,{}_{ {}_\bullet }}}{#1}}}%
\def\trait#1{{\;{}^{\underline{\phantom{#1}}}\;}}%
\def\ftrait{\;{{}^{\underline{\phantom{\longleftrightarrow}}}}\;}%
\def\bdf{{\leftarrow\!\mapsto}}%
\let\nsubset\subsetn%
\def\supsetn{{\supseteq \hspace{-1,2em} {}_{\hbox{ \Tiny $/$}}}}%
\let\nsupset\supsetn
\def\bulletO{{\stackrel{0}{\bullet}}}%
\def\bulletun{{\stackrel{1}{\bullet}}}%
\def\a{\alpha}%
\def\wa{\widetilde{a}}%
\def\A{\Bbb A}%
\let\aA=\dA%
\def\cA{{\check \A}}%
\def\hgC{\widehat{\goth{C}}}%
\def\b{\beta}%
\def\hb{\widehat{\beta}}%
\def\wb{\widetilde{b}}%
\def\hB{\widehat{\cal{B}}}%
\def\C{\Bbb C}%
\def\CC{\cal{C}}%
\def\bC{{\boldsymbol C}}%
\def\d{\delta}%
\def\hd{\widehat{\delta}}%
\def\D{\cal{D}}%
\def\wD{\widetilde{\D}}%
\def\e{{\varepsilon}}%
\def\he{\widehat{\varepsilon}}%
\def\hE{{\widehat{\cal{E}}}}%
\def\wE{{\widetilde{E}}}%
\def\f{\phi}%
\def\wf{\widetilde{\f}}%
\def\hff{\widehat{f}}%
\def\F{\cal{F}}%
\def\F{\cal{F}}%
\def\wF{{\widetilde{\F}}}%
\def\hF{\widehat{\F}}%
\def\sF{\underline{\F}}%
\def\uh{\underline{h}}%
\def\g{\gamma}%
\def\wg{\widetilde{\gamma}}
\def\G{\cal{G}}%
\def\hG{\widehat{\G}}%
\def\hhG{\widehat{\hG}}%
\def\bG{\mbox{\boldmath $G$}}%
\def\H{\cal{H}}%
\def\wh{\widetilde{h}}%
\def\hH{\widehat{\cal H}}%
\def\ui{\upsilon}%
\def\H{\cal{H}}%
\def\hH{\widehat{\cal{H}}}%
\def\l{\lambda}%c'est L-minuscule et pas i-majuscule
\def\J{{\cal{J}}}%
\def\wJ{\widetilde{\J}}%
\def\L{\Lambda}%
\def\hL{\widehat{\Lambda}}%
\def\M{\cal{M}}%
\def\wM{\widetilde{\M}}%
\def\n{\eta}%
\def\wnu{\widetilde{\nu}}%
\def\hgN{\widehat{\goth{N}}}%
\def\N{\Bbb N}%
\def\O{\cal{O}}%
\def\wO{\widetilde{\O}}%
\def\hO{\widehat{\O}}%
\def\p{\pi}%
\def\wp{\widetilde{\pi}}%
\def\Q{\Bbb Q}%
\def\r{\rho}%
\def\hr{\widehat{\r}}%
\def\hgR{\widehat{\goth{R}}}%
\def\s{\sigma}%
\def\S{\Sigma}%
\def\wS{\widetilde{\S}}%
\def\SS{\cal{S}}%
\def\gS{{\goth{S}}}%
\def\bS{{\Bbb S}}%
\def\hsep{\widehat{Sep}}%
\def\t{\tau}%
\def\htau{\widehat{\tau}}%
\def\wt{\widetilde{\t}}%
\def\hT{\widehat{\cal{T}}}%
\def\U{\cal{U}}%
\def\V{\cal{V}}%
\def\v#1{\relax\ifmmode{{\nu}}\else{$\check{\hbox{\sc #1}}$}\fi\relax}%
\def\X{\cal{X}}%
\def\hX{\widehat{\X}}%
\def\w{\omega}%
\def\ww{{\widetilde{\omega}}}%
\def\uw{\overline{\omega}}%
\def\wz{\widetilde{z}}%
\def\z{\zeta}%
\def\Z{\Bbb Z}%
\def\cZ{\cal{Z}}%
\def\stw{a(x, y)\,dx + b(x, y)\,dy}%
\def\std{{A(x, y; \,t)\,dx + B(x, y; \,t)\,dy}}%
\def\hSL{{{$\widehat{\scriptstyle SL}\,$}}}%
\let\hsl=\hSL%
\def\heSL{{\sim_{\widehat{\scriptscriptstyle SL}}}}%
\def\1f{{\L_{\C^2, 0}^1}}
\def\fa{{\f_{\a}}}%
\def\fao{{\f_{\a}^{(0)}}}%
\def\fau{{\f_{\a}^{(1)}}}%
\def\fb{{\f_{\b}}}%
\def\fbo{{\f_{\b}^{(0)}}}%
\def\fbu{{\f_{\b}^{(1)}}}%
\def\fab{{\f_{\a\b}}}%
\def\fabp{{\f_{\a\b}'}}%
\def\fabu{{\f_{\a\b}^{(1)}}}%
\def\fabd{{\f_{\a\b}^{(2)}}}%
\def\xab{\xi_{\a\b}}%\a\b}%
\def\limproj{\mathop{\oalign{lim\cr\hidewidth$\longleftarrow$\hidewidth\cr}}}%
\def\limind{\mathop{\oalign{lim\cr\hidewidth$\longrightarrow$\hidewidth\cr}}}%
\author[J.-F. Mattei]{Jean-Fran\c cois Mattei}
\address{Laboratoire de Math\'ematiques Emile Picard,
UMR CNRS 5580\\
Universit\'e Paul Sabatier\\
118, route de Narbonne\\
31062 TOULOUSE Cedex 4 France\\
fax : 33 (0)5 61 55 82 00
}
\email{mattei@picard.ups-tlse.fr}
\author[E. Salem]{Eliane Salem}
\address{
Institut de Math\'ematiques,\\
Universit\'e de Paris 6\\
175 rue du Chevaleret, 75013 Paris Cedex. Et~:
Department of Mathematics\\
Ben Gurion University of the Negev\\
POB 653\\
Beer Sheva 84105, Israel\\
}
\email{salem@math.jussieu.fr}
\title[Modules formels de feuilletages holomorphes]{Modules formels locaux de feuilletages holomorphes}
\alttitle{Formal local modulus for holomorphic singular foliations}
\NumberTheoremsIn{subsection}\NoSwapTheoremNumbers
\begin{document}
\begin{abstract}
Nous donnons une liste compl\`ete d'invariants formels locaux d'une large classe de 1-formes diff\'e\-rentielles
formelles $\w \in \Bbb C [[ x, y]]dx + \Bbb C [[ x, y]]dy$. \\
\indent Une d\'eformation
$\widehat{\mbox{SL}}$-\'equisingu\-li\`ere est une d\'eformation \'equir\'eductible qui, apr\`es r\'eduction, laisse
invariant chaque type singulier formel local ainsi que  la repr\'esentation d'holo\-nomie de chaque composante du diviseur
exceptionnel. Nous caract\'erisons les 1-formes de types formel fini (t.f.f.),
i.e. celles qui poss\`edent une d\'eformation
$\widehat{\mbox{SL}}$-\'equisinguli\`ere verselle, en donnant un crit\`ere combinatoire explicite de finitude. \\
\indent Les
formes t.f.f. contiennent un ouvert dense pour la topologie de Krull dans l'ouvert des 1- formes de deuxi\`eme esp\`ece.
\end{abstract}
\begin{altabstract}\hfill\\
We give a complete list of formal invariants for a large class of formal differential 1-forms
$\w \in \Bbb C [[ x, y]]dx + \Bbb C [[ x, y]]dy$.\\
\indent A $\widehat{\mbox{SL}}$-equisingular deformation is an equireducible deformation which leaves invariant both the local
formal types and the holonomy representation of the components of the exceptional divisor. We characterize the 1-forms Êwith
finite formal type (t.f.f), i.e. those which admit a semi-universal $\widehat{\mbox{SL}}$-equisingular
Êdeformation, and we give an explicit combinatorial criterion of finiteness. \\
\indent The set of 1-forms with finite formal type
contains a dense open set (in the sense of Krull's topology)in the set of 1-forms of the second kind.
\end{altabstract}

%%\dedicatory{A Fr\'ed\'eric Pham}

\maketitle

\tableofcontents

\section*{Introduction}
La classification formelle des germes d'\'equations diff\'erentielles
ordinaires \`a l'origine de $\C^2$ dont la partie lin\'eaire est non-nulle et
diagonalisable est bien connue depuis H. Poincar\'e \cite{Poin} et la th\`ese de H. Dulac
\cite{Dulac}. De mani\`ere g\'en\'erale il s'agit de comparer des 1-formes diff\'erentielles \`a coefficients
des s\'eries formelles
\begin{equation}\label{formsta}
\w = \stw\,,
\end{equation}
$$a(x , y) = \sum_{i,j} a_{ij}x^iy^j \,, \quad b(x, y) = \sum_{i,j}
b_{ij}x^iy^j \in  \C [[ x , y ]]
$$
\noindent par la relation de {\it conjugaison orbitale formelle}~: les deux 1-formes $\w_j = a_j(x, y)
dx + b_j(x, y) dy\,$, $j = 1, 2$, sont formellement conju\-gu\'ees, $\w_1 \sim_{for} \w_2\,$, s'il
existe une {\it unit\'e formelle} $u(x, y) \in
\C[[x, y]] \,$,
$u(0, 0) \not=0$ et un {\it diff\'eomorphisme formel}
$
\phi(x, y) := \left( \phi_1(x, y), \phi_2 (x, y)\right)\,$, $\phi (0, 0) = (0,0)\,$, $ det\left(
D_0\phi\right)
\not= 0\,,
$ tels que~:
\begin{equation}\label{cojint}
\w_2\, = \,u\cdot \phi^\ast \w_1 \, := \, U \cdot \left[\, a_1\left(\phi_1  , \phi_2
\right)\, d\phi_1 + b_1\left(\phi_1  ,
\phi_2  \right)\,d\phi_2 \,\right]\,,
\end{equation} o\`u $D_0 \phi$ d\'esigne la matrice jacobienne de $\phi$ \`a
l'origine. Classiquement, pour les 1-formes les plus simples, les 1-formes r\'eduites,
une liste compl\`ete des classes de conjugaison formelle est obtenue par la
construction de mod\`eles, les "formes normales", cf. \cite{Dulac}\cite{Huk}. Elle consiste \`a
choisir $\phi$ et $u$ qui annulent le plus possible de coefficients du d\'evelopement en s\'erie de la
1-forme. En fait, comme l'a mis en \'evidence J. Martinet
\cite{martbour}, elle correspond \`a une "jordanisation" de l'application adjointe agissant sur les jets
infinis de champs de vecteurs, associ\'ee au champ dual de la 1-forme. \\

Rappelons que la 1-forme (\ref{formsta}) est dite {\it r\'eduite} si la matrice
\begin{equation}\label{unjetdew}
\left(
\begin{array}{cc} a_{10}&a_{01}\\
- b_{10}&- b_{01}
\end{array}
\right)
\end{equation} est diagonalisable avec des valeurs propres $\l_1$, $\l_2$ qui
v\'erifient~: $\l_1 \not= 0\,$, $\l_2 /\l_1 \notin
\Bbb Q_{> 0}\,$. La liste des formes normales, qui sont polynomiales en les variables $x, y$ et
lin\'eaires en les param\`etres, cf. (\ref{sec.def.equising}), fait aparaitre qu'une fois fix\'e le rapport des
valeurs propres
$\l_2/\l_1$, deux invariants num\'eriques suffisent \`a la classification formelle~: un invariant
discret (un \'el\'ement de $ \N$) qui, lorsque la 1-forme $\w$ est analytique, donne la
classification topologique \cite{Ctop} et un invariant continu (un \'el\'ement de
$\C$). \\

La classe la plus simple de formes non-r\'eduites est celle des {\it 1-formes nilpotentes}~: la matrice (\ref{unjetdew})
est nilpotente. La classification de ces 1-formes a \'et\'e abord\'ee par la recherche de
formes normales, cf. F. Takens
\cite{Tak} \cite{lorpnorm} \cite{StrZol}. Le probl\`eme rencontr\'e est l'obtention d'une forme normale qui soit \`a la
fois canonique
\cite{Loray} \cite{Paulzol} et convergente. Le caract\`ere g\'en\'eriquement divergent-Gevrey de la forme normale
canonique \cite{Mir} ne permet plus de consid\'erer cette expression comme une liste de
1-formes-mod\`eles.
L'autre approche du probl\`eme, initialis\'e par D. Cerveau et R. Moussu dans
\cite{Cer-Mou} puis compl\'et\'e par R. Meziani \cite{Meziani} et d'autres auteurs \cite{BMS} consiste
\`a mettre en
\'evidence un syst\`eme complet et non-redondant d'invariants formels. L'invariant consid\'er\'e, qui
recouvre presque tous les cas, est une  repr\'esentation du groupe libre
\`a deux g\'en\'erateurs dans le groupe des diff\'eomorphismes formels d'une variable
$$
\hH_\w : \Z\a_1\star \Z\a_2 \fle \,\widehat{Diff}(\C , \, 0)  := \left\{ \left. h(z)
\in \CÊ[[\,z
\,]]\,\;\;\right|
\;  h(0)= 0
\, ,
\; h'(0)
\not= 0\right\}\,.
$$
Il s'interpr\`ete comme rendant compte de la structure transverse du feuilletage et est \'evidement de dimension
infini. La classification est donc obtenue pour cette classe de 1-formes, mais les liens entre les deux approches
demeure \`a ce jour encore tr\`es mys\-t\'erieux.\\

Nous abordons ici le cas g\'en\'eral. Notre approche consiste \`a localiser le probl\`eme "au voisinage d'une
1-forme", ensuite d\'efinir un groupe d'invariants qui g\'en\'eralisent l'invariant $\hH_\w$ du cas nilpotent
et enfin fixer ce premier groupe d'invariants et chercher une famille compl\`etant la liste. \\

La localisaton se fait en ne s'int\'eressant qu'aux
germes de familles analytiques de 1-formes formelles. Plus pr\'ecisemment, une {\it
d\'e\-for\-mation\footnote{Il serait plus correct de dire~: un germe de d\'eformation.} transversalement formelle d'une
1-forme formelle
$\w$ de base le germe
$P :=\left(
\C^p ,
\, 0
\right)$} est la donn\'ee d'une 1- forme qui s'\'ecrit
$$
\n = \std\, , \quad
\n_0 := \n|_{t=0} =  \w,
$$
$$A(x, y; t) = \sum_{i, j}A_{ij}(t)x^iy^j  \,, \quad B(x, y; t) = \sum_{i, j}B_{ij}(t)x^iy^j \in \C\{ t\} [[ x, y ]]\,,
$$
les coefficients $A_{ij}(t)$ et $B_{ij}(t) \in \C\{ t\}$ \'etant tous convergents sur un m\^eme polydisque
ouvert. Deux d\'eformations $\n$ et $\n'$ sont dites {\it t.f.-conjugu\'ees},  s'il existe dans
$\C^p$ un polydisque ouvert
$K$  cent\'e \`a l'origine et des s\'eries $\phi^1$, $\phi^2$,
$u$, $c \in \O_{\C^p} \left( K \right) [[x, y]]$, $u_0( 0, 0, 0) \not= 0\,$,  telles que~:
$$
\phi^\ast \,\n  = \, u \cdot \n' + c dt
\, ,
\quad
\phi|_{t=0} = (x, y)\,,\quad  \quad \phi  := \left(\phi^1 , \, \phi^2
\right)\,.
$$ D\'efinisons un {\it $d$-invariant formel} comme un objet $\cal I(\w)$
attach\'e \`a chaque 1-forme for\-melle $\w$, tel que
pour toute d\'efor\-ma\-tion $\n$
t.f.-conjugu\'ee
\`a la {\it d\'efor\-ma\-tion constante}
$\n^{cst} \equiv \w\,$, on a~:
$\cal I(\n|_{t=t_0}) =\cal I(\n|_{t=0})$ pour tout $t_0$ petit.\\

Avant de d\'efinir les groupes d'invariants rappelons l'existence d'une  {\it r\'educ\-tion des
singularit\'es de
$\w$} \cite{Seid} \cite{M-M}, c'est \`a dire d'une application holomorphe canonique not\'ee
$E_\w : \M_\w \fle \C^2$, obtenue comme composition finie d'applications d'\'eclatements simultan\'e de
points, telle qu'en chaque point $m$ du {\it diviseur de r\'educ\-tion}
$\D_\w:= E_\w^{-1}(0)$
l'image r\'eci\-pro\-que
$E_\w^\ast (\w)$ peut s'\'ecrire, dans des
coordonn\'ees locales $u$, $v$,  appropri\'ees~:
${E_\w}^\ast (\w) = \ u^pv^q \ww_{\,m}\,$, avec $\D_\w = \{ u^\e v
= 0 \}$,
$p, q \in \N\,$, $\e = 0$ ou
$1$, et
$\ww_{\,m}$ est une 1-forme formelle en $m$, non-singuli\`ere ou bien
singuli\`ere \`a singularit\'e isol\'ee et r\'eduite. Lorsqu'en un point $m$ de $\D_\w\,$, la
restriction du germe
$\ww_m$ au diviseur $\D_\w$ n'est pas identiquement nulle, $\w$ poss\`ede un nombre
infini de {\it s\'eparatrices formelles},  c'est
\`a dire descourbes irr\'eductibles formelles \`a l'origine de $\C^2$ dont la param\`etrisation
$\gamma (t)
\in
\C [[ t ]]^2$ satisfait $\gamma^\ast \w \equiv 0$. On dit alors que
$\w$ est {\it dicritique.} Lorsque $\w$ n'est pas dicritique, le nombre de s\'eparatrices
est fini $\geq 1$,  d'apr\`es C. Camacho- P. Sad
\cite{C-S}.\\

{\it \noindent Dans ce tout ce travail nous supposons
$\w$ non-dicritique}\\

Nous distinguerons trois types de $d$-invariants. Le premier sera constitu\'e d'invariants
discrets, ce sont les plus \'el\'ementaires; ils ne d\'ependent en g\'en\'eral que des s\'eparatrices formelles. Le
deuxi\`eme rendra compte des singularit\'es obtenues apr\`es r\'eduction; ces invariants sont en g\'en\'eral contenu dans un
espace de dimension finie. Le troisi\`eme type d'invariant est plus complexe. L'objet de ce travail est de les mettre en
\'evidence. Nous montrons en particulier qu'ils constituent g\'en\'eriquement, dans un sens pr\'ecis\'e au chapitre
\ref{tffgen} des familles de dimensions finies. Plus pr\'ecis\'ement distingons~:
\\

$\bullet$ {\bf les $d$-invariants formels globaux du processus la r\'eduction.}
Ils d\'ecrivent
seulement le processus de r\'e\-duc\-tion, en oubiant la position exacte, sur
chaque composante de $\D_\w$,  des centres d'\'ecla\-tement et des points du {\it lieu singulier
de r\'eduction}
$$
\S_\w :=\left\{ m \in \D_\w \; | \; \ww_m (m) = 0\right\}\,.
$$ Classique\-ment on construit un graphe connexe pond\'er\'e et fl\`ech\'e, {\it
l'arbre dual de r\'eduction} $\aA[\w]
$ de la mani\`ere suivante~: on se donne biunivoquement un sommet de $\aA[\w]
$ pour chaque composante irr\'eductible $D$ de
$\D_\w$ puis on le pond\`ere par l'auto-intersection de $D$ et on lui attache autant
de fl\`eches que de points de
$\S_\w\cap \left( D- Sing(\D_\w)\right) $, enfin on relie deux som\-mets par une ar\`ete chaque fois
que les composantes irr\'eductibles correspondantes s'inter\-sectent. Pour une tr\`es large classe de 1-formes, les {\it
formes de deuxi\`eme esp\`ece} d\'ecrites en (\ref{subs.red.des.sing.}), la donn\'ee de
$\aA[\w]
$
\'equivaut
\`a la donn\'ee du {\it type d'\'equisingularit\'e}
\footnote{ Le {\it type d'\'equisingularit\'e} d'une courbe formelle d'\'equation
$f(x, y) :=  \sum_{i, j = 1}^{\infty} f_{ij}x^iy^j = 0$ peut \^etre d\'efini comme le
type topologique commun \`a tous les germes de courbes analytiques, \`a l'origine de
$\C^2$, d'\'equations
$\sum_{i, j = 1}^N f_{ij}x^iy^j = 0\,$, pour $N\in \N$ assez grand.} du germe de
courbe formelle
$Sep(\w)\,$, form\'ee de toutes les s\'eparatrices formelles de $\w$.
\\

$\bullet$ {\bf Les $d$-invariants formels locaux issus de la r\'eduction.}  Remarquons d'abord que
 $E_\w^\ast (\w)$ d\'efinit une 1-forme formelle en chaque point du diviseur $\D_\w$.
De mani\`ere pr\'ecise,
$E_\w^\ast (\w)$ est une 1-forme {\it transversalement formelle}~: elle s'\'ecrit
localement $A(u, v)\, du + B(u, v)\,dv\,$, les coefficients
$A$, $B$ \'etant des sections locales du faisceau $\hO$, de base $\D_\w$, obtenu en
compl\'etant, pour la topologie
$\goth I_{\D_\w}$-adique, la restriction \`a $\D_\w$ du faisceau
$\O_{\M_\w}$ des fonctions holomorphes sur $\M_\w\,$, $\goth I_{\D_\w}\subset
\O_{\M_\w}\,$ d\'esignant le faisceau d'id\'eaux des fonctions nulles sur
$\D_\w$. Des expressions pr\'ecises de $ A$ et de $B$ sont donn\'ee en
(\ref{explfora}) et
(\ref{explforb}). On doit en retenir que $E_\w^\ast (\w)$ peut \^etre vu comme un champ,
le long de $\D_\w$, de 1-formes formelles dont les coefficients varient
holomorphiquement le long de $\D_\w$. Les sections locales de $\hO$ s'appellent {\it
fonctions transversalement for\-mel\-les}. On d\'efinit de la m\^eme mani\`ere
\footnote { ou encore \`a l'aide de l'extension des scalaire
$\,\iota^{-1}\left(\O_{\M_\w}|\D_\w\right) \fle \hO$, o\`u
$\iota :
\D_\w \ifle \M_\w$ d\'esigne l'in\-jec\-tion canonique. } les faisceaux des {\it
champs de vecteurs, formes diff\'erentielles, diff\'eomorphismes locaux...
transversale\-ment formels}. Nous \'ecrirons {\it en abr\'eg\'e t.f.} pour
transversa\-lement for\-mel. \\

La collection $\left( \ww_m \right)_{m \in \D_\w}$ d\'efini un sous-module
not\'e $\wF_\w$ du faisceau de base $\D_\w$ des 1-formes t.f., que nous appelons {\it feuilletage r\'eduit associ\'e \`a
$\w$}.  La {\it classe de conjugaison t.f. du feuilletage
$\wF_\w$ en un point
$m
\in
\D_\w\,$} est la classe de con\-ju\-gaison du germe $\ww_m$ pour la relation d'\'equivalence d\'efinie par les
\'egalit\'es (\ref{cojint}), mais avec ici
$\phi$ et
$u$ transversale\-ment formels. Nous notons $\left[
\wF_{\w , \,m} \right]_{tf}$ cette classe. La collection
$$
\widehat{L} ( \w )  :=  \left( \left[ \, \wF_{\w , \,m} \right]_{tf} \right)_{m \in\S_\w},
$$
est clairement un $d$-invariant formel de $\w$. Nous l'appelons {\it
$d$-invariant local t.f. de $\w$}. Remarquons que lorsque
$\w$ est {\it deuxi\`eme esp\`ece} (\ref{especes}) le $d$-invariant formel
$L(\w)$ est "de dimension finie~: il peut \^etre contenu dans une union finie
d'exemplaires de $\C $. Cela ressort tout simplement de l'expression des formes normales pour la
conjugaison t.f. qui, dans ces cas, sont les m\^emes que les formes normales de la conjugaison formelle
(\ref{subs.red.des.sing.}).\\

$\bullet$ {\bf Les $d$-invariants formels semi-locaux issus de la r\'eduction.} En
choisissant une courbe lisse formelle en un point
$m_0\notin \S_\w$ d'une composante irr\'eductible
$D$ de $\D_\w$, ainsi qu'une coordonn\'ee formelle sur cette courbe, on d\'efinit de
la m\^eme mani\`ere qu'en
\cite{M-M} la {\it repr\'esentation (formelle) d'holonomie de $\wF_\w$} le long
cette composante~:
$$
\H_{\w, D}  : \pi_1 ( D - \S_\w ; \, {m_0}) \fle \widehat{Diff}(\C, 0)\,.
$$ La classe $\left[ \H_{\w, D}  \right]_{for}$ de ce morphisme pour la composition
\`a gauche avec les automorphisme int\'erieurs de
$\widehat{Diff}(\C, 0)$ ne d\'epend d'aucun choix. La collection de ces classes est
visiblement un $d$-invariant formel de $\w$. \\

Finalement, en regroupant ces trois types d'invariants, nous obtenons la collection~:
$$
\widehat{SL}(\w) \; := \; \left( \; \aA[\w]  \; , \; \widehat{\cal L}(\w) \;  , \;
\left( \left[ \H_{\w, D} \right]_{tf} \right)_{D \in comp( \D_\w ) } \right)
$$

o\`u $comp( \D_\w)$ d\'esigne l'ensemble des composantes irr\'eductibles de $\D_\w$. Appelons cet invariant {\it
invariant complet d'\'equisingularite semi-locale formelle} ou encore {\it de \hsl-\'equisingularit\'e associ\'e
\`a
$\w$}.
\noindent Dans \cite{MStop} \cite{MSmr} nous montrons sous une hypoth\`ese
g\'en\'e\-rique tr\`es faible, que lorsque
$\w$ est holomorphe ces $d$-invariants formels sont des $d$-invariants
\footnote{ Nous montrons pr\'ecis\'ement que l'arbre dual de r\'eduction, le type analytique
de $\wF$ en chaque point $m \in
\S_m$, ainsi que le type analytique de la repr\'esentation d'holonomie de chaque
composante de $\D_\w$, reste constant dans chaque d\'eformation topologiquement
triviale du germe de feuilletage singulier en l'origine d\'efini par $\w$. }  topologiques.
Il est maintenant naturel de rechercher une classification \`a "type \hsl\, fix\'e". \\

Une d\'eformation t.f. $\n := \std$ de
$\w$ de base $P := \left( \C^p,\,0 \right)$ sera dite {\it \'equir\'eductible} s'il
existe une "r\'eduction en famille". Cela signifie
pr\'ecis\'ement l'existence d'une vari\'et\'e $\M_\n$ de dimension
$2 + p$,
munie d'un diviseur $\D_\n
\subset \M_\n$ \`a croisements normaux, d'un sous espace lisse $\S_\n \subset \D_\n$ de dimension
$p$, d'une application holomorphe propre $ E_\n \, : \, \M_\n  \fle\, \C^2 \times P\,$ obtenue par une succession
d'\'eclatements de courbes lisses,  avec $E_\n^{-1}(0\times P) =
\D_\n\,$, telle que la restriction de $\pi_\n :=pr_2 \circ
E_\n$ \`a $\S_\n$ est \'etale et que, pour chaque $t
\in P$ assez petit, la {\it fibre de $E_\n$ au dessus de
$t$}
$$ E_{\n, t} \; :\;  \M_{\n , t} :=   \pi_\n^{-1}(t)\, \fle \, \C^2 \times \{t\}
$$ est l'application de r\'eduction de la 1-forme $\n_t$ obtenue en restreignant $\n$
\`a $\C^2 \times \{ t\}$. De plus  la fibre du diviseur est le diviseur de r\'eduction de
$\n_t$ et la fibre de $\S_\n$  est le lieu singulier de r\'eduction de
$\n_t$~:
$$
\D_\n \cap \pi^{-1}(t) \, = \, \D_{\n_t}\,,\qquad
\S_\n \cap \pi^{-1}(t) = \S_{\n_t}\,.
$$ On voit facilement que l'espace $ \M_\n$ est un produit $C^\infty$ au voisinage de
$\D_\n$~: il exis\-te un germe de diff\'eomorphisme
$C^\infty$ au desus de $P$,
$
 \underline{\Psi} : \left(\M_\n , \,\D_\w \right) \, \iso \, \left(\M_\w \times P  ,
\,\D_\w \times
\{0\}Ê\right)\,,$ $\underline{\Psi} = \left( \Psi , \pi_\n\right)$, induisant une
famille de diff\'eomorphismes
\begin{equation}\label{diffibrr}
\Psi_t : \left( \M_{\n_t}, \D_{\n_t} \right) \fle \left( \M_\w \times \{t\}, \D_\w \times \{t\} \right)\,.
\end{equation}
De plus on a~: $\Psi_t(\S_{\n_t}) = \S_\w\times \{ t\} $.
Nous pouvons maintenant d\'efinir pr\'ecisement la notion de  {\it d\'eformation
\`a type \hsl constant}, ou encore de {\it d\'eformation \hsl-\'equisinguli\`ere}, par les
propri\'et\'es suivantes~:
\begin{enumerate}
\item $\n$ est une d\'eformation \'equir\'eductible de $\w$,
\item $\wF_{\n}$ est {\it t.f.-triviale} le long de chaque composante irr\'eductible de $ \S_\n$,
c'est
\`a dire qu'en chaque point $c$ de $\S_\w
\subset \S_\n$ il existe un diff\'eo\-mor\-phismes $\phi_{c}$ du germe
$\left( \M_{\n}, c\right)$, sur le germe $\left( \M_\w \times P, (c,0)\right)$,
fibr\'e au dessus de $P$, t.f. le long de
$\D_\n$ et tel que
$\phi_c^\ast\, \wF_{\w^{cste}}Ê\,= \wF_{\n}
$,
\item pour chaque composante irr\'eductible $D$ de $\D_\w = \pi_\n^{-1}(0) \cap
D_\n$, il existe une famille analytique de diff\'eo\-mor\-phismes formels d'une variable
$\phi_t(z)$, $\phi_t(0) = 0\,$, telle
que pour
$t$ assez petit on ait~:
$$
\tau_{\phi_t} \, \circ \hH_{\n_t \,,\, D} \,= \, \hH_{\w \, , \, D}\,,\circ \, {
\Psi_t\,}_\ast \,, \qquad
$$
avec
$\tau_{\phi_t}(g) :=   \phi_t
\circ g
\circ\phi_t^{-1}\,$, et ${\Psi_t\,}_\ast$ d\'esignant l'automorphisme induit par (\ref{diffibrr}) sur
les groupes fondamentaux consid\'er\'es.
\end{enumerate}
La "classification locale en $\w$ " sera achev\'ee si l'on obtient une d\'eformation
 $\n$ de $\w$ qui est {\it \hsl-verselle} dans le sens suivant~:
\begin{itemize}
\item $\n$ \it est \hsl-\'equisinguli\`ere
et toute d\'eformation \hsl-\'equisinguli\`ere $\n'$ est t.f. conjugu\'ee \`a une d\'eformation qui
s'\'ecrit $\underline{\lambda}^\ast \n$ o\`u $\underline{\lambda}(x, y, t)  = (x, y, \lambda (t)$ et $\lambda$ est un ger\-me
d'application holomorphe de l'espace des param\`etre de $\n$ dans celui de $(\n'$.
\end{itemize}

\vspace{1em}

$\bullet$  {\bf  Les singularit\'es de type formel fini (t.f.f.).}   Consid\'erons d'abord la situation infinit\'esimale.
Le faisceau (de base $\D_\w$) essentiel, d\'ecrit en  (\ref{sytr} ), est le quotient
$$
\hT_{\wF_\w} = \left.{\hB_{\wF_\w}}\right/{\hX_{\wF_\w}}
$$
du faisceau des germes aux points de $\D_\w$ des champs de vecteurs t.f. basiques pour $\F_\w$ par le
sous-faisceau des champs tangents \`a $\wF_\w$. On appellera {\it d\'efor\-mation infinit\'esimale de
$\w$} tout \'element de l'espace $H^1\left(\U \,;\, \hB_{\wF_\w} \right)$, o\`u $\U$ d\'esigne un
{\it recouvrement distingu\'e} du diviseur $\D_\w$, i.e. ses ouverts
sont de deux types~: l'intersection de $\D_\w$ avec un petit polydisque de coordonn\'ees centr\'e en
un point singulier de
$\wF_\w\,$, ou bien une composante irr\'edictible de $\D_\w$ \'epoint\'ee de toutes les
singularit\'es de $\wF_\w$ qu'elle portait. Les conditions 2. et 3. ci-dessus donnent, pour toute
d\'eformation \hsl-\'equisinguli\`ere $\n$ de $\w$ d'espace de param\`eres $(\C^p , 0)$, des germes $\phi_U$ de
diff\'eomorphismes t.f. le long de chaque $U \in \U$, qui trivialisent $\wF_\n$. La non-trivialit\'e au premier ordre de
$\n$ est d\'ecrite par "l'application d\'eriv\'ee \`a l'origine" suivante, explicit\'ee en (\ref{derin})~:
\begin{equation}\label{abcd}
\left[ \DD{\sF_\n}{t\phantom{\sF}} \right]_{t = 0} \, : \, T_0 \C^p \fle
H^1 \left( \U \, ; \,\hB_{\wF_\w} \right) \, , \quad
Z \efle  \left[\left(
\DD{\phi_V\circ \phi_U^{-1}}{t\phantom{\circ \phi_U^{-1}aa"}}\right)_{t = 0}\cdot\, Z\right]\,.
\end{equation}
\noindent Nous montrons en (\ref{readefin}) que toute application lin\'eaire de $T_0\C^p$ dans $H^1 \left( \U \, ;
\,\hB_{\wF_\w}
\right)$ peut \^etre r\'ealis\'ee comme la d\'eriv\'ee (\ref{abcd}) d'une d\'eformation \hsl-\'equisinguli\`ere $\n$ d'espace de param\`etre
$(\C^p , 0)$. D\'efinissons~:
\begin{itemize}
\item {\it $\w$ est de type formel fini, t.f.f. en abr\'eg\'e, si $\hb(\w) :=
dim_\C H^1
\left( \U \, ; \, \hB_{\wF_\w} \right)$ est fini}.
\end{itemize}

\noindent Nous devrons souvent nous limiter \`a la classe des 1-formes $\w$ qui satisfont la propri\'et\'e suivante~: si
tous les {\it groupes d'holonomies} $H_{D, \w} := Im(\cal H_{D, \w)}$
sont finis, alors il existe un point singulier $c \in Sing(\wF_\w)$ en lequel $\wF_\w$ ne poss\`ede pas de germe
d'int\'egrale premi\`ere t.f.. Une telle 1-forme sera dite {\it bonne}. L'\'enonc\'e suivant regroupe les th\'eor\`emes
(\ref{readefin}) et (\ref{thver})
\\

\noindent{\bf Th\'eor\`eme principal 1. }{\it Soit $\w$ bonne. Alors~:
\begin{enumerate}
\item $\w$ est t.f.f. si et seulement si $\w$ poss\`ede une
d\'eformation
\hsl-verselle;
\item une d\'eformation
$\n$ de
$\w$ est \hsl-verselle si et seulement si l'application d\'eriv\'ee $\left[ \DD{\sF_\n}{t\phantom{\sF}} \right]_{t = 0}$ est
surjective.
\end{enumerate}}

La condition "$\w$ bonne" peut vraissemblablement \^etre affaiblie, au prix de d\'evelop\-pements plus longs. Mais le
principal int\'et\^et de cette condition est de permettre un calcul simple de $\hb(\w)$. De mani\`ere pr\'ecise, en
(\ref{subsenerf}) nous associons \`a $\w$ un arbre  bicolor\'e, partiellement orient\'e et pond\'e\-r\'e $\goth
N^\ast(\w)$, appel\'e {\it nerf complet de } $\w$. Il se construit de mani\`ere simple  \`a partir seulement de l'invariant
\hsl($\w$). Notons~:
$$\widehat{\tau}(\w) := dim_\C H^1 \left( \U \,;\, \hT_{\wF_\w}\right)\,\quad \hbox{et} \quad \hd(\w) := dim_\C H^1 \left(
\U
\,;\,
\hX_{\wF_\w}\right)$$
et d\'esignons par $\hgC (\w)\,$
l'espace topologique obtenu
\`a partir de $\hgN^\ast(\w)$ en contractant et en identifiant \`a un seul et m\^eme point touts les sommets et
ar\^etes de poids 0.
Nous montrons~:\\

\noindent{\bf Th\'eor\`eme principal 2. }{\it Soit $\w$ bonne.  Alors~:
\begin{enumerate}
\item $\w$ est t.f.f. si et seulement si la partie rouge de
$\goth N^\ast(\w)$ est connexe et r\'epulsive;

\item  $\hb(\w) \leq \hd (\w) + \widehat{\tau}(\w)$ et l'\'egalit\'e est r\'ealis\'ee d\`es que
$\w$ est de deuxi\`eme esp\`ece  (\ref{especes}) sans facteur int\'egrant formel.

\item $\hd\,(\w) \; = \; \sum_{c} \frac{(\nu_c - 1) \, (\nu_c - 2) }{2}\, $
o\`u $c$ parcourt l'ensemble de tous les points singuliers
apparaissant  dans la r\'eduction (\ref{ard}) de $\w\,$ et
$\nu_c$ d\'esigne la
multiplicit\'e de l'\'eclat\'e divis\'e de $\w$ au point $c$.

\item $\widehat{\tau}(\w) = rang_\Z H_1 \left( \hgC (\w)\,;\, \Z \right)$.
\end{enumerate} }

\noindent Lorsque $\w$ est de deuxi\`eme esp\`ece
$\hd (\F )$  est \'egal \`a la dimension de la strate \`a
$\mu$-cons\-tant de l'union des s\'eparatrices formelles de $\w$, c.f.
(\ref{minmult}).\\

Ce proc\'ed\'e combinatoire permet aussi d'expliciter une base de $H^1 \left( \U \,;\, \hT_{\wF_\w}\right)$. Cette
pr\'ecision permet de voir que pour toute d\'eformation \hsl-\'equisinguli\`ere $\n$ d'une 1-forme bonne $\w$ de base
quelconque
$(\C^p, 0)$, "les familles de d\'eformations infinit\'esimales" forment un $\hO_{\C^p, 0}$-module de type fini. Plus
pr\'ecis\'ement d\'e\-signons par  $\hB_{\wF_\n}^v$ le faisceau de base $\D_\w \subset \D_\n$ des germes aux points de
$\D_\w$ de champs de vecteurs de $\M_\n$ qui sont verticaux et basiques pour $\wF_\n$.\\

{\bf Th\'eor\`eme de pr\'eparation.} {\it Si $\w$ est bonne et t.f.f., alors $H^1 \left( \U \,;\, \hB_{\wF_\n}^v\right)$ un
$\hO_{\C^p, 0}$-module de type fini, qui est libre lorsque $\w$ ne poss\`ede pas de facteur
int\'egrant formel.}\\

\noindent Ce th\'eor\`eme enonc\'e en (\ref{enimdir}) et montr\'e en (\ref{primdir}) est un ingr\'edient
indispensable
\`a la d\'emonstration du th\'eor\`eme principal 1. Sa d\'emonstration ne peut pas se faire par des arguments de
g\'eom\'etrie analytique, car le faisceau
$\hB_{\wF_\n}^v$ n'est pas un faisceau de modules sur le faisceau des fonctions t.f.. Elle utilise le caract\`ere
constructif de la d\'emonstration du th\'eor\`eme principal 2, qui permetd'expliciter des bases de d\'eformations infinit\'esimales.
Ainsi le sh\'ema logique de ce papier est~:
$$
\hbox{\it Th\'eor\`eme principal 2}\; \Longrightarrow \;\hbox{\it Th\'eor\`eme de pr\'eparation}
\; \Longrightarrow \; \hbox{\it Th\'eor\`eme principal 1}\;.
$$

Nous examinons aussi la g\'en\'eralit\'e de ces r\'esultats. Nous d\'ecrivons d'abord l'ensemble $\widehat{\goth E}^{2}$ des 1-formes
formelles de deuxi\`eme esp\`ece, en le stratifiant par "strates
d'\'equisingularit\'e". Chaque srate est  pro-constructible de codimension finie, cf. th\'eo\-r\`eme (\ref{prstrpred}). Nous
montrons ensuite le th\'eor\`eme suivant qui est pr\'ecis\'e en (\ref{genexpl})~:\\

\noindent{\bf Th\'eor\`eme principal 3. }{\it L'ensemble $\widehat{\goth E}^2$  et le sous ensemble $\widehat{\goth
E}_{tff}$ des $\w \in \widehat{\goth E}$ qui sont t.f.f
sont des ouverts de la topologie de Krull. De plus $\widehat{\goth E}^2_{tff}$ est Krull-dense dans $\widehat{\goth E}^2$.}\\

\centerline{\bf  Plan de l'article}

\vspace{1,5em}

{\bf Chapitre \ref{sec.arbres}.} Les techniques de construction de d\'eformations, de type
Kodaira-Spencer, vont consister \`a decouper
$\wF_\w$ le long des ouverts d'un recouvrement adapt\'e $\U$,  puis \`a effectuer des recollements t.f. le long des
intersections de ces ouverts. La premi\`ere obstruction rencontr\'ee vient du caract\`ere formel de la vari\'et\'e ainsi
obtenue. L'objet de ce chapitre est de lever cette obstruction~: les {\it th\'eor\`emes de stabilit\'e (\ref{mth.co.})} et
de {\it d\'etermination finie} (\ref{d.q.fin.})
 expriment, "avec contr\^ole des jets", la propri\'et\'e suivante
\begin{itemize}
\item \it toute vari\'et\'e formelle obtenue \`a partir de $\M_\w$
par recollement  t.f. est t.f. iso\-morphe \`a une vari\'et\'e holomorphe, qui est la c\^ime d'un arbre de m\^eme arbre
dual que
$\w$.
\end{itemize}
Ce chapitre est technique. Nous sugg\`erons de limiter une
premi\`ere lecture au vocabulaire n\'ec\'essaire  \`a l'ennonc\'e des
th\'eor\`emes (\ref{mth.co.}) et (\ref{d.q.fin.}).\\

{\bf Chapitre \ref{sec.def.feuil}.}  Nous d\'eveloppons une technique tr\`es g\'en\'erale de construction de
d\'efor\-ma\-tions
\hsl-\'equisinguli\`eres. La notion cl\'e est celle de {\it syst\`eme semi-local coh\'erent} (\ref{coherent}). Elle permet
d'\'etablir un {\it th\'eor\`eme de r\'ealisation} de d\'eformation \hsl-\'equi\-singuli\`eres (\ref{realis.}), avec controle
des jets de la 1-forme d\'efinissant la d\'eformation.\\

{\bf Chapitre \ref{sec.def.equising}.}  Nous d\'egageons l'importante classe des 1-formes {\it de deuxi\`eme
esp\`ece}, pour lesquelles l'equir\'eduction d'une d\'eformation est \'equivalente \`a la constance de
l'arbre dual pond\'er\'e par les auto-inter\-sec\-tions, cf. le
th\'eor\`eme (\ref{deuxequ}). Le th\'eor\`eme (\ref{cardeuxesp}) donne des condition \'equivalentes tr\`es simples qui
caract\'erise cette classe.\\

{\bf Chapitre \ref{sec.e.sing}.}  Nous introduisons la notion de \hsl-\'equisingularit\'e, de d\'eformation infinit\'esimale et,
apr\`es avoir donn\'e des techniques de construction de d\'eformations, nous prouvons le th\'eor\`eme principal 1 en admettant le
th\'eor\`eme de pr\'epa\-ration.\\

{\bf Chapitre \ref{carsintff}.} Nous prouvons le th\'eor\`eme principal 2, puis le th\'eor\`eme de pr\'epa\-ration. \\

{\bf Chapitre \ref{tffgen}.} Nous stratifions l'espace de formes de deuxi\`eme esp\`ece et nous montrons le th\'eor\`eme
principal 3.

\section{Pr\'eliminaires sur les arbres d'\'eclatements}\label{sec.arbres}
L'objet de ce chapitre est de prouver des propri\'et\'es de "stabilit\'e" et de
"d\'eter\-mi\-na\-tion
finie" pour des espaces  obtenus \`a partir de $\C^2$
par une succession d'\'eclatements ponctuels au dessus de 0,  Th\'eor\`emes
(\ref{d.q.fin.}) et
(\ref{mth.co.}).
\subsection{Arbres et arbres duaux} \label{A.A.dual.}
\addcontentsline{toc}{section}{\hspace{0,8em} {}\thesubsection .  Arbres et arbres duaux}
Nous appelons ici {\it arbre au dessus d'une vari\'et\'e holomorphe connexe
$Q$ de dimension $p$} la donn\'ee
d'un diagramme commutatif
$\A_Q$~:
\begin{equation}\label{diag.arbre}
\begin{array}{ccccccccc}
\M^h & \fle \,\cdots \, \fle & \M^j & \stackrel{E^j}{\fle} & \M^{j-1} &
\fle \cdots
\stackrel{E^1}{\fle} & \M^0 & \stackrel{\pi}{\fle} & Q \\
\bigcup & & \bigcup &  &\bigcup & & \bigcup & & \\
\S^h & \fle \,\cdots \, \fle & \S^j & \fle & \S^{j-1} & \fle \,\cdots \,
\fle &\S^0 & & \\
\bigcup & & \bigcup &  &\bigcup & & \bigcup & & \\
S^h & \fle \,\cdots \, \fle & S^j & \fle & S^{j-1} & \fle \,\cdots \, \fle
&S^0 & &
\end{array}
\end{equation}
o\`u, pour chaque $j = 0,\ldots ,h\,$ : $\M^j$ est une vari\'et\'e
analytique complexe lisse de dimension $2 + p\,$ appel\'ee  $\hbox{\it
j}^{\hbox{\it \`eme}}$
{\it espace \'eclat\'e}, $\S^j$ est un sous-ensemble analytique ferm\'e de
$\M^j$ de dimension $p$ appel\'e
$\hbox{\it j}^{\hbox{\it \`eme}}$ {\it lieu de singularit\'es}, $S^j$ est
une sous-vari\'et\'e analytique lisse ferm\'ee (non n\'ecessairement connexe) de
$\S^j$ de dimension $p$
appel\'e $\hbox{\it j}^{\hbox{\it
\`eme}}$ {\it centre d'\'eclatement}, \'eventuellement $S^h$ peut \^etre vide, et
tels que, en notant
$$E_j := E^0 \circ \cdots \circ E^j\,,\quad\pi_j
:= \pi \circ E_j\,
, \quad \D^j := E_j^{-1}(S^0),,$$
avec $j = 0, \ldots ,h$, on ait~:
\begin{enumerate}
\item[0.] $\pi$ est une submersion,
\item chaque $E^{j+1}$ est l'application d'\'eclatement de centre $S^j\,$,
\item chaque $S^j$ est une union de composantes irr\'eductibles de $\S^j\,$
et $\S^0 =
S^0\,$,
\item chaque $\S^j$ est contenu dans $\D^j$ et la restriction de $\pi_j$ \`a
$\S^j$ est propre \`a fibres
finies,
\item pour chaque $j = 1, \ldots, h\,$ la restriction de $\pi_j$ \`a chaque
com\-po\-sante con\-nexe de
$S^j$ est  un biholomorphisme sur $Q\,$.
\end{enumerate}
\begin{defi}\label{arreg}
\noindent Nous dirons que {\it l'arbre $\A_Q$ est  r\'egulier} si chaque
$\S^j$ est lisse et la
restriction de $\pi_j$ \`a chaque com\-po\-sante con\-nexe de $\S^j$ est  un
biholomorphisme sur
$Q\,$.
\end{defi}

\noindent Nous noterons
\begin{equation}\label{not.arbre}
\A_Q = (\M^j\,, E^j\,, \S^j\,, S^j\,, \pi_j\,, \D^j\,)_{j = 0,\ldots ,h}\;
.
\end{equation}
L'entier $h$
s'appelle la {\it hauteur de l'arbre}, la vari\'et\'e $\M^0\,$  son  {\it socle},
$\M^h\,$ sa {\it cime}
et
$\D^j\,$ le {\it $\hbox{\it j}^{\hbox{\it \`eme}}$ diviseur exceptionnel}.
Pour
simplifier nous d\'esignerons les {\it donn\'ees de cime} de la
mani\`ere suivante~:
\begin{equation}\label{cime}
\wM := \M^h\,,\ \ \ \wS := \S^h\,,\ \ \ \wD := \D^h\,,\ \ \ \wE := E^h\,,\
\ \ \wp :=
\p^h\,.
\end{equation}

On d\'efinit {\it l'image inverse}
de $\A_Q$ par une application holomorphe $\l$ d'une vari\'et\'e holomorphe
$Q'$ dans $Q\,$, comme le diagramme
$\l^{\star}\A_Q$ obtenu de mani\`ere naturelle \`a partir de
(\ref{diag.arbre}) par produits fibr\'es de $\l$
et des $\pi_j\,$. C'est un arbre au dessus de
$Q'$. En particulier lorsque
$Q'$ se r\'eduit \`a un seul point $t \in Q$ cette op\'eration donne {\it la
fibre de $\A_Q$ au dessus de $t$},
c'est \`a dire l'arbre not\'e~:
\begin{equation}\label{arbrfibr}
\A_Q(t) = (\M^j_t\,, E^j_t\,, \S^j_t\,, S^j_t\,, \pi_{j, t}\,,
\D^j_t\,)_{j = 0,\ldots ,h}\; ,
\end{equation}
qui s'obtient en restreignant les espaces et les applications du diagramme
(\ref{diag.arbre}) aux fibres
$\pi_j^{-1}(t)\,$. En particulier $\S^0_t = S^0_t$ est r\'eduit \`a un point.\\

Fixons un point $t_0$ de $Q\,$. Nous appellerons {\it germe de $\A_Q$ au
dessus de $\left( Q, t_0
\right)$} le diagramme obtenu \`a partir de (\ref{diag.arbre}) en rempla\c cant
pour chaque
$j = 1, \ldots, h\,$
les es\-paces et les
applications par leurs germes le long des diviseurs exceptionnels
$\D^j_{t_0} \subset \D^j$ de
$\A_Q(t_0)\,$
et, pour $j=0\,$, en rempla\c cant $\M^0$ et $\S^0 = S^0$ et par leurs germes
au point
$\left\{ m_0Ê\right\} := \S^0_{t_0} = S^0_{t_0}\,$.
\\

Deux germes d'arbres au dessus de $P := \left( \C^p , 0 \right)\,$ de m\^eme
hauteur $h\,$
$$\A_P = (\M^j\,, E^j\,, \S^j\,, S^j\,, \pi_j\,, \D^j\,)_j\,,\quad
    \A'_P = ( {\M'}^j\,,  {E'}^j\,,  {\S'}^j\,, {S'}^j\,,  {\pi'}_j\,,  {\D'}^j\,)_j$$
sont dits {\it isomorphes}, resp. {\it hom\'eomorphes}, s'il existe le long
de chaque diviseur
$\D^j_0$ de $\A_P(0)$ des germes de biholomorphismes, resp.
d'hom\'eomorphismes
$$\f_j : (\M^j, \D^j_0\,)\,\fle\,( {\M'}^j,  {\D'}^j_0\,)\,,\ \ j = 0,\ldots ,h$$
qui envoient lieux singuliers sur lieux singuliers, centres sur centres et
commutent aux
applications d'\'eclatement et respectent les fibrations sur $P\,$, i.e.
$$\f_j(\S^j) =  {\S'}^j\,,\quad\f_j(S^j) = {S'}j\,,\quad {E'}^j\circ \f_j =
\f_{j-1}\circ E^j\,,\quad
\p'_j\circ\f_j = \p_j\,,$$
avec $j=0,\ldots,h\,$. Nous noterons~:
$$\f_{_{\bullet}}\; :\; \A_P \stackrel{\sim}{\rightarrow} {\A_P'}\, , \qquad
\f_{_{\bullet}}\; := \; (\f_j)_{j=0,\ldots,h}\,.$$

Un {\it arbre dual} $\,\A^{\ast}\,$ est un graphe fini connexe, sans
cycle, pond\'er\'e en associant
un nombre entier \`a chaque sommet et fl\'ech\'e en associant des fl\`eches en nombre fini
\`a certains sommets. Le nombre $v({\sf s})$ de fl\`eches et
d'ar\^etes attach\'ees \`a un sommet
${\sf s}$ s'appelle la {\it valence de $\,\sf  s\,$} et l'ensemble
$Ad({\sf s})$ de ces fl\`eches et ar\^etes
s'appelle {\it ensemble des \'el\'ements adjacents \`a $\,\sf s\,$}. De m\^eme,
l'ensemble (constitu\'e de 2, resp.
1 \'el\'ements) des sommets sur lesquels s'attache une ar\^ete, resp. une
fl\`eche, not\'ee $\,\sf b\,$, se note
$\,Ad( {\sf b})\,$ et est appel\'e {\it ensemble des \'el\'ements adjacents \`a
$\,{\sf  b}\,$}.

\begin{defi}\label{adu}
Un arbre dual $\A^\ast$ est dit {\it associ\'e \`a un germe d'arbre $\A_P$} s'il
existe une correspondance
biunivoque entre les composantes irr\'eductibles du diviseur $\wD$ de la
cime $\wM$ de $\A_P\,$ et les
sommets de
$\A^\ast\,$,
$${\sf D} \rightarrow s({\sf D}) \, \quad \hbox{ou} \quad D({\sf s})
\leftarrow {\sf s}\, ,$$
telle que~:
\begin{itemize}
\item deux sommets $\sf s$ et $\sf s'$ de $\A^{\ast}\,$ sont reli\'es par
une ar\^ete si
$D ({\sf s}) \cap D ({\sf s'}) \not= \emptyset\,$,
\item le nombre de fl\`eches attach\'ees \`a un sommet $\sf s$ de $\A^\ast$ est
\'egal au nombre de
composantes irr\'eductibles de $D({\sf s}) \cap \wS$ qui ne sont pas des
composantes irr\'eductibles du
lieu singulier de
$\wD\,$,
\item chaque sommet $\sf s$ de $\A_P^{\ast}\,$ est pond\'er\'e par la
premi\`ere
classe de Chern
$e({\sf s}) \in H^2(\,D({\sf s})\,, \Z) \simeq \Z\,$ du fibr\'e normal \`a $D
({\sf s})\,$ dans $\wM\,$.
\end{itemize}
\end{defi}

Un tel arbre dual est unique. Il se notera $\A^\ast_P\,$.
Rappelons sans d\'emonstration quelques propri\'et\'es bien connues.
Soient $\A_P$, et $\A'_P$ deux germes d'arbres et notons $\wM$, resp.
$\wM'$ leurs cimes, $\wD$, resp. $\wD'$ leurs diviseurs de cime,  $\wS$,
resp. $\wS'$ et $\widetilde{S}$, resp. $\widetilde{S'}$ leurs lieux de
singularit\'es et leurs centres de cime. D\'esignons aussi par
$\wD_0\,$, resp. par
$\wD'_0$ les diviseurs de cimes de
$\A_P (0)$, resp. de
${\A_P'} (0)\,$.
\begin{prop}\label{isom.arbr}
Soit $\wf : (\wM ,\, \wD_0)
\stackrel{\sim}{\rightarrow} (\wM' ,\,\wD'_0)\,$ un  germe de
biholomorphisme, resp. d'hom\'eomorphisme au
dessus de $P\,$ v\'erifiant~:
\begin{equation}\label{isocime}
\wf ( \wD ) = \wD' \, , \qquad \wf\,(\wS) = \wS' \, \qquad
\wf(\,\widetilde{S}) =
\widetilde{S'}\,.
\end{equation}
Alors $\A_P$ et ${\A_P'}$ ont m\^eme hauteur $\,h\,$ et
$\wf$ induit un isomorphisme, resp. un hom\'e\-omor\-phis\-me
de germes d'arbres  $\f_{_{\bullet}} : \A_P \iso \A'_P\,$ tel que $\f_h =
\wf\,$.
\end{prop}

\begin{prop}\label{iso.arb.}
Supposons $\A_P\,$ et ${\A_P'}\,$ r\'eguliers. Les as\-sertions suivantes
sont \'equivalentes~:
\begin{enumerate}
\item  $\A_P$ et $\A_P'$ sont hom\'eomorphes,
\item  Il existe un germe
$\wf : (\wM ,\, \wD_0) \stackrel{\sim}{\rightarrow} (\wM' ,\,\wD'_0)\,$
d'hom\'eomorphisme au dessus de $P\,$ qui satisfait
(\ref{isocime}),
\item  $\A_P^{\ast} = {\A_P'}^{\ast}\,$.
\end{enumerate}
\end{prop}

%%%%%%%%%%%%%%%%%%%%%%%

\subsection{Arbres infinit\'esimaux et rel\`evements}\label{arbr.inf}
\addcontentsline{toc}{section}{\hspace{0,8em} {}\thesubsection .  Arbres infinit\'esimaux et rel\`evements}
Arbres et arbres duaux
Une pond\'eration \'equi\-va\-lente des sommets de l'arbre dual d'un germe d'arbre
$\A_P$ de hauteur $h\,$ au
dessus de $P\,$, not\'e encore (\ref{diag.arbre}), est la {\it "pond\'eration
alg\'ebrique"} $D \efle m(D)$. Avant de la d\'efinir, fixons
quelques notations.  Pour chaque niveau $j = 0,\ldots ,h\,$
d\'esignons par~:
\begin{quotation}
\item[$\O_{(j)}$] la restriction au j-\`eme diviseur
exceptionnel $\D^j_0 \subset \D^j \subset \M^j\,$ de l'arbre $\A_P(0)$, du
faisceau $\O_{\M^j}$
des germes de
fonctions holomorphes de $\M^j\,$, i.e.
$\O_{(j)}\,= \iota_j^{-1}(\,\O_{\M^j}\,)\,$ o\`u  $\iota_j : \D^j_0 \ifle
\M^j\,$ est l'application d'inclusion; pour $j =0\,$, $\pi^{-1}(0) \cap
\S^0$ est r\'eduit \`a un point
$m_0$ et
$\O_{(0)}$ est l'anneau $\O_{\M^0, \, m_0}$ des germes de fonctions
holomorphes en ce point,
\item[$I_j$ ] le faisceau d'id\'eaux de $\O_{(j)}$ form\'e des germes qui
s'annulent sur le diviseur $\D^j\,$, $j \geq 1\,$,
\item[$I_0$ ] l'id\'eal de l'anneau $\O_{(0)}$ form\'e des germes de fonctions
qui s'annulent sur $\S^0 =
S^0\,$,
\item[$\J_{(j)}$] le faisceau d'id\'eaux ${E_j}^{\ast} (I_0)
\subset \O_{(j)}\,$,  qui est localement engendr\'e par les fonctions du
type $f \circ E_j\,$,
$f\,$ appartenant \`a  l'id\'eal $I_0 \subset \O_{\M^0}\,$, $\,\J_{(0)} =
I_0\,$.
\end{quotation}
\noindent Pour une r\'eunion finie $ {\D'}\,$ de composantes irr\'eductibles de
$\D^j\,$, on notera
\begin{quotation}
\item[$I_{ {\D'}}$] le sous-faisceau d'id\'eaux de $\O_{(j)}\,$ form\'e des
germes qui
s'annulent sur  $ {\D'}\,$.
\end{quotation}

\noindent Enfin, pour simplifier, nous noterons~:
\begin{equation}\label{id.arbres}
\wJ := \J_{(h)}\,,\qquad \wO := \O_{(h)}\,\qquad
\widetilde{I} := I_{\wD}\,.
\end{equation}

\noindent Visiblement $\J_{(j)}$ se d\'ecompose de la mani\`ere suivante~:
\begin{equation}\label{decomp.id.}
\J_{(j)} =: \prod_{D\in comp(\D^j)} {I_{D}}^{m\,(D)}\,,\qquad j = 1,\ldots
, h\,,
\end{equation}
 o\`u $comp(\D^j)\,$ d\'esigne la collection des composantes irr\'eductibles de
$\D^j\,$ et les
exposants $m\,(D)\,$ se calculent par le proc\'ed\'e tr\`es simple suivant qui
ne d\'epend que de
l'arbre dual associ\'e \`a $\A_P\,$~:
\begin{itemize}
\item  $m\,(\D^1) = 1\,$,
\item  si $D' \subset \D^j\,$ est le transform\'e strict d'une composante
$D\,$ de
$\D^{j-1}\,$, alors $m\,(D') = m\,(D)\,$,
\item  si $D'\,$ est cr\'e\'e par un \'eclatement de centre connexe $C \subset
S^{j-1}\,$, alors
$m\,(D')\,$ est la somme des "multiplicit\'es" $m\,(D)\,$ des composantes
irr\'eductibles
$D\,$ de $\D^{j-1}\,$ qui contiennent $C\,$.
\end{itemize}

\begin{lemm}\label{cohom.ecl.}
Pour tout $j = 0,\ldots ,h$, les propri\'et\'es suivantes sont v\'erifi\'ees~:
\begin{enumerate}
\item les applications $E_j\,$ induisent les isomorphismes~:
$$H^0(\D^j_0;\: \J_{\,(j)}) \simeq  I_0 \simeq H^0(\D^j_0;\: I_{\D^j})
\,.$$
\item a) $H^1(\D^j_0; \: \O_{(j)}) = 0\,$,\phantom{MMM}  b) $H^1(\D^j_0;
\: \J_{\,(j)}) = 0\,$,
\item a) $H^2(\D^j_0; \: \O_{(j)}) = 0\,$,\phantom{MMM}  b)  $H^2(\D^j_0;
\: \J_{\,(j)}) = 0\,$,
\end{enumerate}
\end{lemm}
\begin{proof}
Les propri\'et\'es 1) d\'ecoulent directement du th\'eor\`eme classique d'Hartogs.
La d\'emonstration de 2.a) se fait ais\'ement par r\'ecurrence sur la hauteur $h\,$ de
l'arbre; au
premier cran c'est un calcul simple de s\'eries de Laurent. L'induction se
fait \`a l'aide de la
suite exacte de Mayer-Vietoris.
La propri\'et\'e 2.b) d\'ecoule de 2.a) car les faisceaux $\J_{(j)}\,$
sont engendr\'es par leurs sections globales.
Les propri\'et\'es 3.a) et 3.b) r\'esultent de l'existence de recouvrements de
Stein
sans intersection trois \`a trois.
\end{proof}

Consid\'erons pour chaque $j = 0, \ldots , h\,$ et $k \in \N$ les espaces
annel\'es suivants. Ils
s'interpr\`etent en consid\'erant les $k$-i\`emes voisinages infinit\'esimaux
\footnote{
Remarquons que la notion de voisinage infinit\'esimal est ici relative \`a la
filtration
par les puissances du faisceau d'id\'eaux
$\J_{(j)}$ de $\wO_{(j)}\,$.
}
de chaque diviseur
$\D^j$ dans
$\M^j\,$, puis en germifiant le long du $j$-i\`eme diviseur $\D^j_0 \subset
\D^j$ de $\A_P(0)\,$.

\begin{equation}\label{ar.inf.def.}
{\M^{j}}^{[k]} := \left(\, \left| {\M^{j}}^{[k]}\right| := \D^j_0\,;\quad
\O_{{\M^{j}}^{[k]}} := \O_{(j)} /\J_{\,(j)}^{k + 1}\,\right)\,.
\end{equation}
\noindent D'apr\`es ce qui pr\'ec\`ede les germes des applications d'\'eclatement
$E^j$ ainsi que
$\p\,$ se factorisent
en des morphismes d'espaces annel\'es~:
$${E^j}^{[k]} : {\M^{j}}^{[k]} \fle {\M^{j}}^{[k-1]} \,  \quad
\hbox{et}\quad
\pi^{[k]} : {\M^{0}}^{[k]} \fle P\, .$$
\noindent On obtient ainsi un
{\it germe d'arbre infinit\'esimal au dessus de $P\,$}, c'est \`a dire le
diagramme commutatif~:
$$
\begin{array}{ccccccccccc}
\wM^{[k]} & := & {\M^{h}}^{[k]} & \fle \,\cdots \, \fle & {\M^{j}}^{[k]} &
\stackrel{{E^j}^{[k]}}{\fle} & {\M^{j-1}}^{[k]} & \fle \cdots
\stackrel{{E^1}^{[k]}}{\fle} &
{\M^{0}}^{[k]} & \stackrel{\pi^{[k]}}{\fle} & P \\
 & & \bigcup & & \bigcup &  &\bigcup & & \bigcup & & \\
 & & \S^h & \fle \,\cdots \, \fle & \S^j & \fle & \S^{j-1} & \fle \,\cdots
\, \fle &\S^0 & & \\
& & \bigcup & & \bigcup &  &\bigcup & & \bigcup & & \\
& & S^h & \fle \,\cdots \, \fle & S^j & \fle & S^{j-1} & \fle \,\cdots \,
\fle &S^0 & &
\end{array}
$$
Nous
d\'esignons par $\A_P^{[k]}\,$ ce diagramme et notons~:
$\,\pi_j^{[k]} := \pi^{[k]} \circ {E_j}^{[k]}\,$.

\begin{lemm}\label{descente}
Consid\'erons une composante irr\'eductible $D\,$ de
$\D^j\,$ cr\'e\'ee par l'\'ecla\-te\-ment d'une composante irr\'eductible $S$
de $S^{j-1}\,$ et notons $D_0 := D \cap
\D_0^j\,$ et $\{ s_0 \} := S \cap \D_0^{j-1}\,$. Soit
$\tau : \left({\M^{j}}^{[k]}, \, D_0\right) \rightarrow \left(\C^2\times P
,\, 0\right)$
un germe le long de
$D_0\,$ de morphisme au dessus de $P$ dont l'application sous-jacente
envoie $D$ sur $\{ 0 \} \times
P\,$.  Alors $\tau\,$ se factorise \`a travers
${E^j}^{[k]}\,$ en un germe de morphisme
$$\tau^{\flat} : \left({\M^{j - 1}}^{[k]}, \,
s_0\right) \fle \left( \C^2\times P,\, 0 \right)$$
au dessus de $P$ qui induit un germe de
biholomorphisme de $\left(S, \, s_0 \right)\,$ sur  $\left( \{0\} \times
P, \, 0
\right)\,$.
\end{lemm}
\begin{proof}
Le morphisme $\tau\,$ est enti\`erement d\'etermin\'e par les images
$$f := \tau^{\ast}(x)\,,\ g := \tau^{\ast}(y)\ \in H^0\left( D_0 \,;\,
\O_{{\M^{j}}^{[k]}}\right)\,$$
des deux premi\`eres coordonn\'ees de $\C^2\times P\,$ par le co-morphisme
$$ \tau^{\ast} : \O_{\C^2\times P , \, 0} \fle
\tau_{\ast}\left( \O_{{\M^{j}}^{[k]}}\right)\,.$$
On d\'eduit facilement des \'egalit\'es 2.b) et 3.b) de (\ref{cohom.ecl.}) que
$f$ et $g$ sont induites par
des sections globales $F\,$, $G \in H^0\left(D_0\,;\,\O_{\M^{j}}\right)\,$
nulles sur $D\,$. Par le
th\'eor\`eme d'Hartogs elles se factorisent en des germes
$F^{\flat}\,$, $G^{\flat} \in \O_{\M^{j-1}, \, s_0}\,$ nuls sur
$S\,$. On d\'efinit alors $\tau^{\flat}\,$ en prenant respectivement pour
${\tau^{\flat}}^{\ast}(x)\,$ et
${\tau^{\flat}}^{\ast}(y)\,$ les restrictions de $F^{\flat}\,$ et de
$G^{\flat}\,$ \`a
${\M^{j-1}}^{[k]}\,$. \end{proof}

La proposition suivante peut \^etre consid\'er\'ee comme une version
infini\-t\'e\-simale de la pro\-po\-sition
(\ref{isom.arbr}).

\begin{prop}\label{desc.morph.}
Soit
$\widetilde{\phi} : \wM^{[k]} \fle {\wM'}{}^{[k]}\,$ un morphisme entre les
cimes de deux germes
d'arbres infinit\'esimaux $\A_P^{[k]}$ et ${\A'}_P^{[k]}\,$. Supposons que
$\widetilde{\phi}$
induise un germe
le long de $\wD_0\,$, not\'e $| \wf |\,$, de biholomorphisme entre les
diviseurs de cime
$\wD\,$ et
$\wD'\,$ et v\'erifie (\ref{isocime}). Il existe alors un isomorphisme de
germes d'arbres infinit\'esimaux
${\phi}_{\bullet} : \A_P^{[k]} \fle {\A'_P}^{[k]}\,$
tel que ${\phi}_h = \widetilde{\phi}\,$, o\`u $h$ d\'esigne la
hauteur (commune) de ces arbres.
\end{prop}

\begin{proof}
L'application $|\widetilde{\phi} |\,$ transforme les composantes de $\wD$
de poids $-1$, au sens de
(\ref{A.A.dual.}),  en les composantes de poids $-1$ de $\wD'\,$, car le
voisinage infinit\'esimal
d'une composante d\'etermine son fibr\'e normal. Ainsi
$\Psi := { {E'}^h}^{[k]}\circ |\widetilde{\phi}|\,$ envoie sur ${S'}{h-1}\,$
la r\'eunion  $\wD" :=
\left({E^h}\right)^{-1}(S^{h-1})\,$ des derni\`eres composantes cr\'e\'ees de
$\wD\,$. Comme
$ {\M'}^{h-1}\,$ est analytiquement trivial le long de chaque composante $S'$
de ${S'}{h-1}\,$,
$$( {\M'}^{h-1}, \,S') \stackrel{\sim}{\fle} (\C^2\times P, \:0\times P)\,,$$
d'apr\`es le th\'eor\`eme de contraction de Grauert \cite{B-S},
on peut appliquer le lemme pr\'ec\'edent \`a $\Psi\,$ et, au voisinage de $\D"$,
factoriser
$\widetilde{\phi}$ dans ${\M^{h-1}}^{[k]}$ par ${E^h}^{[k]}$. Sur $\wM -
\wD"$ le
morphisme $E^h$ est un isomorphisme et la factorisation de
$\widetilde{\phi}$ existe trivialement.
Elle est unique et se recolle ainsi avec celle construite au voisinage de
$\D"\,$.
D'o\`u la conclusion par r\'ecurrence sur $h\,$.
\end{proof}

Les diff\'eomorphismes au dessus de $P$ entre les socles $\M^0$ et ${ {\M'}}^0$
de deux germes d'arbres
$\A_P$ et $\A'_P$ ne se rel\`event pas en g\'en\'eral aux espaces \'eclat\'es $\M^j$
et
${ {\M'}}^j\,$. Cependant, lorsqu'un relev\'e existe, il est unique. Cette
derni\`ere propri\'et\'e est en d\'efaut
lorsqu'il s'agit de relever un isomorphisme ${\M^0}^{[k]} \iso
{{ {\M'}}^0}^{[k]}\,$ de voisinages
infinit\'esimaux.
\\

De mani\`ere g\'en\'erale consid\'erons un biholomorphisme $f$ entre deux vari\'et\'es
lisses au
dessus d'une vari\'et\'e  $Q$,
$$f  : M \iso M'\, ,\quad \pi : M\fle Q\, , \quad   {\pi'} : M' \fle Q\,.$$
Supposons que $\pi$ et $ {\pi'}$ sont des submersions et que l'on a
$ {\pi'} \circ f = \pi \,$. Un calcul direct en coordonn\'ees permet de voir
que $f$ se rel\`eve de mani\`ere
unique \`a travers des \'eclatements
$$E : \widetilde{M} \fle M \, , \qquad  {E'} : \widetilde{M'} \fle M'$$
de centres lisses connexes, respectivement $C \subset M$ et $C' \subset
M'\,$, d\`es que $C$ et $C'$ sont
isomorphes \`a
$Q$ via $\pi$ et $ {\pi'}$  et que $f(C) = C'\,$. De plus le rel\`evement
$\widetilde{f} :
\widetilde{M} \fle \widetilde{M'}$ est unique et est un biholomorphisme.
On voit facilement que le
$k$-jet, $j^k_D(\widetilde{f})\,$, de $\widetilde{f}$ le long de $D :=
E^{-1}(C)$ ne d\'epend que du
$(k+1)$-jet,
$j^{k+1}_C(f)\,$, de $f$ le
long de
$C\,$. Cependant $j^k_C(f)$ ne suffit pas \`a determiner
$j^k_D(\widetilde{f})$.
Ainsi l'application $f^{[k]} : M^{[k]} \fle {M'}^{[k]}$ induite par $f$
sur les voisinages
infinit\'esimaux d'ordre $k$ de $C$ et $C'$ admet une infinit\'e de relev\'es
$\widetilde{f}^{[k]} : \widetilde{M}^{[k]} \fle \widetilde{M'}^{[k]}\,$
sur les
voisinages infit\'esimaux d'ordre $k$ de $D$ et $D':= { {\pi'}}^{-1}(C')\,$
dont les restrictions
\`a $\widetilde{M}^{[ k - 1 ]}$ co\"{\i}ncident.
\\

A l'aide
de ces remarques g\'en\'erales il est ais\'e de prouver le lemme suivant.

\begin{lemm}\label{jets-dif.arbres}
Consid\'erons deux germes d'arbres $\A_P$ et $\A'_P$ au dessus de $P$
et  $\A_P^{[k]}\,$, ${\A'_P}^{[k]}$ les germes d'arbres infinit\'esimaux
respectifs induits.
Soit $\Phi_{{}_{\bullet}}\,,\ \Psi_{{}_{\bullet}} : \A_P^{[k]}
\stackrel{\sim}{\rightarrow}
{\A'_P}^{[k]}$ deux isomorphismes au dessus de $P\,$. Pour $l < k\,$
notons $\Phi_{{}_{\bullet}}^{[l]}$
et
$\Psi_{{}_{\bullet}}^{[l]}\,$ leurs restrictions \`a $\A_P^{[l]}$ et
${\A'_P}^{[l]}\,$.
Si, pour un $r< k\,$ et un $j_0 < k-r\,$, l'\'egalit\'e $\Phi_{j_0} =
\Psi_{j_0}\,$ est v\'erifi\'ee, alors
on a~: $\Phi_{j_0+r}^{[k-r]} = \Psi_{j_0+r}^{[k-r]}\,\,$.
\end{lemm}

\begin{prop}\label{mont.morph.}
Soit $\Phi_0 : \M^0 \stackrel{\sim}{\rightarrow}  {\M'}^0\,$ un
diff\'eomorphisme au
dessus de $P$ entre les socles de deux germes d'arbres $\A_P$ et $\A'_P$
de m\^eme hauteur
$h\,$, avec  $\Phi_0(S^0) = {S'}^0\,$. Si l'isomorphisme induit
${\Phi_0}^{[h]} : {\M^0}^{[h]} \stackrel{\sim}{\rightarrow} { {\M'}^0}^{[h]}$
se rel\`eve en
un isomorphisme d'arbres infinit\'esimaux
$\Psi_{{}_{\bullet}} : \A_P^{[h]} \stackrel{\sim}{\rightarrow}
{\A'_P}^{[h]}\,$,
${\Psi_0}^{[h]} = {\Phi_0}^{[h]}\,$,
alors $\Phi_0$ se rel\`eve aussi en un unique isomorphisme de germes d'arbres
$\Phi_{{}_{\bullet}} :\A_P \stackrel{\sim}{\rightarrow} \A'_P$ et
${\Phi_{{}_{\bullet}}}^{[h]} = \Psi_{{}_{\bullet}}\,$.
\end{prop}

\begin{proof}
Raisonnons par r\'ecurrence sur la hauteur $h$ des arbres. On a vu que
$\Phi_0 $ se
rel\`eve en un diff\'eomorphisme $\Phi_1 :\M^1 \stackrel{\sim}{\rightarrow}
 {\M'}^1\,$, puisque $\Phi_0(S^0) = {S'}^0\,$. D'apr\`es
(\ref{jets-dif.arbres})  on a
${\Phi_1}^{[h-1]} = {\Psi_1}^{[h-1]}\,$. En particulier au niveau
ensembliste,
la restriction de $\Phi_1$ au diviseur $\D^1$ est \'egale \`a
l'application
$|\Psi_1| \,$ induite par $\Psi_1$ sur $\D^1\,$. Ainsi, puisque
$\Psi_{{}_{\bullet}}$ est un isomorphisme d'arbres infinit\'esimaux, on a
aussi~:
$\Phi_1(\S^1) =  {\S'}^1$ et $\Phi_1(S^1) = {S'}^1\,$. Ceci prouve la
proposition pour $h=1\,$.
Lorsque $h > 1$, on applique l'hypoth\`ese de r\'ecurrence \`a chaque germe de
$\Phi_1$ le long
d'une composante irr\'eductible de $S^1\,$.
\end{proof}

\subsection{Les faisceaux de diff\'eomorphismes transversalement formels}\label{gro.dif.t}
\addcontentsline{toc}{section}{\hspace{0,8em} {}\thesubsection .  Les faisceaux de diff\'eomorphismes t. f.}
D\'esi\-gnons encore
par $\A_P$ un germe d'arbre au dessus  de $P$ et conservons les
notations (\ref{diag.arbre}), (\ref{not.arbre}), (\ref{cime}),
(\ref{arbrfibr}), (\ref{id.arbres}) et
(\ref{ar.inf.def.}). \\

En chaque point $m$ du diviseur de cime $\wD_0 \subset \wD$ de l'arbre
$\A_P(0)$ consid\'erons les germes
d'automorphismes holomorphes de $\wM\,$, resp. $\wM^{[k]}\,$,
qui commutent avec $\widetilde{\pi}\,$
resp. $\widetilde{\pi}^{[k]}\,$,
valent l'identit\'e en restriction \`a la cime $\wM_0 \subset \wM$ de
$\A_P(0)\,$
et laissent invariant le germe de $\wD$ en $m$. On obtient ainsi des
faisceaux de groupes de base $\wD_0\,$. Nous les noterons~:
\begin{equation}\label{gr.dif.}
\G_{\A_P}\quad \hbox{ou encore}\quad  \G \,, \qquad
\hbox{resp.}\quad\G_{\A_P}^{[k]}\quad
\hbox{ou encore}\quad  \G^{[k]}\,.
\end{equation}
Les sections de ces faisceaux respectent d'apr\`es (\ref{decomp.id.})
le faisceau d'id\'eaux
$\wJ := \J_{\,(h)}\,$ d\'efini en (\ref{id.arbres}). On dispose ainsi
pour $0\leq l \leq k $, de morphismes de restrictions
$$\r_l : \G \fle \G^{[l]}\, , \qquad \r_l^k : \G^{[k]} \fle \G^{[l]}\,,
\qquad
0 \leq lÊ\leq k\,.$$
Les faisceaux de base $\wD_0$ des {\it germes d'automorphismes l-tangents
\`a l'identit\'e le long de
$\wD$} d\'e\-fi\-nis par
\begin{equation}\label{gr.dif.t}
\G_l := ker(\r_l)\qquad \hbox{et} \qquad \G_l^{[k]} := ker(\r_l^k)
\end{equation}
forment des faisceaux de sous-groupes distingu\'es de $\G\,$ et de
$\G^{[k]}$ . En parti\-culier $\G_0$ est form\'e des
\'el\'ements de $\G$ qui valent l'identit\'e en restriction \`a $\wD\,$.
\\

On appelle {\it faisceau des germes d'automorphismes transversalement
formels le long de $\wD$},
respectivement {\it transversalement formels et l-tangents \`a l'identit\'e le
long de $\wD$} les limites
projectives
\begin{equation}\label{defautf}
\hG := \limproj_k{\G^{[k]}}\, ,\qquad
\hG_l := \limproj_{k \geq l}{\G_l^{[k]}}\,.
\end{equation}
On dispose encore de projections canoniques $\hr_l : \hG_l \fle
\hG_l^{[k]}\,$.

\begin{rema}\label{rem1} On aurait aussi pu filtrer $\wO$ par les
puissances du faisceau d'id\'eaux
$I_{\wD}$ d\'efini en (\ref{arbr.inf}) et poser
$$\wM^{[[k]]} := \left( \wD_0 \,; \:
\O_{\wM}/{\widetilde{I}}^{k+1}\right)\, ,
\qquad
\hhG := \limproj_{k}{\G^{[[k]]}}\,,$$
o\`u $\G^{[[k]]}$ est le faisceau des germes, aux points de $\wD_0\,$,
d'automorphismes de
$\wM^{[[k]]}$ qui valent l'identit\'e sur $\wM_0\,$ et commutent aux
projections sur $P\,$. Comme
d'apr\`es (\ref{decomp.id.})
ces deux filtrations
sont \'equivalentes, les faisceaux $\hhG\,$ et $\hG\,$ sont canoniquement
isomorphes. De m\^eme les sous-faisceaux $\hhG_l \subset \hhG\,$ d\'efinis de
mani\`ere similaire
sont canoniquement isomorphes \`a $\hG_l\,$.
\end{rema}

\begin{rema}\label{rem2}
En un point r\'egulier $m$ de $\wD_0\,$, dans un syst\`eme de
coordonn\'ees locales  $(z_1, z_2 ;\, t_1, \ldots , t_p)$ dans lesquelles
$\wD = \{z_1 = 0\,\}\,$, $\wp =
t := (t_1,\ldots ,t_p)\,$ et  $\wJ_{, m} = ( z_1^r )\,$, les \'el\'ements de
la fibre $\hG_{l,\, m}$ de
$\hG_l$ s'identifient aux couples $\left( F_1, F_2 \right)$ de s\'eries
formelles en
$z_1$
$$F_j = z_j + z_1^{r\,(l+1)}\sum_{k=0}^{\infty} \,A_{j,\,k}(z_2;
t)\,z_1^k\,,\quad j = 1, 2\,,
$$
dont les coefficients $A_{j,\,k}\,$ sont des fonctions
holomorphes sur un voisinage commun de $m$
dans $\wD\,$. Lorsque $m$ est un point singulier de $\wD_0$ avec $\wD = \{z_1z_2
= 0\,\}\,$ et
$\wJ_{,m} = (z_1^rz_2^s )\,$, alors les \'el\'ements de $\hG_{l,\,m}$ s'identifient
aux couples $\left( F_1, F_2 \right)$ avec
$$F_j = z_j + (z_1^rz_2^s)^{l+1}
\sum_{k=0}^{\infty}\,\left(A_{j,\,k}^1(z_1; t) + A_{j,\,k}^2(z_2;
t)\right)\!  (z_1z_2)^k \,,\quad j= 1, 2\,, $$
les coefficients $A_{j,\,k}^r$
de ces s\'eries \'etant encore holomorphes sur un m\^eme voisinage de $m$ dans
$\wD\,$.  Cette derni\`ere
\'ecriture est bien s\^ ur unique si l'on exige que $A_{1,\,k}^2(0; t)
\equiv  A_{2,\,k}^2(0; t)  \equiv 0\,$.
On d\'eduit imm\'ediatement de ces \'ecritures que, si l'on se donne pour
chaque
$l
\in
\N$  un \'el\'ement
$\Phi^l$ de $\hG_{l, m}$, les suites
$$\Phi_n := \Phi^n\circ\cdots\circ\Phi^0\quad \hbox{et}\quad
\Psi_n := \Phi^0\circ\cdots\circ\Phi^n$$
convergent dans $\hG_{, m}\,$, pour la topologie $I_{\wD}$-adique. D'autre
part tout \'el\'ement
de $\hG_{, m}$ est limite de telles suites et, $\Phi^0 ,\ldots,\Phi^s$
\'etant construits,
$(\Phi^{s+1})$ est unique modulo $I_{\wD}^{s+2}\,$.
\end{rema}

%%%%%%%%%%%%%%%%%%%%%%%

\subsection{Th\'eor\`emes de stabilit\'e et de d\'etermination finie}\label{th.st.det}
\addcontentsline{toc}{section}{\hspace{0,8em} {}\thesubsection .  Th\'eor\`emes de stabilit\'e et de d\'etermination
finie}
Nous conservons les notations du paragraphe pr\'ec\'edent.

\begin{defi}\label{doncrit}
Nous appelons {\it donn\'ees critiques d'un germe d'arbre $\A_P$ au dessus
de
$P$}  l'ensemble $\bC\,(\A_P)$ dont
les \'el\'ements, appel\'es {\it ensembles critiques} sont :
\begin{enumerate}
\item les points de l'ensemble singulier de cime $\wS_0$ de l'arbre
$\A_P(0)\,$,
\item les composantes connexes du compl\'ementaire de $\wS_0$ dans le
diviseur excep\-tion\-nel de cime
$\wD_0$ de $\A_P(0)\,$;
\end{enumerate}
\end{defi}
On appelle {\it donn\'ees critiques de
dimension 0, resp. de dimension 1}, et on note $\bC_0\,(\A_P)$, resp.
$\bC_1\,(\A_P)$  les ensembles des donn\'ees
critiques du type 1. ou 2..  Deux ensembles critiques seront
dits {\it adjacents} si leurs  adh\'erences s'intersectent; on note~:
\begin{equation}\label{adjac}
\begin{array}{rl}
\phantom{aaa} & Ad(K) := \left\{ \,K' \in \bC\,(\A_P)
\/\,;\/\,\overline{K} \cap
\overline{K'} \not=
\emptyset \, \right\}\\
\phantom{aaa} & Ad_i(K) := \left\{\, K' \in  Ad(K) \/\,;\/\, dim(K') = i\,
\right\}\, ,\quad i = 0 \; \hbox{ou} \;1\,.
\end{array}
\end{equation}

\begin{defi}\label{rec.dist.}
Nous appelons {\it recouvrement distingu\'e de $\wD_0\,$} tout
recouvrement ouvert
$\U := \left(U_{\a}\right)_{\a \in \bC\,(\A_P)}\,$ du type suivant~:
\begin{enumerate}
\item si $dim(\a) = 1\,$, on pose $U_{\a} := \a\,$,
\item si $dim(\a) = 0\,$, $U_{\a}$ est la trace sur $\wD_0$ d'un petit
polydisque
$\{ |z_1| < \e \; , \; |z_2| < \e \}\,$, o\`u $z_1\,$, $z_2$ sont des
coordonn\'ees de $\wM_0$ dans
lesquelles $\wD_0$ est monomial et $z_1(\a) = z_2(\a) = 0\,$.
\end{enumerate}
\end{defi}
Les propri\'et\'es remarquables d'un tel recouvrement, que nous utiliserons
constamment, sont~: a) les
ouverts de $\U$ sont de Stein et admettent un syst\`eme fondamental de
voisinages de Stein dans
$\wM\,$ et, b) toute intersection 3 \`a 3 d'ouverts de $\U$ est triviale.

\begin{rema}\label{rem3} En calculant \`a partir des expressions de la
remarque (\ref{rem2}), on voit que
si $\U$ est un recouvrement distingu\'e de $\wD_0\,$,
pour $\a \in \bC(\A_P)$ et $0 \leq l \leq k\,$ les applications canoniques
$$\G\,(U_{\a}) \fle \G^{\,[k]}\,(U_{\a})\, ,\qquad
\hG\,(U_{\a}) \fle \hG^{\,[k]}\,(U_{\a})\ \,,$$
$$\G_l\,(U_{\a}) \fle \G_l^{\,[k]}\,(U_{\a})\, ,\qquad
\hG_l\,(U_{\a}) \fle \hG_l^{\,[k]}\,(U_{\a})\ \,,$$
sont surjectives.
\end{rema}

Rappelons qu'\'etant donn\'e un faisceau de groupes $\bG$ et un recouvrement
ouvert
localement fini $\cal W := (W_i)_{{}_{i\in I}}\,$ de la base de $\bG$,
on appelle {\it 1-cocycle \`a valeurs dans $\bG\,$} et on note
$\left(g_{ij}\right)
\in Z^1\left(\cal W; \bG\right)$ toute famille
$$\left(g_{ij}\right)_{(i, j) \in I\check{\times}I}\,,\quad g_{ij} \in \bG
(W_i\cap W_j)\,,\quad
I\check{\times}I := \{(i, j) \in I^2\,;\
 W_i \cap W_j \not= \emptyset \}\,, $$
telle que
$g_{ij} = g_{ji}^{-1}$ et
$g_{ij} = g_{ik}g_{kj}$ chaque fois que $W_i \cap W_j \cap W_k \not=
\emptyset \,$.
Deux 1-cocycle ${\cal C} := (g_{ij})\, $ et ${\cal C}' := (g'_{ij})$ de
$Z^1\left(\cal W; \bG\right)$ sont
dits {\it cohomologues} s'il existe
un \'el\'ement $\cal C^0 := \left(g_{i}\right)$ de $Z^0\left(\cal W; \bG\right) :=
\prod_{i \in
I}\bG\,(W_i)\,)\,$
tel que
\begin{equation}
g'_{ij} = g_i \cdot g_{ij} \cdot g_j^{-1}\,.
\end{equation}
Nous noterons alors ${\cal C} \approx_{\bG}{\cal C}'\,$ ou plus pr\'ecis\'ement $\cal C' = \cal C^0 \star
\cal C$ et nous  d\'esignerons par $[{\cal C}\,]_{{}_{\bG}}\,$
ou par $[g_{ij}]_{{}_{\bG}}\,$ la classe d'\'equivalence de ${\cal C}$;
lorsqu'aucune confusion n'est possible nous
omettrons l'indice $\bG\,$ dans ces notations.
\\

Soit $\U$ un recouvrement distingu\'e de $\wD_0\,$.
Nous pouvons
associer
\`a un cocycle ${\cal C} := \left(\f_{\a\b}\right) \in Z^1\left(\U;
\G_l\right)\,$, $l \geq  0$ le
germe, le long d'une courbe $\D_{0, \,{\cal C}} \simeq \wD_0\,$, de la vari\'et\'e
holomorphe $\M_{{\cal C}}\,$ obtenue
par recollement en quotientant l'union disjointe $\sqcup_\a (V_\a,\, U_\a)$ des germes, le long des
$U_\a$, de voisinages $V_\a$ de $U_\a$
dans la cime $\wM$ de $\A_P\,$, par la relation d'\'equivalence
\begin{equation}\label{recvar}
m \sim_{{\cal C}} m' \ssi \hbox{il existe}\ (\a, \b) \in
{\bC}\,(\A_P)\check{\times}{\bC}\,(\A_P)\quad \hbox{tel que}
\end{equation}
$$m \in V_{\a}\,,\quad m' \in V_{\b} \quad \hbox{et}\quad m'
=\f_{\b\a}(m)\,.$$
\noindent Les propri\'et\'es des $\f_{\a\b}$ font que cette vari\'et\'e est
naturellement munie
d'un diviseur $\D_{{\cal C}} \supset \D_{0,\, \cal C}\,$,
d'un germe le long de
$\D_{0, {\cal C}}$ de submersion holomorphe $\p_{{\cal C}}$ sur $P\,$ et
d'un germe le long de $\wD_0$ de
plongement
$$
\r_{{\cal C}} : \left( \wM_0 \, ,\, \wD_0 \right) \hookrightarrow  \left(
\M_{{\cal C}}\,, \,
\D_{0,\, \cal C} \right)
$$
qui induit un isomorphisme entre les germes
$\left( \wM_0 , \, \wD_0 \right)$ et $\left(  \pi_{\cal C}^{-1}(0) , \,
\D_{0,\, \cal C} \right)\,$.
Lorsque
${\cal C} \in Z^1\left(\U; \G_l\right)\,$, $l \geq  0$ on dispose aussi
d'un plongement
\begin{equation}\label{iddiv}
\ui_{{\cal C}}^{[l]} : \wM^{[l]} \hookrightarrow \M_{{\cal C}}\,,\quad
\hbox{avec}\quad
\p_{{\cal C}}\circ \ui_{{\cal C}}^{[l]} = \wp^{\,[l]}\,,
\end{equation}
qui induit un isomorphisme entre les diviseurs $\wD \simeq \wM^{[0]}$ et
$\wD_{{\cal C}}\,$. En restriction \`a $\wD_{0 ,\,{\cal C}}\,$ cet
isomorphisme co\"{\i}ncide avec $\r_{{\cal C}}\,$.

\begin{rema}\label{coho.dif.} Deux cocycles ${\cal C} :=
\left(\f_{\a\b}\right)\,$ et
${\cal C}' := \left(\f'_{\a\b}\right)\,\in Z^1\left(\U; \G\right)\,$ sont
cohomologues, $\cal C' = \cal C^0 \star\cal C$ avec  ${\cal
C}^{\, 0} := \left(\f_{\a}\right)\, \in
Z^0\left(\U; \G\right)\,$, si et seulement s'il existe un germe de
diff\'eomorphisme
$\Phi_{{\cal C}^{\, 0}} : \M_{{\cal C}} \stackrel{\sim}{\fle} \M_{{\cal
C}'}\,$
qui satisfait les relations de
commutation~:
$$\p_{{\cal C}'} \circ \Phi_{{\cal C}^{\, 0}}Ê\, = \p_{{\cal C}}\,,\qquad
\r_{{\cal C'}} = \Phi_{{\cal C}^0} \circ \r_{{\cal C}}\,.$$
Le diff\'eomorphisme $\Phi_{{\cal C}^{\, 0}}\,$ se construit par recollement
des $\f_{\a}\,$~: la relation
de cohomologie exprime la compatibilit\'e avec les "changements de cartes"
$\f_{\a\b}\,$. Lorsque
${\cal C}\,$,
${\cal C}'\in Z^1\left(\U; \G_l\right)\,$, pour exprimer la relation de
cohomologie dans le faisceau $\G_l\,$, il faut
rajouter la condition suppl\'ementaire~: $\Phi_{{\cal C}^{\, 0}} \circ
\ui_{{\cal C}}^{[l]} = \ui_{{\cal C}'}^{[l]}\,$.
\end{rema}

Donnons nous un 1-cocycle ${\cal C} \in Z^1\left(\U; \G\right)\,$. Pour $k
\in \N$,
d\'efinissons le germe, le long de $\D_{0,\, \cal C}\,$, du voisinage
infinit\'esimal $\M_{{\cal
C}}^{\,[k]}\,$ de
$\D_{{\cal C}}\,$ dans
$\M_{{\cal C}}$ en posant~:
$$\left| \M_{\cal C}^{[k]}\right| := \D_{0,\, \cal C} \, , \qquad
\O_{\M_{{\cal C}}^{\,[k]}} :=
i_{0,\, \cal C}^{-1} \left(
\O_{\M_{{\cal C}}}\left/
\prod_{{}_{D\in comp(\D_{{\cal C}})}} \left(I_{{\cal C},\,
D}^{m\,(D)}\right)^{k +
1}\right.\right)$$
o\`u $i_{0,\, \cal C} : \D_{0, {\cal C}} \hookrightarrow \M_{{\cal C}}\,$
est l'application
d'inclusion,
$I_{{\cal C}, D}$ est le faisceau sur
$\D_{{\cal C}}\,$ for\-m\'e des germes de fonctions nulles sur  une
composante irr\'eductible $D\,$
de $\D_{{\cal C}}$ et les $m\,(D)\,$  sont d\'efinis par la d\'ecomposition
(\ref{decomp.id.}) de
$\wJ\,$.
Lorsque ${\cal C} \in Z^1\left(\U; \G_l\right)\,$, $l < k\,$, le morphisme
$\ui_{{\cal C}}^{[l]}\,$ se factorise clairement en un plongement
$\underline{\ui}_{{\cal C}}^{\,[l]} : \wM^{\,[l]} \hookrightarrow
\M_{{\cal C}}^{\,[k]}\,$, (i.e. le
co-morphisme  associ\'e \`a $\underline{\ui}_{{\cal C}}^{\,[l]}\,$ est
surjectif). Un calcul direct en
coordonn\'ees permet de voir que $\underline{\ui}_{{\cal C}}^{\,[l]}\,$ se
factorise en un
isomorphisme
$\underline{\ui}_{{\cal C}}^{\,[l]} : \wM^{\,[l]} \iso \M_{{\cal
C}}^{\,[l]} \hookrightarrow
\M_{{\cal C}}^{\,[k]}\,$. \\

En reprenant, pour un germe d'arbre
quelconque, la d\'emonstration du lemme (1.3.1) de \cite{Minv} qui ne concerne que
les arbres produits, on
obtient la premi\`ere affirmation du lemme suivant. La seconde affirmation
r\'esulte
directement du fait que les exposants $m\,(D)\,$ ne d\'ependent que de
l'arbre dual.

\begin{lemm}\label{stab.cimes}
Le germe de  $\M_{{\cal C}}\,$ le long de $\D_{0, {\cal C}}\,$ est
biholomorphe au germe
$(\wM'\,,\wD_0'\,)$ de la cime d'un germe d'arbre $\A'_P$ au dessus de
$P\,$, qui poss\`ede le
m\^eme arbre dual que $\A_P\,$,
$$\begin{array}{ccc}
(\M_{{\cal C}}, \D_{0, {\cal C}})\!\!\!\!\!
&\stackrel{\s}{\stackrel{\sim}{\fle}}&\!\!\!\!\!(\wM', \wD_0') \\
& {}_{\p_{{\cal C}}}\!\!\!\searrow \phantom{\, P \,} \swarrow\!\!\!
{}_{\wp'}& \\
& P &
\end{array}$$
De plus le biholomorphisme $\s$ induit pour tout $k\,$ un isomorphisme
$\s^{\,[k]}\,$ entre
$\M_{{\cal C}}^{\,[k]}\,$ et
${\wM'}{}^{\,[k]}\,$.
\end{lemm}

Fixons un cocycle ${\cal C} := \left(\xi_{\a\b}\right)\, \in Z^1\left(\U;
\G\right) $. Comme en (\ref{gr.dif.t}) notons $\G_{{\cal C}}$ le faisceau de base
$\D_{0,{\cal
C}}\,$ des germes de biholomorphismes de $\M_{\cal C}\,$ qui commutent avec
$\pi_{\cal C}$, valent l'identit\'e en restriction \`a
$ \wM_{0,\, \cal C}$ et par $\G_{{\cal C},\, k} \subset \G_{{\cal C}}$ le
sous-faisceau des germes
qui valent aussi l'identit\'e en restriction \`a
$\M_{{\cal C}}^{\,[k]}\,$. Consid\'erons l'isomorphisme
$\wD_0 \simeq \D_{0,\, \cal C}$ comme une identification et sur $\wD_0$
comparons les cohomologies
du recouvrement $\U$ \`a valeurs dans les faisceaux $\G_{{\cal C},\, k}\,$
et
$\G_k\,$.
Pour chaque ouvert $U_{\a}$ de $\U$ nous disposons par construction d'un
germe
de biholomorphisme $u_{\a} : (\wM, U_{\a}) \iso (\M_{{\cal C}}, U_\a)\,$;
ce qui donne
les identifications~:
$$ \kappa^0 : Z^0\left(\U; \G\right) \iso Z^0\left(\U; \G_{{\cal
C}}\right)\, , \qquad
 \kappa^0_k : Z^0\left(\U; \G_k\right) \iso Z^0\left(\U; \G_{{\cal C},\,
k}\right) \, ,$$
$$\left(\Psi _{\a}\right) \mapsto  \left(\Psi_{{\cal C}, \a}\right) :=
\left({u_{\a}}_{\ast}(\Psi_
{\a})\right) := \left(u_{\a}\circ\Psi _{\a}\circ u_{\a}^{-1}\right)$$
Pour tout
$(\a, \b) \in {\bC}\,(\A_P)\check{\times}{\bC}\,(\A_P)\,$ nous disposons
de deux germes de
dif\-f\'eo\-mor\-phis\-mes  $u_{\a}\circ \iota_{\a\b} \,$ et $u_{\b}\circ
\iota_{\a\b}\,$, o\`u
$\iota_{\a\b}\,$ d\'esigne le germe d'inclusion  $(\wM, U_{\a}\cap U_{\b})
\ifle (\wM, \wD_0)\,$.  Le
recouvrement $\U\,$ n'a pas d'intersection d'ouverts trois \`a trois
non-triviale. Ainsi, pour tout
faisceau de groupes $\bG\,$ sur $\wD_0\,$, le groupe $Z^1(\U ;\, \bG)\,$
s'identifie au groupe
$$
\overline{Z}^1(\U ;\, \bG) :=\prod_{( \a,Ê\b) \in \cal A} \bG(U_{\a}\cap U_{\b})\,,
$$
$$
\mbox{avec} \quad \cal A := \left\{\left. (\a, \b) \in
{\bC}\,(\A_P)\check{\times}{\bC}\,(\A_P) \;  \right/  \; dim U_{\a} < dim U_{\b}\:
\right\}\,.
$$

Choisissons d'identifier les ensembles de 1-cocycles de la mani\`ere
suivante~
$$\kappa^1 : \overline{Z}^1\left(\U; \G\right) \iso \overline{Z}^1\left(\U;
\G_{{\cal C}}\right) \, , \qquad
\kappa_k^1 : \overline{Z}^1\left(\U; \G_k\right) \iso
\overline{Z}^1\left(\U; \G_{{\cal C},\,
k}\right)\, ,$$
$$\left(\Psi _{\a\b}\right) \mapsto  \left(\Psi_{{\cal C},\, \a\b}\right)
:=
\left(({u_{\a}\circ \iota_{\a\b}})_{\ast}(\Psi _{\a\b})\right)\,.$$
On voit facilement que les relations de cohomologie s'expriment dans les
espaces
$\overline{Z}^1$ par~:
\begin{equation}\label{tr.co.tor.}
\left(\Psi_{{\cal C},\, \a\b}\right) \approx_{\scriptstyle \G_{\cal C}}
\left(\Psi_{{\cal C},
\,\a\b}'\right)
\ssi \left(\Psi_{\a\b}\circ\xi_{\a\b}\right) \approx_{\scriptstyle \G}
\left(\Psi_{\a\b}'\circ\xi_{\a\b}\right)
\end{equation}
$$
\left(\Psi_{{\cal C},\, \a\b}\right) \approx_{\scriptstyle \G_{\cal C, k}}
\left(\Psi_{{\cal C},
\,\a\b}'\right)
\ssi \left(\Psi_{\a\b}\circ\xi_{\a\b}\right) \approx_{\scriptstyle \G_k}
\left(\Psi_{\a\b}'\circ\xi_{\a\b}\right)
$$
\\

\begin{enonce}{Th\'eor\`eme de Stabilit\'e}\label{d.q.fin.}
Soient $\A_P\,$ un germe d'arbre au dessus de $P$, de hauteur $h\,$, $\U$
un recouvrement
distingu\'e du diviseur de cime $\wD_0$ de $\A_P(0)$ et soit $\G\,$,
res\-pec\-tivement $\G_k\,$, $k \in
\N\,$, les faisceaux introduits en (\ref{gro.dif.t}).
Alors
pour tout $k \in \N\,$ les applications suivantes
$$\phantom{\hbox{resp.}} H^1\left(\U; \G_{h}\right) \fle H^1\left(\U;
\G\right)\ \,,\qquad
\left[\f_{\a\b}\right]_{\G_{h}}
\mapsto \left[\f_{\a\b}\right]_{\G}\, , $$
$$\hbox{resp.} \quad H^1\left(\U; \G_{k+h}\right) \fle H^1\left(\U;
\G_k\right)\ \,,\qquad
\left[\f_{\a\b}\right]_{\G_{k+h}}
\mapsto \left[\f_{\a\b}\right]_{\G_k}$$
sont les applications constantes $ \left[\f_{\a\b}\right]_{\G_{h}} \mapsto
\left[Id_{\a\b}\right]_{\G}\,$,
resp. $ \left[\f_{\a\b}\right]_{\G_{k+h}} \mapsto
\left[Id_{\a\b}\right]_{\G_k}\,$, o\`u
$Id_{\a\b}$ d\'esigne le germe de l'application identit\'e le long de
$U_{\a}\cap U_{\b}\,$.
\end{enonce}

\begin{proof} Donnons nous un 1-cocycle ${\cal C} :=
\left( \f_{\a\b} \right) \in Z^1(\U; \G_{k+h})\,$ D'apr\`es la
proposition (\ref{stab.cimes}), $\M_{{\cal C}}$ est diff\'eomorphe \`a la
cime d'un germe d'arbre $\A'_P$ au dessus de $P\,$, $\s : \M_{{\cal C}}
\iso \wM'\,$.
Le diff\'eomorphisme $\wt := \s^{[k+h]} \circ \ui_{{\cal C}}^{[k+h]} :
\wM^{[k+h]} \iso {\wM'}{}^{[k+h]}\,$
induit gr\^ace
\`a (\ref{desc.morph.}) un isomorphisme d'arbres infinit\'esimaux
$$\t_{{}_{\bullet}} : \A_P^{[k+h]} \iso {\A'}_P^{[k+h]}\ \,,\qquad
\hbox{avec}\quad \t_{h} = \wt\,.$$
Soit $T : \M^0 \iso  {\M'}^0$ un diff\'eomorphisme entre les socles de $\A_P$
et de $\A'_P$
induisant $\t_0\,$, i.e. $T^{[k+h]} = \t_0\,$. Relevons le
(\ref{mont.morph.}) en un (unique)
isomorphisme $T_{{}_{\bullet}}: \A_P \iso \A'_P\,$,  $T_0 =T\,$. En
comparant
$T_{{}_{\bullet}}$ \`a $\t_{{}_{\bullet}}$ \`a l'aide de
(\ref{jets-dif.arbres}) on voit que pour
tout $j \leq h\,$, on a $T_j^{[k+h-j]} = \t_j^{[k+h-j]}\,$. En particulier
$T_{{}_{\bullet}}^{[k]} = \t_{{}_{\bullet}}^{[k]}\,$ et $H :=
T_h^{-1}\circ \s :  \M_{{\cal C}} \iso
\wM$ est un diff\'eomorphisme qui vaut l'identit\'e en restriction \`a
$\wM^{[k]}\,$, c'est \`a dire
$(H \circ \ui_{{\cal C}})^{[k]} = \hbox{identit\'e}_{\wM^{[k]}}\,$. On
conclut alors par
(\ref{coho.dif.}) que l'on a bien $[\f_{\a\b}]_{\G_k} =
[Id_{\a\b}]_{\G_k}\,$.\\

\indent La m\^eme d\'emonstration s'applique aux faisceaux $\G_h$ et $\G\,$.
\end{proof}

\begin{coro}\label{d.q.fin.tord.}
Avec les m\^emes notations que dans le th\'eor\`eme ci-dessus, fixons $k \in
\N\,$ et un 1-cocycle ${\cal
C} := \left(\xi_{\a\b}\right) \in Z^1\left(\U; \G\right)\,$. Alors pour
tout
${\cal C}' := \left(\Psi_{\a\b}\right) \in Z^1\left(\U; \G_{h}\right)\,$,
resp.
${\cal C}' := \left(\Psi_{\a\b}\right) \in Z^1\left(\U;
\G_{k+h}\right)\,$, il existe
$\left(\Psi_{\a}\right) \in Z^0\left(\U;
\G \right)\,$, resp.
$\left(\Psi_{\a}\right) \in Z^0\left(\U;
\G_k\right)\,$ tel que~:   $\Psi_{\a\b} \circ \xi_{\a\b} = \Psi_{\a} \circ
\xi_{\a\b} \circ
\Psi_{\b}^{-1}\,$.
\end{coro}
\begin{proof} Conservons les notations de la d\'emonstration du th\'eor\`eme.
Puis\-que le dif\-f\'eo\-mor\-phis\-me $\s$ donn\'e par (\ref{stab.cimes})
induit des diff\'eomorphismes
$\s^{[k+h]} : \M_{{\cal C}}^{[k+h]} \is {\wM'}{}^{[k+h]}\,$, il induit aussi
des isomorphismes entre les
faisceaux de groupes $\G_{{\cal C}}$ et  $\G_{\A'_P}\,$, resp.
$\G_{{\cal C},\, k}$ et $\G_{\A'_P, k}\,$, d\'efinis en (\ref{gr.dif.}).
Appliquons maintenant le th\'eor\`eme pr\'ec\'edent \`a $\A'_P\,$. Les conclusions
de ce th\'eor\`eme se
"traduisent" sur le germe d'arbre $\A_P\,$ \`a l'aide de
(\ref{tr.co.tor.}). On obtient imm\'ediatement
le r\'esultat.
\end{proof}

\begin{enonce}{Th\'eor\`eme de D\'etermination Finie}\label{mth.co.}
Soient $\A_P\,$ un germe d'ar\-bre au dessus de $P$ de
hauteur
$h\,$,
$\U$ un recouvrement distingu\'e du diviseur de c\^{\i}me $\wD_0$ de $\A_P(0)\,$,
$k \in\N\,$ et
$\G\,$, $\hG\,$, resp. $\G_k\,$, $\hG_k\,$, les faisceaux (\ref{gr.dif.t})
des germes de
diff\'eomorphismes holomorphes, transversalement formels, resp.
holomorphes, transversalement formels
$k$-tangents \`a l'identit\'e, d\'efinis au chapitre (\ref{gro.dif.t}). Alors
pour tout $k \in \N\,$
les applications canoniques suivantes~:
\begin{enumerate}
\item $\chi : H^1(\U; \G ) \fle H^1 (\U;
\G^{[h]} )\,,\quad
 [\f_{\a\b} ]_{\G} \mapsto
 [\f_{\a\b}^{[h]} ]_{\G^{[h]}}\,$ ,
\item $\chi_k : H^1 (\U; \G_k ) \fle H^1 (\U;
\G_k^{[k+h]} )\,,\quad
 [\f_{\a\b} ]_{\G_k} \mapsto
 [\f_{\a\b}^{[k+h]} ]_{\G_k^{[k+h]}}\,$ ,
\item $\widehat{\chi} : H^1 (\U; \hG ) \fle H^1 (\U;
\G^{[h]} )\,,
\quad  [\f_{\a\b} ]_{\hG}
\mapsto  [\f_{\a\b}^{[h]} ]_{\G_k^{[h]}}\,$ ,
\item $\widehat{\chi}_k : H^1(\U; \hG_k) \fle H^1(\U;
\G_k^{[k+h]})\,,
\quad  [\f_{\a\b} ]_{\hG_k}
\mapsto  [\f_{\a\b}^{[k+h]} ]_{\G_k^{[k+h]}}\,$ ,
\item $\underline{\chi} : H^1 (\U; \G ) \fle H^1 (\U;
\hG ) \, \is
H^1 (\U;
\hhG )\,, \quad   [\f_{\a\b} ]_{\G} \mapsto
 [\f_{\a\b} ]_{\hG}\,$ ,
\item $\underline{\chi}_k : H^1 (\U; \G_k ) \fle H^1 (\U;
\hG_k ) \, \is
H^1 (\U;
\hhG_k) \,, \quad   [\f_{\a\b} ]_{\G_k} \mapsto
 [\f_{\a\b} ]_{\hG_k}\,$ ,
\end{enumerate}
sont bijectives.
\end{enonce}

\begin{proof}
Les assertions 5. et 6. d\'ecoulent directement des assertions
pr\'ec\'edentes. La surjectivit\'e de
$\chi\,$, $\chi_k\,$, et de $\widehat{\chi}\,$, $\widehat{\chi}_k$ r\'esulte
de (\ref{rem3}).

Montrons l'injectivit\'e de $\chi_k\,$.
Donnons nous des \'el\'ements
$\left(\fab \right)$ et $\left(\fabp \right)$ de $Z^1(\U; \G_k)$ qui
induisent des \'el\'ements $\G_k^{[k+h]}$-cohomologues dans $Z^1(\U;
\G_k^{[k+h]})\,$
$$\fab^{[k+h]} =
\varphi_{\a}\circ\fabp^{[k+h]}\circ\varphi_{\b}^{-1}\,\,,\qquad
\left(\varphi_{\a}\right) \in Z^0(\U; \G_k^{[k+h]})\,.$$
Le 0-cocycle $(\varphi_{\a})$ provient d'un cocycle $(\fa) \in Z^0(\U;
\G_k)\,$,
$$\fab^{[k+h]}  = \left(\fa\circ\fabp\circ\fb^{-1}\right)^{[k+h]}\,,\qquad
\fa^{[k+h]} =
\varphi_{\a}\,.$$ Pour obtenir le r\'esultat il suffit d'appliquer le
corollaire
(\ref{d.q.fin.tord.}) avec
$${\cal C} := \left(\xi_{\a\b}\right) :=
\left(\fa\circ\fabp\circ\fb^{-1}\right)\,, \qquad
 {\cal C}' := \left(\fab\circ\xi_{\a\b}^{-1}\right)\,.$$

Prouvons maintenant l'injectivit\'e de $\widehat{\chi}_k\,$. Soient
$(\fab)\,$ et
$(\fabp ) \in Z^1(\U; \hG_k)$ qui induisent la m\^eme classe dans $H^1(\U;
\G_k^{[k+h]})\,$.
Comme dans B) nous avons un \'el\'ement $(\fa)$ de $Z^0(\U; \hG_k)$ tel que~:
$$\left( \fab^{[k + h]}\right) := \left(\fa\circ\fabp
\circ\fb^{-1}\right)^{[k + h]}\,.$$
Donnons nous $(\xab) \in Z^1(\U; \G_k)$ tel que
$\xab^{[k+h]} = \fab^{[k+h]}\,$.
Nous allons mon\-trer successivement les \'egalit\'es, dans $H^1(\U; \hG_k)$,
des classes
$$\left[\fab\right]_{\hG_k} = \left[\xab\right]_{\hG_k}\,, \qquad
\left[\fa\circ\fabp \circ\fb^{-1}\right]_{\hG_k} =
\left[\xab\right]_{\hG_k}\,;
$$
Ce qui permet de
conclure.\\

D\'ecomposons $\fab$ de la mani\`ere suivante~:
$$\fab = R_{\a\b}^{(2)}\circ\fabu\circ\xab\,, \qquad (\fabu) \in Z^1(\U;
\G_{k+h})
\,,\qquad (R_{\a\b}^{(2)}) \in Z^1(\U; \hG_{k+2h})\,.$$
En utilisant le corollaire (\ref{d.q.fin.tord.}) on obtient~:
$(\fao) \in Z^0(\U; \G_k)$ v\'erifiant
$$\fabu\circ \xab =
\fao\circ\xab\circ(\fbo)^{-1}\,.$$
Ainsi on a~:
$$\fab = \fao\circ\widetilde{R}_{\a\b}^{(2)}\circ\xab\circ
(\fbo)^{-1}\,,$$
avec
$$\left(\widetilde{R}_{\a\b}^{(2)}\right) := \left((\fao)^{-1}\circ
R_{\a\b}^{(2)}\circ
\fao\right) \in Z^1(\U; \hG_{k+2h})\,,$$
puisque $\hG_{k+2h}$ est un sous-groupe distingu\'e de $\hG_k\,$.
D\'ecomposons alors
 $\widetilde{R}_{\a\b}^{(2)}$ de la mani\`ere suivante~:
$$\widetilde{R}_{\a\b}^{(2)} =
R_{\a\b}^{(3)}\circ\fabd\,,\qquad
(\fabd) \in  Z^1(\U;\G_{k+2h})\,,  \qquad(R_{\a\b}^{(3)}) \in
Z^1(\U;
\hG_{k+3h})\,.$$
On obtient de nouveau \`a l'aide du corollaire (\ref{d.q.fin.tord.}) une
d\'ecomposition~:
$$\fabd\circ\xab = \fau\circ\xab\circ \left(\fbu\right)^{-1}\,,\qquad
\left(\fau\right) \in Z^0(\U; \G_{k+h})\,;$$
et
$$\widetilde{R}_{\a\b}^{(2)}\circ\xab = \fau\circ
\widetilde{R}_{\a\b}^{(3)}\circ \xab \circ
(\fbu)^{-1}\,,$$
$$(\widetilde{R}_{\a\b}^{(3)}) := \left(
(\fau)^{-1}\circ R_{\a\b}^{(3)}\circ\fau
\right)\in Z^1(\U; \hG_{k+3h})\,.$$
Ce qui donne~:
$$\fab = \left(\fao\circ\fau\right) \circ \widetilde{R}_{\a\b}^{(3)}\circ
\xab \circ
\left(\fbo\circ\fbu\right)^{-1}\,.$$
De cette mani\`ere, en it\'erant cette construction, on obtient des
suites  $$\left(\fa^{(j)}\right) \in Z^0(\U; \G_{k+jh})\,, \qquad
\left(R_{\a\b}^{(j)}\right) \in
Z^1(\U; \hG_{k+jh})\,,$$
qui, pour tout $n \in \N\,$,  v\'erifient l'\'egalit\'e~:
$$\fab = \left(\fao\circ\cdots\circ\fa^{(n)}\right) \circ
R_{\a\b}^{(n+2)} \circ \xab  \circ
\left(\fbo\circ\cdots\circ\fb^{(n)}\right)^{-1}\,.$$
D'apr\`es la remarque (\ref{rem2}) la suite
$\left(\fao\circ\cdots\circ\fa^{(n)}\right)_{n\in\N}$ admet, pour la topologie de krull de $\hG_k(U_{\a}) \,$ une limite
$\Phi_{\a}$ et $R_{\a\b}^{n+2}$ tend vers l'identit\'e. Par continuit\'e on
obtient la d\'ecom\-po\-sition
cherch\'ee~:
$$\fab = \Phi_{\a} \circ  \xab \circ  \Phi_{\b}^{-1}\,, \qquad
\qquad (\Phi_{\a}) \in
Z^0(\U; \hG_k)\,.$$
En appliquant cette m\^eme m\'ethode au cocycle
$ \left(\fa\circ\fabp \circ\fb^{-1}\right)\,$, on obtient l'\'egalit\'e
$$\left[\fa\circ\fabp \circ\fb^{-1}\right]_{\hG_k} =
\left[\xab\right]_{\hG_k}\,.$$
Et donc~: $\left[\fa\circ\fabp \circ\fb^{-1}\right]_{\hG_k} =
[\fab]_{\hG_k}\,$.\\

Les d\'emonstrations de l'injectivit\'e de $\chi$ et de $\widehat{\chi}$ sont
similaires; ce qui ach\`eve la
d\'emonstration.
\end{proof}

%%%%%%%%%%%%%%%%%%%%%%%%%%%%%%%%

\section{D\'eformations de feuilletages}\label{sec.def.feuil}  Soit M une vari\'et\'e
holomorphe connexe de dimension~$n$, notons $\O_M$ le faisceau des germes
de fonctions
holomorphes, $\L_{M}$ le faisceau des germes de 1-formes  diff\'erentielles
holomorphes et $\X_M$ le
faisceau des germes de champs de vecteurs holomorphes sur $M\,$. On
d\'esigne  toujours par $P$ le
germe $\left(
\Bbb{C}^p , \, 0 \right)\,$.\\

{\bf Dans ce chapitre et dans tous les suivants nous emploierons le terme "arbre"
pour "germe d'arbre".}

\subsection{G\'en\'eralit\'es sur les feuilletages}\label{gen.feuil}
\addcontentsline{toc}{section}{\hspace{0,8em} {}\thesubsection .  G\'en\'eralit\'es sur les feuilletages}
D\'efinissons un {\it feuilletage holomorphe de codimension q
(\'eventuellement singulier) sur $M$}
comme la donn\'ee d'un faisceau de sous-modules $\Lambda _{\cal F}$ de
$\Lambda _M\,$ localement
libres  de rang
$q$, v\'erifiant en chaque point $m$ les relations d'int\'egrabilit\'e
suivantes~:
\begin{equation}\label{ecr.def}
\w^{1}\wedge\cdots\wedge \w^{q}\wedge  d\w^{i} = 0 \ , \qquad \,\forall i
= 1,\ldots, q\quad
\forall \ \w^{1},\ldots , \w^{q} \in\Lambda _{{\cal F}, m}
\end{equation}
Cette d\'efinition n'est pas la plus g\'en\'erale
\footnote{On aurait pu seulement imposer \`a $\Lambda_{\F}$ d'\^etre  coh\'erent
de {\it rang g\'en\'erique}
$q$, c'est \`a dire de rang $q$ presque partout.}  mais nous n'aurons \`a
faire ici qu'\`a deux types de
feuilletages, les feuilletages de codimension 1,
\begin{equation}\label{rel.int.}
\Lambda _{\F,\, m} = (\w_{m})\  , \qquad  \w_{m}\wedge d\w_{m}= 0 \,,
\end{equation} ou bien leurs "d\'eformations \`a $p$ param\`etres" que nous
d\'efinirons pr\'ecis\'ement en
(\ref{def.equired.})~:
 $$\Lambda_{\F_P} = (dt_1,\ldots , dt_p, \Omega)\ ,\ \qquad dt_1 \wedge
\ldots \wedge dt_p \wedge
\Omega \wedge d\Omega= 0\, , \qquad
\Omega\mid {}_{\{t = 0\}} = \w\,,$$
 avec $t_1, \ldots , t_p$ les coordonn\'ees canoniques sur $P := \left(\Bbb
{C}^p, \, 0
\right)\,$.

\vspace{1 em}

En dehors du {\it lieu singulier de $\F$}, i.e. le sous-ensemble
analytique $Sing(\F)$ d\'efini par
le radical du faisceau d'id\'eaux
$$ I_{\F } := (\buildrel q\over {\wedge}\!\Lambda
_{\cal F})\cdot(\buildrel q\over \otimes\cal X_M) \,,
$$
on a un feuilletage r\'egulier de codimension $q$ not\'e $\F^{reg}\,$.
Lorsque
$Sing(\F )$ n'a pas de composante de codimension 1, on voit facilement
qu'aucun point singulier
$m$ n'est "\'eliminable": on ne peut pas prolonger $\F^{reg}$ au voisinage
de $m$ en un feuilletage
r\'egulier. En divisant localement les g\'en\'erateurs de  $\Lambda _{{\cal
F},m}$ par le $p.g.c.d.$ de
leurs  coefficients, on construit un unique feuilletage $sat(\F)$ appel\'e
{\it satur\'e de $\F$}, qui
co\"{\i}ncide avec $\F$ sur la partie r\'eguli\`ere $M-Sing(\F )\,$ de $M$ et qui
est maximal pour cette
derni\`ere propri\'et\'e. On a l'\'egalit\'e~:

$$\Lambda_{sat(\F),\, m} = \left\{\,
\n \in \Lambda_{M,\, m} \ ;Ê\ \n\wedge \w^1\wedge \cdots \wedge \w^q = 0 \
,
 \ \ \forall  \ \w^1,\ldots , \w^q \in \Lambda _{\cal F,\, m}\, \right\} $$

{\it L'image r\'eciproque $f^{-1}\F$ de $\F$} par une application holomorphe
$f$ d'une vari\'et\'e
holomorphe $M'$ dans $M$  est le feuilletage de codimension $q$ localement
d\'efini par les images
r\'eciproques
$f^{\ast}\w\,$, $\w \in \Lambda _{\F , \, m}$ lorsque celles-ci engendrent
un faisceau  localement
libre de rang $q$. Quand $f$ est un plongement et les $f^{\ast}\w\,$, $\w
\in
\Lambda _{\F , \, m}$ tous identiquement nuls, nous dirons que $M'$
est une {\it vari\'et\'e
int\'egrale} ou {\it invariante} de
$\F\,$. Lorsque $M'$ est la partie r\'eguli\`ere d'un sous-ensemble analytique
ferm\'e $S$ de $M$, nous
dirons que $S$ est un {\it ensemble int\'egral} ou encore {\it invariant} de
$\F\,$.

\begin{defi}\label{trans.strict}  Nous  appelons {\it transform\'e strict de
$\F$ par $f$} et
notons  $f^{\ast}\F$ le satur\'e de
$f^{-1}\F\,$.
\end{defi}

Lorsque f est un biholomorphisme, les deux feuilletages $\F$ et $\cal H :=
f^{-1}\F\,$ sont dits
{\it holomorphiquement  conjugu\'es}; et l'on note $\F \sim _{hol} \cal
H\,$. Signalons la propri\'et\'e
d'image directe suivante qui s'obtient facilement \`a l'aide du th\'eor\`eme
d'extension de Levi des
fonctions m\'eromorphes.

\begin{prop}\label{im.dir} Soient $h : M \fle M'$ une application propre
entre deux vari\'et\'es
holomorphes de m\^eme dimension et
$\F$ un feuilletage holomorphe (singulier) sur $M\,$. Supposons que  la
codimension de
$h(Sing(\F))$ soit $\geq 2\,$ et que la restriction de $h$ au
compl\'ementaire de $Sing(\F)$ soit un
diff\'eomorphisme sur son image. Alors il existe  sur $M'$ un feuilletage
unique not\'e $h_{\ast}\F\,$
et appel\'e image directe de $\F$ par $h$ tel que
$h^{\ast}h_{\ast}\F = \,sat(\F)\,$.
\end{prop}

Nous aurons aussi \`a consid\'erer la position relative  d'un feuilletage par
rapport \`a une
hypersurface \`a croisement normal. Cela nous am\`ene \`a poser~:

\begin{defi}\label{singcouples} Soit $\F$ un  feuilletage holomorphe
satur\'e de codimension $q$ au
voisinage d'un ensemble analytique $S$. On appelle {\it lieu singulier du
couple $(\F, S)\,$}
l'ensemble analytique $Sing(\F, S)\,$ d\'efini par le radical du faisceau
d'id\'eaux
$${\cal J}_{\F,\, S} = sat \,\left( I_{\F, \, S} \right) \, , \qquad
I_{\F,\, S} := (\buildrel
q\over {\wedge}\!\Lambda _{\F} \,)\cdot\, (\buildrel q\over
\otimes\X_{M,\, S})\,,$$ o\`u : $\X_{M, S} \subset \X_M$  est le
sous-faisceau des germes de champs
de vecteurs holomorphes sur $M\,$ tangents \`a $S\,$  (i.e. \`a la partie
lisse de $S$) et o\`u, pour un id\'eal $I := (f_1,\ldots,f_r)$ de  $\O_{M,m}$, on
d\'esigne par
$sat
\,(I)$ l'id\'eal
$\left( {f_1/g} , \ldots ,{f_r/g} \right)\,$, $g := p.g.c.d.(f_1,\ldots,f_r)$

\end{defi}

Visiblement $Sing(\F)$ et $Sing(\F, S)$ co\"{\i}ncident en dehors de $S$.  On
voit facilement~:

\begin{prop}\label{critsing} Soit $S$ une hypersurface de $M$
\`a croisements normaux. Un
point $m\in S$ appartient \`a
$M - Sing\,(\F, S)$ si et seulement si en ce point $\F$ est r\'egulier et
chaque composante
irr\'eductible locale de
$S$ est soit invariante soit transverse \`a $\F\,$ en tout point.
\end{prop}

Si l'on n'avait pas effectu\'e dans (\ref{singcouples}) l'op\'eration de
saturation de $I_{\F,\, S}$
on aurait de plus comme points singuliers ceux des composantes de $S$ qui
sont vari\'et\'es int\'egrales
de $\F$.\\

Consid\'erons
maintenant un sous espace analytique $S = (\mid\!S\!\mid\,  ,\, {\cal
O}_M/I_S)$ de $M\,$, non
n\'ecessairement r\'eduit. D\'esignons par
$\widehat{M}^S$ l'espace annel\'e donn\'e par

\begin{equation}\label{esptrform}
 |\widehat{M}^{S}| \, := |S| , \, \hO_M^S \,: =\, \limproj_{\,k \in
\N}\left({ {\iota}^{-1}(\O_M)
}\over {\,{\iota}^{-1}({I_S}^{k+1})} \right) \,,
\end{equation} o\`u $\iota : S \ifle M$ d\'esigne l'application d'inclusion.
Nous dirons que
$\widehat{M}^{S}$ est un {\it espace transversalement formel} et   les
\'el\'ements de $\hO_M^S$
seront appel\'es {\it germes de  fonctions transversalement formelles le
long de $S\,$}.

\vspace{1 em}

Classiquement l'espace transversalement formel $\widehat{M}^S$ s'identifie
\`a celui obtenu de la
m\^eme mani\`ere \`a partir de l'espace analytique r\'eduit associ\'e \`a $S$.

\vspace{1 em}

Nous n'aurons \`a consid\'erer ici que des sous-espaces de dimension 0, (o\`u
l'on retrouve la notion
usuelle de s\'erie formelle), et des hypersurfaces localement d\'efinies dans
de bonnes  coordonn\'ees
$(z_1,\ldots ,z_n)$ par une \'equation du type
$z_1 = 0\,$, resp. $z_1z_2 = 0\,$. Les \'el\'ements de  $\hO_{M,m}^S$
s'expriment alors
respectivement comme les s\'eries~:
\begin{equation}\label{explfora}
\sum_{k=0}^{\infty} \,A_k(z_2, \ldots ,z_n)\,z_1^k\, ,
\end{equation}
{ou bien}
\begin{equation}\label{explforb}
\sum_{k=0}^{\infty}\,\left(A_k^0(z_3 \ldots , z_n) +  A_k^1(z_1 ; \, z_3,
\ldots z_n) +  A_k^2(z_2
;\, z_3 \ldots , z_n)\right)\!(z_1\,z_2)^k\, ,
\end{equation}
\noindent les rayons de convergence des coefficients $A_k\,$, $A_k^j$
\'etant tous minor\'es par un
m\^eme r\'eel $>0\,$.

\vspace{1 em}

Soit $\cal J$ un faisceau coh\'erent d'id\'eaux de $\hO_M^S$ et $S_{\cal
J}:=\, supp(\hO_M^S/{\cal J})
\subset S$. On appelle {\it sous-espace de $M$ transversalement formel le
long de $S\,$} l'espace
annel\'e $Z_{\cal J}= (S_{\cal J}\, ,\hO_M^S/{\cal J})$. Lorsque l'id\'eal
$\cal J$ est \'egal \`a son
radical, on dit que $Z_{\cal J}$ est un {\it sous-ensemble de $M$
transversalement formel le long
de $S\,$}. L'{\it intersection de deux sous-ensembles transversalement
formels} est d\'efini par le
radical du produit de leurs id\'eaux. Lorsque $I$ est un faisceau d'id\'eaux
de $\O_M$ (avec $I$ \'egal
\`a son radical), nous notons
$\widehat{I}^S\,=\,I
\otimes_{\O_M}
\hO_M^S\,$. Dans ce cas $S_{\widehat{I}^S} = Z \cap S$ o\`u $Z$ est le
sous-ensemble analytique de
$M$ d\'efini par $I$.
\\

Soit $S'$ un sous-espace analytique d'une vari\'et\'e $M'$. Une {\it
application transversalement
formelle} $f :{\widehat{M'}}^{S'} \rightarrow \widehat{M}^S\,$, que l'on
notera aussi $f :(M',S')
\rightarrow (M,S)$ est un morphisme d'espaces annel\'es donn\'e par
$ |f| : |S'| \rightarrow |S|$ et $f^{\ast} : \hO_M^S \rightarrow
\hO_{M'}^{S'}$.

\vspace{1 em}

Pour tout sous-ensemble quelconque $K'$ de $S'$, on appelle {\it germe de
$\widehat{M'}^{S'}$ le
long de $K'$}, le sous-espace annel\'e $(K', i^{-1}(\hO_{M'}^{S'}))$ o\`u $i :
K'
\hookrightarrow S'$ d\'esigne l'application d'inclusion. Le  {\it germe  de
$f$ le long d'un
sous-ensemble $K' \subset S'$}, not\'e $f :(M',K')
\rightarrow (M,S)$, est un morphisme d'espace annel\'e entre les espaces
annel\'es $(K',
i^{-1}(\hO_{M'}^{S'}))$ et ${\widehat M}^{S}$.

\vspace{1 em}

 Via l'extension des scalaires
$\iota^{-1}(\O_M) \hookrightarrow {\hO}_M^S\,$ on obtient les faisceaux
(de base $S$) des  {\it
1-formes diff\'erentielles} et des {\it champs de vecteurs transversalement
formels le long de $S\,$}
en posant~:
$$\widehat{\L}_M^S := \iota^{-1}(
\L_{M})\otimes_{\iota^{-1}(\O_M)}\hO_M^S\,,\qquad
\hX_M^S := \iota^{-1}(\X_M)\otimes_{\iota^{-1}(\O_M)}{\hO}_M^S \,.$$
Lorsque $S$ est r\'eduit \`a un
point, $S=\left\{m\right\}$, on retrouve les notions usuelles de s\'eries,
1-formes diff\'erentielles
et champs de vecteurs formels et nous notons simplement ces espaces par
${\hO}_{M,\,m}\,$,
$\widehat{\L}_{M,\, m}\,$,
$\hX_{M,\, m}\,$.

\begin{defi}\label{germefeuilltrf}  Un {\it feuilletage $\hF$ de $M$ de
codimension $q$
transversalement formel le long d'un sous-espace analytique
$S\,\subset M $} est la donn\'ee
 d'un faisceau $\L_{\hF}$ de sous-modules localement libres de rang $q$ du
faisceau
$\widehat{\L}_M^S$ qui satisfait en chaque point $m$ de $S$ les relations
d'int\'egrabilit\'e
(\ref{rel.int.}).
Le {\it germe de $\hF$ le long
d'un sous-ensemble $K$}
de $S$, not\'e $(\hF,K)$ est la donn\'ee de la restriction $\L_{\hF,K} = i
^{-1} (\L_{\hF})$, $i : K \hookrightarrow S$ d\'esignant l'application d'inclusion.
\end{defi}

\begin{rema}\label{im.dir.for} Avec ces d\'efinitions, on \'etend sans peine \`a
un feuilletage
transversalement formel $\hF$ le long de
$S$, les notions de lieux singulier, de saturation, ainsi que les notions
d'image r\'eciproque
$f^{-1}(\hF)$ et de transform\'e strict $f^{\ast}(\hF)$ par une application
transversalement
formelle  $f : {\widehat{M'}} ^{S'}
\rightarrow
\widehat M ^{S}$. On dispose aussi d'une notion de sous-espace
transversalement formel invariant
par $\hF$. On peut aussi d\'efinir la  relation de conjugaison formelle
entre deux feuilletages
transversalement formels $\hF$ et $\hF'$ par un isomorphisme
transversalement formel, que l'on
notera par $\hF \sim _{for} \hF'\,$.
\\

Remarquons aussi que toute application holomorphe $h$ d'une vari\'et\'e
holomorphe $M'$ dans une
vari\'et\'e holomorphe $M$ induit une application transversalement formelle
${\widehat M'}{}^{S'}
\rightarrow
\widehat{M}^{S}$ avec ici $S' = h^{-1}(S)$ i.e. $|S'|:=h^{-1}(|S|)$,
$\O_{S'} :=\O_{M'}/I_{S'}$ et
$I_{S'}:= h^{\ast}(I_S) \O_{M'}$. En particulier, l'image r\'eciproque par
$h$ d'un feuilletage
transversalement formel le long de $S$ est un feuilletage transversalement
formel le long de
$h^{-1}(S)$.  D'autre part l'op\'eration d'image directe commute avec la
compl\'etion, d'apr\`es les
th\'eor\`emes de comparaison  de Grauert
\cite[page 269]{B-S}. La proposition de type Hartogs (\ref{im.dir}) s'\'etend ainsi \`a
l'image directe d'un
feuilletage transversalement formel  par une application holomorphe. Ainsi lorsqu'on
consid\`ere un feuilletage formel \`a
l'origine de
$\Bbb {C}^2$  d\'efini par une 1-forme
formelle
\begin{equation}\label{fdfstan}
\w := \stw\;\; \in \;\; \hL_{\Bbb {C}^2 ,\, 0}
\end{equation} son transform\'e strict par une succession finie $E$
d'\'eclatements en des points sera
un feuilletage transversalement formel le long de $E^{-1}(0)$.
R\'eciproquement tout feuilletage
transversalement formel le long de $E^{-1}(0)$ provient d'un feuilletage
formel satur\'e \`a l'origine
de $\Bbb {C}^2$.
\end{rema}

\begin{rema}\label{flowbox} Supposons \`a pr\'esent que $S$ est une
sous-vari\'et\'e analytique de $M$, de
codimension \'egale \`a la codimension
$q$ du feuilletage $\hF$, et invariante par $\hF$. Fixons un point $m
\in S - (Sing(\hF) \cap
S)$. On peut voir qu'il existe des coordonn\'ees $z_1 , z_2 , \ldots , z_n$
transversalement
formelles  de $\hO_{M,m}^S$ telles que $\L_{\hF,m}$ soit engendr\'e par
$dz_1 , dz_2 , \ldots , dz_q
$. C'est la version transversalement formelle du Th\'eor\`eme de Frobenius qui
se montre de mani\`ere
identique. En effet, \'etant donn\'es une submersion $R : (M,m) \rightarrow (S,m)$
et un germe en
$m$ de champ transversalement formel $X$ qui se projette en un champ
holomorphe sur $S$ via $R$, on
voit facilement qu'il existe une application $\Phi : R^{-1}(m) \times \Bbb
{C} \rightarrow
(\widehat{M}^{S},m)$ transversalement formelle le long de $m \times \Bbb
{C} $ telle que
$\frac {\partial \Phi}{\partial t} = X \circ \Phi$, o\`u $t$ d\'esigne la
variable sur $\Bbb {C}$
\cite[page 505]{M-M}.

Ainsi on d\'efinit comme d'habitude {\it l'holonomie de $\hF$ le
long de $S - (Sing(\hF)
\cap S)$} comme une repr\'esentation  $\rho_S$ de $\pi_1(S - (Sing(\hF) \cap
S), m_0)$ dans
$\widehat{Diff}(T,m_0)$ o\`u $(T,m_0)$ est un germe de sous-vari\'et\'e
holomorphe de dimension $q$
transverse \`a $S$ en un point fix\'e $m_0 \in S$ et $\widehat{Diff}(T,m_0)$
d\'esigne le groupe des
diff\'eomorphismes formels de $(T,m_0)$. L'image  $H_S$ de $\rho_S$
s'appelle le {\it groupe
d'holonomie de $S$ par rapport \`a $\hF$}.

Comme dans le cas holomorphe $\rho_S$ classifie le
germe de $\hF$ le long de l'ensemble invariant $S$. Plus
pr\'ecis\'ement, donnons nous $R :
(M,S)
\rightarrow S$ un germe de submersion le long de $S$, avec $R^{-1}(m_0) =
T$, la restriction de
$R$ \`a $S$ \'etant l'identit\'e. Alors les repr\'esentations d'holonomie
respectives $\rho_S$ et $\rho'_S$
de deux feuilletages $\hF$ et $\hF'$ de $M$ transversalement formels le
long de $S$ sont
formellement conjugu\'ees
$\rho'_S = \phi \circ \rho_S \circ \phi^{-1}$ avec $\phi \in
\widehat{Diff}(T,m_0)$ si et
seulement s'il existe un isomorphisme transversalement formel $\Phi : (M,
S) \rightarrow (M,S)$
tel que $\Phi^{\ast}\hF = \hF'$.
\end{rema}

\subsection{D\'eformations formelles d'un feuilletage}\label{def.equired.}
\addcontentsline{toc}{section}{\hspace{0,8em} {}\thesubsection .  D\'eformations formelles d'un feuilletage}
Soit
$\F$ un feuilletage formel \`a l'origine de $\Bbb{C}^2$  d\'efini par une 1-forme formelle
\begin{equation}\label{fdfst}
\w := \stw\;\; \in \;\; \hL_{\Bbb{C}^2 ,\, 0}
\end{equation} et soit $P$ le germe $\left( \Bbb{C}^p , \,0 \right)\,$.
Une {\it d\'eformation
formelle de $\F$ le long de $0 \times P\,$}, de para\-m\`etres
$P\,$,  est la donn\'ee \`a unit\'e multiplicative pr\`es, d'une forme
diff\'erentielle formelle du type
\begin{equation}\label{standdef}
\n := \std
\end{equation}
telle que
$A (x, y; t) =
\sum_{j, k = 0}^{\infty} A_{jk}(t)\,x^j y^k  \,$, $B(x, y ; t) =
\sum_{j, k = 0}^{\infty} B_{jk}(t)\,x^j y^k \in
\Bbb{C}\{\!t\!\}[[x,\,y]]\,,$  avec $t := (t_1,
\ldots ,t_p) \, ,$ les rayons de convergence des $A_{jk}(t)\,$,
$B_{jk}(t)\,$ \'etant minor\'es par un
r\'eel
$> 0\,$ commun et
$A(x, y ; 0) = a(x, y) \,$, $B(x, y ; 0) = b(x, y) \,$, $A(0,
0 ; t) =  B(0, 0 ; t) =
0\,$.

Plus g\'en\'eralement, consid\'erons un feuilletage $\F$ de codimension 1
transversalement formel le
long d'un sous-ensemble analytique $\,S\,$ (non n\'ecessairement ferm\'e)
d'une vari\'et\'e holomorphe
$\,M\,\,$.

\begin{defi}\label{deftrform} Une {\it d\'eformation transversalement formelle de
$(M,\, S, \,\F)$
para\-m\`etr\'ee par
$P\,$}, est la donn\'ee
\begin{equation}\label{equ: deformation}
\sF_P = \left(\M_P, \,\SS_P,\, \pi,\, \s ; \,\F_P\right) \,,
\end{equation}
\noindent {\bf 1) }d'une {\it d\'eformation $\left(\M_P, \,\SS_P,\, \pi
\right)$ de $(M, S)$ de base
$P$}, c'est \`a dire~:
\begin{enumerate}
\item [a)] d'un germe le long de $S$ de plongement holomorphe $\sigma$ de
$(M, S)$ dans une
vari\'et\'e holomorphe $\M_P$ munie d'un sous-ensemble analytique ${\cal S}_P$
tel que
$\sigma (S) \subset {\cal S}_P\,$,
\item [b)]   d'un germe de submersion  $\pi : ({\cal M}_P, \sigma (S))
\rightarrow P$ tel que la
restriction de  $\pi$ \`a ${\cal S}_P$ soit plate et v\'erifie~: $\pi^{-1}(0)
= \sigma (M)\,$,
$\pi^{-1}(0)\cap {\cal S}_P = \sigma (S)\,$,
\end{enumerate}
\noindent {\bf 2) }du germe $\F_P$ le long de $\sigma(S)$ d'un feuilletage
transversalement formel
le long de
${\cal S}_P\,$, de m\^eme dimension que $\F$ qui v\'erifie~:
\begin{enumerate}
\item[a)]  les feuilles de $\F_P$ sont contenues dans les fibres de
$\pi\,$,
\item[b)]  $\s^{\ast}\F_P = \F\,$.
\end{enumerate}
\end{defi}

La condition 2.a. signifie que tout germe en un point de $\sigma(S)$ d'un
champ de vecteurs $X$
transversalement formel le long de ${\cal S}_P$ qui annule $\L_{\F_P}$
est  {\it vertical}, i.e.
$T\pi\! \cdot\! X \equiv 0\,$.
\\

Remarquons que lieu singulier $Sing\,(\F_P, \SS_P)$, que l'on d\'efini comme en  (\ref{singcouples}), est aussi
le sous-ensemble de $M$
transversalement formel le long de $S$ donn\'e par le radical du faisceau
d'ideaux
\begin{equation}\label{idealredef} sat\left({\Lambda _{\cal F_P}}\!\cdot
\hX_{\!\M_{P/\pi},
\SS_P}\right)\,,
\end{equation}
o\`u $\,\hX_{\!\M_{P/\pi}, \SS_P}\,$ d\'esigne le faisceau des germes de
champs de vecteurs verticaux
transversalement formels le long de $\SS_P$ et tangents \`a $\SS_P\,$.
Nous appellerons cet ensemble transversalement formel {\it le lieu singulier
relatif de la d\'eformation $\sF_P$} et le notons
$Sing_P\,(\sF_P)\,$.

\begin{rema}\label{defeno} Lorsque la restriction de $\pi$ \`a $\SS_P$ est
propre on peut, pour des
petites valeurs du param\`etre $t \in P\,$, consid\'erer les fibres des
donn\'ees de (\ref{equ:
deformation}) au dessus de $t\,$. On obtient un espace analytique compact
r\'eduit
$\SS_P(t) := \SS_P \cap \pi^{-1}(t)\,$, un germe de vari\'et\'e $\M_P(t) :=
\left(  \pi^{-1}(t) , \,
\SS_P(t) \right)\,$  et un feuilletage $\F_P(t)\,$ sur $\M_P(t)$
transversalement formel le long
de $\SS_P(t)$, obtenu en restreignant $\F_P$ \`a $
\pi^{-1}(t)\,$.  La propret\'e de $\pi$ permet donc de consid\'erer $\sF_P$
comme le germe en $0 \in
P$ d'une "famille" de feuilletages transversalement formels

$$\sF_P(t) := \left( \M_P(t) , \, \SS_P(t) , \, \F_P(t) \right)\,.$$
D'apr\`es (\ref{idealredef}) le
lieu singulier relatif de $\sF_P$ est la "famille" des lieux singuliers
des $\F_P(t)$.
\end{rema}

\begin{rema}\label{defforme} Lorsque $\F$ est un feuilletage formel \`a
l'origine de $\Bbb {C}^2$, i.e. $M
= \Bbb {C}^2$ et $S = \{ 0
\}$, la condition de platitude de la restriction de $\pi$ \`a
$\SS_P$ lui impose d'\^etre \'etale au dessus de $P$ et l'on peut supposer que
$\SS_P = 0 \times P\,$.
On retrouve alors la notion de d\'eformation formelle le long de $0 \times
P\,$ d\'efinie par
(\ref{standdef}). Le lieu singulier relatif de $\sF_P$ est d\'efini par les
\'equations $A(x, y ; \,t)
=  B(x, y ; \,t) = 0\,$ et le feuilletage formel $\F_P(t)$ est alors donn\'e
par la 1-forme formelle
$\n_t \in \hL_{\Bbb{C}^2Ê\times t , \, (0 , t)}$ obtenue en fixant la
valeur du  param\`etre $t \in
P$ dans (\ref{standdef}).
\end{rema}

Donn\'e un germe d'application holomorphe $\,\r : Q := \left(\Bbb{C}^q, \,0
\right) \fle P\,$ on
appelle  {\it d\'eformation obtenue \`a partir de $\,\sF_P\,$ par le
changement de param\`etres
$\,\r\,$} la d\'eformation transversalement formelle
$$\r^{\ast}\!\sF_P := (\r^{\ast}\!\M_P ,\, \r^{\ast}\!\SS_P,\,
\r^{\ast}\!\pi,\, \r^{\ast}\!\s
;\,\r^{\ast}\!\F_P)$$  dont les composantes sont d\'efinies par les
diagrammes cart\'esiens
$$\hbox{$ {\begin{array}{ccc} {\r^{\ast}\!\M_P} &
{\!\!\!\buildrel{i_{\r}}\over \longrightarrow} &
{\M_P} \\ {\r^{\ast}\pi}{\downarrow}\phantom{\r^{\ast}\pi}
&{\!\!\!\Box}&{\pi\downarrow
\phantom{\pi}}\\  {Q}&{\buildrel{\r}\over \longrightarrow} &{P}
\end{array}} {\qquad } {\begin{array}{ccc} {\r^{\ast}\!\SS_P} &
{\!\!\!\buildrel{\!\!\!i_{\r}}\over \longrightarrow} & {\SS_P} \\
{\r^{\ast}\pi}{\downarrow}\phantom{\r^{\ast}\pi}
&{\!\!\!\Box}&{\pi\downarrow
\phantom{\pi}}\\  {Q}&{\!\!\!\buildrel{\r}\over \longrightarrow} &{P}
\end{array}}
$}$$ avec $\r^{\ast}\F_P := {i_{\r}}^{\ast}\F_P$.  Comme d'habitude cette
op\'eration est
contravariante.
\\

Deux d\'eformations transversalement formelles de $(M,\, S, \,\F)$ de m\^eme
espace $P$ de
pa\-ra\-m\`et\-res
$$
\sF_P := (\M_P,\, \SS_P, \,\pi,\, \s; \,\,\F_P)\quad \mbox{et} \quad\sF '_P := ( {\M'}_P,
\,\SS'_P,
\, {\pi'}, \, {\S'}; \,\F'_P)
$$ sont dites {\it formellement conjugu\'ees} s'il
existe un germe
d'isomorphisme $\Phi$ transversalement formel qui conjugue
$\F_P\,$ \`a $\F'_P\,$, fait commuter les submersions $\pi$ et $ {\pi'}\,$  et
vaut l'identit\'e au
dessus de l'origine :
$$\Phi : (\M_P, \s(S)) \iso ( {\M'}_P,  {\S'}(S))\, \quad
 {\pi'} \circ \Phi = \pi\,, \quad
\Phi\circ \s =  {\S'}\,, \quad \Phi^{\ast}\F'_P=\F_P\,.$$
On notera alors
$\Phi^{\ast}\sF'_P = \sF_P\,$, ou encore
$\sF'_P {\sim_{P, \,for}} \sF_P\,$. On dit que $\sF_P$ est {\it formellement triviale} s'il existe
une conjugaison formelle
entre
$\sF_P$ et la {\it d\'eformation constante}~:
\begin{equation}\label{defcte}
\sF_P^{cst}:= \left(M \times P,\, S \times P,\, \pi,\, i; \,\F_P^{cst}
\right),
\end{equation} o\`u
$\pi$ d\'esigne la projection canonique de $M\times P\,$ sur $P\,$, $i\,$
l'inclusion de $M \cong M
\times {0}\,$ dans $M\times P\,$ et $\F_P^{cst}$ est le {\it feuilletage
transversalement formel
constant} d\'efini comme $pr^{-1}\F$, avec $pr : M \times P
\rightarrow M$ la projection canonique.

\subsection{Construction de d\'eformations formelles par collage}
\addcontentsline{toc}{section}{\hspace{0,8em} {}\thesubsection .  Construction de d\'eformations formelles par collage}
Fixons un
feuilletage formel $\F$ \`a l'origine de $\Bbb{C}^2\,$, et un germe d'arbre
non n\'ecessairement
r\'egulier
\begin{equation}\label{arbnot}
\A_P = \left( \M^j\,, E^j\,, \S^j\,, S^j\,, \pi_j\,, \D^j\,\right)_{j =
0,\ldots ,h}
\end{equation} au dessus de $P := \left( \Bbb{C}^p , \, 0 \right)\,$. Pour
les donn\'ees de cimes de
$\A_P$ et celles du germe d'arbre $\A_P(0)$ d\'efini en (\ref{arbrfibr}),
nous conservons les
notations (\ref{cime}). Le transform\'e strict $\F'= \wE_0^{\ast}\F\,$ de
$\F$ par $\wE_0\,$, cf.
(\ref{trans.strict}) est un feuilletage transversalement formel le long du
diviseur de cime $\wD_0
\subset \wD\,$ de
$\A_P(0)\,$. Comme dans (\ref{germefeuilltrf}) on d\'esigne par $\left( \F' ,Ê\, W \right)$
le
germe de $\F'\,$ le long d'un sous-ensemble quelconque
$W \subset \wD_0$.
Nous supposons dans tout ce qui suit que $\wD_0$ est un ensemble int\'egral
de $\F'\,$.

\begin{defi}\label{comp.arb.fe} Une d\'eformation transversalement formelle
$\sF_P$ de $\F$ le long
de $0 \times P\,$, de param\`etres
$P$ est dite {\it compatible avec l'arbre $\A_P$} si
\begin{enumerate}
\item $\M^0 = \Bbb{C}^2 \times P\,$, $\S^0 = S^0 = 0 \times P$ et pour
tout $j = 1,
\ldots h$ le lieu singulier relatif  du transform\'e strict de $\F_P$ sur
$\M^j$ intersecte $\D^j$
suivant $\S^j\,$.
\item le diviseur de cime $\wD$ de $\A_P$ est un ensemble int\'egral du
feuilletage $\F'_P :=
\wE^{\ast}\F_P\,$,
\item la restriction de la projection $\pi : \M^h \rightarrow P$ au lieu
singulier relatif de
$\sF'_P$ est finie (c'est \`a dire $\pi_{\ast}(\hO_{Sing_P(\sF_P)})$ est un
$\O_P$-module fini),
\item $Sing_P\left(\sF'_P\right) \cap \wD_0 = Sing(\F')\,$.
\end{enumerate}
\end{defi}

\vspace{1 em}

Nous allons donner des proc\'ed\'es de construction de telles d\'eformations.
pour cela nous faisons les deux
hypoth\`eses suivantes :
\begin{itemize}
\item l'ensemble $\wS_0$ des singularit\'es de la
cime $\A_P(0)$ est \'egal au
lieu singulier de $\F'$,
\item le diviseur $\wD_0$ est un ensemble
int\'egral de $\F'$.
\end{itemize}
Reprenons le vocabulaire et les notations introduites dans la definition (\ref{doncrit}).

\begin{defi}\label{bon.syst.lac} Nous appelons {\it bon syst\`eme de lacets}
associ\'e \`a un \'el\'ement
critique $K$ toute collection $\Gamma_K$ de lacets telle que ~:
\begin{enumerate}
\item Si $K \in \bC_0(\A_P) $ alors $\Gamma_K =(\g_{K, L})_{L\in Ad_1(K)}$ o\`u $\g_{K, L}$ est un lacet simple trac\'e sur
$L$ (cf. (\ref{adjac})) bordant un disque conforme $W_{K,\, L} \subset \overline L$ v\'erifiant
$\overline{W_{K, \,L}} \cap Sing(\F') = \{ K \}$;
\item Si $K \in \bC_1(\A_P) $ alors  $\Gamma_K$ est une collection $(\g_{K,\, m})_{m\in Ad_0(K)}$ de
lacets  simples trac\'es sur $K$ de m\^eme origine $m_0$ et bordant des disques conformes $V_{K, m}\subset \overline{K}$ tels
que
$\overline{V_{K,m}}\cap Sing(\F') = \{m\}$ et $\overline{V_{K,m}}\cap \overline{V_{K,m'}} = \{m_0\}$ pour $m \not= m'$.
\end{enumerate}
\end{defi}
Fixons  un germe de courbe lisse $(T, m_0)$ transverse \`a $K \in \bC_1(\A_P)$ en  l'origine commune d'un bon syst\`eme de
lacet $\Gamma_K$ $m_0$  ainsi qu'un bon syst\`eme de
lacets. Ordonnons et orientons les lacets de ce syst\`eme~: $\Gamma_K = \left( \gamma_j\right)_{j = 1, \ldots , v(K)}$ de mani\`ere
que leurs classes dans $\pi_1( K, m_0)$ v\'eri\-fient~: $\ptsur\gamma_{v(K)}^{-1} = \ptsur\gamma_1 \circ \cdots \circ
\ptsur\gamma_{v(K)-1}$. Notons $\underline{h}_j \in \hO_{T_, m_0}$ l'holonomie du feuilletage $\F'$ le long de $\gamma_j$, $j =
1,
\ldots , v(K)$.
\begin{prop}\label{constrcol} Pour chaque $j = 1, \ldots , v(K)$, donnons-nous un germe de
famille de
diff\'eomorphismes formels
$h_{j, \,t} \in \hO_{T, \, m_0 }\,$ d\'ependant analytiquement du param\`etre
$t \in P\,$  et tels que
$$h_{j, \, 0} = \underline{h}_j\,,\qquad \mbox{et} \qquad h_{v(K), \,t}^{-1} = h_{1, \,t}\circ\cdots h_{v(K) - 1, \,t}\,. $$
Il existe alors une d\'eformation transversalement
formelle
$\sF'_{P, \, K}$ de $\left( \wM_0,\, K , \,\F' \right)$ de param\`etres $P$ r\'ealis\'ee
sur le produit
$\left(\wM_0 \times P,\, K \times P\right)\,$ telle que pour chaque $j = 1, \ldots , v(K)$, l'holo\-no\-mie  du
feuilletage (de
dimension 1) $\F'_{P, \, K} \,$  le long du lacet
$\gamma_j$,  calcul\'ee sur la transversale $\left( T \times P , \, (m_0,
0)\right)\,$ est \'egale \`a
$h_j(u, t) := h_{j, \,t}(u)\,$.
\end{prop}

\begin{proof} Nous raisonnons par induction sur $ v=v(K)$. Pour $v(K) =
1\,$, $K$ est
conform\'ement un disque et le r\'esutat est imm\'ediat~: d'apr\`es la remarque
(\ref{flowbox}), les
feuilletages de feuille $K$ et leurs d\'eformations sont tous triviaux le
long de $K\,$.
\\

Pour $v := v(K) > 2$, notons $Ad_0(K)=\{m_1,\ldots ,m_v\} \subset
\overline{K} \subset \Bbb P^1$.
Alors $\overline{K} - \{m_v\}$  est biholomorphe \`a $\Bbb C$ et peut \^etre
recouvert par deux
disques conformes
$U_1$ et
$U_2$ d'intersection connexe et simplement connexe, et tels que
$$ U_1 \cap Ad_0(K)=\{m_1\},\quad U_2
\cap Ad_0(K)=\{m_2,\ldots ,m_{v-1}\},\quad \hbox{avec}$$
$$\g_1 \subset U_1,\quad \hbox{et} \quad \{\g_2,\ldots ,\g_{v-1}\} \subset
U_2/,.$$
Supposons donn\'ee
pour chaque $k = 1, 2$ une d\'eformation transversalement formelle
$\sF'_{P, \, U_k}$ du germe de $\F'$ le long de $U_k$. Par la remarque
(\ref{flowbox}), ces
d\'eformations sont formellement triviales  le long de $U_1 \cap U_2\,$ et
donc formellement
conjugu\'ees.  Soit $\,\hG\,$ le faisceau (\ref{defautf}) de base $\wD_0\,$ d'automorphismes transversalement formels introduit
en (\ref{gro.dif.t}). Donnons
nous  une section $\Phi_{12}\,$ de
$\,\hG\,$ au dessus de
$U_1
\cap U_2,\,$ telle que
$\Phi_{12}^{\ast} \, \sF'_{P , \, U_2}\, =\,  \sF'_{P , \, U_1}\,$.

On utilise alors le r\'esultat suivant, dont une br\`eve preuve sera donn\'ee
plus bas
\begin{equation}\label{sl} H^1(\U,\hG) = 0 \, , \qquad \U := \left( U_1 ,
\, U_2 \right)\,.
\end{equation} Ainsi $\Phi_{12}$ se d\'ecompose en
\begin{equation}\label{equhom}
\Phi_{12} = \Phi_2 \circ \Phi_1^{-1} \, , \qquad \Phi_k \in \hG(U_k),\quad
k=1,2.
\end{equation} Et les d\'eformations $\Phi_1^{\ast} \, \sF'_{P, \, U_1} $ et
$\Phi_2^{\ast} \, \sF'_{P, \, U_2} $ se recollent le long de
$U_1 \cap U_2\,$ en la d\'eformation cherch\'ee.
\\

Pour d\'emarrer l'induction, il reste \`a traiter le cas $v(K) = 2\,$.  Il se
fait comme pr\'ec\'edemment.
L'ensemble $K$ est conform\'ement le plan complexe \'epoint\'e
$\Bbb{C}^{\ast}$ que l'on recouvre par trois secteurs $U_k\,$, $k= 1, 2,
3$ sans intersection 3 \`a 3
et d'intersections 2 \`a 2 connexes et simplement connexes. L'holonomie de
$\F'$ est donn\'ee par un
seul \'el\'ement $h \in \widehat{Diff}(T,m_0)$ o\`u $(T,m_0)$ est une
transversale \`a $K$ en $m_0 \in U_1
\cap U_2\,$. Consid\'erons la d\'eformation $h_t$ de $h$ donn\'ee par
l'hypoth\`ese, et posons
$\widetilde{h}_t= h_t \circ h^{-1}$. Le long de $U_1
\cap U_2,\,$ le feuilletage
$\F'$ est trivial d'apr\`es (\ref{flowbox}). Donnons nous une section $\Phi$
de $\hG$ sur $U_1
\cap U_2,\,$ qui est \'egale \`a $\widetilde{h}_t$ sur $(T,m_0)$. Consid\'erons
pour $ \U = (U_1, U_2,
U_3)$ le cocycle
$(\Phi_{ik}) \in Z^1 \left(\U ; \, \hG \right)$ donn\'e par $\Phi_{12} =
\Phi\,$,$\,\,\Phi_{23} = id$
et
$\Phi_{31} = id$. On applique \`a nouveau l'\'egalit\'e (\ref{sl}) qui est
encore vraie pour le
recouvrement $ \U = (U_1, U_2, U_3)$. Comme pr\'ec\'edemment, on obtient une
d\'eformation
transversalement formelle $\sF'_P$ du germe de $\F'$ le long de $K$. On
v\'erifie que l'holonomie de
cette d\'eformation le long de $K$ est \'egale \`a $h_t$.
\end{proof}

\begin{preuvede}{{de \rm (\ref{sl})}} Cette \'egalit\'e peut se montrer
directement en r\'esol\-vant
l'\'equation cohomologique (\ref{equhom}). Pour cela on explicite cette
\'equation en d\'eve\-loppant en
s\'erie par rapport \`a la coordonn\'ee $z$ qui correspond \`a une
submersion-\'equation  de $K$. En
proc\'edant par induction suivant les coefficients de $z^k$, on est ramen\'e \`a
r\'esoudre l'\'equation
lin\'eaire
$A_{12} = A_2 - A_1 $, $A_{12} \in \O^2_{K}(U_1 \cap U_2)$ d'inconnues
$A_k \in \O^2_{K}(U_k)$.
Ceci est toujours possible puisque $H^1(\U,\O^2_{K})=0$. La d\'emonstration est
similaire dans le cas du
recouvrement de
$\Bbb{C}^{\ast}$ par trois secteurs.
\end{preuvede}

Donnons nous maintenant un arbre $\A_0$ au dessus de $\left\{ \! 0
\! \right\}\,$ compatible
(\ref{comp.arb.fe}) avec un feuilletage formel $\F$ \`a l'origine de
$\Bbb{C}^2\,$ et pour chaque \'el\'ement critique
$K \in \bC\,(\A_0)$ de $\A_0$ un syst\`eme $\Gamma_K := \left( \gamma_{K, L}
\right) \,$ de bons
lacets (\ref{bon.syst.lac}).  Notons encore $\F'$ le feuilletage
transversalement formel le long du
diviseur de cime
$\wD_0$ de
$\A_0\,$, transform\'e strict de $\F\,$ et $(\F', K) \,$ le germe
(\ref{germefeuilltrf}) de $\F'$ le
long de $K \in \bC\,(\A_0)\,$

\begin{defi}\label{ssloc} Une collection de d\'eformations
transver\-sa\-lement formelles des germes $(\F', K)
\,$, ${K \in \bC\,(\A_0)}$,
\begin{equation}\label{sytsl}
\bS := \left(\sF'_{P, K}\right)_{K \in \bC\,(\A_0)}\,,Ê\quad \sF'_{P, K} := \left(W_{P, K},\, \SS_{P, K},\,  {\pi'}_{P, K},\,
{\S'}_{P,  K};\,\F'_{P, K}\right)\,,
\end{equation}
telle que $\SS_{P, K}$ soit une vari\'et\'e int\'egrale de
$\F'_{P, K}$ sera appel\'ee
{\it syst\`eme semi-local de d\'eformations formelles de
$\F'\,$}. Le syt\`eme $\bS$ est dit {\it r\'ealis\'e dans une d\'eformation $\left(\M_P,
\,\SS_P,\, \pi, \, \s
\right)$ de $(\wM_0, S_0)$} si les
$W_{P,
\,K}$ sont des voisinages ouverts de $\s (K)$ dans $\M_P$ qui intersectent
$\SS_P$ suivant
$\SS_{P, K}$, et si de plus $ {\pi'}_{P, K}\,$ et $ {\S'}_{P, K}$ sont des
restrictions de $\pi$ et
$\s\,$.
\end{defi}
Une {\it conjugaison formelle de deux syst\`emes semi-locaux} $\left( \Phi_K
\right)_{K \in \bC\,(\A_0)}\,$ est  une collection  de germes le long de
chaque ensemble critique, de conjugaisons formelles des
d\'eformations constituant les syst\`emes semi-locaux.
\\

Donn\'e un syst\`eme semi-local $\bS = \left(\sF'_{P, K}\right)_{K \in
\bC\,(\A_0)}\,$  de
d\'eformations transversalement formelles de $(\F', K),$ consid\'erons pour
chaque couple d'ensembles
critiques $(m, K)\,$, $dim(K) = 1\,$, $m \in Ad_0(K)$, les  deux familles
de diff\'eomorphismes
formels d\'ependant holomorphiquement de  $t \in P\,$
\begin{itemize}
\item l'holonomie $h_{m,\, K ;\, t}$ du feuilletage $\F'_{P, m}$ le long de
$\g_{m, K}$ et \\
\item l'holonomie $h_{K,\, m ;\, t}$ de $\F'_{P, K}\,$ le long de $\g_{K,
m}$.
\end{itemize} Ces familles sont bien d\'efinies modulo la relation de
conjugaison   par une famille
de diff\'eomorphismes formels d\'ependant holomorphiquement de  $t\,$.

\begin{defi}\label{coherent} Le syst\`eme $\bS$ sera dit {\it coh\'erent} si
pour chaque
$(m, K)\,$ avec ${K \in \bC\,(\A_0)},$ $dim(K) = 1,\,$ et $m \in Ad_0(K)$,
il existe une  une
famille de diff\'eomorphismes formels $\phi_{m, \, K;\,  t}$ d\'ependant
holomorphiquement de
$t\,$, et qui vaut l'identit\'e pour $t=0\,$,  telle que
\begin{equation}\label{systemscoh}
h_{K,\, m;\, t}\, = \, \phi_{m, \, K;\,  t}\circ h_{m,Ê\, K;\, t} \circ
(\phi_{m,\, K; \,
t})^{-1}\,.
\end{equation}
\end{defi}

Consid\'erons une d\'eformation transversalement formelle $\sF_P$ de $\F$
compatible avec un arbre
$\A_P$ au dessus de
$P$ tel que $\A_P(0) = \A_0\,$. Les germes $\left( \F'_P , \, K \right)\,$
de la transform\'ee
stricte
$\F'_P$ de $\F_P$ sur la cime de $\A_P\,$,  le long de chaque ensemble
critique  $K \in
\bC\,(\A_0) = \bC\,(\A_P)\,$ forment visiblement un syst\`eme coh\'erent not\'e
$\bS(\sF'_P)\,$.  R\'eciproquement on a~:

\begin{enonce}{Th\'eor\`eme de r\'ealisation}\label{realis.} Soient $\F$ un feuilletage
formel \`a l'origine de
$\Bbb{C}^2\,$ compatible avec un arbre $\A_0$ au dessus de $\{0\}\,$,  et
$\bS$ un syst\`eme
semi-local coh\'erent de d\'eformations transaversalement formelles du
transform\'e strict
$\F'$ de
$\F$ sur la cime de $\A_0\,$. Il existe alors un germe d'arbre $\A_P$ au
dessus de $P$ tel que
$\A_P(0) = \A_0$ et une d\'eformation $\sF_P$ de $\F$ de param\`etres $P$
compatible avec $\A_P$ telle
que $\bS= \bS(\sF'_P)$.
\end{enonce}

\begin{proof} Notons (\ref{sytsl}) le syst\`eme $\bS\,$.  Le long de chaque
ensemble critique $K$
les germes d'espaces $W_{P, \, K}$ et $\SS_{P, K}$ sont holomorphiquement
triviaux, cf.
(\cite{Minv} Lemme 1.2.2  page 304).  On supposera donc que les
d\'eformations $\sF'_{P, \, K}$ sont
r\'ealis\'ees sur :
$$W_{P,\, K} = \left( \wM_0 \times P\,, K \times 0 \right)\, ,  \qquad
\SS_{P ,\,K} = \left( \wD_0
\times P  \, , K \times 0
\right)\, ,$$  o\`u $\wM_0$ et $\wD_0$ sont la cime et le diviseur de cime
de $\A_0\,$. Identifions
$\wM_0$ et $\wD_0$ \`a
$\wM_0 \times 0$ et $\wD_0 \times 0 \,$. Quitte \`a restreindre, les traces
$U_K $ des $W_{P, \,K}$
sur $\wD_0$ forment un recouvrement distingu\'e $\,\U$ de $\wD_0$. Notons
$\sF_{P ,Ê\, U_K}$ le
germe de $\sF_{P ,Ê\, K}$ le long de
$U_K\,$.
\\

Deux ouverts $U\,$, $V \in \U$ s'intersectent s'ils correspondent \`a des
\'el\'ements critiques
adjacents. La condition de coh\'erence permet alors de construire \`a l'aide
de (\ref{flowbox}) une
conjugaison formelle de d\'eformations
$\widehat{\Phi}_{UV}$ entre les  germes de $\sF_{P ,Ê\, U}$ et de $\sF_{P
,Ê\, V}$ le long de $U
\cap V\,$, valant l'identit\'e en restriction \`a $\wD_0\times P\,$.
Appliquons le th\'eor\`eme
(\ref{mth.co.}) au cocycle $\widehat{\CC} :=
\left(\widehat{\Phi}_{UV}\right)
\in Z^1\left(\U ;\,\hG\right)\,$.  On obtient un cocycle holomorphe
${\CC} :=
\left({\Phi}_{UV}\right) \in Z^1 \left( \U ; \, \G
\right) \,$ et une cochaine formelle
$\left(\widehat{\Phi}_{U}\right) \in Z^0 \left( \U ; \, \hG \right) \,$
telle que
$\widehat{\Phi}_{UV} = \widehat{\Phi}_U \circ \Phi_{UV} \circ
\widehat{\Phi}_V^{-1} \,$. La
relation de congugaison
${\widehat{\Phi}_{UV}}^{\ast} \sF_{P ,Ê\, U} = \sF_{P ,Ê\, V}$ donne~:
$${{\Phi}_{UV}}^{\ast}{\widehat{\Phi}_{U}}^{\ast} \sF_{P ,Ê\, U}  =
{\widehat{\Phi}_{V}}^{\ast}
\sF_{P ,Ê\, V}\,.$$ Cette \'egalit\'e signifie que lorsqu'on construit comme
en (\ref{recvar}) une
vari\'et\'e holomorphe $\M_\CC\,$  en recollant par les biholomorphismes
$\Phi_{UV}\,$  des voisinages
des
$U \in \U$ dans $\wM_0 \times P\,$, les feuilletages transversalement
formels
${\widehat{\Phi}_{U}}^{\ast} \F_{P ,Ê\, U} \,$ se recollent aussi en un
feuilletage
transversalement formel le long d'un diviseur (isomorphe \`a $\wD_0 \times
P\,$). D'apr\`es le lemme
(\ref{stab.cimes}) $\,\M_\CC\,$ est biholomorphe \`a la cime d'un arbre au
dessus de
$P$ et l'on conclut par image directe cf.(\ref{im.dir.for}) .
\end{proof}

Nous aurons besoin au chapitre g\'en\'ericit\'e d'un r\'esultat plus pr\'ecis
permettant de construire des
d\'eformations avec contr\^ole des jets. Le th\'eor\`eme ci-dessous n'est pas
optimal pour ce contr\^ole
mais suffira \`a notre usage.\\

De mani\`ere g\'en\'erale nous dirons que deux d\'eformations transversalement
formelles d'un m\^eme
feuilletage r\'ealis\'ees sur le m\^eme espace
$$\sF'_P = (\M_P, \SS_P, \pi, \s; \F'_P)\ \ \ \hbox{et}\ \  \ \sF''_P =
(\M_P, \SS_P, \pi, \s;
\F''_P)$$  sont {\it tangentes \`a l'ordre $k\in \N\, $ le long d'un
sous-espace analytique} non
n\'ecessairement r\'eduit de $\M_P\,$ :
$$Z: = \left(\left |Z \right| \, ,\,\- \O_Z := \O_{\M_P}/{\goth
A}\right)\,
\hbox{avec} \left| Z \right| \subset \SS_P,\,$$  si en chaque point de
$\SS_P\,$ les feuilletages $\F'_P$ et $\F"_P$ peuvent \^etre d\'efinis par des
formes diff\'erentielles
dont les coefficients ne diff\`erent que d'un \'el\'ement de l'id\'eal
$$\widehat{\goth A}^{k+1}\, := \, {\goth A}^{k+1} \otimes_{\O_{\M_P}}
\hO_{\M_P}^Z \, \subset \,
\hO_{\M_P}^Z \,.$$
En d'autres termes on a l'\'egalit\'e des {\it
restrictions}  des feuilletages
$\F'_P$ et $\F''_P$ au {\it k-i\`eme voisinage  infinit\'esimal} de $Z$ dans
$\M_P\,$, c'est \`a dire l'\'egalit\'e des faisceaux~:
\begin{equation}\label{jet.feuil}
\L_{\F'_P} \otimes_{{\hO^Z}_{\M_P}} \left(\left.{\hO^Z}_{\M_P} \right/
\goth A^{k+1} \right) =
\L_{\F"_P} \otimes_{\hO^Z_{\M_P}} \left(\left.{\hO^Z}_{\M_P} \right/ \goth
A^{k+1} \right)
\,.
\end{equation}

\vspace{1 em}

Fixons de nouveau un arbre $\A_P$ de hauteur $h$ au dessus de $P$, not\'e
(\ref{arbnot}). D\'esignons
par
\begin{equation}\label{divanal}
\Bbb{D} := \left( \wD \,; \; \left. \O_{\wM} \right/  \wJ \right)\textsf{}
\end{equation} le sous-espace analytique de la cime $\wM$  de $\A_P$
d\'efini par le faisceau
d'id\'eaux  $\wJ$ image  r\'eciproque sur $\wM$ du faisceau d'id\'eaux $I_P
\subset \O_{\Bbb{C}^2\times
P}$ des fonctions nulles sur
$0\times P\,$ cf. (\ref{id.arbres}).

\begin{enonce}{Th\'eor\`eme de r\'ealisation tangente}\label{de.real.tg.} Soient $\F$ un
feuilletage formel \`a
l'origine de $\Bbb{C}^2,\,$\,\,$\A_P$ un germe d'arbre au dessus de $P$
tel que $\F$ est
compatible avec $\A_P(0)$ au dessus de $\{0\}$. Notons $\F'$ la
transform\'ee stricte de $\F$ sur la
cime $\wM_0$ de $\A_P(0)$. On se donne une d\'eformation  transversalement
formelle
$\sF_P$ de
$\F$ de param\`etres $P$ compatible avec
$\A_P\,$, et un syst\`eme semi-local
$\bS := \left(\sF"_{P,\, K}\right)_{K \in \bC\,(\A_P(0))}$ de d\'eformations
formelles de $\F'$
r\'ealis\'ees sur la cime $\wM\,$ de $\A_P$. Supposons qu'au voisinage de
chaque ensemble critique $K$
le transform\'e strict
$\sF'_P$ de
$\sF_P$ sur $\wM$ est tangent (\ref{sytsl}) \`a l'ordre $k + 2h\,$ \`a
$\sF"_{P, K}\,$, le long de
$\Bbb{D}\,$, $k
\in
\N\,$. Alors il existe un germe de d\'eformation
$\underline{\cal H}_P$ de
$\F$ transversalement formelle le long de $0 \times P\,$, $k$-tangente \`a
$\sF_P$ le long de $0
\times P$ (consid\'er\'e comme espace r\'eduit), compatible avec $\A_P$ et dont
le transform\'e strict
$\underline{\cal H}'_P$ sur $\wM$ induit un syst\`eme semi-local
$\bS(\underline{\cal H}'_P)$ formellement conjugu\'e \`a $\bS\,$.
\end{enonce}

\begin{proof}  On proc\`ede de la m\^eme mani\`ere que pour la d\'emonstration du
th\'eor\`eme  de r\'ealisation. Les
hypoth\`eses de tangence donnent maintenant, toujours gr\^ace au Th\`eor\`eme
\ref{mth.co.}, 3. des
cocycles v\'erifiant~:
$$\left(\Phi_{UV}\right) \in Z^1(\U; \G_{k+2h})\, ,\quad
\left( \widehat{\Phi}_UÊ\right) \in Z^0(\U; \hG_{k+2h}) ,\quad
{{\Phi}_{UV}}^{\ast}
{\widehat{\Phi}_{U}}^{\ast} \F"_{P,\, U}  =
{\widehat{\Phi}_{V}}^{\ast}\F"_{P, V}\,$$ o\`u $\U$
d\'esigne encore le recouvrement distingu\'e de $\wD_0$ form\'e des
intersections
$U_K := W_{K, \,P} \cap \wD_0$, et
$\F"_{P,\, U_K}$ d\'esigne le germe de $\F"_{P,\, K}$ le long de $U_K\,$.
D'apr\`es le th\'eor\`eme de
stabilit\'e (\ref{d.q.fin.}) le cocycle $\left( \Phi_{UV }\right)$ vaut
l'identit\'e dans
$H^1(\U; \G_{k+h})\,$. Ainsi
$$\Phi_{UV} = {\Phi'}_U \circ {\Phi'}_V^{-1} \, , \qquad \left(\Phi'_U
\right) \in Z^0 \left( \U ;
\, \G_{k+h}Ê\right) \,, $$ et on obtient
\footnote{On aurait aussi pu appliquer l'isomorphisme 2. du th\'eor\`eme
\ref{mth.co.} pour obtenir
directement les conjugaisons
$\Psi_U\,$.}
$$\widehat{\Psi}_U^{\ast}Ê \F"_{P,\, U}Ê = \widehat{\Psi}_V^{\ast}Ê
\F"_{P,\, V}Ê \, , \quad U\cap
V
\not=
\emptyset\, ,Ê\quad
\widehat{\Psi}_U := \widehat{\Phi}_U \circ {\Phi'}_U    \in \hG_{k+h}
\left( U \right)\, .$$ Ces
relations d\'efinisent un germe le long de $\wD_0$ de feuilletage global
${\cal H'}_P\,$,
transversalement formel le long de
$\wD$.  Consid\'erons la d\'eformation ${\cal H}_P$ de $\F\,$,
transversalement formelle le long de $0
\times P\,$, obtenue par image directe (\ref{im.dir}) de
${\cal H'}_P$  par l'application $\wE : \wM \fle \Bbb{C}^2 \times P$
compos\'ee des \'eclatements de
$\A_P\,$. La  transform\'ee stricte de ${\cal H}_P\,$ sur $\wM$ est ${\cal
H'}_P\,$. Par
construction on a
\begin{equation}\label{tan}
\hL_{\cal H'_P} = \hL_{\F'_P} \quad \mbox{modulo}\quad  \wJ^{k+h+1}
\hL_{\wM}\,, \qquad
\hL_{\wM} := \L_{\wM}\otimes \hO_{\wM}^{\Bbb{D}}\,.
\end{equation}  Soit $\n$ et $\xi$ les 1-formes diff\'erentielles formelles
\`a l'origine de
$\Bbb{C}^2 \times P$  d\'efinissant  $\F_P$ et ${\cal H}_P$ respectivement.
Il reste \`a montrer que
les coefficients de $\n - \xi$ appartiennent \`a l'ideal
$\widehat{I}_P^{k+1}\,$ o\`u
$\widehat{I}_P:= (x, y)$ d\'esigne l'id\'eal de
${\hO_{\Bbb{C}^2 \times P,0}}^{0 \times P}$ engendr\'e par les coordonn\'ees
$x$ et $y$.\\

Les images r\'eciproques de $\n$ et $\xi$ sur $\wM$ d\'efinissent des
faisceaux de sous-modules de
$\hL_{\wM}$ qui se d\'ecomposent en
$$\left( \wE^{\ast}(\n) \right) = \goth I_\n \cdot \L_{\F'_P} \, , \quad
\left( \wE^{\ast}(\xi) \right) = \goth I_\xi \cdot \L_{\cal H'_P} \,
,\qquad
\goth I_\n \,, \; \;\goth I_\xi \subset \hO_{\wM}^{\Bbb{D}}\,.
$$  Il d\'ecoule de (\ref{tan}) que $\wE^{\ast}(\n) - \wE^{\ast}(\xi)$ est
une section globale de
$\wJ^{k+h+1} \hL_{\wM}\,$. On voit facilement par r\'ecurrence sur $h$ que
les champs de vecteurs
$x^h \,\dd x$ et $x^h \,\dd y$ sur
$\Bbb{C}^2 \times P$ se rel\`event  sur $\wM\,$ en des champs de vecteurs
holomorphes globaux $X :=
\wE^{\ast}(x^h \,\dd x)$ et $Y :=\wE^{\ast}(x^h\,\dd y)\,$. Ainsi en
posant $\widetilde{x} = x
\circ \widetilde{E}$ et $\widetilde{y} = y
\circ \widetilde{E}$, on voit que
$$ f := \widetilde{x}^{-h}\,\left(\wE^{\ast}(\n) \cdot X  -
\wE^{\ast}(\xi) \cdot X \right)
 \quad \hbox{et} \quad g :=  \widetilde{x}^{-h}\,\left(\wE^{\ast}(\n)
\cdot Y  - \wE^{\ast}(\xi)
\cdot Y
\right)\,$$ sont des sections globales de $\wJ^{k+1},\,$ et que ce
faisceau est engendr\'e par les
sections globales
$\widetilde{x}^\a \widetilde{y}^\b \,$,
$\a + \b = k + 1\,$. Comme $H^1 \left( \wD_0 ; \, \hO_{\wM}^{\Bbb{D}}
\right) = 0$ d'apr\`es
(\ref{cohom.ecl.}) on peut d\'ecomposer
$$f = \sum_{\a + \b = k + 1} f_{\a\b} \,\widetilde{x}^\a
\widetilde{y}^\b\,, \quad  g = \sum_{\a +
\b = k + 1} g_{\a\b} \,\widetilde{x}^\a \widetilde{y}^\b\,, \quad
\mbox{avec} \quad  f_{\a\b} , \, g_{\a\b} \in  H^0 \left( \wD_0 ; \,
\hO_{\wM}^{\Bbb{D}}
\right)\,.$$ Les coefficients $f_{\a\b}$ et $g_{\a\b}$ se redescendent
(\ref{im.dir.for}) par
$\wE$ en des \'el\'ement $f_{\a\b}^{\flat}(x, y ; \,t)$ et
$g_{\a\b}^{\flat}(x, y ; \,t)$ de
${\hO_{\Bbb{C}^2\times P,\, 0}}^{0\times P}$.  Ainsi
$$\n \cdot \dd x - \xi \cdot \dd x = \sum_{\a + \b = k + 1}
f_{\a\b}^{\flat}(x, y ; \,t) \, x^\a
y^\b $$ et $$ \n \cdot \dd y - \xi \cdot \dd y = \sum_{\a + \b = k + 1}
g_{\a\b}^{\flat}(x, y ;
\,t) \,  x^\a y^\b $$ appartiennent \`a
$\widehat{I}_P^{k+1}\, ,$ ce qui ach\`eve la d\'emonstration.
\end{proof}

\section{Equir\'eduction}\label{sec.def.equising}
\subsection{Singularit\'es de deuxi\`eme esp\`ece}\label{subs.red.des.sing.}
\addcontentsline{toc}{section}{\hspace{0,8em} {}\thesubsection .  Singularit\'es de deuxi\`eme esp\`ece}
Nous rappelons et pr\'ecisons ici quelques notions classiques sur la
r\'eduction des singularit\'es de
feuilletages. Pour plus de d\'etails nous renvoyons le lecteur \`a \cite{C-M}
\cite{M-M} ou encore \`a
\cite{martbour}. Donnons-nous une vari\'et\'e
holomorphe $M$ de dimension deux,  une courbe \`a croisements normaux $D
\subset M\,$ et un feuilletage
formel $\F$ en un point $m_0 \in D\,$.

\begin{defi} Nous disons que  {\it $\F$ est r\'eduit en $m_0$} si $\F$ est~:
soit r\'egulier, soit d\'efini
par une 1-forme diff\'erentielle formelle $\w$ qui, dans des coordonn\'ees
locales appropri\'ees
$(x , \,y)\,$, $x\,(m_0) = y\,(m_0) = 0\,$, s'\'ecrit
\begin{equation}
\w = ( \l \,x + \cdots)\,dy  + (\mu \,y + \cdots )\,dx \, , \qquad \l \, ,
\mu \, \in \C \, , \quad
\l \not= 0\,,
\quad \frac{\mu}{\l} \notin \Q_{< 0}\, ,
\end{equation} o\`u $+ \cdots$ d\'esigne des termes d'ordre sup\'erieur i.e. des
\'el\'ements de l'id\'eal $(x,
\, y )^2 \subset \hO_{M, \,m_0}\,$.\\
\indent Nous disons que {\it le couple $\left( \F ,\, D \right)$ est
r\'eduit en $m_0$ } si
$\left( \F , \, D \right)$ est r\'egulier en $m_0\,$ cf.
(\ref{singcouples}), ou bien si $\F$ est
r\'eduit et chaque germe $D_{m_0}'$ de composante irr\'eductible locale de $D$
en $m_0$ est une courbe
invariante de $\F\,$ i.e. la restriction de $\w$ \`a $D_{m_0}'$  est
identiquement nulle.
\end{defi}

La classification formelle des feuilletages singuliers r\'eduits est bien
connue~: il existe des
coordonn\'ees formelles
$z_1\,$, $z_2 \in \hO_{M ,\,m_0}\,$, $z_1 (m_0) = z_2 (m_0) = 0\,$  et $u
\in \hO_{M ,\,m_0}\,$,
$u(m_0) \not= 0$  tels que $\w = u \, \ww\,$  et $\ww$ est l'une des
1-formes suivantes,  dites {\it
formes normales}~:
\begin{enumerate}
\item $\ww := \l_1z_1\,dz_2 +\l_2z_2\,dz_1\,$ avec $\l_1 \l_2 \not=
0\,$,$\l_1 + \l_2 = 1\,$,
$\l_2/\l_1 \in
\C - \Q\,$,
\item $\ww := q\,z_1\,dz_2 + p\,z_2\,dz_1\,$ avec $p, q  \in \N^{\ast}\,$,
$(p,q)=1\,$,
\item $\ww :=
q\,z_1(1+\z\,(z_1^pz_2^q)^k)\,dz_2+p\,z_2(1+(\z-1)\,(z_1^pz_2^q)^k)\,dz_1$
avec
$p, q, k  \in \N^{\ast}\,$, $(p,q)=1\,$, $\z \in \C\,$,
\item $\ww := (\z z_2^p - p)z_1\,dz_2 +z_2^{p+1}\,dz_1 \,$, avec $p \in
\N^{\ast},\,\,\,$$\z \in
\C\,$.
\end{enumerate} On dit dans le premier cas que $\ww$ est {\it
non-r\'esonnant}. Dans le cas 2. on dit
que $\ww$ est {\it r\'esonnant lin\'earisable} et {\it r\'esonnant
non-lin\'earisable} dans le cas 3.. Enfin
dans le dernier cas on dit que $\ww$ est un {\it selle-n\oe ud}. De plus
ces \'ecritures sont uniques~:
deux feuilletages d\'efinis par deux 1-formes distinctes du type 1. \`a 4. ne
sont pas formellement
conjugu\'es. \\

Remarquons que dans tous les cas les courbes formelles
$z_j = 0\,$,
$j = 1, 2$ sont des courbes formelles invariantes (\ref{sec.def.feuil}) de
$\F\,$. Ce sont en fait les seules. Si $\F$ est un selle-n\oe ud le germe
de courbe $\{ z_2 = 0 \}$
est appel\'ee {\it vari\'et\'e forte} et le germe
$\{ z_1 \} = 0$ est appel\'e {\it vari\'et\'e faible} de $\F \,$.

\begin{defi}\label{sntagent}
Supposons que $\F$ est un selle-n\oe ud et que $\left( \F , \, D
\right)$ est r\'eduit.  Nous disons que $\F$ est un {\it selle-n\oe ud transverse
\`a $D$ en $m_0$} si $D$ est
lisse au point $m_0$ et si son germe est une vari\'et\'e forte de $\F\,$. Dans
le cas contraire nous
dirons que $\F$ est {\it un selle-n\oe ud tangent \`a $D$ en $m_0\,$}. Lorsque de
plus  l'holonomie locale de
la vari\'et\'e faible de $\F$ est p\'eriodique nous dirons que {\it $\F$ est un
selle-n\oe ud tangent
r\'esonnant}.
\end{defi}

Nous avons suppos\'e $\F$ formel au point $m_0\,$. Lorsque $\F$ est
transversalement formel le long de
$D\,$ il est naturel de chercher \`a obtenir la forme normale de $\F$ dans
des coordonn\'ees
transversalement formelles le long de $D\,$.

\begin{prop}\label{forntf}\cite{Ra-Ma.s.n.}, \cite{Ra-Ma}
 Soit $\F$ un feuilletage transversalement formel le long de $D$ admettant
$m_0$ comme singularit\'e
isol\'ee. Si $(\F , \,D)$ est r\'eduit en
$m_0$ et n'est pas en ce point un selle-n\oe ud tangent \`a $D\,$, il existe
alors des coordonn\'ees
$z_1\,$, $z_2 \in \hO^D_{M ,\,m_0}\,$, $z_1 (m_0) = z_2 (m_0) = 0\,$
transversalement formelles le
long de $D$ et $u \in \hO^D_{M ,\,m_0}\,$, $u(m_0) \not= 0$  telles que
$\w = u \, \ww\,$  et $\ww$
est l'une des 1-formes 1. \`a 4. ci-dessus.
\end{prop}

Lorsqu'un
feuilletage formel $\F$ \`a l'origine de $\C^2$ d\'efini par une 1-forme {\it
\`a singularit\'e isol\'ee},
$Sing (\w ) = \{0 \}\,$
\footnote{i.e. $ dim_\C \left(\left. \O_{\C^2 , \, 0} \right/ (a, b )
\right) \, < \infty$ }
$$\w := \stw $$
 n'est pas r\'eduit on construit un arbre au dessus de
$\left\{0\right\}\,$,  dans le sens du chapitre 1, mais de hauteur \`a
priori infinie, not\'e
\begin{equation}\label{ard}
\left(\,M_\w^j, E_\w^j, \,\S_\w^j, \, C_\w^j,\,  \pi_{\w ,Ê\,j},
\,\D_\w^j\,\right)_j\, ,\quad
\hbox{ou} \quad
\left(\,M_\F^j, E_\F^j, \,\S_\F^j, \, C_\F^j,\,  \pi_{\F ,Ê\,j},
\,\D_\F^j\,\right)_j
\end{equation}  en d\'efinissant par induction la succession d'\'eclatements
de la mani\`ere suivante~:
\begin{enumerate}
\item $M_\w^0 := \C^2\,,\qquad \S_\w^0 :=\left \{0\right\} =: {C}_\w^0\,$,
\item $\S_\w^j := Sing(\F^j,  \D_\w^j)\,$ o\`u $\F^j\,$ d\'esigne le
transform\'e strict de $\F$ par
l'application $E_{\w , \,j} := E_\w^1 \circ \cdots \circ E_\w^j\,$
compos\'ee des \'eclatements de centres
${C}_\w^k,\,$
$k = 0,\ldots,j-1\,$  et $\D_\w^j := E_{\w , \,j}^{-1}(0)\,$,
\item $C_\w^j \subset \S_\w^j$ est l'ensemble des points de $\D_\w^j\,$ o\`u
le couple $\,\left(
\F^j,\, \D_\w^j \right)\,$ n'est pas r\'eduit.
\end{enumerate}

Un th\'eor\`eme classique de Bendixon-Seidenberg \cite{bendixon} \cite{Seid} \cite{M-M}  affirme que l'arbre
de r\'e\-duc\-tion d'un
feuilletage formel \`a l'origine de
$\C^2$  est de hauteur finie : avec les notations ci-dessus, il existe un
entier not\'e $h_\w$ ou
$h_\F$ pour lequel $C_\w^{h_\w} = \emptyset\,$.

Dans toute la suite du texte nous d\'esignerons indiff\'eremment par
$\A[\F]\,$ ou $\A[\w]\,$ l'arbre
(\ref{ard}) que nous appelons {\it arbre de r\'eduction de $\F$ ou de $\w$}.
Nous noterons aussi
\begin{equation}
\;\; E_{\F} :=  E_\w := E_{\w , \,h_\w}\, , \quad   M_\w := M_\w^{h_\w}\,,
\quad \D_\F  := \D_\w :=
E_\w ^{-1}(0)\, ,
\quad \wF := E_\w^{\ast}(\F)\,.
\end{equation}

\begin{defi}\label{especes} Nous dirons que~:\\
\indent  $\F$ est {\it non-dicritique} si $\D_\F$ est un ensemble invariant de
$\wF\,$,\\
%%%\indent  $\F$ est {\it de premi\`ere esp\`ece} s'il est
%%%non-dicritique et si aucun
%%%point singulier de $\wF$
%%%n'est de type selle-n\oe ud tangent r\'esonnant,\\
\indent  $\F$ est {\it de deuxi\`eme esp\`ece} s'il est non-dicritique et si aucun
point singulier de $\wF$
n'est de type selle-n\oe ud tangent.\\
\indent $\F$ est {\it semi-hyperbolique}
\footnote{Ces feuilletages sont aussi appel\'es {\it courbes g\'en\'eralis\'ees} dans \cite{C-L-S}}
s'il est non-dicritique et si aucun
point singulier de $\wF$ n'est
de type selle-n\oe ud.
\end{defi}
%%%%%%%%%%%%%%%%%%%%%%%%%%%%%%%%
%%j'ai supprim\'e premi\`ere esp\`ece%%
%%%%%%%%%%%%%%%%%%%%%%%%%%%%%%%%
L'arbre dual fl\'ech\'e et pond\'er\'e $\A^{\ast}[\F]$ de $\A[\F]$, cf.
(\ref{adu})  poss\`ede une
seconde pond\'eration, la {\it pond\'e\-ra\-tion par multiplicit\'es} que nous
allons d\'efinir maintenant. \\

De mani\`ere g\'en\'erale, donnons-nous une application holomorphe $E$ d'une
surface lisse $M$ sur $\C^2$
telle que $\D := E^{-1} (0)$ soit une hypersurface et une 1-forme formelle
$\n \in \hL_{\C^2 , \,0}$
non n\'ecessairement \`a singularit\'e isol\'ee. Le sous-faisceau de modules
$\left( E^{\ast}(\n)
\right)$ de
$\hL^\D_M$ engendr\'e par
$E^{\ast}(\n)$ se d\'ecompose de mani\`ere unique en
\begin{equation}\label{pondmul}
\left( E^{\ast}(\n) \right) = \widehat{J} \cdot \prod_{D \in comp(\D)}
\widehat{I}_D^{\,m_D(\n)}\,
\cdot \,
\hL_{\wF_\n}\, , \qquad m_D\,(\n) \in \N\,,
\end{equation}  o\`u $comp(\D)$ d\'esigne l'ensemble des composantes
irr\'eductibles de $\D\,$,
$\widehat{I}_D$ le faisceau des fonctions transversalement formelles le
long de $\D$ nulles sur
$D\,$, $\F_\n\,$ le feuilletage d\'efini par $\n\,$,
$E^\ast \F_\n$ son transform\'e strict (\ref{trans.strict}) et
$\widehat{J}$ un faisceau d'id\'eaux dont les z\'eros ne contiennent aucune
composante irr\'eductible de
$\D\,$.  Lorsque $\n$ est \`a singularit\'e isol\'ee en 0, on a $\widehat{J} =
(1)\,$.
\begin{defi}\label{multfdiv} L'entier $m_D (\n )\,$ d\'efini par
(\ref{pondmul}) s'appelle la
{\it multiplicit\'e de $\n$ suivant la composante irr\'eductible $D\,$}. Il se note
aussi $m_D \left( \F_\n
\right)\,$.
\end{defi}
\noindent Ce  nombre s'interpr\`ete comme le plus grand entier $k$ tel qu'en
chaque point la $k$-i\`eme
puissance d'une \'equation r\'eduite locale de $D$ divise
$ E^{\ast}\,(\n)\,$.\\

Enrichissons la donn\'ee de $\A^{\ast}[\F]$ en associant \`a chaque sommet $s$
correspondant \`a une composante
irr\'eductible
$D$ de $\D_\w\,$, d'une part l'auto-intersection $e ( s ) := \langle D, D \rangle$ et d'autre part  {\it
la multiplicit\'e} $m({ s}) :=
m_D(\F) \in \N \,$. On obtient un arbre doublement pond\'er\'e. Nous y
rajoutons une donn\'ee
suppl\'ementaire en munissant du symbole "fl\`eche double" chaque sommet dont
le diviseur correspondant,
appel\'e {\it composante dicritique}, n'est pas invariant par $\wF\,$.
L'arbre obtenu est not\'e
$\cA[\F]\,$.

\begin{defi}\label{adudoub}
L'arbre $\cA[\F]\,$ s'appelle {\it arbre dual doublement pond\'er\'e associ\'e \`a
$\F\,$}.
\end{defi}

Dans \cite{C-L-S} on trouvera des formules liant ces multiplicit\'es avec le
{\it nombre de Milnor de} $\F$
d\'efini par
\begin{equation}\label{nbrMilnor}
\mu_0\, (\F) := dim_\C \left( \left. \hO_{\C^2 , \,0} \right/ (a , \, b)
\right)\,,
\end{equation}
et avec la {\it multiplicit\'e alg\'ebrique du feuilletage \`a
l'origine} c'est \`a dire l'ordre du
premier terme du d\'eveloppement en s\'erie de
$\w$ que l'on note
$\nu_0(\w)$ ou $\nu_0(\F)\,$. Les auteurs d\'efinissent aussi le {\it nombre
de Milnor
$\mu_0\,(\F , Z)$ de $\F$ suivant une courbe invariante lisse $Z$} comme
la multiplicit\'e \`a l'origine
de la restriction \`a $Z$ d'un champ de vecteurs formel \`a singularit\'e isol\'ee
d\'efinissant $\F\,$.  Ils prouvent la formule
suivante, lorsque $\D_\w$ est sans
composante dicritique~:

\begin{equation}\label{formcsn}
\nu_0(\F) + 1 = \sum_{D \in comp(\D_\w)}\, \sum_{m \in Sing(\wF \cap D)}
m\,(D)
\left(\mu_m(\wF,\, D) - \e(m) + 1 \right)
\end{equation}
\noindent o\`u $\e(m)$ est le nombre 1 ou 2 de branches de $\D_\w$ au point
$m$ et $m\,(D)$ est d\'efini
comme en (\ref{arbr.inf}) par l'\'egalit\'e de faisceaux d'id\'eaux
$$ {E_\w}^{\ast}\, \goth M  = \prod_{D \in comp(\D_\w)} \widehat{I}_D^{\;
m\,(D)}\, ,
\qquad \goth M := (x, \,y ) \subset \hO_{\C^2 , \,0}\, .$$
\begin{rema}
$\F$ est un selle-n\oe ud tangent \`a $Z$ si et seulement si $\mu_0 (\F , Z
) > 1\,$.
\end{rema}

La formule (\ref{formcsn}) permet de voir par induction que, pour $\F$
non-dicritique,   les
multiplicit\'es $m_D\,(\w)$ s'expriment comme des fonctions affines des
nombres de Milnor $\mu_m\,(\wF,
\, D)\,$, $D \in comp (\D_\w)\,$, $ m \in Sing (\wF )\,$, dont les
coefficients sont des entiers
$\geq 0$ qui ne d\'ependent que de
$\A^{\ast}\,[\w]\,$. Le pas de l'induction repose sur les \'egalit\'es~:
\begin{equation}\label{formmult}
\;Ê\; m_{D_c}\,(\w) = \nu_c\,(\wF^{j}) + \sum_{D \in Ad (c)} m_{D}\,(\w)
\, , \quad  Ad (c) :=
\left\{ \left. D \in comp(\D_\w^{j}) \; \right/ \; c \in D \right\}\, ,
\end{equation} o\`u $c \in C_\w^{j}$ et $D_c\,$ est la composante de $\D_\w$
cr\'e\'ee par l'\'eclatement de
$c\,$. Les feuilletages de deuxi\`eme esp\`ece sont caract\'eris\'es par les
\'egalit\'es
$\mu_m\,(\wF , \, D) = 1\,$. On en d\'eduit la
\begin{rema}\label{minmult}
Si $\F$ est de deuxi\`eme esp\`ece, alors
\begin{enumerate}
\item les
 multiplicit\'es $m_D\,(\w)$ sont donn\'ees par $\A^{\ast}[\w]\,$,
\item dans l'ensemble
$\left\{ w'\in \hL_{\C^2, \,0} \; \left/ \; Sing\,(w') = \{0\} \;
\mbox{et} \;\; \A^{\ast}[w'] \equiv
\A^{\ast}[w] \right. \right\}\,$ la pond\'e\-ra\-tion par multiplicit\'es de
$\cA[w]\,$ r\'ealise pour
chaque sommet un minimum.
\end{enumerate}
\end{rema}
\indent Sous l'hypoth\`ese de non-dicriticit\'e, $\F$ ne poss\`ede qu'un nombre
fini de courbes formelles
invariantes \`a l'origine appel\'ees {\it s\'eparatrices formelles de $\F$} par
certains auteurs. Soit
$\hff \in \hO_{\C^2 , \,0}$ une \'equation r\'eduite de l'union de ces courbes
que l'on notera
$\hsep\,(\F)\,$ ou $\hsep\,(\w)\,$.
Dans \cite{C-L-S} les auteurs prouvent~:
\begin{equation}\label{inegmudi}
\begin{array}{l}
\nu_0\,(\w) \geq \nu_0 \,(d\hff)\, ,\qquad \mu_0 \,(\w ) \geq
\mu_0\,(d\hff)\, , \\
\\
m_D\,(\w) \geq m_D\,(d\hff )\,, \quad D \in comp ( \D_\w ) \, ,
\end{array}
\end{equation}
ainsi que l'\'equivalence
\begin{equation}
 \F \;\mbox{  est semi-hyperbolique}\; \ssi \; \mu_0 \,(\w ) =
\mu_0\,(d\hff) \, .
\end{equation}
\indent Nous allons donner plusieurs caract\'erisations des singularit\'es de
deuxi\`eme esp\`ece. La plus
int\'eressante fera intervenir les faisceaux de base $\D_\w$ suivants~: le
faisceau
$\hX_\wF\,$  form\'e des germes aux points de
$\D_\w$ des champs de vecteurs $X$ transversalement formels le long de
$\D_\w$ qui sont tangents \`a
$\wF\,$ i.e.
$\L_\wF \cdot X = 0\,$ et le faisceau
$\hX_{\widehat{S'}}$ form\'e des  champs $X$ transversalement formels
tangents \`a
$\widehat{S}' := E_\w^{-1}\left(
\hsep(\F)Ê\right)\,$, i.e. $X \cdot I_{\widehat{S'}} \subset
I_{\widehat{S'}}\,$, o\`u
$I_{\widehat{S'}}$ est le radical du faisceau d'id\'eaux $\left( \hffÊ\circ
E_\w \right)\,$. D\'esignons
par
$\a$ le morphisme $\hX_{\widehat{S'}} \fle \hO_{\M_\w}^{\, \D_\w}\,$ qui \`a
un champ $X$ associe
${E_\w}^{\ast}(\w) \cdot X\,$.
Le th\'eor\`eme suivant, annonc\'e dans \cite{Mdouady}, donne une
interpr\'etation g\'eom\'etrique de la propri\'et\'e prouv\'ee dans le lemme cl\'e (3.2) de \cite{Mdouady}.

\begin{theo}\label{cardeuxesp} Soit $\F$ un feuilletage formel
non-dicritique \`a l'origine de $\C^2$
d\'e\-fini par une 1-forme formelle $\w$ \`a singularit\'e isol\'ee et $\hff$ une
\'equation r\'eduite de
$\hsep(\F)\,$. Les propri\'et\'es suivantes sont \'equivalentes~:
\begin{enumerate}
\item $\F$ est de deuxi\`eme esp\`ece,
\item $\nu_0\,(\w) = \nu_0 \,(d\hff)\,$,
\item $\cA[\w]$ = $\cA[d\hff]\,$,
\item pour chaque $D \in comp(\D_\w)\,$ on a l'\'egalit\'e $m_D\,(\w) =
m_D\,(d\hff )\,$,
\item il existe un $D \in comp(\D_\w)\,$ pour lequel $m_D\,(\w) =
m_D\,(d\hff )\,$,
\item on a la suite exacte~:
\begin{equation}\label{exgeom}
0 \fle \hX_\wF \fle \hX_{\widehat{S'}} \stackrel{\a}{\fle} \left(
\hffÊ\circ E_\w \right) \fle 0\, .
\end{equation}
\end{enumerate}
\end{theo}

\begin{proof} Dans \cite{C-L-S} la d\'emonstration de la formule (\ref{formcsn})
reste valable pour une
succession quelconque $E$ d'\'eclatements. Appliquons-la \`a $\w\,$ et \`a
$d\hff\,$ pour $E = E_\w\,$.
Comme le lieu singulier de $E^\ast (\F_{d\hff})\,$ est contenu dans celui
de $\wF\,$, on obtient
l'in\'egalit\'e
$\nu_0\,(\w )
\geq
\nu_0\,(d\hff )\,$. Aux points singuliers $m$ de $\wF\,$, $\mu_m\, (E^\ast
(\F_{d\hff}) , \, D) \,$
est \'egal \`a 1, ou \`a 0 si
$E^\ast (\F_{d\hff})$ est r\'eguli\`ere en $m\,$. On en d\'eduit imm\'ediatement
que l'\'egalit\'e
$\nu_0\,(\w ) =
\nu_0\,(d\hff )\,$ est r\'ealis\'ee si et seulement si
$Sing \left( E^\ast (\F_{d\hff}) \right) = Sing (\wF ) $ et
$\mu_m\, (\wF , \, D) = 1\,$ pour chaque point singulier
$\,m\,$ et chaque composante irr\'eductible $D\,$ de
$\,\D_\w\,$; c'est \`a dire si et seulement si $\,\w\,$ est de deuxi\`eme
esp\`ece, ce qui prouve $1. \ssi
2.\,$.
\\

Trivialement on a $\, 3. \Longrightarrow 4.\,$ et $\, 4. \Longrightarrow
5.\,$.  Pour l'implication
$\, 5. \Longrightarrow 2.\,$ raisonnons par contrapos\'ee et supposons que
$\nu_0\,(\w) > \nu_0\,(d\hff
)\,$. La composante $D_0$ de $\D_\w$ cr\'e\'ee par l'\'eclatement de l'origine
satisfait $m_{D_0}\,(\w) >
m_{D_0}\,(d\hff)\,$.  A l'aide de (\ref{inegmudi}) appliqu\'ee en chaque
point singulier des
transform\'es stricts de $\w$ et de $d \hff$ sur $D_0$, on en d\'eduit que
$m_D\,(\w) > m_D\,(d\hff)$ pour tout
$D
\in comp\,(\D_\w)\,$.\\

Pour montrer $\,1. \Longrightarrow 3.\,$ raisonnons par r\'ecurrence sur la
hauteur de
$\A[\w]\,$. Il suffit d'apr\`es la propri\'et\'e ($\star$) ci-dessus de montrer
l'\'egalit\'e
$\A^{\ast}[w] = \A^{\ast}[d\hff]\,$.  Lorsque $h_\w = 1$ les s\'eparatrices
formelles de $\w$ sont des
courbes lisses deux \`a deux transverses et il suffit de voir que leur
nombre $p$ est $ \geq 3\,$,
c'est \`a dire que $\hff = 0$ n'est pas d\'ej\`a une courbe \`a croisement normal.
Sinon $\nu_0\,(d\hff) = 1$
et d'apr\`es
$\, 1.
\Longrightarrow 2.\,$ on a aussi $\nu_0\,(\w) = 1\,$. Dans ce cas les
seuls feuilletages qui se
r\'eduisent par un seul \'eclatement sont du type
$\, \w = (y + \cdots )\,dx + (-x + \e y + \cdots ) \,dy\,$ avec $\e = 0$
ou $ 1\,$. Ces deux
possibilit\'es sont exclues car $\w$ est dicritique si $\e = 0$ et on
obtient un selle-n\oe ud tangent
pour $\wF$ lorsque $\e = 1Ê\,$. La d\'emonstration du pas de r\'ecurrence est
triviale.
\\

Calculons maintenant les fibres du co-noyau de $\a\,$ aux points $\,c\,$
de $\D_\w\,$. Nous
distinguerons plusieurs cas. Supposons d'abord que $\wF$ est r\'egulier et
$\D_\w$ est lisse en $c\,$.
Fixons des coordonn\'ees transversalement formelles locales $z_1, \, z_2\,$
dans lesquelles $\hff \circ
E_\w =  z_2^{m_D\,(d\hff ) + 1}\,$. La composante $D$ de $\D_\w$ portant
$c$ est $\left\{ z_2 = 0
\right\}\,$. Au point $c$ la fibre de $\hX_{\widehat{S}'}$ est engendr\'ee
par les champs $\dd z_1\,$
et
$z_2\, \dd z_2Ê\,$ et le germe de $\wF\,$ est d\'efini par une 1-forme
transversalement formelle qui
s'\'ecrit
$\ww = z_2^{m_D\,(\w)}\,\left( A(z_1, z_2)\, z_2\,dz_1 + B(z_1, z_2) \,
dz_2 \right)\,$ avec $B(0,
\,0) \not= 0\,$. Ainsi
\begin{equation}\label{relmuex}
\left( \hff \circ E_\w \right)_c = \left( z_2^{m_D\,(d\hff ) + 1} \right)
\quad \mbox{et} \quad
coker(\a)_{\, c} = \left( z_2^{m_D\,(\w ) + 1} \right)\,.
\end{equation}

On en d\'eduit l'implication $\,6. \Longrightarrow 4.\,$.\\

Pour l'implication r\'eciproque il reste \`a prouver l'exactitude de la suite
(\ref{exgeom}). aux points
singuliers de
$\wF\,$ car (\ref{relmuex}) donne l'exactitude aux points r\'eguliers.
Consid\'erons d'abord un point
singulier $\,c\,$ sur une composante
$D$ de
$\D_\w$ o\`u
$\D_\w$ est lisse.  On sait que
$\mu_c\,(\wF,Ê\, D) = 1\,$ car
$\,4. \Longrightarrow 1.\,$. Ainsi dans de bonnes coordonn\'ees locales en
$c$ on peut \'ecrire
$$\ww = z_2^q\,\left( \, A(z_1, z_2) z_2\,dz_1 + ( \l + B(z_1, z_2) )
z_1\, dz_2 \,\right)\,, \quad
B(0, \,0) = 0\,, \quad \l \not= 0\, ,$$
$$\hbox{et} \quad \hff \circ E_\w = z_1 z_2^{q + 1}\, ,Ê\quad D = \left\{
z_2 = 0 \right\}\, , \quad
q := m_D\,(\w) = m_D\,(d\hff)\,.$$
\noindent Maintenant $\hX_{\widehat{S}'}$ est engendr\'ee par $z_1\,\dd
z_1\,$ et
$z_2\, \dd z_2Ê\,$ ce qui donne bien
$\left( \hff \circ E_\w \right)_c =  coker(\a)_{\, c} = \left( z_1\,z_2^{q
+ 1} \right)\,$.\\

Enfin lorsque $\,c\,$ est l'intersection de deux composantes irr\'eductibles
$D'$ et $D"$ de $\D_\w\,$,
la singularit\'e $c$ est une selle, toujours gr\^ace \`a $\,4. \Longrightarrow
1.\,$. En posant
$p:= m_{D"}\,(\w) = m_{D"}\,(d\hff)\,$ et $q := m_{D'}\,(\w) =
m_{D'}\,(d\hff)\,$, quitte \`a bien
choisir les coordonn\'ees, on a
$$\ww = z_1^p z_2^q\,\left( \,( \l_2  + \cdots) z_2\,dz_1 + ( \l_1 +
\cdots ) z_1\, dz_2 \,\right)\,,
\quad
\l_1 \l_2 \not= 0\, ,$$
$$\hff \circ E_\w = z_1^{p+1} z_2^{q + 1}\, ,Ê\quad  D' = \left\{ z_2 = 0
\right\}\, , \quad D" =
\left\{ z_1 = 0 \right\} \, .$$  De nouveau $\hX_{\widehat{S}'}$ est
engendr\'e par $z_1\,\dd z_1\,$ et
$z_2\, \dd z_2Ê\,$, ce qui donne bien
$\left( \hff \circ E_\w \right)_c =  coker(\a)_{\, c} = \left(
z_1^{p+1}\,z_2^{q + 1} \right)\,$.
\end{proof}

\begin{coro}\label{redsecesp}
Si $\F_\w$ est de deuxi\`eme esp\`ece, alors la r\'eduction de $\F_\w$ est identique \`a la r\'eduction de ses s\'eparatrices
formelles, c'est \`a dire~: $\A[\F_\w] \equiv \A[d\hat{f}]$, ou encore~: $M_\w = M_{d\hat{f}}$ et $E_\w =
E_{d\hat{f}}$.
\end{coro}

\subsection{Equir\'eductibilit\'e}\label{defer}
\addcontentsline{toc}{section}{\hspace{0,8em} {}\thesubsection .  Equir\'eductibilit\'e}
Dans ce qui suit, consid\'erons un feuilletage
formel $\F$ \`a l'origine de
$\C^2\,$ d\'efini par une 1-forme \`a singularit\'e isol\'ee $\w\,$ et
une d\'eformation $\sF_P$ de $\F$ transversalement formelle le long de
$0 \times P$ de param\`etres $P := (\C^p , \, 0)\,$, cf.
(\ref{standdef}). Conservons les
notations (\ref{ard}),
(\ref{defeno}) et (\ref{defforme}).
La notion d'\'equir\'eductibilit\'e correspondra, comme dans la th\'eorie
d'\'equisingularit\'e de courbes, \`a la
possibilit\'e d'effectuer une "r\'eduction en famille".

\begin{defi}\label{defequ} Supposons $\F$ non-r\'eduit. Nous disons que
$\sF_P$ est {\it \'equir\'eductible} s'il
existe un (germe d') arbre
$\A_P = \left( \M^j\,, E^j\,, \S^j\,, S^j\,, \pi_j\,, \D^j\, \right)_{j =
0,\ldots ,h}$ au dessus de
$P\,$, r\'egulier au sens de (\ref{arreg}) et tel que
\begin{enumerate}
\item $Sing_P(\sF_P) = \S^0 = 0 \times P\,$ et $\M^0 = \C^2 \times P\,$,
\item pour chaque $j = 1,\ldots ,h\,$, $\S^j$ est le lieu singulier
$Sing\,(\F_P^j , \, \D^j )$ du couple form\'e du transform\'e strict $\F_P^j$
de $\F_P$ sur $\M^j$ et du
$j$-i\`eme diviseur exceptionnel, cf. (\ref{singcouples}) et
(\ref{idealredef}),
\item pour toute valeur assez petite de $t \in P$ l'arbre $\A_P(t)\,$ est
exactement l'arbre de
r\'eduction du feuilletage formel $\F_P(t)\,$.
\end{enumerate} Si $\F$ est d\'ej\`a r\'eduit nous disons que $\F_P$ est
\'equir\'eductible si~: $Sing_P(\sF_P)
= 0 \times P$ et $\F_P(t)$ est r\'eduit pour
$t$ assez petit.
\end{defi}

Il est clair que l'arbre $\A_P$ lorsqu'il existe est uniquement d\'etermin\'e
par $\sF_P\,$. Nous
noterons cet arbre
$$\A[\sF_P] = \left( \M_{\sF_P}^j\,, E_{\sF_P}^j\,, \S_{\sF_P}^j\,,
S_{\sF_P}^j\,, \pi_{\sF_P,
\,j}\,, \D_{\sF_P}^j\,
\right)_{j = 0,\ldots ,h}\, ,\quad \mbox{ou encore}$$
$$
\A[\n] = \left( \M_{\n}^j\,, E_{\n}^j\,, \S_{\n}^j\,, S_{\n}^j\,, \pi_{\n,
\,j}\,, \D_{\n}^j\,
\right)_{j = 0,\ldots ,h}$$ lorque $\sF_P$ est d\'efinie par la 1-forme
$\n\,$, comme en
(\ref{standdef}). L'\'equir\'eductibilit\'e donne donc~:
$$
\A[\sF_P](t) = \A[\F_P(t)]\, .
$$

Remarquons les deux exemples int\'eressants de d\'eformations
non-\'equir\'eductibles donn\'ees par les
1-formes suivantes, o\`u $t \in \left( \C , \, 0 \right)$ est le param\`etre
de d\'eformation et
$\l \in (\Bbb R -Ê\Q )_{>0}\,$ est fix\'e~:
\begin{equation}\label{exneq}
\n_1 = x\,dy + y ( y - t ) \,dx \quad \mbox{et}
\quad \n_2 := (t - \l )x \, dy + y \, dx \,  .
\end{equation} Dans le premier la condition 3. est r\'ealis\'ee sans que la
condition 1. le soit. Le
deuxi\`eme exemple montre que l'\'equir\'eductibilit\'e n'est pas une condition
analytique~: $\F$ est r\'eduit,
la condition 1. est r\'ealis\'ee mais pour
$t \in \l + \Q_{<0}$ le feuilletage $\F_P(t)$ est dicritique et donc
non-r\'eduit. On peut aussi
construire des exemples faisant appara\^{\i}tre ces situations apr\`es
\'eclatements. Cependant on voit par un
calcul imm\'ediat~:

\begin{prop}\label{rempre} Supposons la condition 1. de (\ref{defequ})
satisfaite et $\F_P(t)$
non-dicritique pour $t$ assez petit. Alors le lieu singulier relatif
(\ref{idealredef}) du transform\'e
strict $\F_P^1$ de $\F_P$ par l'application d'\'eclatement de $0 \times P$
est un sous-ensemble
analytique de codimension pure 1 du diviseur exceptionnel. Il est fini au
dessus de $P$ si et
seulement si les conditions \'equivalentes suivantes, dites
d'\'equimultiplicit\'e, sont satisfaites~:
\begin{equation}\label{equm} (i) \quad \nu_{(0 ,\, t)}\left(
\F_P(t)\right) =  \nu_0\left( \F
\right)\,,
\qquad (ii) \quad  m_{\D^1(t)}\left( \F_P(t)\right) = m_{\D^1(0)}\left(
\F\right)\, ,
\end{equation} o\`u $\D^1(t)$  d\'esigne le diviseur cr\'e\'e par l'\'eclatement de
l'origine dans $\C^2 \times
\{ t\}\,$. De plus, lorsque $\F$ est r\'eduit chaque
$\F_P(t)$ est aussi r\'eduit et lorsque $\F$ est un selle-n\oe ud, chaque
$\F_P(t)$ est aussi un
selle-n\oe ud.
\end{prop}

La condition
$Sing_P(\sF_P) =  0 \times P\,$ peut s'exprimer en disant que le nombre de
Milnor
$\mu_{(0 , t)}\left( \F_P(t) \right)$ est ind\'ependant de $t\,$, cf.
(\ref{nbrMilnor}). Ainsi,
toujours sous la condition
$\F_P(t)$ non-dicritique, les conditions 1. et 2. de (\ref{defequ})
\'equivalent aux conditions~:
constance du nombre de Milnor et \'equimultiplicit\'e \`a chaque \'etape du
processus de r\'eduction. Le th\'eor\`eme suivant \'etend aux feuilletages de deuxi\`eme esp\`ece un crit\`ere d'\'equir\'eductibilit\'e
prouv\'e dans \cite{Mdouady}.

\begin{theo}\label{deuxequ}
Soit $\F$ un feuilletage formel \`a l'origine de $\C^2\,$ \`a singularit\'e
isol\'ee non-r\'eduite
et non-dicritique et soit
$\sF_P$ une d\'eformation de $\F$ formelle le long de $0 \times P\,$ cf
(\ref{standdef}), de lieu
singulier relatif $0 \times P\,$. Alors~:
\begin{enumerate}
\item $\sF_P$ est \'equir\'eductible si et seulement si pour chaque
$t
\in P$ assez petit, les arbres fl\`ech\'es et doublement pond\'er\'es (par
auto-intersection et par
multiplicit\'es) $\cA[\F_P(t)]\,$ et
$\cA[\F]\,$ d\'efinis en (\ref{subs.red.des.sing.}) sont \'egaux
\footnote {L'\'egalit\'e signifie l'isomorphisme, en un sens clair, d'arbres
fl\`ech\'es et doublement
pond\'er\'es, puisque ces arbres, comme graphes dans
$\Bbb R^2\,$, ne sont  en fait d\'efinis qu'\`a isomophisme pr\`es. }.
\item Si de plus $\F$ est de deuxi\`eme esp\`ece, $\sF_P$ est \'equir\'eductible
si et seulement si pour
chaque $t \in P$ assez petit, les arbres duaux de r\'eduction (fl\`ech\'es et
pond\'er\'es par
auto-intersection)
$\aA[\F_P(t)]\,$ et $\aA[\F]\,$ sont \'egaux.
\end{enumerate}
\end{theo}

\begin{rema}
La condition $Sing_P (\sF_P ) = 0 \times P $ est une hypoth\`ese essentielle
dans ce th\'eor\`eme. Dans
l'exemple de la d\'eformation $\sF_\C$ donn\'ee par $\n := y \, dx + x (x-t)
\, dy \,$ cette
hypoth\`ese n'est pas satisfaite bien qu'en chaque point $( 0 , t)$ le
feuilletage $\F_\C (t)$ soit
r\'eduit et de multiplicit\'e alg\'ebrique un. On construit facilement des
d\'eformations satisfaisant
$\cA[\F_P(t)] = \cA[\F]\,$
o\`u ce ph\'enom\`ene se produit \`a une \'etape du processus
de r\'eduction. Ces d\'eformations ne sont pas \'equir\'eductibles car
$Sing_P (\sF_P ) \not= 0 \times P $.
\end{rema}

\begin{preuvede}{du th\'eor\`eme} Montrons d'abord l'\'equivalence 1. Une des
implication est triviale. Soit
$\n := A(x, y;
\,t) \,dx + B(x, y; \, t)\,dy$ une 1-forme d\'efinissant $\sF_P\,$.  Notons
$\n_t$ la restriction de
$\n$ \`a
$\left( \C^2\times P , \, (0 ; \,t) \right) \,$ et
$$
\left( \M_{\n_t}^j , E_{\n_t}^j , \S_{\n_t}^j , S_{\n_t}^j , \pi_{{\n_t},
\,j} ,\D_{\nu_t}^j
\right)_{j = 0 ,
\ldots , h_{\n_t}}
$$
l'arbre de r\'eduction $\A[\F_P(t)]\,$. D\'esignons par
$\n_t^j$ la 1-forme
diff\'erentielle sur $\M_{\n_t}^j\,$, image r\'eciproque de $\n_t$ par la
compos\'ee $ E_{\n_t}^1 \circ
\cdots \circ  E_{\n_t}^j\,$ et par $\F_P^j(t)$ le transform\'e strict de
$\F_P(t)$ sur
$\M_{\n_t}^j\,$. Supposons satisfaite l'\'egalit\'e
\begin{equation}\label{hypar}
\cA[ \sF_P (t) ] = \cA[ \F ]\, .
\end{equation} Nous allons construire par induction un arbre
d'\'equir\'eduction de $\sF_P\,$, mais
auparavant faisons quelque remarques utiles.\\

La difficult\'e principale est qu'on ne dispose pas de notion de
"continuit\'e" de la correspondance $ t
\efle \cA[\F_P(t) ]\,$. En effet lorsque $\cA[\F]$ poss\`ede des sym\'etries
non-triviales on ne peut pas en g\'en\'eral, \`a partir seulement de la donn\'ee de
$\cA[\F_P(t)]$,
associer de fa\c con canonique une composante
irr\'eductible du diviseur de r\'eduction
$\D_{\n_t}$ de $\F_P(t)$ \`a chaque sommet de $\cA[\F_P(t)]\,$. La condition
(\ref{hypar}) donne
cependant quelques informations globales. On sait en effet retrouver \`a
partir de la pond\'eration par
auto-intersection "l'ordre de cr\'eation des composantes". Plus pr\'ecis\'ement
effectuons sur
$\cA[\F_P(t)]\,$ les op\'erations qui correspondent aux contractions des
composantes de $\D_{\n_t}$
d'auto-intersection -1, c'est \`a dire les "mouvements"  consistant \`a
remplacer
$$
{}^{ {}^{ \displaystyle \cdots} }
 \stackrel{( {e'},\, m')} {\stackrel{\bullet}{\swarrow \!\!{}_{\goth f'_1} \cdots \,\,{}_{\goth f'_{r'}}\!\!\searrow}}
\; { {}^{ \overline{\phantom{(AAA)}}} }\;
 \stackrel{(-1,\, m)} {\stackrel{\bullet}{\swarrow \!\!{}_{\goth f_1} \cdots \,\,{}_{\goth f_{r}}\!\!\searrow}}
\; { {}^{ \overline{\phantom{(AAA)}}} }\;
 \stackrel{(e",\, m")} {\stackrel{\bullet}{\swarrow \!\!{}_{\goth f"_1} \cdots \,\,{}_{\goth f"_{r"}}\!\!\searrow}}
\; {}^{ {}^{ \displaystyle  \cdots} }
$$
par
$$
{}^{ {}^{ \displaystyle  \cdots} }
 \stackrel{( {e'} + 1,\, m')} {\stackrel{\bullet}{\swarrow \!\!{}_{\goth f'_1} \cdots \,\,{}_{\goth f'_{r'}}\!\!\searrow}}
\; { {}^{ \overline{\phantom{(AAA)}}} }\;
 \stackrel{(e" + 1,\, m")} {\stackrel{\bullet}{\swarrow \!\!{}_{\goth f"_1} \cdots \,\,{}_{\goth f"_{r"}}\!\!\searrow}}
\; {}^{ {}^{ \displaystyle  \cdots} }
$$
%$$ {}^{ {}^{\mbox{
et
$$
{}^{ {}^{ \displaystyle \cdots} }
 \stackrel{( {e'},\, m')} {\stackrel{\bullet}{\swarrow \!\!{}_{\goth f_1} \cdots \,\,{}_{\goth f_{r}}\!\!\searrow}}
\; { {}^{ \overline{\phantom{(AAA)}}} }\;
 \stackrel{(-1,\, m)} {\stackrel{\bullet}{\swarrow \!\!{}_{\goth f'_1} \cdots \,\,{}_{\goth f'_{r'}}\!\!\searrow}}
\; {}^{ {}^{ \displaystyle  \cdots} }
$$
par
$$
{}^{ {}^{\displaystyle\cdots} }
 \stackrel{( {e'} + 1,\, m')} {\stackrel{\bullet}{\swarrow \!\!{}_{\goth f_1} \cdots \,\,{}_{\goth f_{r+1}}\!\!\searrow}}
%\!\!\!\! {}^{ {}^{ \displaystyle \longrightarrow}{}_{\goth f_{r+1}} }
\; {}^{ {}^{ \displaystyle  \cdots} }
$$
On obtient un nouveau graphe fl\`ech\'e et doublement pond\'er\'e not\'e $\cA^{h
- 1}[\F_P(t)]\,$, qui est
en fait "l'arbre dual fl\'ech\'e et doublement pond\'er\'e associ\'e aux diviseur
$\D_{\n_t}^{h - 1}\,$". En r\'ep\'etant cette op\'eration on obtient tous les
arbres duaux interm\'edaires
$\cA^j[\F_P(t)]\,$ associ\'es au diviseurs $\D^j_{\n_t}\,$, et cela de
mani\`ere uniquement combinatoire
\`a partir de $\cA[\F_P(t)]\,$. Ainsi la condition (\ref{hypar}) donne les
\'egalit\'es
\begin{equation}\label{crarint}
\cA^j[\F_P(t)] = \cA^j[\F] \, , \qquad j = 1 \ldots h \, ,
\end{equation}
et la constance du nombre
$$ m_j(t) := \sum_{D \in comp\left( \D_{\n_t}^j \right)} m_D (\n_t)\,,
$$
o\`u $comp\left( \D_{\n_t}^j \right)$ d\'esigne l'ensemble des composantes irr\'eductibles.
Mais $\F_P(t)$ est non-dicritique, car $\cA[\F_P(t)]$ n'a pas de symbole "double fl\`eche" et la
multiplicit\'e $m_D(\n_t)\,$
d'une composante $D \subset \D_{\n_t}^{j + 1}$ cr\'e\'ee par l'\'eclatement d'un
point
$c$ de $S_{\n_t}^j$ est \'egale \`a la multiplicit\'e alg\'ebrique de  $\n_t^j\,$
en ce point. Le nombre
$m_{j+1}(t) - m_j(t)$ est donc \'egal \`a la somme des multiplicit\'es
$\nu_c \left( \n_t^j \right)$ de $\n_t^j$ aux points $c \in
S_{\n_t}^j\,$.  D'autre part $\nu_c
\left( \n_t^j \right) - \nu_c \left( \F_P^j(t) \right)$ est la somme des
multiplicit\'es $m_D(\n^j_t)$
des composantes irr\'eductibles $D$ de $\D_{\n_t}^j$ qui passent par $c\,$.
Ainsi le nombre
$$\Delta_j(t) := \sum_{c \in S_{\n_t}^j}\nu_c \left( \n_t^j \right) -
\nu_c \left( \F_P^j(t)
\right)$$ se calcule  \`a partir de $\cA^j[\F_P(t)]$ et $\cA^{j +
1}[\F_P(t)]\,$. Il est constant
d'apr\`es (\ref{crarint}).  On en d\'eduit que
\begin{equation}\label{ctmu}
\widetilde{m}_j(t) := \sum_{c \in S_{\n_t}^j} \nu_c \left( \F_P^j(t)
\right) =  m_{j + 1}(t) -
m_j(t) - \Delta_j(t)
\end{equation} est aussi constant. Nous sommes maintenant en mesure de
d\'ecrire l'induction.\\

Consid\'erons le transform\'e strict $\F^1_P$ du feuilletage $\F_P$ par
l'application
$E^1 : \M^1 \fle \C^2 \times P $ d'\'eclatement de centre $S^0 := 0 \times
P\,$.  D'apr\`es
(\ref{rempre}) son lieu singulier $\S^1$ est une courbe, contenue dans le
diviseur exceptionnel $\D^1
:= (E^1)^{-1}(0 \times P)$ car $Sing_P(\sF_P) = 0 \times P\,$ qui est fini
au dessus de $P$ via
$\pi_1 := \pi_0 \circ E^1\,$,
$\pi_0(x, y ; \, t) := t\,$ d'apr\`es (\ref{rempre}). La restriction
$\F^1_P(t)$ de $\F^1_P$ \`a $\M^1(t)
:= \pi_1^{-1}(t)$ est le transform\'e strict de
$\F_P(t)$ par l'\'eclatement de l'origine et $\S^1(t) := \S^1 \cap \M^1(t)$
est le lieu singulier de
$\F^1_P(t)\,$. Ainsi
$\# \S^1(t) $ est le nombre de fl\`eches et d'ar\^etes port\'ees par le sommet
de $\cA[\F_P(t)]$
correspondant au premier diviseur cr\'e\'e. Ce nombre est constant d'apr\`es
(\ref{crarint}). On en d\'eduit
que $\S^1$ est \'etale au dessus de $P$ via $\pi_1\,$.  Toujours d'apr\`es
(\ref{rempre}) les composantes
de $\S^1$ sont de deux types~: ou bien $\F_P^1(t)$ est r\'eduit en chaque
point, ou bien $\F_P^1(t)$
n'est r\'eduit en aucun point.  Notons $S^1$ l'union des composantes du
second type.\\

\indent Supposons maintenant que l'on ait it\'er\'e $n$ fois cette
construction, c'est \`a dire que l'on
dispose d'un arbre au dessus de $P$ de hauteur $n$
$$\A_n := \left( \M^j , E^j , \S^j , S^j , \pi_j , \D^j \right)_{j = 0 ,
\ldots , n}$$ qui satisfait
les propri\'et\'es suivantes~:  notons comme d'habitude~: $E_j := E^1 \circ
\cdots \circ E^j\,$, $\pi_j := \pi_0 \circ E_j\,$ et  d\'esignons par
$\F_P^j$ le transform\'e strict de
$\F_P$ par $E_j\,$ et par $\F_P^j(t)$ sa restriction \`a $ \M^j(t) :=
\pi_j^{-1}(t)\,$, alors
\begin{enumerate}
\item[(i)] $\S^j$ est le lieu singulier de $\F_P^j\,$,
\item[(ii)] $S^j$ est l'union des composantes irr\'eductibles de $\S^j$
form\'ee des points $m'$ o\`u le
germe de $\F_P^j(\pi_j(m'))$ n'est pas r\'eduit.
\end{enumerate} Param\`etrisons les composantes $S^{\,n}_1 , \ldots ,
S^{\,n}_k$ de $S^n\,$ par des
sections $c_1(t) , \ldots , c_k(t)$ de $\pi_n\,$  et notons $\wnu_r(t)$ la
multiplicit\'e
alg\'ebrique $\nu_{c_r(t)}
\left( \F_P^n(t) \right)$ de
$\ \F_P^n(t)\,$ au point $c_r(t)\,$. Le nombre
$$\wnu^n(t) := \sum_{r = 1}^k \wnu_r(t)$$ est \'egal au nombre
$\widetilde{m}_n(t)$ d\'efinit en
(\ref{ctmu}), car $\A_n(t)$ est une partie de l'arbre de r\'eduction de
$\F_P(t)\,$. Il est constant. Comme les applications
$t
\efle
\wnu_r(t)$ sont localement d\'ecroissantes, chaque multiplicit\'e
$\nu_{c_r(t)} \left( \F_P^n(t) \right)$ est aussi constante. Ainsi le lieu
singulier $\S^{n + 1}$ du
transform\'e strict de
$\F_P$ par l'application $E^{n+1}$ d'\'eclatement de  $S^{\,n}$ est fini au
dessus de $P$
via $\pi_{n + 1} := \pi_n \circ E^{n+1}\,$. Le nombre de points de $\S^{n
+ 1}(t) := \S^{n + 1} \cap
\pi_{n + 1}^{- 1}(t)$ est \'egal au nombre total de fl\`eches et d'ar\^etes de
$\cA^{n+1}[\F_P(t)]\,$. Ce
nombre est donc constant d'apr\`es (\ref{crarint}) et $\S^{n + 1}$ est lisse
\'etale au dessus de $P\,$.
On note $S^{\, n + 1}$ l'union des composantes de $\S^{n + 1}$ contenant
un point non-r\'eduit de
$\F^{n + 1}_P(0)\,$. Les conditions (i) et (ii) ci-dessus sont satisfaites
et l'on peut continuer
l'induction jusqu'\`a r\'eduction compl\`ete.\\

Montrons maintenant 2. Nous reprenons la m\^eme induction. La seule
diff\'erence dans cette
construction est qu'il faut maintenant d\'eduire la constance des
multiplicit\'es $\wnu_r(t)$ de
l'hypoth\`ese~: "$\F$ est de deuxi\`eme esp\`ece". Conservons les m\^emes
notations. Puisquel'on a
$\aA[\F_P(t)] = \aA[\F]\,$, il existe pour chaque $t$ une bijection
$\r_t : \left\{ 1 , \ldots k \right\} \fle \left\{ 1 , \ldots k \right\}$
telle qu'en d\'esignant par
$\left( \F_{P}^n (t) ,\, m \right)\,$ le germe de $\F_P^n(t)$ en $m \in
\pi_n^{-1}(t)\,$, on ait
l'\'egalit\'e des arbres duaux de r\'eduction~:
$$\aA\left [ \,\left( \, \F_{P}^n (t) ,\, c_r(t)\, \right) \,\right]  =
\aA\left[\, \left( \,\F_{P}^n
(0) ,\, c_{\r_t(r)}(0)\, \right)\, \right]\,.$$ On voit facilement que
l'hypoth\`ese~: "$\F$ est de
deuxi\`eme esp\`ece" implique que chaque
$\left(\, \F_{P}^n (0) ,\, c_{r}(0)\, \right)\,$ est aussi de deuxi\`eme
esp\`ece. La propri\'et\'e de
minimalit\'e $(\ast \ast)$ de (\ref{minmult}) et la d\'ecroissance des $t
\efle \wnu_r(t)$ donnent
$$\wnu_{\r_t(r)}(0) \; \leq \; \wnu_r(t)Ê\;  \leq \; \wnu_r(0)\,.$$ En
sommant ces in\'egalit\'es on
obtient $\wnu^n(t) = \wnu^n(0)\,$. De la m\^eme mani\`ere que pr\'ec\'edemment on
conclut \`a la constance de
chaque $\wnu_r(t)\,$. \end{preuvede}

\subsection{G\'en\'ericit\'e de l'\'equir\'eduction}\label{genequred}
\addcontentsline{toc}{section}{\hspace{0,8em} {}\thesubsection .  G\'en\'ericit\'e de l'\'equir\'eduction}
Supposons la d\'eformation
$\sF_P$ fix\'ee au d\'ebut de (\ref{defer}) \`a nombre de Milnor constant, ou
bien, ce qui revient au m\^eme~: $Sing_P\left(
\sF_P
\right) = 0 \times P\,$. Nous nous proposons de d\'ecrire
l'ensemble $NR\left( \sF_P \right)$ des points $t_0 \in P$ tels que o\`u le germe de
$\sF_P$ en $(0, t_0)$ n'est pas une d\'eformation
\'equir\'eductible de $\F_P(t_0)$
\begin{equation}\label{enseqr}
NR\left( \sF_P \right) :=
\left\{ t_0 \in P\left/  \sF_P \phantom{a}\mbox{non-\'equir\'eductible en}\phantom{a} t_0 \right.\right\}.
\end{equation}
L'exemple $\n_2$ de (\ref{exneq}) montre que $NR\left(\sF_P \right)$ n'est
pas
analytique. L'obstruction y r\'eside dans la dicriticit\'e de certains
feuilletages $\sF_P (t)\,$. Pour obtenir un ensemble analytique
ferm\'e, affaiblissons la notion de "forme r\'eduite".\\

\begin{defi}\label{preredu}
Nous disons qu'une 1-forme formelle
$
\w' \, := \, ( \a x + \b y ) \, dx + (\zeta x + \xi y ) \,dy + \cdots\,$,
avec $\a\,$, $\b\,$, $\zeta\,$, $\xi
\, \in \C\,$
est {\it pr\'e-r\'eduite} si les valeurs propres de la matrice
$$
\left(
\begin{array}{cc}
- \zeta & - \xi \\
\a & \b
\end{array}
\right)
$$
sont distinctes. Un feuilletage formel en un point $\cal F\,$
sera dit {\it pr\'e-r\'eduit} si son satur\'e est
d\'efini en ce point par une forme pr\'e-r\'eduite.
\end{defi}
\noindent Un calcul imm\'ediat permet de voir que
$\w'$ est pr\'e-r\'eduite si et seulement si le polyn\^ome homo\-g\`ene de plus bas
degr\'e du d\'eveloppement de $\w' \cdot \left( x\, \dd x + y \, \dd y
\right)\,$ est de degr\'e deux, avec deux racines simples (dans
$\Bbb P^1\,$). Ces points sont alors les singularit\'es du transform\'e strict
de
$\F_{\w'}$ par l'\'eclatement de l'origine. Il est clair, par cette
interpt\'etation g\'eom\'etrique,
que~:
\begin{itemize}
\item[($\star$ )] l'ensemble des valeurs du param\`etre $t_0 \in P$ o\`u
$\F_P(t_0)$ n'est pas pr\'e-r\'eduite est analytique ferm\'e.
\end{itemize}

\noindent On d\'efinit de mani\`ere naturelle l'{\it arbre de pr\'e-r\'eduction}
d'un feuilletage formel $\F$ par la succession
d'\'eclatements d'\'eclatements qui consiste, \`a chaque \'etape, \`a \'eclater
simultan\'ement les points singuliers non-pr\'e-r\'eduits. Cet arbre
est fini, puisque toute forme r\'eduite est pr\'e-r\'eduite. On obtient aussi la notion de {\it
d\'eformation \'equi-pr\'e-r\'eductible} en remplacant dans (\ref{defequ}) la
condition 3. par la condition~:

\begin{enumerate}\it
\item[3'.] pour toute valeur assez petite de $t \in P$ l'arbre $\A_P(t)\,$
est l'arbre de
pr\'e-r\'eduction du feuilletage formel $\F_P(t)\,$.
\end{enumerate}

\noindent Consid\'erons l'ensemble $NPR\left( \sF_P \right) $ des $t_0 \in P$ tels que le germe de $\sF_P$ en $(0, t)$ n'est
pas une d\'eformation \'equi-pr\'e-r\'eductible de $\sF_P(t_0)$,
$$NPR \left( \sF_P \right) :=
\left\{\, t_0 \in P \;\left/\;  \sF_P \phantom{a}\mbox{non-\'equi-pr\'e-r\'eductible en}
\phantom{a} t_0\right.\,
\right\}\, $$
Par des arguments standard et \`a l'aide de
($\ast$), on montre~:

\begin{prop}\label{genprer}
$NPR \left( \sF_P \right)$ d\'efinit un germe \`a l'origine de sous-ensemble analytique ferm\'e de $P$, de
dimension strictement inf\'erieure \`a la dimension de $P$,
\'eventuellement vide.
\end{prop}
\noindent Pour d\'ecrire $NR\left(\sF_P \right)$ il reste \`a analyser comment
s'effectue le passage de la
pr\'e-r\'eduction \`a la r\'eduction. Remarquons d'abord que le diviseur de c\^{\i}me de l'arbre de
pr\'e-r\'eduction
peut poss\`eder des {\it composantes dicritiques}, c'est \`a dire des
composantes irr\'eductibles
non-invariantes. On appellera ces composantes,
{\it composantes dicritiques du premier type}. On montre facilement, par
r\'ecurence sur la hauteur de l'arbre de pr\'e-r\'eduction de
$\sF_P(0)\,$, le lemme suivant.

\begin{lemm}\label{dicpreman}
Supposons la deformation $\sF_P$ pr\'e-equir\'eductible. Alors l'ensemble
$\hbox{\rm Dic}^{(1)}\left( \sF_P \right)$ des $t \in P$ tels que
la r\'eduction de
$\F_P(t)$ admette une composante dicritique du premier type est analytique
ferm\'e, \'eventuellement vide ou \'egale \`a $P$.
\end{lemm}

\indent  Sur le diviseur de c\^{\i}me $\cal D'$
de l'arbre de pr\'e-r\'eduction
de $\F$, consid\'erons un point $m$ o\`u le germe de la transform\'ee stricte $\wF'$ de $\cal F$ est pr\'e-r\'eduit mais non-r\'eduit.
Dans de bonnes coordonn\'ees $z_1\,$,
$z_2\,$, $\F'$ est donn\'e par une 1-forme du type $\a z_2 \,d z_1 - \b z_1 \,d z_2 \,+
\cdots\,$ avec
$\b /
\a
\in
\Q_{>0}\,$ et $\beta/\alpha \not= 1$.
Nous laissons au lecteur le
soin de v\'erifier les affirmations ci-dessous, en examinant chaque \'eventualit\'e~: $m$ est un point
r\'egulier de $ {\D'}$, ou bien $m$ est un point singulier de
$\cal D'\,$~:\\

\indent 1) $\;\;\b / \a \notin \N^{\ast}
\cup
\frac{\displaystyle 1}{\displaystyle \phantom{{}^{\ast}}\N^{\ast}}\;$ . Le
feuilletage est alors t.f. lin\'earisable dans des coordonn\'ees transversalement formelles,
il admet une int\'egrale premi\`ere transversalement-formelle-m\'eromorphe
\footnote
{
C'est \`a dire un \'el\'ement du corps des fractions de l'anneau des germes, au
point $m$, de fonctions transversalement formelles.
}et une ou plusieurs composantes dicritiques apparaissent dans la
succession d'\'eclatements suppl\'ementaires que l'on
effectue, au dessus du point $m$, pour aboutir \`a des singularit\'es
r\'eduites. Nous appellerons ces composantes, {\it composantes
dicritiques du deuxi\`eme type}.\\

\indent 2)  $\;\; \b / \a \in  \N^{\ast}
\cup
\frac{\displaystyle 1}{\displaystyle \phantom{{}^{\ast}}\N^{\ast}}\;$. Il
existe alors des coordonn\'ees "normalisantes"
transversalement formelles
dans lesquelles le feuilletage est donn\'e par la 1-forme
$$
\;( n \,z_1 + \l \,z_2^n ) \,d z_2\, -\, z_2 dz_1\,, \qquad \l \in \C\,, \quad n \geq 2.
$$
\noindent Lorsque $\l$ est non-nul, $\{z_2 = 0\}$ est l'unique
s\'eparatrice. Cette courbe
correspond au diviseur, qui est lisse au point
$m$. La r\'eduction se fait par une cha\^{\i}ne de $n$ \'eclatements
suppl\'ementaires. Elle cr\'ee
$n-1$ singularit\'es de type selle et une singularit\'e de type selle-n\oe ud,
situ\'ee au point d'intersection de l'avant
derni\`ere composante et de la derni\`ere composante cr\'e\'e, dont c'est l'unique
singularit\'e. L'holonomie de cette
derni\`ere composante est l'identit\'e et ce selle-n\oe ud est
tangent-r\'esonant, c.f. (\ref{sntagent}). Par contre, si $\l = 0$, le
feuilletage se r\'eduit encore par une cha\^{\i}ne de $n$ \'eclatements faisant
appara\^{\i}tre $n-1$ singularit\'es de type selle, mais le dernier
diviseur cr\'ee est dicritique, partout transverse au transform\'e strict du
feuilletage. Nous appellerons une telle composante, {\it
composante dicritique du troisi\`eme type}. \\

\indent Signalons enfin que le coefficient $\l$ et les $n$-jets au point
$m$ de $z_1$ et $z_2$ sont donn\'es alg\'ebriquement par le
$n$-jet d'une 1-forme d\'efinissant le feuilletage au point $m$. Lorque $\l$
est non-nul, on peut le choisir \'egal \`a 1, quitte \`a
effectuer une homothetie. En particulier, la distinction entre les deux
types de formes normales~: $\l = 1 $ ou bien $\l = 0\,$, est
donn\'ee par le
$n$-jet du feuilletage le long du diviseur
\footnote{
Deux formes formelles le long d'une courbe $\wD'\,$ {\it ont m\^eme $l$-jet
le long de $\wD'\,$} si elles diff\`erent d'une forme \`a
coefficients dans $\cal I_{\wD'}^{l + 1}\,$, o\`u $\cal I_{\wD'}$ d\'esigne le
faisceau des fonctions transversallement formelles le
long de
$\wD'$ qui s'annullent sur $\wD'\,$.}
.\\

Consid\'erons les ensembles $\hbox{\rm Dic}^{(r)}\left( \sF_P \right)$, $r = 1, 2, 3$,
des $t \in P$
tels que la r\'eduction de $\F_P(t)$ admette une composante
dicritique du type $r$, et notons $\hbox{\rm Dic}\left( \sF_P \right) := \cup_{r=1}^3 \hbox{\rm Dic}^{(r)}\left( \sF_P
\right)$. On d\'eduit sans trop de peine de ce qui pr\'ec\`ede, les
propri\'et\'es suivantes~:

\begin{prop}\label{preddeux}
Pour les feuilletages (non-dicritiques) de deuxi\`eme esp\`ece, la
pr\'e-r\'eduction des singularit\'es est \'egale \`a la r\'eduction des
singularit\'es.
\end{prop}

\begin{theo}\label{geredd}
Supposons la d\'eformation $\sF_P$ \`a nombre de Milnor constant et supposons
aussi $\hbox{\rm Dic}\left( \sF_P
\right) \not= P$. Alors~:
\begin{enumerate}
\item $\hbox{\rm Dic}^{(2)}\left( \sF_P \right) \cup \hbox{\rm Dic}^{(3)}\left( \sF_P \right)
\subset NR\left(\sF_P \right) - NPR\left( \sF_P \right)$
\item $\hbox{\rm Dic}^{(3)}\left( \sF_P \right) \subset \overline{\hbox{\rm Dic}^{(2)}\left( \sF_P \right)}
= {NR\left(\sF_P \right) - NPR\left( \sF_P \right)}$
\item $NR\left(\sF_P \right) - NPR\left( \sF_P \right)$ est ou bien vide, ou bien une union d\'enombrable
de sous ensembles analytiques ferm\'es de $P$ de dimension $< dim P$ et $\overline{NR\left(\sF_P \right) -
NPR\left(
\sF_P
\right)}$ est un sous-ensemble semi-analytique r\'eel $\not= P$.
\end{enumerate}
Supposons, de plus, que la d\'eformation $\sF_P(t)$ d\'epende alg\'ebriquement
de $t$, i.e. existe
une forme diff\'erentielle
$\Omega := A (x, y ; t) \, dx + B (x, y ; t) \,dy \,$, avec $A (x, y ; t)$, $B (x, y ; t) \in \C[t] \, [[ x ,
y]]$ qui d\'efinit $\sF_P (t)$ chaque fois que l'on fixe $t \in P$
proche de
l'origine. Alors, dans les propri\'et\'es 1. et 2. ci-dessus, on peut
remplacer les termes~: "analytique" et "semi-analytique r\'eel" par
les termes "alg\'ebrique" et "semi-alg\'ebrique r\'eel".
\end{theo}
%%%%%%%%%%%%%%%%%%%%%%%%%%%%%%%%%%%%%%%%%

%%%%%%%%%%%%%%%%%%%%%%%%%%%%%%%%%%%%%%%%%
\section{\'Equisingularit\'e semi-locale}\label{sec.e.sing}  Fixons un
feuilletage formel $\F$
\`a l'origine de $\C^2\,$, non-dicritique, donn\'e par une 1-forme formelle \`a
singularit\'e isol\'ee
$\w := \stw\,$. Notons encore $E_\F : \M_\F \fle \C^2$ l'application de
r\'eduction de $\F\,$ ainsi
que $\D_\F$ le diviseur exceptionnel et
$\wF$ le transform\'e strict de $\F$ sur $\M_\F\,$. Donn\'ee une d\'eformation
\'equir\'eductible
$\sF_P$ de $\F$ de param\`etres $P := \left( \C^p, \, 0 \right)\,$, nous
noterons  $E_{\sF_P} :
\M_{\sF_P}
\fle
\C^2
\times P$ l'application d'\'equir\'eduction de $\sF_P\,$,
$\pi_{\sF_P}$ le compos\'e de $E_{\sF_P}$ avec la projection canonique sur
$P\,$,
 $\D_{\sF_P}\,$ le diviseur exceptionnel $E_{\sF_P}^{-1}\left( 0 \times P
\right)$ et
$\wF_P\,$ le transform\'e strict de
$\F_P$ sur $\M_{\sF_P}\,$. Nous identifions encore $\M_\F$ \`a la fibre de
$\pi_{\sF_P}^{-1}(0)$ et $\D_\F$ \`a l'intersection de cette fibre avec
$\D_{\sF_P}\,$. Enfin fixons
un recouvrement distingu\'e $\U$ de $\D_\F\,$, c.f. (\ref{rec.dist.}).

\subsection{Equivalence semi-locale de d\'eformations}\label{def.e.s.}
\addcontentsline{toc}{section}{\hspace{0,8em} {}\thesubsection .  Equivalence semi-locale de d\'eformations}
Consid\'erons deux d\'efor\-mations $\sF_P$ et
$\sF'_P$ de
$\F$ transversalement formelles le long de $0 \times P\,$.

\begin{defi} Nous disons que $\sF_P$ et $\sF'_P$ sont {\it semi-localement
\'equi\-valentes}, en
abr\'eg\'e {\it
\hSL -\'equivalentes} et nous notons $\sF_P \heSL \, \sF'_P\,$, si $\sF_P$
et $\sF'_P$ sont
\'equir\'eductibles et si elles induisent
des syst\`emes semi-locaux (\ref{ssloc}) formellement conjugu\'es.
\end{defi}

Explicitons cette d\'efinition.  La relation $\sF_P \heSL \,\sF'_P\,$
signifie l'existence pour
chaque $U \in \U\,$, d'un germe le long de $U$ de diff\'eomorphisme
$\Psi_U : \left( \M_{\sF'_P}, \, U \right) \fle \left( \M_{\sF_P}, \, U
\right)$  transversalement
formel le long de $\D_{\sF'_P}\,$, qui commute aux projections
$\pi_{\sF'_P}\,$ et
$\pi_{\sF_P}\,$, vaut l'identit\'e sur $\M_\F$ et conjugue les germes le
long de $U$ de
$\wF'_P\,$ et $\wF_P\,$.\\

Fixons une d\'eformation \'equir\'eductible $\sF_P\,$. Nous pouvons exprimer
comme un espace de
cohomologie l'ensemble des classes formelles de d\'eformations
\'equir\'eductibles de $\F\,$,
transversalement formelles le long de
$0
\times P\,$,
 qui sont {\it \hSL -model\'ees} sur $\sF_P\,$~:
$$\widehat{Mod} \left( \sF_P \right) := \frac{ \left\{ \, \sF'_P \,
\hbox{d\'eformation
\'equir\'eductible de $\F$} \; \left/ \;
\sF_P \heSL \,\sF'_P\; \right. \right\}}{\sim_{P, \,for}} \; , $$
\noindent o\`u $\sim_{P, \,for}$ d\'esigne comme en (\ref{def.equired.}) la relation
d'\'equivalence formelle
entre d\'eformations. La classe de conjugaison de
$\sF'_P$ sera not\'ee $\left\{ \sF'_P \right\} \,$.\\

\indent Consid\'erons le faisceau $\hhG \simeq \hG$ d\'efini en (\ref{rem1})
des germes aux points de $\D_\F$ de diff\'e\-omor\-phismes
transversalement formels le long de $\D_{\sF_P}$, qui commutent avec
$\pi_{\sF_P}$ et valent
l'identit\'e sur
$\M_\F\,$. Soit
$\widehat{Aut}Ê\left(Ê\sF_P \right)$ le sous-faisceau de $\hG$ form\'e des
\'el\'ements qui
laissent invariant $\wF_P\,$. Supposons l'\'equivalence $\sF_P \heSL
\,\sF'_P\,$ satisfaite.
Avec les notations
pr\'ec\'edentes, la collection $\Phi_{UV} := \Psi_U \circ
\Psi_V^{-1}\,$,
$ U , V
\in \U\,$, $U \cap V \not= \emptyset\,$, est une 1-cochaine \`a valeurs dans
$\widehat{Aut}Ê\left(Ê\sF_P
\right)\,$. Elle d\'efinit un 1-cocycle
$\left[ \Psi_U \circ \Psi_V^{-1} \right] \in H^1 \left( \U , \,
\widehat{Aut}Ê\left(Ê\sF_P
\right) \right)\,$, not\'e
$\left[ \sF'_P \right]\,$. Soit $\sF"_P$ une autre d\'eformation \hSL
-\'equivalente \`a
$\sF_P\,$. En interpr\'etant la relation de cohomologie comme une relation
de compatibilit\'e de
conjugaisons locales, on voit imm\'ediatement~:
$$\sF'_P \; \sim_{P, \,for} \;\sF"_P\quad \Longleftrightarrow \quad
\,\left[ \, \sF'_P
\,\right] = \left[\, \sF"_P\, \right]\,.$$
\noindent On a ainsi une injection
\begin{equation}\label{idco}
\widehat{Mod} \left( \sF_P \right) \longrightarrow H^1 \left( \U , \,
\widehat{Aut}Ê\left(Ê\sF_P \right) \right)\, , \qquad
\left\{ \sF'_P \right\} \efle \left[ \sF'_P \right]\, .
\end{equation}

\begin{theo}[de r\'ealisation]\label{intco} L'application (\ref{idco}) ci-dessus est
bijective.
\end{theo}

\begin{proof} La m\'ethode de construction est similaire \`a celle d\'ej\`a
utilis\'ee pour montrer le
th\'eor\`eme de r\'ealisation  (\ref{realis.}).  Soit $\CC := \left( \Phi_{UV}
\right) \in Z^1
\left( \U , \, \widehat{Aut}Ê\left(Ê\sF_P \right) \right)\,$.
Appliquons l'assertion 5. du th\'eor\`eme
de d\'etermination finie (\ref{mth.co.}) \`a l'arbre d'\'equi\-r\'eduction $\A [
\sF_P ]\,$ et au cocycle
induit par $\CC$ dans $Z^1 \left( \U \, ; \, \hhG \right) \,$.
Il existe le long de chaque $ U \in \U$  un germe $\Phi_U$
de diff\'eomorphisme transversalement formel et des germes de
biholomorphismes $
\widetilde{\Phi}_{UV}$ le long des intersections non-vides $U \cap V$, $\,
U , \, V \in
\U\,$, qui satisfont les relations de commutation avec
$\pi_{\sF_P}\,$, valent l'identit\'e au dessus de 0 et tels que~:
${\Phi}_{UV} = {\Phi}_U \circ \widetilde{\Phi}_{UV} \circ {\Phi}_V^{-1}
\,$. Soit
$\M_\CC\,$ la vari\'et\'e holomorphe obtenue  en recollant par les
biholomorphismes
$\widetilde{\Phi}_{UV}\,$  des voisinages, dans $\M_{\sF_P}\,$, des
ouverts $U \in \U$. Elle est
naturellement munie d'un diviseur $\D_\CC\,$, d'une submersion $\pi_\CC$
sur $P$ et d'un
plongement $\s_\CC : \M_\F \ifle \M_\CC\,$. D'apr\`es le lemme
(\ref{stab.cimes}) elle est
biholomorphe (au dessus de $P$) \`a la cime $ {\M'}$ d'un arbre $\A'_P$ au
dessus de
$P$; de plus
$\A'_P(0)$ est l'arbre de r\'eduction de $\F\,$. L'invariance de $\sF_P$ par
${\Phi}_{UV}$ donne
$${\widetilde{\Phi}_{UV}}^{\ast}\left( {{\Phi}_{U}}^{\ast}\left( \sF_{P}
\right) \right) =
{{\Phi}_{V}}^{\ast}\left( \sF_{P}\right)\,.$$
\noindent Cette \'egalit\'e signifie que $\M_\CC$ est aussi munie d'un
feuilletage $\F_\CC\,$,
transversalement formel le long de
$\D_\CC\,$ tel que $\s_\CC^{\ast}\left( \F_\CC \right) =\wF\,$. On obtient
ainsi sur $ {\M'}$ un
feuilletage transversalement formel le long du diviseur exceptionnel qui
est une d\'eformation de
$\wF\,$. Son image directe, cf. (\ref{im.dir}), (\ref{im.dir.for}), sur
$\C^2 \times P\,$ est, par
construction, une d\'eformation de $\F$ qui est \hSL -\'equivalente \`a
$\sF_P\,$.
\end{proof}

\begin{rema}\label{indepsd} A priori tout ce qui pr\'ec\`ede semble d\'ependre du
choix du recouvrement
distingu\'e $\U$ de $\D_\F\,$. Cela n'est pas g\^enant~: tout d'abord, pour ce
qui nous int\'eresse dans
ce travail, nous pouvons fixer $\U\,$; d'autre part nous aurions pu
d\'efinir la notion de
\hsl -\'equivalence ainsi que l'espace de cohomologie associ\'e en passant \`a
la limite sur l'ensemble
des recouvrements distingu\'es de
$\D_\F$ (qui forme clairement un syst\`eme inductif). Mais cette pr\'ecaution
est superflue. Il est en
effet possible de d\'emontrer que, si les transform\'es stricts de deux
d\'eformations \'equir\'eductibles
de
$\F$ sont formellement conjugu\'es au voisinage d'un \'el\'ement
$U$ de $\U\,$, alors ils le sont aussi au voisinage de tout \'el\'ement $U'\,$
d'un recouvrement
distingu\'e $\U'$ tel que $U \subset U'\,$. Il en est de m\^eme pour des
intersections d'ouverts
de recouvrements distingu\'es.
\end{rema}

\subsection{D\'eformations \hSL-\'equisinguli\`eres}
\addcontentsline{toc}{section}{\hspace{0,8em} {}\thesubsection .  D\'eformations \hSL-\'equisinguli\`eres}
Notons
$\sF_P^{cst}$ la d\'eformation constante de $\F\,$ de param\`etres $P\,$
d\'efinie en (\ref{defcte}),
c'est \`a dire la d\'eformation d\'efinie par
$\stw\,$, consid\'er\'ee comme 1-forme transversalement formelle le long de
$0 \times P\,$ dans
$\left( \C^2 \times P , \, 0 \right)\,$. C'est \'evidemment une
d\'e\-for\-ma\-tion \'equir\'eductible et
l'application d'\'equir\'eduction $E_{\sF_P^{cst}}$ est le produit $E_\F
\times Id_P : \M_\F \times P \fle
\C^2
\times P$ de l'application de r\'eduction de $\F$ par l'identit\'e de $P\,$.
Rappelons qu'une d\'eformation est
dite {\it formellement triviale} si elle est formellement conjugu\'ee \`a
$\sF_P^{cst}\,$.

\begin{defi} Une d\'eformation $\sF_P$ transversalement formelle le long de
$0 \times
P\,$ est dite {\it \hSL -\'equisinguli\`ere} si elle est \'equir\'eductible et
\hSL -\'equivalente \`a
$\sF_P^{cst}\,$.
\end{defi}

La proposition (\ref{intco}) donne l'identification.
\begin{equation}\label{cosla}
\widehat{Mod} \left( \sF_P^{cst} \right) \, \iso  \,
H^1\left( \U ; \,\widehat{Aut}Ê\left(Ê\sF_P^{cst} \right)
\right)\, ,
\end{equation}

\indent Un germe de $\widehat{Aut}Ê\left(Ê\sF_P^{cst}\right)$ s'interpr\`ete
comme un (germe
de) famille
$\Phi_t$ de diff\'eomorphismes de $\M_\F\,$, transversalement formels le
long de $\D_\F\,$, qui
laissent invariant le feuilletage r\'eduit $\wF$ et telle que $\Phi_0$ soit
l'identit\'e. Aussi nous
adopterons souvent la notation
$\left( \Phi_{W, \, t}\right)_t$ pour les sections de
$\widehat{Aut}Ê\left(Ê\sF_P^{cst}\right)$ sur un ouvert  $W$
de  $\D_\F\,$.

\begin{lemm}\label{lemcroix} Soient $\sF_{P_1}$ et $\sF_{P_2}$ deux
d\'eformations \hsl
-\'equisinguli\`eres de $\F\,$ de param\`etres $P_j := \left( \C^{p_j} , \, 0
\right)\,$, $j = 1 , \,2\,$. Il existe alors une d\'eformation \hsl
-\'equisinguli\`ere $\sF_{P_1 \times
P_2}$ de
$\F$ de param\`etres
${P_1 \times P_2}$ qui est conjugu\'ee \`a $\sF_{P_1}\,$, respectivement \`a
$\sF_{P_2}\,$ lorsqu'on
restreint l'espace des param\`etres \`a $P_1 \times 0\,$,  respectivement \`a $0
\times P_2\,$.
\end{lemm}

\begin{proof} Notons $\left( \Phi^j_{UV , \, t_j} \right) \in
Z^1\left( \U ; \,\widehat{Aut}Ê\left(Ê\sF_{P_j}^{cst} \right) \right)
\,$ un cocycle associ\'e \`a
$\sF_{P_j}\,$,  par la bijection (\ref{cosla}), $j = 1\,$, $2\,$.
Le recouvrement $\U$ n'ayant pas d'intersection 3 \`a 3 non-triviale, les
"compos\'es"
$\Phi^1_{UV , \, t_1} \circ \Phi^2_{UV , \, t_2}\,$, $U \cap V \not=
\emptyset\,$, d\'efinissent un
\'el\'ement de
$Z^1 \left( \U , \, \widehat{Aut}Ê\left(Ê\sF_{P_1 \times P_2}^{cst}
\right) \right)\,$. La
proposition (\ref{intco}) donne la d\'eformation cherch\'ee.
\end{proof}

\noindent Ce lemme, qui peut s'interpr\'eter comme une propri\'et\'e
d'irr\'eductibilit\'e de l'espace
des d\'eformations
\hsl -\'equisinguli\`eres, est essentiel dans la d\'emonstration du th\'eor\`eme de
versalit\'e (\ref{thver}).

\subsection{Champs basiques, champs transverses}
\addcontentsline{toc}{section}{\hspace{0,8em} {}\thesubsection .  Champs basiques, champs transverses}
Un champ de vecteurs de $\M_{\sF_P}$ transversalement formel le
long d'un ouvert de
$\D_{\sF_P}$ sera appel\'e ici  {\it basique pour $\sF_P\,$} s'il est
tangent \`a $\D_{\sF_P}$ et si
son  flot laisse le feuilletage $\wF_P\,$ invariant.  Les germes de tels
champs forment un
faisceau de base $\D_{\sF_P}$ que nous notons $\underline{\hB}_{\wF_P}\,$.
D\'esignons par
$\underline{\hB}^{\,v}_{\wF_P}$ le sous-faisceau de
$\underline{\hB}_{\wF_P}$ des germes de champs de vecteurs {\it basiques
verticaux}, c'est \`a
dire basiques et qui annulent l'application tangente
$T\pi_{\sF_P}\,$. D\'esignons aussi par $\underline{\hX}^{\,v}_{\wF_P}$ le
sous-faisceau de
$\underline{\hB}^{\,v}_{\wF_P}$ des champs verticaux tangents \`a $\wF_P\,$.
Ce dernier est un faisceau de
modules coh\'erenr, en fait localement libre de rang 1, sur le faisceau
d'anneaux
$\underline{\hO}^{\D_{\sF_P}}_{\M_{\sF_P}}\,$, de base $\D_{\sF_P}$, des
germes de fonctions
transversalement formelles le long de $\D_{\sF_P}\,$, cf.
(\ref{esptrform}). Ce n'est plus le cas
de
$\underline{\hB}^{\,v}_{\wF_P}$ qui poss\`ede seulement une structure de
faisceau de $\pi_{\sF_P}^{\ast}\left(
\O_P
\right)$-modules.

\begin{defi}\label{sytr}
Nous appelons {\it champ transverse}
\footnote{
Appel\'e aussi {\it sym\'etrie transverse} dans \cite{C-M}.
}
de $\wF_P$ toute section du
faisceau $\underline{\hT}_{\wF_P}$ qui est d\'efini par la suite exacte
\begin{equation}\label{exchtr}
0 \fle  \underline{\hX}^{\,v}_{\wF_P} \fle \underline{\hB}^{\,v}_{\wF_P}
\fle
\underline{\hT}_{\wF_P} \fle 0 \, .
\end{equation}
\noindent Nous notons $\{ X \}$ la section de $\underline{\hT}_{\wF_P}$
d\'efinie par une section
$X$ de $\underline{\hB}^{\,v}_{\wF_P}\,$.
\end{defi}

\indent Restreignons la base de ces faisceaux \`a $\D_\F \subset
\D_{\sF_P}\,$, en d'autres termes
consid\'erons les faisceaux
$${\hB}^{\,v}_{\wF_P} := i^{-1}\left( \underline{\hB}^{\,v}_{\wF_P}
\right) \,, \qquad
\hX^{\,v}_{\wF_P} := i^{-1}\left( \underline{\hX}^{\,v}_{\wF_P} \right) ,
$$
$$
\hT_{\wF_P} := i^{-1}\left( \underline{\hT}_{\wF_P} \right)\, , \qquad
{\hO}_{\M_{\sF_P}} :=  i^{-1}\left(
\underline{\hO}^{\D_{\sF_P}}_{\M_{\sF_P}}\right) \, ,$$
\noindent o\`u $i : \D_\F \ifle \D_{\sF_P}\,$ est l'application
d'inclusion.\\

Explicitons les
fibres de ces faisceaux en un point $m \in \D_\F$ o\`u le feuilletage
$\wF_P$ est r\'egulier. Fixons
des  coordonn\'ees transversalement formelles
$z_1\,$, $z_2\,$, $t_1, \ldots , t_p\,$ dans lesquelles $\wF_P$ est d\'efini
par le champ
$\dd{z_1}\,$, $\D_{\sF_P}$ par $z_2 = 0\,$ et
$\pi_{\sF_P} = (t_1 , \ldots , t_p)\,$. Alors les \'el\'ements de $\hB^{\,
v}_{\wF_P , \, m}$ sont les
germes de champs transversalement formels (en la variable $z_2\,$) qui
s'\'ecrivent
\begin{equation}\label{fca}
\a (z_1, z_2, t_1, \ldots , t_p ) \, \dd z_1 + \b (z_2, t_1, \ldots , t_p
)\,z_2 \dd z_2\, ,
\end{equation}
et $\hX^{\, v}_{\wF_P , \, m}$ est le ${\hO}_{\M_{\sF_P}, \,
m}$-module libre de
base
$\dd{z_1}$. On voit ainsi que, si la coordonn\'ee $z_2$ est d\'efinie sur un
ouvert $W$ de $\,\D_\F -
Sing(\wF)\,$, i.e.
$z_2 \in {\hO}_{\M_{\sF_P}}(W)\,$, alors l'application de
restriction
$\hT_{\wF_P}(W) \fle \hT_{\wF_P , \,m}$ est surjective. D'o\`u~:

\begin{prop}\label{prsec} Sur l'ouvert $\D_\F - Sing \left( \wF \right)\,$
le faisceau
$\hT_{\wF_P}$ est localement constant. En particulier, donn\'es deux ouverts
$U$ et $W$ de
$\D_\F\,$, $W \subset U\,$, dont les groupes fondamentaux en un point de
base $m \in W$ satisfont
$\pi_1(W ;\, m) = \pi_1(UÊ;\, m)\,$ et tels que $( U - W ) \cap Sing\left(
\wF \right) =
\emptyset\,$, alors,   toute section de $\hT_{\wF_P}$ sur $W$ se prolonge
de mani\`ere unique  en
une section de $\hT_{\wF_P}$ sur $U\,$.
\end{prop}

\begin{lemm}\label{surjtau}
Soit $W$ un ouvert d'une composante irr\'eductible $D$ de $\D_\F\,$, $W
\not= D\,$. Alors
l'application canonique ${\hB}^{\,v}_{\wF_P} ( W ) \fle {\hT}_{\wF_P}( W
)$ est surjective.
\end{lemm}

\begin{proof} Visiblement $\hX^{\,v}_{\wF_P}\,$ est un faisceau de
${\hO}_{\M_{\sF_P}}$-modules localement libre de rang 1. A
l'aide de la suite exacte
de cohomologie associ\'ee \`a (\ref{exchtr}) il suffit de prouver que
$H^1 \left( W \, ; \, {\hO}_{\M_{\sF_P}} \right) $ est nul.
La restriction de $\underline{\hO}^{\D_{\sF_P}}_{\M_{\sF_P}}$ \`a tout
voisinage de
Stein
$V $ de $W$ dans $\D_{\sF_P}\,$ est isomorphe \`a
${\O_V}^{\N}\,$ d'apr\`es les \'ecritures (\ref{explfora}) et
(\ref{explforb}). La conclusion r\'esulte
de la coh\'erence de $\O_V$ et de l'existence d'un syst\`eme fondamental de
voisinages de Stein de $W$
dans
$\D_{\sF_P}\,$.
\end{proof}

\begin{prop}\label{dlousski}
La suite exacte (\ref{exchtr}) associ\'ee au recouvrement distingu\'e $\U$ de
$\D_\F$ induit une suite
exacte longue~
$$
\begin{array}{ccccccccc}
0 & \fle & H^0\left(\U \, ; \,\hX^{\,v}_{\wF_P}  \right) & \fle &
H^0\left(\U \, ; \,
{\hB}^{\,v}_{\wF_P} \right) &
\fle &  H^0\left(\U
\, ; \, \hT_{\wF_P} \right) &  & \\
 &\stackrel{\delta_0}{\fle} & H^1\left(\U \, ; \,\hX^{\,v}_{\wF_P}
\right) & \fle &  H^1\left(\U
\, ; \, {\hB}^{\,v}_{\wF_P} \right) &
\fle &  H^1\left(\U
\, ; \, \hT_{\wF_P} \right) & \stackrel{\delta_1}{\fle} & 0
\end{array}
$$
\end{prop}

\begin{proof}
Consid\'erons le diagramme suivant dont les lignes sont exactes d'apr\`es le
lemme ci dessus~:
$$
\begin{array}{ccccccccc}
0 & \fle &\displaystyle \prod_{U\in \U} \hX^{\,v}_{\wF_P}( U ) & \fle &
\displaystyle \prod_{U\in
\U}
\hB^{\,v}_{\wF_P}( U ) & \fle & \displaystyle \prod_{U\in \U}
\hT^{\,v}_{\wF_P}( U ) & \fle & 0\\
& & \downarrow & & \downarrow & & \downarrow & &      \\
0 & \fle & \displaystyle \prod_{U , V\in \U} \hX^{\,v}_{\wF_P}( U \cap V )
& \fle &  \displaystyle
\prod_{U , V
\in \U}
\hB^{\,v}_{\wF_P}( U \cap V) &
\fle & \displaystyle \prod_{U , V \in \U} \hT^{\,v}_{\wF_P}( U \cap V) &
\fle & 0
\end{array}\; .
$$
Comme le recouvrement $\U$ est sans intersection trois \`a trois, la suite
du serpent de ce
diagramme correspond \`a la suite longue de cohomologie de la proposition.
\end{proof}

\vspace{1,5em}

\indent Pour simplifier le texte d\'esignons ici par $\cal N_P$ l'un des
trois faisceaux
\begin{equation}\label{notfcan} {\hB}^{\,v}_{\wF_P}\, , \qquad
\hX^{\,v}_{\wF_P}\ , \qquad
\hT_{\wF_P}\,.
\end{equation}  Faisons varier $\cal W$ dans le syst\`eme inductif de tous
les recouvrements
distingu\'es de $\D_\F\,$.

\begin{prop}\label{liminco}
Chaque application canonique de
$H^1 \left( \U\, ;\, \cal N_P \right)$ dans $\limind_{\cal W} H^1 \left(
\cal W\, ; \, \cal N_P
\right)$ est bijective.
\end{prop}

\begin{proof} Consid\'erons d'abord le cas $\cal N_P =
\hX^{\,v}_{\wF_P}\,$. On a vu dans la
d\'emonstration du lemme (\ref{surjtau}) que tout recouvrement distingu\'e est
acyclique
pour
ce faisceau. On conclut par le th\'eor\`eme de Leray. Le cas $\cal N_P =
\hT_{\wF_P}$ est une cons\'equence directe de la proposition (\ref{prsec}).
Enfin le cas
$\cal N_P =
\hB^{\,v}_{\wF_P}$ s'obtient par passage \`a la limite inductive des suites
exactes
longues (\ref{dlousski}).
\end{proof}

Lorsque $P$ est r\'eduit \`a $\{ 0 \}\,$ on a $\sF_P = \F$ et les faisceaux
(\ref{notfcan}) sont
constitu\'es de champs de vecteurs transversalement formels sur
$\M_\F$. Nous les notons
\begin{equation}\label{notfcano}
\hB_{\wF}\,, \qquad \hX_\wF\, ,  \qquad \hT_\wF\,,
\end{equation}
\noindent et d\'esignons par $\cal N_0$ l'un quelconque de ces faisceaux.
Pour une d\'eformation
\hsl-\'equisinguli\`ere
$\sF_P\,$ avec
$P$ quelconque, la restriction
$X|\M_\F$ d'une section $X$ d'un des faisceaux $\cal N_P$ de
(\ref{notfcan}) \`a
$\M_\F$ est une section du faisceau $\cal N_0$ correspondant.
La trivialit\'e locale de $\wF_P$
implique que, pour tout ouvert $W$ contenu dans un ouvert $U$ du
recouvrement distingu\'e $\U\,$, les
applications
$\cal N_P(W) \fle \cal N_0(W)\,$, $X \efle X|\M_\F$ sont surjectives.
Ainsi,
$(t)$ d\'esignant l'id\'eal maximal de $\O_P\,$, les morphismes de faisceaux
et les morphismes de
$\O_P$-modules
\begin{equation}\label{isotens}
\cal N_P \otimes_{\O_P} \O_P / (t) \fle \cal N_0 \, , \qquad
\cal N_P(W) \otimes_{\O_P} \O_P / (t) \fle \cal N_0(W) \, ,
\end{equation}
$W \subset U \in \U\,$, donn\'es par $X \otimes \stackrel{\cdot}{f(t)} \,
\efle
\,f(0)\cdot \left. X
\right|\M_\F\,$  sont des isomorphismes.

\begin{prop}\label{tens} Les applications canoniques
$$H^1 \left( \U \, ; \, \cal N_P\otimes_{\O_P} \O_P /(t) \right) \fle H^1
\left( \U \, ; \,
\cal N_P
\right)
\otimes_{\O_P} {\O_P}/{(t)}$$ sont des isomorphismes.
\end{prop}
\begin{proof} Comme $\U$ est sans intersection non-triviale 3 \`a 3, on a
la suite exacte
$$\prod_{U \in \U} \cal N_P(U)  \fle \prod_{\stackrel{\scriptstyle U, V
\in \U}{U \cap V
\not= \emptyset}} \cal N_P (U \cap V)
\fle H^1 \left( \U \, ; \, \cal N_P \right) \fle 0 \, .$$ La conclusion
r\'esulte de l'exactitude \`a
droite de $\otimes_{\O_P} \O_P / (t)$  et des isomorphismes
(\ref{isotens}).
\end{proof}

\noindent Au chapitre suivant nous \'etudierons en d\'etail les faisceaux
(\ref{notfcano}) et leurs
cohomologies.

\subsection{Vitesses de d\'eformation}
\addcontentsline{toc}{section}{\hspace{0,8em} {}\thesubsection .  Vitesses de d\'eformation}
Soit $\cZ$ un \'el\'ement du $\O_P$-module $\X_P$ des germes de
champs de vecteurs sur
$P := \left( \C^p , \, 0 \right)\,$ et $\sF_P\,$ une d\'eformation \hsl
-\'equisinguli\`ere de
$\F\,$. Cette d\'eformation  \'etant triviale sur des voisinages dans
$\M_{\sF_P}$ de chaque ouvert $U \in \U\,$, le champ $\cZ$ se rel\`eve
suivant $\pi_{\sF_P}$ sur ces
voisinages, en des champs  basiques $\cZ_U\,\in \hB_{\wF_P} ( U )$, i.e.
$T \pi_{\sF_P} \cdot \cZ_U
= \cZ \circ
\pi_{\sF_P}  \,$. Les diff\'erences
$\cZ_{UV} := \cZ_V - \cZ_U\,$, $U \cap V \not= \emptyset\,$, sont des
champs basiques verticaux.
Le cocycle
$\left( \cZ_{UV} \right) \in  Z^1 \left( \U \, ; \, \hB^{\, v}_{\wF_P}
\right)$ d\'epend du choix
des rel\`evements, mais sa classe de cohomologie n'en d\'epend pas.

\begin{defi} Nous appelons {\it application de Kodaira-Spencer de $\sF_P$}
l'application
$$\left[ \DD{\sF_P}{t\phantom{\sF}} \right] : \X_P \fle H^1 \left( \U \, ;
\, \hB^{\, v}_{\wF_P}
\right) \, , \quad
\cZ \efle \left[ \cZ_{UV} \right]\, .$$
\end{defi}
\noindent Un calcul direct montre que si $\left[ \Phi_{UV , \, t}Ê\right]
\in  H^1 \left(
\U , \, \widehat{Aut}Ê\left(Ê\sF_P^{cst}
\right)
\right)\,$ est la classe du cocycle associ\'e \`a $\sF_P$ par (\ref{intco}),
on a l'\'egalit\'e
\begin{equation}\label{cokodsp}
\left[ \cZ_{UV} \right] = \left[
\DD{\Phi_{UV , \, t}}{t \phantom{\Phi_{UV}}}Ê\circ {\Phi^{-1}_{UV , \, t}}
\cdot \cZ
\right]\,.
\end{equation}
\noindent Par construction la classe $\left[ \cZ_{UV} \right]$ est
exactement l'obstruction \`a
l'existence d'un rel\`evement global de
$\cZ$ en un champ basique. Plus pr\'ecis\'ement le sous-module
$$\cal K := ker \left( \left[ \DD{\sF_P}{t\phantom{\sF}} \right] \right)
\; \subset \; \X_P$$ est
exactement constitu\'e des champs de vecteurs sur $P$ qui poss\`edent un relev\'e
global
$\widetilde{\cZ} \in H^0 \left( \D_\F ; \, \hB_{\wF_P} \right)$. On en
d\'eduit im\'ediatement~:

\begin{prop} Une d\'eformation \hsl -\'equisinguli\`ere $\sF_P$ est formellement
triviale si et
seulement si $\left[ \DD{\sF_P}{t\phantom{\sF}}
\right]$ est l'application identiquement nulle.
\end{prop}

Nous allons raffiner cette proposition et obtenir la trivialit\'e de $\sF_P$
suivant les feuilles
d'un feuilletage de $P\,$. Remarquons d'abord que
$\cal K$ est un sous-module {\it involutif} de $\X_P\,$, i.e. stable par
crochet de Lie. En effet
le crochet de Lie commute aux images directes des champs de vecteurs
(lorsqu'elles sont d\'efinies)
et le crochet de deux champs basiques est aussi un champ basique. En
g\'en\'eral $\cal K $ n'est pas
un sous-module libre de $\X_PÊ\,$ et ses g\'en\'erateurs peuvent \^etre
singuliers \`a l'origine.
Consid\'erons le sous-espace vectoriel $\cal K(0)$ de l'espace tangent
$T_0P\,$ de
$P$ \`a l'origine, constitu\'e des valeurs
\footnote{ Notons qu'en g\'en\'eral $\cal K(0) \not= \cal K \otimes _{\O_P}
\O_P / (t)\,$, mais
$\cal K(0)$ est l'image de
$\cal K \otimes _{\O_P} \O_P / (t) \fle \X_P \otimes _{\O_P} \O_P / (t)  =
T_0P\,$. } \`a l'origine
des \'el\'ements de $\cal K\,$. Un r\'esultat de D. Cerveau \cite{C} donne
l'existence  d'un sous-module
libre involutif de
$\cal K\,$, non-unique, qui poss\`ede une base constitu\'ee de champs de
vecteurs r\'eguliers dont les
valeurs \`a l'origine forment une base de
$\cal K(0)\,$. Il d\'efinit un germe de feuilletage holomorphe r\'egulier
$\cal H\,$ sur $P\,$.

\begin{lemm}\label{retrtriv}
Soit $P' \subset P$ un germe de sous-vari\'et\'e strictement transverse \`a
$\cal K(0)\,$, i.e. $T_0 P = T_0 P' \oplus \cal K(0)\,$ et soit $\iota  :
P'
\ifle P$ l'application d'inclusion. Notons $\sF_{P'}$ la d\'eformation $
\iota^{\ast} \sF_P$ obtenue
en restreignant les param\`etres (\ref{def.equired.}) \`a
$P'\,$. Alors il existe un germe de submersion
$\l : P \fle P'$ telle que $\sF_P$ soit formellement conjugu\'e \`a
$\l^{\ast}\sF_{P'}\,$.
\end{lemm}

\begin{proof} Le feuilletage $\cal H$ pr\'ec\'edemment d\'efini peut \^etre
engendr\'e par des champs
$\cZ_j\,$, $j = 1 , \ldots , r\,$, $r := dim \cal K (0) \,$, qui commutent
entre eux. Ces champs
admettent, sur un voisinage de
$\D_{\sF_P}\,$ dans
$\M_{\sF_P}\,$, des rel\`evements globaux en des champs basiques
$\widetilde{\cZ_j}$. L'int\'egration
successive des $\cZ_j$ avec conditions initiales sur $P'\,$ donne une
rectification
$P \iso P' \times \C^r\,$ de $\cal H\,$. Elle induit une submersion $\l\,$
de $P$ sur $P'\,$ qui
consiste \`a suivre les feuilles de $\cal H\,$. L'int\'egration successive des
relev\'es
$\widetilde{\cZ_j}$  donne alors la conjugaison cherch\'ee.
\end{proof}

\subsection{Codimension formelle de $\F\,$}
\addcontentsline{toc}{section}{\hspace{0,8em} {}\thesubsection .  Codimension formelle de $\F\,$}
En d\'erivant un
cocycle $\left[ \Phi_{UV, \, t} \right] \in H^1 \left( \U , \,
\widehat{Aut}Ê\left(Ê\sF^{cst}_P \right)\right)$ associ\'e par (\ref{idco})
\`a une d\'eformation
\hsl -\'equisinguli\`ere $\sF_P$ de $\F$  on obtient des cocycles
$$\left( \left. \DD{\Phi_{UV, \, t}}{t_j\phantom{\Phi_{U,}}}\right |_{t =
0}Ê\right)
 \; \in \; Z^1 \left( \U ; \, \hB_\wF \right) , \qquad j = 1 , \ldots , p$$
\noindent \`a valeurs dans le faisceau $\hB_\wF$ des champs basiques de
$\wF\,$. Leurs classes de
cohomologies not\'ees
$$\left[ \DD{\sF_P}{t_j\phantom{\sF}} \right]_{t = 0}\;
 \in \; H^1 \left( \U ; \, \hB_\wF \right)
 , \quad j = 1 , \ldots , p$$ ne d\'ependent que de $\sF_P\,$. Elles
s'interpr\`etent comme les
"vitesses initiales de d\'eformation" de $\sF_P\,$. Gr\^ace \`a (\ref{cokodsp})
et \`a l'identification
$H^1 \left( \U \, ; \, \hB_\wF \right) \simeq H^1 \left( \U \, ; \,
\hB_{\wF_P}
\right) \otimes_{\O_P} \O_P /(t)$ donn\'ee par (\ref{tens}), l'application
lin\'eaire
\begin{equation}\label{derin}
\left[ \DD{\sF_P}{t\phantom{\sF}} \right]_{t = 0} \, : \, T_0 P \fle H^1
\left( \U \, ; \,
\hB_\wF \right) \, , \quad
\sum_{j = 1}^p \left. \a_j \dd{t_j} \right|_{t=0} \efle \sum_{j = 1}^p
\a_j \left[
\DD{\sF_P}{t_j\phantom{\sF}} \right]_{t = 0}
\end{equation}
\noindent peut \^etre vue comme le tensoris\'e de l'application de
Kodaira-Spencer de $\sF_P\,$, c'est
\`a dire~:
\begin{equation}\label{tenskss}
\left[ \DD{\sF_P}{t\phantom{\sF}} \right]_{t = 0} = \left[
\DD{\sF_P}{t\phantom{\sF}}
\right]\otimes_{\O_P} \O_P /(t)\,.
\end{equation}
\noindent Remarquons enfin que pour tout $k \in \N$ et tout germe
d'application holomorphe
$\l : \left( \C^k , \, 0 \right) \fle P\,$, on a~:
\begin{equation}\label{compder}
\left[ \DD{\l^\ast \sF_P}{s\phantom{\l^\ast \sF}} \right]_{s = 0} =
\left[ \DD{\sF_P}{t\phantom{\sF}} \right]_{t = 0} \circ \left. \DD{\l}{s}
\right|_{s = 0}\,.
\end{equation}

\begin{defi} On appelle {\it d\'eformation \hSL -infinit\'esimale} de $\F$
tout \'el\'ement de $H^1
\left( \U \, ; \, \hB_\wF \right)\,$.
\end{defi}

\noindent On a la propri\'et\'e de r\'ealisation suivante~:

\begin{enonce}{Th\'eor\`eme de r\'ealisation}\label{readefin} Donn\'ees des d\'eformations \hsl
-infinit\'esimales
$v_1
\, , \ldots ,v_q  \in \; H^1 \left( \U \, ; \, \hB_\wF \right)$, il existe
une d\'eformation
\hsl -\'equisinguli\`ere
$\sF_Q$ de $\F$ de param\`etres $Q := \left( \C^q , \, 0 \right)$ telle que
$\left[
\DD{\sF_Q}{t_j\phantom{\sF}} \right]_{t = 0} = v_j\;$,
$j = 1 ,
\ldots , q\,$.
\end{enonce}

\begin{proof} Soient $\left( X_{UV}^j\right) \in Z^1 \left( \U \, ; \,
\hB_\wF \right)$  tels que
$v_j = \left[ X_{UV}^jÊ\right]\,$, $j = 1 , \ldots  , q\,$. Les compos\'es
$\Phi_{UV , \, t} := \Phi_{UV,\, t_1}^1\circ \cdots \circ \Phi_{UV,\,
t_q}^q$ des flots
$\Phi_{UV, \, s}^j$ des $X_{UV}^j$ d\'efinissent un cocycle \`a valeurs dans
$\widehat{Aut}Ê\left(Ê\sF^{cst}_Q\right)\,$. Il d\'efinit par la bijection
(\ref{cosla}) une
d\'eformation \hsl -\'equisinguli\`ere qui convient.
\end{proof}

\begin{defi}\label{stffdef} Un feuilletage formel \`a l'origine de $\C^2$
d\'efini par une 1-forme
diff\'erentielle formelle $\w$  est dit {\it de type formel fini}, en abr\'eg\'e
t.f.f., si $\w$ est
non-dicritique, \`a singularit\'e isol\'ee et
$H^1\left(\U\, ;\,\hB_\wF\right)\,$ est un
$\C$-espace vectoriel de dimension finie. Sa dimension, not\'ee
$\widehat{\b}(\F)$ ou encore
$\widehat{\b}(\w)$ est appel\'ee {\it codimension formelle de $\F\,$, ou de
$\w\,$}.
\end{defi}

Nous caract\'eriserons plus loin les feuilletages t.f.f. par des conditions
combinatoires sur
l'arbre dual de r\'eduction de $\F\,$ que l'on a enrichi d'un nombre fini de
donn\'ees suppl\'ementaires.
Auparavant comparons $\hb \left( \F \right)$ aux nombres
\begin{equation}\label{cotr}
\hd \left( \F \right) := dim_\C \, H^1 \left( \U ; \, \hX_\wF \right)
 \quad \mbox{et} \quad
\htau \left( \F \right) :=  dim_\C \, H^1 \left( \U ; \, \hT_\wF \right)\,.
\end{equation}
\noindent Une version formelle d'un r\'esultat de \cite{Minv} dont la
d\'emonstration se transcrit
lit\-t\'e\-ra\-le\-ment, donne la formule
\begin{equation}\label{fordepinf}
\hd\,(\F) \; = \; \sum_{c} \frac{(\nu_c - 1) \, (\nu_c - 2) }{2}\, ,
\end{equation}
\noindent o\`u dans cette somme, $c$ d\'ecrit l'ensemble
$\sqcup_{j = 0}^h \S_\F^j$  de tous les points singuliers (y compris $0
\in \C^2$) qui
apparaissent  dans le processus de r\'eduction de $\F\,$, cf. (\ref{ard}) et
o\`u $\nu_c$ d\'esigne la
multiplicit\'e alg\'ebrique de transform\'e strict du feuilletage au point
$c\,$. Signalons que cette
formule est valable en toute g\'en\'eralit\'e, m\^eme si
$\F$ est dicritique. \\

\indent Rappelons \cite{C-M} qu'une {\it int\'egrale premi\`ere formelle},
resp. un {\it facteur
int\'egrant formel} de $\F\,$, ou de
$\w\,$, est un \'el\'e\-ment
$f$ de $\hO_{\C^2 ,\, 0}$ tel que
\begin{equation}\label{factint}
df \wedge \w = 0\,, \quad \hbox{resp.} \quad d\,(\w / f) = 0\, .
\end{equation}
Le quotient de deux facteurs int\'egrants est une
{\it int\'egrale premi\`ere m\'eromorphe formelle} de $\w\,$, c'est \`a dire un
\'el\'ement $f$ du corps des
fractions de $\hO_{\C^2 ,\, 0}$ satisfaisant $ df \wedge \w = 0 \,$. Enfin
si $X$ est un champ de
vecteurs formel \`a l'origine de $\C^2$ qui est {\it basique pour $\F\,$},
ie. $L_X
\w \wedge \w = 0\,$, et non-tangent, alors $\w\cdot X\,$ est un facteur
int\'egrant de $\F\,$.\\

Notons $\widehat{Int} \,(\w)\,$ le sous-espace vectoriel de $\hO_{\C^2 ,
\, 0}$ des
facteurs int\'egrants formels qui s'annulent sur les s\'eparatrices formelles
$\hsep\,(\w)\,$ de $\w\,$,
cf. (\ref{minmult}), et
$\hB(\w )
\subset
\hX_{\C^2 ,
\, 0}\,$ le sous-espace vectoriel du module  des champs formels \`a
l'origine de $\C^2$ qui sont
basiques pour
$\w\,$. On a la

\begin{prop}\label{splitcod}
Supposons $\F$ non-dicritique. Alors
$\hb (\F)$ est fini si et seulement si $\htau (\F)$ est fini. De plus on a
l'in\'egalit\'e~:
$$
\hd \left( \F \right) \; + \; \htau \left( \F \right)\;  -
\; \he \left( \F \right) \; \leq \;
\hb \left( \F \right)\;  \leq \;
\hd \left( \F \right) \; + \; \htau \left( \F \right)\;,
$$
avec $\he \left( \F \right) := dim_\C
\left( \left.
\widehat{Int}\,(\w)\; \right/ \left( \w \cdot \hB (\w ) \right) \right)
\,$. Lorsque $\F$
est de deuxi\`eme esp\`ece, l'\'egalit\'e
$$
\hb \left( \F \right) = \hd \left( \F \right)  +  \htau \left( \F \right)
-
 \he \left( \F \right) \, $$
est r\'ealis\'ee et $\he (\F ) = 0$ ou $1\,$.
\end{prop}

\begin{proof} La suite exacte longue (\ref{dlousski}) donne
\begin{equation}\label{sel}
0 \fle N \fle H^1 \left( \U ; \, \hX_\wF \right)  \fle H^1 \left( \U ; \,
\hB_\wF
\right)
\fle H^1 \left( \U ; \, \hT_\wF \right) \fle 0
\end{equation}
avec $N := coker \left( H^0 \left( \U ; \,
\hB_\wF \right)  \fle H^0 \left( \U ; \, \hT_\wF \right)  \right) \,$.
Remarquons que $ H^0
\left( \U ; \, \hB_\wF \right)$ s'iden\-ti\-fie \`a $\hB (\w )\,$. En effet
tout champ de vecteurs
holomorphe global sur un voisinage de $\D_\F$ dans $\M_\F$ est l'image
r\'eciproque d'un champ de
vecteurs de
$(\C^2 , \, 0 ) \,$; la commutation des op\'erations de compl\'etion et
d'image directe
\cite[page 269]{B-S} \'etend ce fait aux champs formels. De m\^eme
l'\'evaluation par $ E_\F^{\ast}
\,(\w)$ induit une injection
\begin{equation}\label{ffint}
0 \fle  H^0 \left( \U ; \, \hT_\wF \right) \fle \widehat{Int}\,(\w)\,.
\end{equation}
La formule d'Euler de
la suite (\ref{sel}) donne l'in\'egalit\'e cherch\'ee.
\\

Supposons maintenant $\F$ de deuxi\`eme esp\`ece. Il suffit de montrer que
l'injection ci-dessus est un
isomorphisme. Pour cela donnons nous un facteur int\'egrant $f \in
\widehat{Int}\,(\w)$ et montrons l'existence en tout point $m \in \D_\F$
d'un germe $X_m \in
\hB_{\wF , \, m}$ tel que $  E_\F^{\ast} \,(\w) \cdot X_m = f \circ
E_\F\,$.\\

Lorsque $\wF$ est r\'egulier en $m\,$, fixons des coordonn\'ees
transversalement formelles $z_1,
\,z_2\,$ appropri\'ees o\`u
$\, E_\F^{\ast} \,(\w) = u(z_1, z_2) \,{z_2}^\a\, d z_2 \,$, avec $\,u(0,
0) \not= 0 \,$ et $\a \in
\N\,$. Les deux facteurs int\'egrants $ \,u(z_1, z_2)\, {z_2}^\a \, $ et $
\, f \circ E_\F \, $ de $\,
E_\F^{\ast} \,(\w) \, $ diff\`erent d'un facteur multiplicatif qui est une
in\-t\'e\-gra\-le premi\`ere
formelle m\'eromorphe, ce qui donne~:
$$f \circ E_\F\, = \, u(z_1, z_2)\, {z_2}^\b\,l\,(z_2) \, , \quad
\mbox{avec} \quad l\,(z_2) \, \in
\,
\C[[z_2]]\, , \quad l\,(0) \not= 0
\, , \quad \b \, \in \, \Z
\, .$$
Par hypoth\`ese $f$ est une \'equation (non-n\'ecessairement r\'eduite) de
$\hsep\,(\w)\,$. La
caract\'erisation (\ref{cardeuxesp}) des feuilletages de deuxi\`eme esp\`ece
donne
$\b \geq \a + 1\,$. Le champ $X_m := l\,(z_2) \,{z_2}^{\b - \a}\,
\dd z_2\,$ est basique et convient. \\

Lorsque $m$ est un point singulier de $\wF\,$ il existe d'apr\`es
(\ref{forntf}) des coordonn\'ees
transversalement formelles $z_1 ,
\, z_2\, $ dans lesquelles $\D_\F$ est d\'efini par un mon\^ome et $
E_\F^{\ast} \,(\w)$ s'\'ecrit $ \,
u(z_1, z_2)\, z_1^\a z_2^\b \,
\ww\,$,
$u(0, 0) \not= 0\,$, o\`u
$\ww$ est l'une des formes normales 1. \`a 4. de (\ref{subs.red.des.sing.}).
On voit facilement sur
ces expressions qu'il existe toujours un germe $Z_m \in \hB_{\wF , \, m}$
non-tangent au
feuilletage $\wF\,$. Cela sera explicit\'e en d\'etail au paragraphe suivant.
Ainsi
$h := u(z_1, z_2)\, z_1^\a z_2^\b \, \ww (Z_m)$ est un facteur int\'egrant
de $ E_\F^{\ast}
\,(\w)\,$. Si $\wF$ ne poss\`ede pas au point
$m$ de germe d'int\'egrale premi\`ere formelle non-constante, on a $f \circ
E_\F = c \, h\,$, $ c \in
\C\,$ et le champ $X_m := c\, Z_m$ convient.  Sinon $\ww$ s'\'ecrit
$q\,z_1\,dz_2 + p\,z_2\,dz_1\,$
avec $p, q  \in \N^{\ast}\,$ et $(p,q)=1\,$. Alors
$\C[[z_1^p z_2^q]]$ est l'anneau des germes en $m$ d'int\'egrales premi\`eres
for\-mel\-les cf.
\cite{M-M}. Visiblement le champ radial
$R_m := z_1
\dd{z_1} + z_2
\dd{z_2}\,$ est basique et non-tangent. Il vient $f \circ E_\F = l\,(z_1^p
z_2^q)\, h(z_1,\,
z_2)\,$, o\`u, comme pr\'ec\'edemment, $l$ est une s\'erie formelle d'une variable
car
$\F$ est de deuxi\`eme esp\`ece. Le champ $X_m := l\,(z_1^p z_2^q)\, h(z_1,\,
z_2)\, R_m$ convient.
\\

Il reste \`a prouver que $\he (\F)$ vaut $0$ ou $1\,$. Visiblement que
$\widehat{Int} \,(\w)\,$ est
un module libre de rang 1 sur l'anneau $\hO_{\F}$ des germes \`a l'origine
de $\C^2\,$, d'int\'egrales
premi\`eres formelles de $\F\,$. De m\^eme $\w \cdot \hB (\w )$ est un sous
$\hO_{\F}$-module. Ainsi,
le r\'esultat est trivial quand $\hO_{\F} = \C\,$.\\

Lorsque $w$ poss\`ede une int\'egrale premi\`ere
formelle  non-constante $f$ dont la d\'ecom\-position en facteurs
irr\'eductibles est ${f_1}^{n_1} \cdots
{f_p}^{n_p}\,$, elle s'\'ecrit
$$\w = u\, f_1 \cdots f_p \left( n_1 \frac{d f_1}{f_1} + \cdots + n_p
\frac{d f_p}{f_p}\right)
\,\quad u \in \hO_{\C^2 , 0} \, , \quad u(0,0) \not= O $$
\noindent et $g := f_1 \cdots f_p$ est une \'equation r\'eduite de
$\widehat{Sep}(\w)\,$.
Ainsi $\widehat{Int} \,(\w)\, = \CÊ[[f
]] \, u g
\,$.
On voit
facilement que $g$ appartient \`a
l'id\'eal $I(\w)$ engendr\'e par les coefficients de $\w\,$ si et seulement si
$f$ est
{\it quasi-homog\`ene}, i.e. $f \in I(df)\,$. Lorsque c'est le cas on
obtient un champ $X$ tel que
$\w \cdot X = u g$ et $\he (\F ) = 0 \,$. Si $f$ n'est pas quasi-homog\`ene,
$u g \notin \w \cdot
\hB ( \w )$ et $\he (\F ) \not= 0\,$. Mais le th\'eor\`eme de Brian\c con
\cite{BS} assure que
l'on peut r\'esoudre l'\'equation $df \cdot X = f^2 \,$. Il vient $\w \cdot X
= u g f\,$ et
$\w \cdot \hB ( \w ) = (f^2) \cdot u g \,$. D'o\`u $\he (\F ) = 1\,$.
\end{proof}

\subsection{D\'eformations \hSL-verselles}\label{defslver}
\addcontentsline{toc}{section}{\hspace{0,8em} {}\thesubsection .  D\'eformations \hSL-verselles}
Dans cette partie nous
montrons, sous une condition g\'en\'erique sur le  feuilletage -$\F$ bon- un th\'eor\`eme de
versalit\'e pour les d\'eformations \hsl-\'equisinguli\`eres.

\begin{defi}\label{fbon} Nous dirons que le feuilletage $\F$ est {\it bon} s'il existe,
soit une  singularit\'e $m$ du transform\'e strict
$\wF$ de $\F$, soit une composante irr\'eductible $D$ du diviseur $\D_\F\,$,  telle que au
voisinage de $m$, ou le long de $D^{\ast} := \left( D - Sing(\wF )\cap D \right)$, toute
int\'egrale premi\`ere transversalement formelle est constante~: $\hO_{\F, m} = \C\,$, ou
$\hO_\F  ( D^{\ast} ) = \C\,$.
\end{defi}
\noindent Remarquons que la condition $\hO_\F (D^{\ast}) = \C\,$ \'equivaut
\`a la non-finitude du groupe d'holonomie $H_D$ de la composante $D\,$, d\'efini en
(\ref{flowbox}).\\

\indent Fixons une d\'eformation $\sF_Q$ de $\F\,$ de param\`etres
$Q :=\left(\C^q ,\, 0\right)\,$,
$q
\in
\N\,$, transversalement formelle le long de $0 \times Q\,$ et
\hsl-\'equisinguli\`ere,  donn\'ee par une forme diff\'erentielle formelle
$$\uw := a(x, y;\, u)\,dx + b(x, y;\, u)\, dy\, .$$

\begin{defi} Nous disons que $\sF_Q$ est {\it \hsl-verselle} si pour tout $p \in
\N$ et toute d\'eformation $\sF_P$ de $\F\,$ de param\`etres $P := \left( \C^p , \, 0
\right)\,$, transversalement formelle le long de
$0
\times P\,$ et \hsl-\'equisinguli\`ere, il existe un germe d'application  holomorphe
$\l : P \fle Q$ tel que $\sF_P$ soit formellement conjugu\'ee \`a $\l^\ast
\sF_Q\,$. Lorsque $Q =
\left\{ 0 \right\}$ et $\sF_Q$ est \hsl-verselle, nous dirons que $\F$ est
{\it \hsl-stable}.
\end{defi}

\noindent En d'autres termes, $\sF_Q$ est \hsl-verselle si toute  d\'eformation formelle
\hsl-\'equisinguli\`ere de
$\F$ est d\'efinie par une 1-forme qui s'\'ecrit $\Phi^{\ast} \left( a(x, y;\,
\l(t) )\,dx + B(x, y;\,
\l (t) )\, dyÊ\right)\,$, o\`u $\Phi (x, y; t) = \left( \underline{\Phi}(x,  y ;\, t) ;\, t
\right)\,$ est un diff\'eomorphisme de $\left( \C^2 \times \C^p , \, 0
\right)$ transversalement formel le long de $0 \times \C^p \,$. En particulier la
\hsl-stabilit\'e  signifie que tout d\'eformation
\hsl-\'equi\-sin\-gu\-li\`ere de $\F$ est formellement triviale.

\begin{enonce}{Th\'eor\`eme de versalit\'e}\label{thver}  Supposons $\F$ bon. Lorsque $Q \not= \{ 0 \}\,$, la
d\'efor\-mation $\sF_Q$  est \hsl-verselle si et seulement si les "vitesses initiales de
d\'eformation"
\begin{equation}\label{vid}
\left[ \DD{\sF_Q}{u_1\phantom{\sF}} \right]_{u = 0} \, , \ldots , \left[
\DD{\sF_Q}{u_q\phantom{\sF}} \right]_{u = 0}\, , \qquad u = (u_1, \ldots ,  u_q)\, ,
\end{equation}
\noindent forment un syst\`eme de g\'en\'erateurs de $H^1 \left( \U ; \, \hB_\wF
\right)\,$. D'autre part
$\F$ est \hsl-stable si et seulement si
$H^1 \left( \U ; \, \hB_\wF \right) = 0 \,$.
\end{enonce}

\noindent Un ingr\'edient essentiel de la d\'emonstration de ce r\'esultat est le  th\'eor\`eme
ci-dessous. Nous le d\'emontrerons au paragraphe (\ref{primdir})
\begin{enonce}{Th\'eor\`eme de pr\'eparation}\label{enimdir} Supposons $\F$ bon et t.f.f.. alors, pour toute
d\'eformation
\hsl-\'equisinguli\`ere $\sF_P$ de $\F\,$, le
$\O_P$-module
$H^1\left( \U ; \hB^v_{\wF_P} \right)$ est de type fini. Si de plus $\F$ ne poss\`ede pas
de facteur int\'egrant formel, alors $H^1\left( \U ; \hB^v_{\wF_P} \right)$ est un
$\O_P$-module libre de rang $\hb(\F)$.
\end{enonce}

\begin{preuvede}{du th\'eor\`eme de versalit\'e.} Supposons que les cocycles  (\ref{vid})
engendrent
$H^1 \left( \U ; \, \hB_\wF \right)\,$. Soit $\sF_P$ une d\'eformation
\hsl-\'equisinguli\`ere de
$\F\,$. Donnons-nous gr\^ace au lemme (\ref{lemcroix}) une d\'eformation
\hsl-\'equisinguli\`ere
$\sF_{P \times Q}$ de param\`etres $P \times Q\,$, qui est conjugu\'ee \`a
$\sF_P$ lorsqu'on restreint les param\`etres \`a $P\times 0\,$ et qui est conjugu\'ee \`a
$\sF_Q$ lorsqu'on  restreint les param\`etres \`a
$0 \times Q\,$.  Pour all\'eger le texte notons~:
$$
\Delta \; := \; \left[ \DD{\sF_{P\times Q}}{(t; u)} \right] \, : \,
\X_{P\times Q} \fle H^1 \left(
\U \, ; \,
\hB^v_{\wF_{P \times Q}} \right) \, \quad \mbox{et}$$

$$\Delta_0 \; := \; \left[ \DD{\sF_{P\times Q}}{(t; u)} \right]_{t = 0\, ,  u = 0} \, : \,
T_{(0;
\,0)} P\times Q \fle H^1 \left( \U \, ; \,
\hB_\wF
\right)\, .\phantom{et}
$$
\noindent Par (\ref{tenskss}), $\Delta_0$ est surjective et son noyau est transverse \`a
$0 \times T_0 Q\,$. On d\'eduit de (\ref{tenskss}) l'\'egalit\'e~:
$\,coker\left(
\Delta \right)\,
\otimes_{\O_{P \times Q}} \left(\left. \O_{P \times Q } \right/ \!\goth m
\right)\; =
\; 0\,$, o\`u $\goth m$ d\'esigne l'id\'eal maximal de $ \O_{P \times Q }\,$.  D'apr\`es
le th\'eor\`eme  (\ref{enimdir}) de pr\'eparation
$\,H^1
\left(\U\, ; \,\hB^v_{\wF_{P\times Q}}\right)\,$ - et \`a fortiori
$\,coker\left(
\Delta \right)\,$, est un
$\O_{P \times Q}$ module de type fini. Le lemme de Nakayama appliqu\'e \`a
$coker \left( \Delta \right)$ donne la surjection de
$\left[ \DD{\sF_{P\times Q}}{(t; u)} \right]$\,. En tensorisant la suite exacte
$$ 0 \fle \cal K \stackrel{\a}{\fle} \X_{P \times Q}
\stackrel{\Delta}{\fle} H^1
\left( \U \, ; \, \hB^v_{\wF_{P \times Q}} \right)
\fle 0\,
$$ par $\otimes_{\O_{P \times Q}} \left(\left. \O_{P \times Q} \right/
\!\goth m\right)\,$, on obtient
$ \cal K ( 0 ) \, = \, ker (\Delta_0 )\,$,
\noindent avec~:
$$\, \cal K(0)\, := \,
\left\{\left. \cal Z( 0 )  \right| \cal Z \in \cal K \right\} \, = \, Im
\left(
\a \otimes_{\O_{P \times Q}}
\left(\left. \O_{P \times Q} \right/ \!\goth m\right) \right)\,.$$   Le lemme
(\ref{retrtriv}) donne une submersion $\L : P \times Q \fle Q$  telle que
$\sF_{P \times Q}$ est formellement conjugu\'e \`a $\L^\ast \sF_Q\,$. On  obtient la
factorisation cherch\'ee en restreignant $\L$ \`a
$P \times 0\,$.\\

\indent R\'eciproquement, supposons que la d\'eformation $\sF_Q$ est
\hsl-verselle et donnons-nous une d\'eformation infinit\'esimale $v \in H^1 \left( \U \, ;
\,
\hB_\wF \right)\,$. D'apr\`es le th\'eor\`eme de r\'ealisation (\ref{readefin}) il  existe
une d\'eformation
\hsl-\'equisinguli\`ere $\sF_\C$ de base $\left( \C , \, 0 \right)$ telle que
$v =
\left[ \DD{\sF_\C}{t\phantom{\sF}} \right]_{t = 0}\,$. a \'equivalence pr\`es  on peut
poser
$\sF_\C = \l^\ast \sF_Q\,$, avec $\l := (\l_1, \ldots , \l_q ) : \left( \C  , \, 0 \right)
\fle Q\,$ holomorphes. Par (\ref{compder}) et (\ref{derin}) on obtient~:
$$v\, =\,  \sum_{j = 1}^q\, \frac{d \l_j}{dt\phantom{{}_j}} (0) \cdot
\left[
\DD{\sF_Q}{u_j\phantom{\sF}} \right]_{u = 0}\, ,$$
\noindent d'o\`u la conclusion.\\

\indent Remarquons enfin que ces d\'emonstrations s'adaptent au cas $Q = \{  0 \}\,$.
Elles prouvent alors que $\F$ est stable si et seulement si $ H^1 \left( \U \, ;
\,\hB_\wF \right) = 0\,$.
\end{preuvede}

Le th\'eor\`eme de versalit\'e, joint au th\'eor\`eme de r\'ealisation  (\ref{readefin})
donne imm\'edia\-te\-ment~:

\begin{theo}\label{exist} Tout feuilletage $\F$ bon et de type formel fini poss\`ede une
d\'eformation
\hsl-verselle $\sF_P$, telle que l'application
\begin{equation}\label{kstens}
\left[ \DD{\sF_P}{t\phantom{\sF}} \right]_{t = 0}\, : \,T_0P
\fle H^1 \left( \U \, ; \,\hB_\wF \right)\,
\end{equation} d\'efinie en (\ref{derin}) est un isomorphisme.
\end{theo}

\section{Caract\'erisation des singularit\'es t.f.f.}\label{carsintff}
\subsection{ Description du faisceau $\hT_\wF\,$}\label{desftrfor}
\addcontentsline{toc}{section}{\hspace{0,8em} {}\thesubsection .  Description du faisceau $\hT_\wF\,$}
Rappelons d'abord quelques r\'esultats de classification formelle des  sous-groupes de
type fini $H$ du groupe
$\widehat{Diff}\left(Ê\C , \, 0 \right)$ des diff\'eomorphismes formels de
$(\C , \, 0 )\,$, cf.
\cite{Paul} ou
\cite{Ilhy}. Nous notons $\hO_H \subset \C[[z]]$ l'anneau des s\'eries  formelles
invariantes par l'action \`a droite $g \star h := h \circ g^{-1}$ de $H$ sur $\hO_{\C,
\, 0} = \C [[z]]\,$. Les
\'el\'ements de
$\hO_H$ sont appel\'es {\it int\'egrales premi\`eres de $H\,$.} De m\^eme un champ  de
vecteurs formel qui s'annule \`a l'origine et qui est invariant par l'action (de
conjugaison) de
$H$ est appel\'e {\it sym\'etrie formelle de $H$}. L'espace de ces champs est not\'e
$\hT_H\,$. On  a~:
\begin{enumerate}
\item $\hT_H$ est soit nul, soit un $\hO_H$-module libre de rang 1;
\item $H$ est fini, $\# H =: p\,$, si et seulement si il existe
$f \in \hO_H$ non-constante. Dans ce cas, il existe une coordonn\'ee  formelle $\wz\,$
telle que $H$  est engendr\'e par la rotation $e^{2i\pi / p} \, {\wz}\,$. On a alors
n\'ecessairement~:
$\hO_H =
\C [[ {\wz}^p ]]\,$,
$\hT_H = \hO_H  \wz \,\dd \w  z\,$;
\item $H$ est commutatif si et seulement si son action laisse invariant un  champ de
vecteurs formel
$Z \in \C [[ z ]] z \, \dd z \,$ non-identiquement nul. Si de plus $ H $  est d'ordre
infini, alors
$\hT_H = \C Z\,$.
\end{enumerate}
\noindent Nous dirons qu'un couple $(f , Z) \in \hO_H \times \hT_H$ est  {\it normalis\'e}
si, avec la convention $\C [[ 1 ]] := \C\,$, il satisfait les relations de dualit\'e~:
\begin{equation}\label{counor}
\hO_H = \C [[ f ]]\, , \quad \hT_H = \C [[ f ]] \, Z \, , \quad
 \hbox{il existe $c \in \C \,$ tel que}\;\; df \cdot Z = c\, f \, .
\end{equation}
\noindent Dans chacun des cas 1., 2., 3. l'existence de tels couples est  claire.\\

\indent Notons $\hO_\wF$ le faisceau de base $\D_\F$ des germes, aux  points de $\D_\F$,
des int\'egrales premi\`eres transversalement formelles de $\wF\,$. Nous allons
expliciter, pour les ouverts $U$ du recouvrement distingu\'e
$\U\,$, les espaces
$\hT_\wF \,(U)$ de sections du faisceau des champs transverses  (\ref{sytr}).

\subsubsection{L'ouvert $U \in \U$ ne contient pas de singularit\'es de
$\wF
\,$}
Fixons un point ${m_0}$ sur $U\,$, ainsi qu'un germe de courbe
$T_{m_0}\,$ analytique ou formelle, lisse et transverse \`a
$U$ en
${m_0}\,$. Nous la munissons d'une coordonn\'ee formelle $z \in
\hO_{T_{m_0}}\,$, $z({m_0}) = 0\,$.   D'apr\`es l'\'ecriture (\ref{fca}) tout champ formel
$Z$ sur $T_{m_0}$  qui s'annule en ${m_0}$ s'\'etend en un germe de champ basique en
${m_0}\,$, d\'efinissant un unique \'el\'ement de $\hT_{\wF , \, {m_0}}\,$. De  plus tout
\'el\'ement de
$\hT_{\wF, {m_0}}$ est de ce type. ainsi
$\hT_{\wF, {m_0}} \simeq \hO_{T_{m_0}} \, z \dd z\,$. La restriction de
$\hT_\wF$ \`a
$U$ est un faisceau localement constant (\ref{prsec}) dont la monodromie  le long d'un
lacet $\gamma$ dans
$U$ d'origine ${m_0}$ s'exprime, avec l'identification pr\'ec\'edente, par $Z
\efle h_\gamma^{\ast}\,(Z)\,$, o\`u $h_\gamma$ est l'image de $[\g ]$ par la
repr\'esentation sur $T_{m_0}$
$$\H_U : \pi_1 ( U ; \, {m_0}) \fle \widehat{Diff}(T_{m_0})\,$$ de  l'holonomie de $\wF$
le long de
$U\,$. Notons $H_U := \H_U \left( \pi_1 ( U ; \, {m_0})Ê\right)$ le {\it  groupe
d'holonomie de
$U\,$}. Toute int\'egrale premi\`ere de $H_U$ se prolonge de mani\`ere unique
\cite{M-M} en une int\'egrale premi\`ere de
$\wF$ le long de $U\,$. On obtient finalement~:
$$\hO_\wF (U ) \simeq \hO_{H_U} \, , \qquad \hT_\wF (U) \simeq \hT_{H_U}\,  .$$ Trois
\'eventualit\'es se pr\'esentent, suivant la nature de
$H_U$.\\

\indent a. {\it $H_U$ est fini, $\# H_U =: p\,$}.  Choisissons pour $z$  une coordonn\'ee
formelle sur
$T_{m_0}$ dans laquelle $H_U$ est un groupe de rotation. alors
$z^p$ se prolonge en  une int\'egrale premi\`ere transversalement formelle
$f_U$ le long de $U\,$ et $\hO_\wF (U) = \C\, [[ f_U]]\,$. Le champ de  vecteurs formel $
Z := z \, \dd z\,$ sur $T_{m_0}\,$, invariant par
$H_U\,$, induit une section $ Z_U$ de $\hT_\wF$  sur
$U\,$ et
$$\hO_\wF (U) = \C\, [[f_U  ]] \simeq \C\, [[z^p  ]] \, ,Ê\quad
\hT_\wF (U) = \C\, [[f_U  ]] \cdot Z_U \simeq \C\, [[z^p  ]] \cdot \,z
\, \dd z \, .
$$

\indent B. {\it $H_U$ est commutatif et infini}. La sym\'etrie formelle  non-identiquement
nulle $ Z$ de $H_U\,$, unique \`a
constante multiplicative pr\`es, s'\'etend encore en  une section
$ Z_U \in \hT_\wF\,(U)\,$ et
$$\hO_\wF\,(U) = \C\,  \qquad \qquad \hT_\wF\,(U) = \C \cdot \, Z_U . $$

\indent C. {\it $H_U$ n'est pas commutatif}. alors
$$\hO_\wF\,(U) = \C  \qquad \mbox{et} \qquad \hT_\wF\,(U) = \{ 0 \} \, .$$  On posera $f_U
:= 1\,$, $ Z_U := 0\,$.

\subsubsection{L'ouvert $U \in \U$ est un voisinage d'une singularit\'e de
$\wF
\,$}Plus pr\'ecis\'ement supposons que $U$ est la trace sur $\D_\F$ d'un  petit
polydisque centr\'e en un point singulier $m\,$. D'apr\`es la proposition (\ref{prsec})
l'application canonique $\hT_\wF
\,(U) \fle \hT_{\wF , m}$ est un isomorphisme. D\'ecrivons $
\hT_{\wF , m}$ lorsqu'il existe en $m$ des coordonn\'ees transversalement  formelles
$z_1\,$, $z_2\,$,
$z_1(m) = z_2(m) = 0\,$, dans lesquelles $\wF$ est d\'efini par une des  1-formes normales
$\ww$ de (\ref{subs.red.des.sing.}). C'est toujours le cas, sauf peut-\^etre lorsque  la
singularit\'e est un selle-n\oe ud tangent. Nous discuterons ce cas s\'epar\'ement.\\

Supposons que $ z_2 = 0 $ est une \'equation d'une branche locale $U_0$ de
$U$ em $m\,$. Notons $h$ l'holonomie de
$\wF$ le long d'un lacet $\gamma$ engendrant $\pi_1 ( U_0^\ast ; \, m_0  )\,$,
$U_0^\ast := U_0 - \{ m \}\,$ et r\'ealis\'ee sur le germe $T_{m_0} := \left\{  z_1 = z_1
( m_0 )
\right\}\,$ de transversale
\footnote{ Ici $T_{m_0}$ est une courbe formelle. } \`a $U_0^\ast\,$ en un  de ses points
$m_0$. Pour une section
$Z$ de $\hB_\wF$ nous d\'esignerons par $\{ Z \}$ sa classe dans
$\hT_\wF\,$. On voit alors que
$$\hT_{\wF , \, m} = \hO_{\wF , \, m} \cdot \{ Z_m \}Ê\, ,$$  avec,  suivant les cas~:
\\

\noindent1.  $\ww := \l_1z_1\,dz_2 +\l_2z_2\,dz_1\,$ avec $\l_1 \l_2 \not=  0\,$,$\l_1 +
\l_2 = 1\,$,  $\l_2/\l_1 \in
\C - \Q\,$,\\
$$\left\{ Z_m \right\}\; =\; \frac{1}{2}\, \left\{ \l_1  z_1\,\dd {z_1} +
\l_2  z_2\,\dd{z_2}\right\} \; = \;
\left\{\l_1  z_1\dd{z_1}\right\} = \left\{\l_2  z_2\dd{z_2}\right\}$$
$$ h(z_2)\; =\; e^{-2i\p\l_2/\l_1}\,z_2\,, \qquad \hO_{\wF , \, m}\; =\;
\C\, .$$

\noindent 2. $\ww := q\,z_1\,dz_2 + p\,z_2\,dz_1\,$ avec $p, q  \in
\N^{\ast}\,$,
$(p,q)=1\,$,
$$\left\{ Z_m \right\}\; =\; \frac{1}{2}\, \left\{ p  z_1\,\dd {z_1} + q
z_2\,\dd{z_2}\right\} \; = \;
\left\{p  z_1\dd{z_1}\right\} = \left\{q  z_2\dd{z_2}\right\}$$
 $$ h(z_2)\; =\; e^{-2i\p p/q}\,z_2\,,\qquad \hO_{\wF,\, m} = \C
\left[\left[z_1^pz_2^q\right]\right]\,.$$

\noindent 3. $\ww :=
q\,z_1(1+\z\,(z_1^pz_2^q)^k)\,dz_2+p\,z_2(1+(\z-1)\,(z_1^pz_2^q)^k)\,dz_1$ avec
$p, q, k  \in \N^{\ast}\,$, $(p,q)=1\,$, $\z \in \C\,$,
$$ \left\{Z_m\right\} =
\left\{-qz_1\dd{z_1} + pz_2\dd{z_2}\right\} =
\left\{\frac{ (z_1^pz_2^q)^k}{ 1 + \z (z_1^pz_2^q)^k}\,z_2\dd{z_2}\right\}  =
\left\{\frac{ (z_1^pz_2^q)^k} { 1 + (\z-1)\,(z_1^pz_2^q)^k}\,z_1\dd{z_1}\right\} \,, $$
$$\hO_{\,\wF, m} = \C\,,\qquad  h(z_2) = e^{\frac{ -2i\p{p}}{  q}}\hbox{exp}
\left(2i\p\frac{ p z_2^{qk + 1}}{ q(1 + \z z_2^{qk})}\dd {z_2}\right)\,.
$$

\noindent 4.  $\;\ww := (\z z_2^{p} - p)z_1\,dz_2 +z_2^{p+1}\,dz_1 \,$,  avec $p \in
\N^{\ast}\,$, $\z \in \C\,$,
$$ \left\{Z_m \right\} = \frac{1}{2}  \left\{\frac { z_2^{p +1}}  {(\z   z_2^p - p) \,}
\dd{z_2} + z_1 \dd{z_1} \right\} = \left\{ \frac{z_2^{p+1}}{(\z z_2^p-p)}
\dd{z_2}\right\} =\left\{z_1 \dd{z_1}\right\}\,,$$

$$\hO_{\wF, m} = \C\,,\qquad
\hbox{$h$ n'est  jamais p\'eriodique \,}.$$

\noindent 5.  $\;\ww$  {\it est un selle-n\oe ud tangent \`a $\{z_2 =  0\}\,$, c.f.}
(\ref{sntagent}).\\
\noindent Dans ce cas, il peut \^etre impossible de r\'ealiser, comme en 4.,  une forme
normale par des coordonn\'ees transversalement formelles \`a $\{ z_2 = 0 \}$. Cependant
L'holonomie
$h$ de cet axe peut quand m\^eme \^etre p\'eriodique. C'est le cas pour la  1-forme
$z_1^{p + 1}\,dz_2 + (\z z_1^p - p)\,z_2\,dz_1\,$, avec
$\z\in \Q\,$, qui admet l'int\'egrale premi\`ere
$$F (z_1, z_2) \; :=\; z_1^n\, z_2^m \, e^{-\efrac{n}{z_1^p}}\, , \qquad \z =:
\frac{m}{n}\, , \quad m , \, n \, \in
\N^{\ast}\, ,
$$
\noindent holomorphe sur $\{z_1 \not= 0\}\,$.

\subsubsection{Extension en un point singulier}Soient $U \,$,
$V\,$ des ouverts de $\U\,$ avec $U
\cap V \not=
\emptyset\,$ et $U_0$ la branche de $U$ contenant $U \cap V \,$. Supposons  que $U$
contient un point singulier
$m$ de
$\wF$ qui n'est pas un selle n\oe ud tangent \`a $U_0\,$.  L'intersection $U
\cap V$ est une couronne.  Nous conservons les notations ci-dessus et supposons que le
point $m_0$ ainsi que le lacet $\gamma$ sont contenus dans $U \cap V\,$. \\

Remarquons que toute section de $ \hO_\wF $ resp. de $ \hT_\wF$ sur $U
\cap V$  s'\'etend en une section de ce faisceau sur $U\,$. En effet le germe en ${m_0}$
d'une  section $X
\in\hT_\wF ( U
\cap V)\,$ peut \^etre repr\'e\-sent\'e par un germe de champ basique tangent \`a
$T_{m_0}$ dont la composante en $\dd{z_1}$ est nulle. La restriction
$X^0$ de ce champ \`a
$T_{m_0}\,$ est invariante par $h\,$. Elle admet une exten\-sion unique en  un champ
basique
$\widetilde{X} \in
\hB_\wF\, (U - \{ m \})\,$. Comme $E_\F^{\ast}(Ê\w ) \cdot \widetilde{X}$  et
$E_\F^{\ast}(Ê\w )
\cdot Z_m$ sont deux facteurs int\'egrants de
$E_\F^{\ast}(Ê\w )\,$, leur quotient $Q$ est une int\'egrale premi\`ere  m\'eromorphe
transversalement formelle. La restriction $Q^0$ de $Q$ \`a
$T_{m_0}$ est invariante par $h\,$, et, quitte \`a remplacer $Q$ par $1 / Q$  on peut
supposer que
$Q^0 \in \hO_{T_{m_0}}\,$.  Dans les cas 1., 3. et 4. $h$ n'est pas  p\'eriodique, ce qui
implique que l'int\'egrale premi\`ere $Q^0$ - et donc aussi $Q$ - est constante $\not=
0\,$. Visiblement $\widetilde{X}  - \mu\, Z_m
\,$, avec $\mu := Q$ ou $1/Q \in \C^{\ast}\,$, est un champ tangent. ainsi
$\{\mu Z_m \}\,$ est une extension de $\widetilde{X}$ sur $U\,$. \\

Il reste \`a examiner le cas 2. Un calcul direct de l'invariance de $X^0$  par $h$ montre
alors que
$X^0 = l\,(z_2^q) z_2 \dd{z_2}\,$, $ l(\xi) \in \CÊ[[\xi]]\,$. ainsi
$l\,(z_1^pz_2^q) z_2 \dd{z_2}\,$ est une extension de $X\,$.\\

La d\'emonstration de l'extension des int\'egrales premi\`eres est similaire. Nous la  laissons au
lecteur.\\

\begin{rema} L'exemple $z_1^{p + 1}\,dz_2 + (\z z_1^p - p)z_2\,dz_1\,$ du selle n\oe ud
tangent donn\'e en 4.b poss\`ede une int\'egrale premi\`ere sur $U - \{m\}\,$, puisque
l'holonomie $h$ est  p\'eriodique, qui ne s'\'etend \'evidemment pas au point $m$. Nous
verrons que la pr\'esence d'une singularit\'e de ce  type pourra emp\^echer le feuilletage
$\F$ d'\^etre de type formel fini.
\end{rema}

\subsubsection{R\'ecapitulation} De m\^eme
qu'en (\ref{counor}) appelons encore {\it couple normalis\'e sur un ouvert}
$W\,$ {\it de}
$\D_\F\,$ tout \'el\'ement
$\left ( f , Z \right)$ de $\hO_\wF ( W ) \times \hT_\wF ( W ) $ qui  v\'erifie,
tou\-jours  avec la convention
$\C [[ 1 ]] := \C\,$, la relation de dualit\'e
\begin{equation}\label{cnor}
\left\{
\begin{array}{l}
\hO_\F ( W) = \C [[ f ]]\, , \quad \hT_\F ( W ) = \C [[ f ]] \, Z \, , \\
\hbox{il existe}\;\; c \in \C \, \hbox{tel que}\,\; df \cdot Z = c\, f \,.
\end{array}
\right.
\end{equation}

\noindent Il est maintenant clair que pour chaque $W \in \U\,$ de tels  couples existent.
L'\'etude cas par cas que l'on vient de faire, compl\'et\'ee par des calculs imm\'ediats,
peut se r\'esumer en la proposition suivante

\begin{prop}\label{passcoin} Soient $U\,$ et $V$ des ouverts de $\U\,$, $U
\cap V \not=
\emptyset\,$ tels que $\wF$ est r\'egulier sur
$V$ et $U$ contient un point singulier $m$ de $\wF$ qui n'est pas un  selle-n\oe ud
tangent. Soit
$T_{m_0}$ un germe de courbe formelle lisse transverse \`a $V$ en un point
$m_0 \in U \cap V\,$. Notons
$H_V\,$, resp. $H_{U^{\ast}}\,$ les sous-groupes de
$
\widehat{Diff}(T_{m_0})$ d'holonomie de $\wF$ le long de $V\,$, resp. le  long de $U^\ast
:= U - \{ m_0
\}\,$. alors les applications naturelles
$$
\hO_\wF (U) \times \hT_\wF (U) \fle \hO_{\wF , \,m} \times \hT_{\wF , \,  m}\, ,\quad
\hO_\wF (U) \times \hT_\wF (U) \fle \hO_\wF (U \cap V) \times \hT_\wF (U
\cap V)\, ,
$$
$$
\hO_\wF (U) \times \hT_\wF (U) \fle \hO_{H_{U^{\ast}}} \times
\hT_{H_{U^{\ast}}} \, , \quad
\hO_\wF (V) \times \hT_\wF (V) \fle \hO_{H_V} \times \hT_{H_V}
$$  sont des isomorphismes respectant les couples normalis\'es. En particulier  toute
section de
$\hO_\wF (V) \times \hT_\wF (V) $ se prolonge en une section de $\hO_\wF  (U \cup V)
\times \hT_\wF (U \cup V)\,$.\\

Supposons de plus que $\hO_\wF (V) \not= \C\,$, et donc aussi $\hO_\wF (U  ) \not= \C\,$.
Fixons des couples normalis\'es
$\left( f_U, Z_U \right)\,$ et $\left( f_V , Z_V \right)\, $ sur $U\,$ et
$V\,$ respectivement. Il existe alors des constantes
$c\,$, $\a \in \C^\ast\,$ et une s\'erie formelle d'une variable
$K(\xi) \,$ telle que
\begin{equation}\label{formex}  f_V =  \a \,f_U^d\, e^{K(f_U)}  \, , \qquad Z_V
=\frac{c}{1 +
\frac{1}{d} \,f_U\, K'(f_U)} \cdot Z_U\, ,
\quad d :=  [ H_V : H_{U^{\ast}}Ê]\,.
\end{equation}

\noindent En particulier on a les identifications~:
$$\hO_\wF (V)\;=\;\C [[ \,f_U^d\, e^{K(f_U)} ]]\; \subset \; \C [[f_U ]]
\; =\; \hO_\wF (U )\, ,$$
$$\hT_\wF (V)\; = \;\C [[ \,f_U^d\, e^{K(f_U)} ]] \, \frac{1}{1 +
\frac{1}{d} \,f_U\, K'(f_U)}
\cdot Z_U\; \subset\; \C [[f_U ]]\cdot Z_U \; = \; \hT_\wF (U )\, .$$
\end{prop}

\subsection{Nerf complet associ\'e \`a $\F\,$}\label{subsenerf}
\addcontentsline{toc}{section}{\hspace{0,8em} {}\thesubsection .  Nerf complet associ\'e \`a $\F\,$}
Consid\'erons {\it le nerf}  $\,\hgN (\F)\,$
du recouvrement distingu\'e
$\,\U\,$~: chaque ouvert
$U$ de $\U$ correspond biunivoquement \`a un sommet $s\,$ de $\,\hgN  (\F)\,$. Nous
noterons alors $U =: U_s\,$. Deux sommets $s $ et
$s'$ sont li\'es par une ar\^ete, not\'ee
$ss'\,$ ou bien $s's\,$, lorsque
$U_{ss'} := U_s\cap U_{s'}\,$ est non-vide. L'ensemble $\gS (\F)$ des sommets de $\hgN
(\F)$ se divise en l'ensemble
$\gS_0 (\F)$ des sommets {\it de type} 0, correspondant aux \'el\'ements de
$\U$ qui sont des voisinages de points singuliers de $\wF$ et en l'ensemble
$\gS_1(\F)$ des sommets {\it de type} 1, correspondant aux \'el\'ements de
$\U$ le long desquels $\wF$ est r\'egulier chaque point. Par construction $\hgN (\F)$ est
un graphe  connexe et simplement connexe.\\

Soit $K$ une partie connexe de $\hgN(\F)$  telle que toutes les ar\^etes de
$K$ joignent deux sommets de $K\,$; nous dirons que $K\,$ est un {\it sous-graphe connexe
de
$\,\hgN (\F)\,$}. appelons {\it valence dans $K$ d'un sommet}
$s$ de $K\,$, le nombre d'ar\^etes de $K$ port\'ees par ce sommet. Les sommets  de valence
1 dans $K$ sont dits {\it extr\'emit\'es} de $K\,$ et leur nombre $v (K) $  s'appellera
{\it valence de $K\,$}.

\subsubsection{Coloriage, pond\'eration et orientation locale de $\hgN
(\F)\,$}
Pour chaque sommet $s$ et chaque ar\^ete $ss'\,$, nous notons~:
$$ E_s := \hT_\wF (U_s) \, , \quad E_{ss'} := \hT_\wF (U_{ss'})\, , \quad   \cal A_s := \hO_\wF
(U_s) \, , \quad \cal A_{ss'} := \hO_\wF (U_{ss'})\, ,
$$ et nous fixons un couple normalis\'e $\left( f_s , \, Z_s \right)\,$ sur $U_s\,$, resp.
$\left( f_{ss'} ,\, Z_{ss'} \right)\,$ sur $U_{ss'}\,$, cf. (\ref{cnor}).  Lorsque $s$ est
de type 0 et la singularit\'e port\'ee par $U_s$ n'est pas un selle-n\oe ud tangent
(\ref{sntagent}), nous prenons pour
$\left( f_{ss'} ,\, Z_{ss'} \right)\,$ la restriction de
$\left( f_{s} ,\, Z_{s} \right)\,$ \`a $U_{ss'}\,$. ainsi, le symbole
$\ast$ d\'esignant un sommet $s$ ou une ar\^ete $ss'\,$, nous avons seulement  trois
possibilit\'es~:
\begin{enumerate}
\item $\cal A_\ast = \C\,$ et $E_\ast = \{ 0 \}\,$,
\item $\cal A_\ast = \C\,$ et $E_\ast = \C \cdot Z_\ast\,$, avec $ Z_\ast \not=  0\,$
\item $\cal A_\ast = \C[[f_\ast]]\,$ et $E_\ast = \C[[f_\ast]] \cdot Z_\ast\,$,  avec
$\cal A_\ast \not= \C \,$, $ Z_\ast \not= 0\,$.
\end{enumerate}  Dans tous les cas, toujours avec la convention
$\C [[ 1 ]] := \C\,$ et $\C \cdot \{ 0\} := \{ 0 \}\,$, nous pouvons
\'ecrire~:
$$ E_\ast = \cal A_\ast \cdot Z_\ast \, \quad \mbox{et} \quad \cal A_\ast = \C[[  f_\ast]] \, .
$$
\noindent Pour chaque ar\^ete $ss'$ nous disposons d'op\'erations de  restrictions qui
s'expriment comme des injections $\C$- lin\'eaires
\begin{equation}\label{aplrstr}
\r_{ss'}^s : E_s \ifle E_{ss'}\, \qquad \hbox{et}\qquad \s_{ss'}^s : \cal A_s
\ifle \cal A_{ss'}
\end{equation}
\noindent et que nous consid\'ererons comme des inclusions. On a ainsi des relations
\begin{equation}\label{simpl} Z_s = a_{s'}^s\,(f_{ss'})\; Z_{ss'} \, ,
\quad f_s = b_{s'}^s\,(f_{ss'})\, ,\quad  a_{s'}^s(\z) \, , \;Ê\; b_{s'}^s(\z) \; \in
\; \C[[\z]]\, .
\end{equation}
\noindent Dans le cas 3. les formules (\ref{formex}) de (\ref{passcoin})  pr\'ecisent les
expressions des s\'eries
$ a_{s'}^s(\z) \,$ et $ b_{s'}^s(\z)\,$. En particulier $a_{s'}^s( 0 )
\not= 0\,$. Dans les autres cas
$a_{s'}^s(\z)$ et $b_{s'}^s(\z)$ sont des constantes, \'eventu\-el\-lement  nulles. \\

Nous allons maintenant orienter les ar\^etes de $\,\hgN (\F)\,$. Nous  convenons que~:
$$
\begin{array}{lcl}
\circ_s \;\rightarrow \; \circ_{s'} & \quad\mbox{signifie~:}\quad &
\mbox{${\r_{ss'}^{s}}\;\;$ n'est pas bijective  et
$\;\;{\r_{ss'}^{s'}}\;\;$ est bijective,}\\
\circ_s \;\leftarrow \; \circ_{s'} & \quad\mbox{signifie~:}\quad &
\mbox{${\r_{ss'}^{s}}\;\;$ est bijective  et
$\;\;{\r_{ss'}^{s'}}\;\;$ n'est pas bijective,}\\
\circ_s \;\leftrightarrow \; \circ_{s'} & \quad\mbox{signifie~:}\quad &
\mbox{${\r_{ss'}^{s}}\;\;$ et
$\;\;{\r_{ss'}^{s'}}\;\;$ sont bijectives,}
\\
\circ_s \;\bdf \; \circ_{s'} & \quad\mbox{signifie~:}\quad &
\mbox{${\r_{ss'}^{s}}\;\;$ et
$\;\;{\r_{ss'}^{s'}}\;\;$  ne sont pas bijectives.}
\end{array}
$$
\noindent Lorsqu'on ne veut pas pr\'eciser $\;\circ_s \;\ftrait \;
\circ_{s'}\;$ signifiera l'un quelconque de ces  cas. \\

Nous allons maintenant colorier en vert ou rouge
$\hgN (\F)\,$ puis pond\'erer les sommets rouges et les ar\^etes vertes qui  joignent deux
sommets verts.\\

\indent Un sommet ou une ar\^ete, not\'e  $\ast$, est {\it colori\'e en rouge} si
$\cal A_\ast = \C\,$ ou, ce qui revient au m\^eme, si
$dim_\C\, E_\ast = 0$ ou 1. Les ar\^etes et les sommets $\ast$ tels que
$\cal A_\ast \not= \C$ sont {\it colori\'es en vert}.  Le sommet
$s$ se dessinera par
$\;\bullet_s\;$ s'il est rouge et par $\;\star_s\;$ s'il est vert. Lorsque  nous ne
d\'esirons pas pr\'eciser la couleur nous le dessinerons par $\;\circ_s\;$.\\

Remarquons que toute ar\^ete joint un sommet de type 1 (correspondant \`a un  ouvert $U
\in \U$  sans singularit\'e du feuilletage) \`a un sommet de type 0 (l'ouvert $U$
correspondant porte une singularit\'e). Un  sommet rouge correspond \`a un ouvert
d'holonomie non-commutative lorqu'il est de type 1.; lorsqu'il est de type 0, il
correspond \`a une singularit\'e sans int\'egrale premi\`ere. Une ar\^ete reli\'ee \`a un
sommet vert est verte. Pour une ar\^ete verte, le sommet de type  0 auquel elle est
reli\'ee est soit vert, soit rouge et correspondant \`a une singularit\'e de type
selle-noeud tangent.  Les seules possibilit\'es, pour une ar\^ete verte, sont~:
\begin{equation}\label{casposs}
\begin{array}{lll} a.\;\;\star_{s'} \;\leftarrow\;  \star_{s}\;,\quad  &  b.\;\;\star_{s'}
\;\leftrightarrow\;
\star_{s}\; ,\quad & \\ c.\;\;\star_{s'} \;\leftarrow\;  \bullet_{s}\;,
\quad & d.\;\;\star_{s'}
\;\bdf\;  \bullet_{s}\; ,\quad  &  e.\;\;\bullet_{s'} \;\bdf\;
\bullet_{s}\; .
\end{array}
\end{equation}
\noindent Le sommet $ s'$ est de type 0 dans le cas a.\,;  Par contre $ s
$ est de type 0 et
$U_{s}$ contient une singularit\'e selle-n\oe ud tangent dans le cas d..  La
configuration c. peut se produire avec
$s'$ de type 0, ou bien avec $s'$ de type 1 et dans ce cas
$U_{s}$ contient un selle-n\oe ud tangent; enfin dans le cas e., soit
$U_s$ soit $U_{s'}$ contient un selle-n\oe ud tangent.\\

On associe \`a un sommet rouge $\bullet_s\,$, {\it son poids}
$$d_s := dim_\C
\,E_s = 0 \; \hbox{ou} \;\; 1\, .$$ Pour toute ar\^ete rouge $ss'$, on a $dim_\C E_{ss'} =
1\,$, puisque  l'holonomie d'une couronne est commutative. Visiblement, pour une ar\^ete
rouge, les seules configurations  possibles sont~:
\begin{equation}\label{posco}
\bulletun\;\leftrightarrow \;\bulletun\; , \qquad\quad
\bulletun \; \leftarrow\;\bulletO\; , \qquad\quad
\bulletO \bdf\;\bulletO\;  ,
\end{equation}
\noindent et dans la derni\`ere configuration, le sommet de type 0  correspond \`a une
singularit\'e de type selle-noeud.\\

\begin{defi}\label{nerfcol} Nous appelons {\it nerf complet de $\F$} et notons $\hgN^\ast
(\F)$ le  nerf colori\'e de $\F\,$,  o\`u chaque sommet rouge est muni de son poids 0 ou 1.
La  partie rouge de $\hgN^\ast
(\F)$ sera d\'esign\'ee par
$\hgR^\ast (\F)\,$. Nous le munissons de la m\'etrique pour
laquelle chaque ar\^ete est isom\'etrique au segment $[0 , \,  1]$ standard.
\end{defi}
\noindent Clairement $\F$ est bon (\ref{fbon}) si et seulement si
$\hgR^\ast (\F)$ est non-vide.

\subsubsection{Parties actives du nerf de $\F$}
Ce sont les
parties de $\hgN^\ast (\F)$ qui  contribueront \`a la dimension de $H^1 \left( \U \, ; \, \hT_\wF
\right)\,$.\\

Nous appelons {\it partie active rouge de $\hgN^\ast (\F ) \,$} tout sous-graphe connexe
de $\hgR^\ast (\F ) \,$ dont les sommets extr\'emit\'e sont de poids 0, tous les autres
sommets \'etant de poids 1. Remarquons que $\hgN^\ast (\F )$ peut ne pas poss\'eder de partie  active rouge.
D'autre part chaque sous-graphe $\bulletO
\bdf\;\bulletO\; $ est une partie active rouge.

\begin{defi}\label{bouquet}  Supposons $\hgR^\ast (\F)$ non-vide et  connexe. Nous
appelons {\it bouquet de cercles assosi\'e \`a
$\F$} l'espace topologique $\hgC (\F)$ obtenu \`a partir de la r\'ealisation
g\'eom\'etrique de $\hgR^\ast (\F)$ en identifiant \`a un seul et m\^eme  point tous les
sommets de poids 0.
\end{defi}
\noindent Cet espace topologique est s\'epar\'e. Il a bien le type d'homo\-topie d'un
bouquet de  cercles, ou d'un point. Le nombre de  cercles se calcule \`a partir des
parties actives $L_1, \ldots , L_r\,$ de $\hgR^\ast (\F ) $ par la  formule~:
\begin{equation}\label{nbrbousph}
\widehat{\s} (\F ) \,:= \, dim_\Z\, H_1\left( \hgC (\F) \, ; \, \Z
\right)\, = \, \sum_{j = 1}^r \,
\left( v(L_j) - 1
\right)\, ,
\end{equation}
\noindent o\`u la valence $v (L_j )$ est le nombre de sommets extr\'emit\'es de
$L_j$.

\subsection{Les crit\`eres de finitude formelle de $\F\,$}\label{thcodsubsec}
\addcontentsline{toc}{section}{\hspace{0,8em} {}\thesubsection .  Les crit\`eres de finitude formelle de $\F\,$}
La finitude
formelle de $\F\,$, en abr\'eg\'e $\F$ est t.f.f., que l'on a  d\'efinie en
(\ref{stffdef}) par $\hb ( \F ) := dim_\C\, H^1
\left( \U \, ;
\, \hT_\wF \right)\, < \infty\,$, \'equivaut par (\ref{splitcod}) \`a celle de
$\htau ( \F ) := dim_\C\, H^1 \left( \U \, ; \, \hB_\wF \right)\,$. La formule
(\ref{fordepinf}) montre que $\hd (\F ) := dim_\C\, H^1 \left(
\U \, ; \, \hX_\wF \right)\,$ se calcule
\`a partir de l'arbre dual doublement pond\'er\'e $\cA [\F ]\,$ introduit en
(\ref{adudoub}).  Lorsque $\F$ ne poss\`ede pas de facteur int\'egrant formel on a $\hb
(\F ) =
\htau (\F ) + \hd (\F )\,$ d'apr\`es (\ref{splitcod}) et nous verrons que $\hb (\F)$ ne
d\'epend que  de
$\hgN^\ast (\F)\,$ si $\hgR^\ast (\F) \not= \emptyset\,$.

Nous disons qu'un
sous-graphe connexe $K$ de $\hgN^\ast (\F)\,$ est {\it r\'epulsif} si toute ar\^ete
attach\'ee \`a un sommet $s \notin K$ est soit  du type simple fl\`eche
$\circ_s \;\rightarrow \; \circ_{s'}$ dirig\'ee vers l'ext\'erieur
\footnote{ Ceci a un sens puisque $K$ est connexe et $| \hgN^\ast (\F)\, |$  simplement
connexe. } de  $K\,$, soit du type
$\circ_s \;\leftrightarrow \; \circ_{s'}\,$. En particulier $\hgN^\ast  (\F)\,$ est
r\'epulsif.

\begin{enonce}{Th\'eor\`eme de codimension}\label{critfr} Soit $\F$ un feuilletage formel  non-dicritique \`a
singularit\'e isol\'ee non-r\'eduite \`a l'origine de $\C^2\,$. Supposons
$\hgR^\ast(\F)\,$ non-vide\footnote{ i.e. la r\'eduction de $\F$ comporte un \'el\'ement
critique $C
\subset \D_\F$ tel que
$\hO_\wF (C) = \C\,$. }. Alors on a~:
\begin{enumerate}
\item $\F$ est t.f.f. si et seulement si $\hgR^\ast (\F)$ est connexe et  r\'epulsif. Dans
ce cas  $\htau (\F )$ est \'egal au nombre
$\widehat{\s}  (\F)\,$ de cercles du bouquet $\hgC (\F)\,$ d\'efini en  (\ref{bouquet}),
et de plus~:
$$\hb \left( \F \right) = \hd \left( \F \right)  +  \htau \left( \F
\right)  -
 \widehat{\e '} \left( \F \right) \,, $$    o\`u $\hd (\F) $ est d\'efini par la formule
(\ref{cotr}) et o\`u $\widehat{\e  '} (\F ) = 0$ ou $1 \,$.
\item Si $\F$ est t.f.f. et de deuxi\`eme esp\`ece (\ref{especes}) alors
$\widehat{\e '} (\F ) = 0\,$ ou 1. Le cas  $\widehat{\e '} (\F ) = 1\,$ se produit
exactement lorsque
$\F$ poss\`ede un facteur int\'egrant (\ref{factint}) formel \`a l'origine de
$\C^2\,$ qui s'annule sur les s\'eparatrices formelles de
$\F$ et il n'existe pas de champ formel \`a l'origine de $\C^2\,$ basique et  non-tangent
pour $\F$. De plus
$\hd (\F )$  est \'egal \`a la dimension de la strate \`a
$\mu$-cons\-tant de l'ensemble des s\'eparatrices formelles $\widehat{Sep}  (\F )\,$, c.f.
(\ref{minmult}).
\end{enumerate}
\end{enonce}

\subsection{D\'emonstration du Th\'eor\`eme de codimension}
\addcontentsline{toc}{section}{\hspace{0,8em} {}\thesubsection .  D\'emonstration du Th\'eor\`eme de codimension}
Une partie des r\'esultats du th\'eor\`eme se d\'eduit de  la proposition (\ref{splitcod}).
Lorsque $\F$ est t.f.f. on a, avec les  notations de (\ref{splitcod})~:
$
\hb \left( \F \right)\; =
\hd \left( \F \right) \; + \; \htau \left( \F \right)\; - \; {\hE'} \left(
\F \right)
$ o\`u
$$ 0 \; \leq \;{\hE'} \left( \F \right) \; \leq \;
 dim_\C
\left( \left.
\widehat{Int}\,(\w)\; \right/ \left( \w \cdot \hB (\w ) \right) \right) \,.
$$
\noindent Si $\F$ est t.f.f. et $\hgR^\ast (\F) \not= \emptyset$, il  n'existe pas
d'int\'egrale premi\`ere globale non-constante le long de
$\D_\F$ et donc tout int\'egrale premi\`ere formelle de $\F$ \`a l'origine de
$\C^2$ est constante. Comme $\F $ est non-dicritique toute int\'egrale premi\`ere
m\'eromorphe formelle \`a l'origine de
$\C^2$ est aussi constante. On en d\'eduit que
$\left.
\widehat{Int}\,(\w)\; \right/ \left( \w \cdot \hB (\w ) \right)$ est de  dimension 0 ou 1,
puisque le quotient de deux facteurs int\'egrants est une int\'egrale premi\`ere
m\'eromorphe. On a donc bien $\;\hE'
\left( \F \right) = 0$ ou 1.\\

\indent Lorque $\F$ est de deuxi\`eme esp\`ece, on a
$$\hE' \left( \F \right) =
 dim_\C
\left( \left.
\widehat{Int}\,(\w)\; \right/ \left( \w \cdot \hB (\w ) \right) \right) =
\hbox{0 ou 1}\, .$$
\noindent Le cas $\hE' \left( \F \right) = 1$ se produit exactement lorsque, \`a l'origine
de $\C^2$, il existe un facteur  int\'egrant formel s'annulant sur les s\'eparatrices
formelles de $\F$, mais tout champ formel basique est tangent.\\

\indent Supposons \`a pr\'esent $\hgR^\ast (\F) $ non-vide et montrons que
$\F$ est t.f.f. si et seulement si $\hgR^\ast (\F) $ est connexe r\'epulsif.\\

\indent Soit $K$ un sous-graphe de
$\hgN^\ast (\F)\,$. Notons
$\gS(K)
\subset
\gS (\F)\,$ l'ensemble des sommets de $K\,$ et
$\gS_i(K)\,$ le sous-ensemble des sommets de $K\,$ de type $i$, avec $i =  0\,$ ou $1\,$.
Nous appelons {\it application de cohomologie associ\'ee \`a
$K$} l'application $\C$-lin\'eaire~:
\begin{equation}\label{eqchomo}
\begin{array}{c} {\displaystyle \Delta_K \; : \;  \prod_{s \in \gS(K)} E_s
\; =: \; Z^0(K)
\;\;\fle Z^1(K)\; :=\; \prod_{(s, s') \in  {\gS}(K)^{\check{2}}}  E_{ss'}}\, ,\\
\\ {\displaystyle \left( X_{s} \right)_s \efle \left( X_{ss'}\right)_{(s,  s')} \, , \quad
\mbox{avec} \quad X_{ss'} := X_{s'} - X_{s}\;} ,
\end{array}
\end{equation} o\`u
\begin{equation}\label{prodchech} {\gS}(K)^{\check{2}} := \left\{ \,\left. (s, s') \in
\gS_0(K)Ê\times
\gS_1(K) \;\right/ \; \hbox{ $ss'$ ar\^ete de $K$ } \; \right\}
\,.
\end{equation}  Nous noterons dans tout ce qui suit~:
$$
\cal H^0\,( K ) := ker \left( \Delta_K \right)\, , \qquad
\cal H^1\,( K ) := coker \left( \Delta_K \right)\, .
$$   Remarquons que deux sous-graphes $K\,$ et $K'\,$ de
$\hgN^\ast (\F)\,$ v\'erifient~:
\begin{equation}\label{rscoh} dim_\C \cal H^1 (K) \geq dim_\C\cal H^1 (K')
\qquad \mbox{d\`es que}
\qquad K' \; \subset K\,.
\end{equation}
\noindent En effet, $\pi$ d\'esignant la projection lin\'eaire $Z^1\,( K )
\fle Z^1\, (K')$ induite par la structure produit, les applications $\Delta_{K'}\,$,
$\pi\circ\Delta_K$ ont m\^eme  image.\\

Visiblement
\begin{equation}\label{conerftr}
\cal H^1\,\left( \hgN^{\ast} (\F)\right) = H^1 \left( \U \, ; \, \hT_\wF
\right)\,.
\end{equation} Pour voir que $\F$ n'est pas t.f.f. il suffit donc de  d\'etecter un
sous-graphe $K$ de $\hgN^{\ast} (\F)$ tel que
$dim_\C \cal H^1 (K) = \infty\,$.\\

\begin{rema}\label{intfais.} En fait $\cal H^1 \left( \hgN^{\ast} (\F )
\right)$ peut s'interpr\'eter comme un espace de cohomologie de Cech~:  munissons la
r\'ealisation  g\'eom\'etrique $\left| \hgN^\ast (\F)
\right|$ de $\hgN^\ast (\F)$ de la topologie dont une base d'ouverts est  la collection
$\cal V$ des r\'ealisations g\'eom\'etriques $\V_s\,$,
$s \in \goth S (\F )\,$,  des parties de $\hgN^{\ast} (\F )$ constitu\'ees  d'un sommet
$s$ et  de toutes les ar\^etes qui y sont attach\'ees. Consid\'erons le faisceau
$\widehat{\Bbb T}_\F$  de base
$\left| \hgN^{\ast} (\F) \right|$ d\'efini par
$$\widehat{\Bbb T}_\F ( \V_s ) := E_s\, , \qquad \widehat{\Bbb T}_\F (\V_s
\cap \V_s') := E_{ss'}\,, \quad (s, s') \in {\gS}(K)^{\check{2}}\, ,$$ les applications
de restriction associ\'ees aux inclusions
$\V_s \cap \V_{s'} \subset \V_s$ et $\V_s \cap \V_{s'} \subset \V_{s'}$
\'etant respectivement les injections
$\,\r_{ss'}^{s}\,$ et $\,\r_{ss'}^{s'}\,$ d\'efinies en (\ref{aplrstr}).   alors
$$
\cal H^j \left(K \right)\, =\, H^j \left( \cal \V_K ; \,\widehat{\Bbb  T}_\F \right)\,
,\qquad j = 0 , \, 1\, ,
$$
\noindent o\`u
$\cal \V_K$ est le recouvrement $\left( \V_s \right)_{s \in \gS(K)}\,$ de
$K\,$.
\end{rema}

\begin{lemm}\label{lemcorou} Soit $K$ une g\'eod\'esique
\footnote{rappelons que nous avons muni $\N^\ast (\F)$ de la  m\'etrique pour laquelle
chaque ar\^ete est isom\'etrique au segment $[0 , \,  1]$ standard. }  de $\hgN^\ast (\F)
\,$ de l'un des types suivants~:
\begin{enumerate}
\item  $\bullet_{s_0} \; \;\ftrait \; \; \star_{s_1} \;
\;\ftrait\cdots\cdots\cdots \ftrait
\;\; \star_{s_n}\;\; \longleftarrow \;\;
\star_{s_{n+1}}\quad$ avec $n \geq 1\,$,
\item  $\bullet_{s_0} \;\; \ftrait\; \; \star_{s_1} \;
\;\ftrait\cdots\cdots\cdots \ftrait
\;\; \star_{s_n}\;\;
\ftrait\; \;\bullet_{s_{n+1}}\quad$ avec $n \geq 1\, ,$
\item  $\bullet_{s_0} \;\; \ftrait \;\; \bullet_{s_1}\;\;$  l'ar\^ete \'etant  verte,
\item  $\bullet_{s_0} \;\;\bdf\;\;\star_{s_1}\;$,
\end{enumerate} les sommets $s_2 , \ldots , s_{n-1}$ \'etant verts et  l'orientation des
ar\^etes ${s_j s_{j+1}}\,$,
$j = 1, \ldots ,n-1$ \'etant quelconque. alors~:
$\;dim_\C\,\cal H^1(K) =
\infty\,$.
\end{lemm}

\begin{proof} Examinons le cas 1. Quitte \`a raccourcir la g\'eod\'esique, on  supposera
gr\^ace \`a (\ref{rscoh}) que toutes les ar\^etes ${s_j s_{j+1}}\,$,
$j = 0, \ldots ,n-1$ sont, soit des doubles fl\`eches
$\star_{s_j}\;\longleftrightarrow\;\star_{s_{j+1}}\,$, soit des fl\`eches  orient\'ees
vers $s_n\,$, i.e.
$\star_{s_j}\; \longrightarrow\;\star_{s_{j+1}}\,$.  En consid\'erant les  op\'erateurs
d'extension et de restriction comme des inclusions, on obtient les identifications
suivantes~:
\begin{equation}\label{idespaces} E_{s_0} \subset E_{s_0s_1} \subset  E_{s_1} \subset
E_{s_1s_2}
\subset
\cdots
\cdots
\subset E_{s_n} \subset \; E_{s_n s_{n+1}} \;\supsetn\; E_{s_{n+1}}\; .
\end{equation}
\noindent qui permettent de consid\'erer tous les espaces comme des sous  espaces de
$E_{s_n s_{n+1}}\,$. Les sommets
$s_n$ et
$s_{n+1}$ sont de type 0 et 1 respectivement.  Fixons un couple normalis\'e
$(f ,  Z) \in  \cal A_{s_n\,s_{n+1}}Ê\times E_{s_n\,s_{n+1}}\,$. Par (\ref{passcoin}) les
identifications ci-dessus  s'\'ecrivent~:
\begin{equation}\label{idchamps}
\C \cdot a(f) \,   Z \;  \subset \;\cdots \cdots\; \subset \; \C[[ f ]]  Z
\;
\supset
 \; \C[[ \wb (f)]] \cdot \wa (f)\,  Z\, .
\end{equation}
\noindent o\`u $a\,$, $\wa\,$, $\wb\,$ sont des s\'eries formelles d'une  variable
satisfaisant~:
$\wb (0) = \wb' (0) = 0$ et $\wa(0) \not= 0\,$. Il suffit de montrer que  le conoyau du
compos\'e
$\widetilde{\Delta}_K$ de l'application de cohomologie $\Delta_K$ d\'efinie  en
(\ref{eqchomo}) avec la surjection lin\'eaire
$$Z^1 (K) \fle E_{s_n , \, s_{n+1}}\,, \qquad \left( X_{j, \, j+1}
\right)_{j= 0 \, , \ldots , n}
\efle \sum_{j=0}^n  X_{j, \, j+1}$$  est de codimension infinie. L'application
$\widetilde{\Delta}_K$ envoie
$\left( X_{j, \, j+1} \right)_{j= 0 \, , \ldots , n}$ sur $X_{n+1} -  X_0\,$. ainsi
l'image de
$\widetilde{\Delta}_K$ s'\'ecrit~:
$$Im\left( \widetilde{\Delta}_K \right) = \left(
\C\left[\left[ \, \wb (f) \right] \right] \wa (f) + \C\, a(f) \right)
\cdot  Z\,,$$ avec $\wb (0) = \wb' (0) = 0\,$; d'o\`u la conclusion. \\

\indent Dans le cas 2. la m\^eme m\'ethode montre que
$Im\left( \widetilde{\Delta}_K
\right)$ est de dimension finie. Le cas 3. est trivial et le cas 4. donne  le m\^eme
r\'esultat que le cas 2..
\end{proof}

Remarquons que, d'apr\`es la liste (\ref{casposs}), $\hgN^{\ast} (\F )\,$ ne  peut
poss\'eder d'ar\^ete du type
$\bullet \; \;\longleftarrow \; \; \star\,$. Il d\'ecoule clairement du  lemme ci-dessus
que si $\F$ est t.f.f., alors $\hgR^\ast (\F)$ est connexe et r\'epulsif, \`a moins qu'il
ne soit vide.\\

\indent Pour obtenir la r\'eciproque d\'efinissons maintenant des "op\'erations
d'\'elagage" d'un sous-graphe de
$\hgN^\ast (\F)$, qui ne changent pas la dimension de l'espace de  cohomologie associ\'e.\\

\indent  Appelons  {\it branche verte s\'ecable} d'un sous-graphe connexe $K$ de
$\hgN^{\ast} (\F )\,$, toute g\'eo\-d\'esique $B \subset K$, reliant deux sommets de $K$,
qui v\'erifie~:
\begin{enumerate}
\item tous les sommets de $B$, autres que ses extr\'emit\'es, sont verts et de  valence 2
dans $K\,$;
\item l'un des sommets extr\'emit\'e de $B\,$, que nous appelons {\it sommet  d'attache de
$B\,$}, est {\it r\'epulsif dans $B$}~: chaque ar\^ete de $B$ est soit une double fl\`eche
$\circ \leftrightarrow \circ\,$,  soit une fl\`eche simple $\circ \rightarrow \circ\,$
orient\'ee vers l'autre sommet extr\'emit\'e de $B\,$, que nous appelons {\it  sommet
libre de $B\,$};
\item le sommet libre de $B$ est vert et de valence 1 dans $K\,$; le  sommet d'attache
est~: soit rouge, soit de valence $\not= 2$, soit vert du type $\longleftarrow\; \star \;
\longrightarrow\,$.
\end{enumerate}
\noindent Remarquons que l'on peut avoir $B = K \,$. Dans ce cas, lorsque  toutes les
fl\`eches sont des doubles fl\`eches, le choix du sommet d'attache est arbitraire.\\

\indent D\'esignons par $El_B(K)$ le  {\it graphe $K$ \'elagu\'e de la branche $B$} c'est
\`a dire le graphe connexe obtenu \`a partir de
$K$ en supprimant toutes les ar\^etes et tous les sommets de $B\,$, sauf le  sommet
d'attache. Lorsque $B = K\,$, $El_B(K)$ est le graphe trivial r\'eduit \`a un seul sommet.
Pour un sous-graphe $L$ de $K\,$,  nous d\'esignons par $\l_L^1$ l'injection de
$Z^1\left(L
\right)$ dans
$Z^1\left( K
\right)$ d\'efinie par
$\l_L^1 \cdot \left( X_{ss'}\right)_{(s, s')} := \left(Y_{ss'}\right)_{(s,  s')}\,$, avec
$Y_{ss'} := X_{ss'}\,$ si
$s\,$, $s'$ sont les sommets d'une ar\^ete de
$L\,$ et $Y_{ss'} = 0$ sinon.

\begin{lemm}\label{secvert} Soit $B$ une branche verte s\'ecable de $K\,$,  l'application
$$\left[ \kappa_{K'}^1 \right] \, : \, \cal H^1\left(K'\right) \fle \cal  H^1\left( K
\right)
\, , \quad
\left[\left( X_{ss'}\right)_{(s, s')}\right] \efle \left[ \kappa_{K'}^1
\cdot \left( X_{ss'}\right)_{(s, s')} \right]\, ,
$$
\noindent avec $ K' := El_B(K)\,$, est un isomorphisme.
\end{lemm}
\begin{proof} Soit
$\circ_{s_0} \; \trait{aa} \;  \star_{s_1} \; \trait{aa}
\cdots\cdots\cdots \trait{aa} \;
\star_{s_{n+1}} \,$ la branche verte s\'ecable $B\,$, $s_0$ d\'esignant le  sommet
d'attache. Consid\'erons l'injection
$\widetilde{\l}^0 : Z^0\left( B - \{s_0\} \right) \fle Z^0\left( K
\right)\,$,  $\left( X_{s}\right)_s
\efle
\left(Y_{s}\right)_s\,$ d\'efinie par $Y_{s} := X_{s}\,$ si $s$ est un  sommet de $B -
\{s_0\}$ et $Y_{s} = 0$ sinon. Le diagramme suivant est commutatif et ses lignes sont
exactes~:
$$
\begin{array}{ccccccccc}  0 & \fle & Z^0\left( B - \{s_0\} \right) &
\stackrel{\widetilde{\l}^0}{\fle} & Z^0\left( K \right) &
\stackrel{\pi^0}{\fle} & Z^0\left(K' \right) &
\fle & 0Ê\\
  &    &\widetilde{\Delta}_B \downarrow \phantom{\Delta_B} & &\Delta_K
\downarrow
\phantom{\Delta_K} &  &
\Delta_{K'} \downarrow \phantom{\Delta_{K'}} & & Ê\\
0 & \fle & Z^1\left( B \right) & \stackrel{\l_B^1}{\fle} & Z^1\left( K
\right) &
\stackrel{\pi^1}{\fle}& Z^1\left(K' \right) &
\fle & 0\; ,Ê\\
\end{array}
$$
\noindent o\`u $\pi^0\,$, $\pi^1$ sont les projections canoniques donn\'ees  par les
structures produits et l'application
$\widetilde{\Delta}_B$ est d\'efinie par la relation de commutativit\'e  qu'elle doit
v\'erifier.  La suite exacte longue associ\'ee donne
$$\cdots \,\fle  \, coker (\widetilde{\Delta}_B)  \fle  \cal H^1\left( K
\right)  \fle  \cal H^1\left(K'
\right)  \fle  0Ê\, .$$
\indent Montrons que $\widetilde{\Delta}_B$ est surjective. Donnons nous  un \'el\'ement
$\left(X_{j, j+1} \right)_{j = 0, \ldots , n}\,$ de
$Z^1\left( B \right)\,$. Faisons les identifications
$$E_{s_0} \subset E_{s_0s_1} \subset E_{s_1} \subset E_{s_1s_2} \subset
\cdots \cdots \subset E_{s_n} \subset E_{s_n s_{n+1}} \subset E_{s_{n+1}}\, ,$$ puisque
$s_0$ est un sommet r\'epulsif de la branche s\'ecable $B\,$. Il est
\'evident que l'on peut r\'esoudre le syst\`eme d'\'equations
$$ X_1 - 0 = X_{0, 1}\,,\quad X_2 - X_1 = X_{1, 2}\, ,\quad \ldots, \quad  X_{n, n+1} =
X_{n+1} - X_n\,,
$$  d'inconnues $X_j \in E_{s_j}\,$.\\
\indent Ainsi $\pi^1$ induit un isomorphisme $\left[ \pi^1 \right] : \cal  H^1\left( K
\right)  \iso  \cal H^1\left(K'
\right)\,$. Comme l'application $\l^1_{K'}$ est une section de $\pi^1$, on  obtient que
$\left[ \l^1_{K'} \right]$ est aussi un isomorphisme.
\end{proof}

\indent Supposons $ \hgR^\ast (\F) $ connexe, non-vide et r\'epulsif et  montrons que $\F$
est t.f.f.. Notons
$K_0 := \hgN^\ast (\F)\,$. aucun sommet vert de $K_0$ n'est du type $\, \longleftarrow
\star  \longrightarrow \,$. Chaque sommet vert
$s_j$ extr\'emit\'e de $ K_0$, $j = 1, \ldots ,r_0\,$ est aussi l'extr\'emit\'e libre
d'une branche verte s\'ecable $B_j\,$ de
$K_0$, le sommet d'attache de $B_j$ \'etant soit un sommet rouge, soit un sommet vert de
valence
$\geq 3\,$. Consid\'erons le sous-graphe $K_1\,$ obtenu en \'elagant $K_0$  (dans un ordre
indiff\'erent) des branches $B_1 ,
\ldots , B_{r_0}\,$. Le graphe $K_1$ est connexe et $ \hgR^\ast (\F) $ est  toujours une
partie connexe r\'epulsive de $K_1\,$. En r\'ep\'etant cette op\'eration on obtient une
filtration  $\hgN^\ast (\F) =: K_0
\, \supset K_1\,\supset\, \cdots\,
\supset K_p$ avec $K_p =  \hgR^\ast (\F) \,$. Le lemme donne~:
\begin{equation}\label{seccorep}
\cal H^1 \left( \hgN^{\ast} (\F) \right) \, = \, \cal H^1 \left( \hgR^\ast  (\F) \right)
\end{equation}

\indent Nous allons maintenant \'elaguer $\hgR^\ast (\F)\,$.\\

\indent Appelons {\it branche rouge s\'ecable} d'un sous-graphe connexe $K$  de $\hgR^\ast
(\F)$ toute g\'eo\-d\'esique
$BÊ\subset K$  qui joint une extr\'emit\'e de $K$ de poids 1 \`a un sommet de
$K$ qui est soit de valence
$\not= 2$ dans $K\,$, soit de valence 2 dans $K$ et de poids 0, et telle  que les autres
sommets sont de valence 2 dans $K$ et de poids 1. Lorsque $B = K$ et tous les sommets de
$B$ sont de poids 1, on choisit arbitrairement l'un des sommets extr\'emit\'e que l'on
appelle {\it sommet d'attache de $B\,$}. Dans  les autres cas on appelle {\it sommet
d'attache de $B$} le sommet extr\'emit\'e qui est de valence $\geq 3$ ou de poids 0.
\\

Chaque sommet extr\'emit\'e de $\hgR^\ast (\F)\,$ qui n'est pas de poids 0 est  aussi
l'extr\'emit\'e libre d'une branche rouge s\'ecable. En it\'erant l'op\'eration qui
consiste \`a \'elaguer successivement toutes ces  branches, on construit de nouveau une
filtration finie d\'ecroissante de premier terme $\hgR^\ast (\F)\,$. Le dernier terme de
cette filtration est soit un point, soit un sous-graphe connexe
$El
\, (\,
\hgR^\ast (\F) \, ) \, \subset \,\hgR^\ast (\F)$ dont les sommets  extr\'emit\'e sont de
poids 0. ainsi $El \, (\, \hgR^\ast (\F)
\, )$ est l'union (non-vide) des parties actives rouges $L_1\,, \ldots ,   L_r$ de
$\hgN^\ast (\F)\,$. On a~:
$$\cal H^1 \left( \hgN^\ast (\F) \right) \; \iso \; \cal H^1 \left(
\hgR^\ast (\F) \right)
 \; \iso \; \cal H^1 \left( El \, \left(\, \hgR^\ast (\F) \, \right)
\right)\,.$$
\noindent  Pour achever la d\'emonstration il reste \`a montrer le

\begin{lemm}\label{cohorouge} On a ~:
$\,dim_\C\,\cal H^1 \left( El \, \left(\, \hgR^\ast (\F) \, \right)
\right) \, = \, dim_\Z\, H_1\,
\left( \hgC (\F)\, ; \, \Z
\right)\,$, o\`u $\hgC (\F)$ d\'esigne le bouquet de cercles associ\'e \`a $\F\,$
d\'efini en (\ref{bouquet}).
\end{lemm}

\begin{proof} Supposons que $El \, (\, \hgR^\ast (\F) \, )\,$ n'est pas  r\'eduit \`a un
point, sinon l'\'egalit\'e est triviale. Visiblement
$$\cal H^1 \left(\, El \,\left( \hgR^\ast (\F)Ê\right)\,
\right)  =  \bigoplus_{j = 1 }^r
\cal H^1 \left(\, L_j\, \right)\, .$$
\noindent On supposera donc que $\hgR^\ast (\F)$ ne poss\`ede qu'une seule  partie active
rouge, not\'ee $K\,$. Soit $v$ la valence de $K$. L'homologie de
$\hgC (\F)$ est celle du complexe simplicial $C$ de dimension 1, obtenu en  rajoutant \`a
$K$ un sommet, not\'e $s_0\,$, et $v$ ar\^etes qui relient $s_0$ aux extr\'emit\'es
$s_1, \ldots , s_v$ de $K$. Consid\'erons maintenant le sous-complexe  simplicial $K'$
obtenu en \^{o}tant \`a $K$ chaque sommet  d'ext\'emit\'e et l'ar\^ete qui lui est
attach\'ee. Toutes les ar\^etes de $K'$ sont  du type
$\;\bulletun \;Ê\leftrightarrow\;
\bulletun\;$ et l'on a~: $E_s \iso E_{ss'} \iso E_{s'}\,$. Ces espaces
\'etant tous de dimension 1, nous les identifions \`a $\C\,$. ainsi $Z^0 (K')\,$,
$Z^1(K')\,$ et $Z^1( K)\,$ peuvent \^etre consid\'er\'es comme des espaces  de  cochaines
simpliciales \`a valeurs dans
$\C\,$~:
$$ Z^0 (K' ) \; \iso \;Z^0 (K' \, ; \C) := \bigoplus_{s \in \gS(K')}\,
\C\cdot {s}\, ,
$$
$$ Z^1 (K') \; \iso \; Z^1(K' \, ; \, \C ) =
\bigoplus _{ss' \in {\gS}(K')^{\check{2}} } \, \C\cdot {ss'} \, ,
$$
$$ Z^1 (K) \; \iso \; Z^1(K \, ; \, \C ) =
\bigoplus _{ss' \in {\gS}(K)^{\check{2}} } \, \C\cdot {ss'} \, ,
$$ o\`u ${\gS}(K')^{\check{2}}\,$ et ${\gS}(K)^{\check{2}}\,$ ont \'et\'e d\'efinis  en
(\ref{prodchech}). D'autre part, $E_s = 0$ lorsque $s$ est un sommet extr\'emit\'e de
$K\,$. ainsi nous pouvons identifier $Z^0 (K')$ \`a $Z^0 (K)\,$ et
$\Delta_K$ peut \^etre vue comme une application de $Z^0 (K' ; \C)$ dans $Z^1 (K ; \C)\,$.
On obtient le  diagramme commutatif suivant dont les lignes sont exactes~:
\begin{equation}\label{ssuuiiex}
\begin{array}{ccccccccc}  0 & \fle & Z^0(K'; \C )    & \stackrel{\iota^0}{\fle} & Z^0 (C
\, ; \, \C  )  & \stackrel{\pi^0}{\fle} &
\bigoplus_{j = 0}^v\,\C \cdot s_j          & \fle & 0  \\
  &      & \Delta_K \downarrow \phantom{\r} &       & \d_C^1\downarrow
\phantom{\d_C^1}          &      & \s \downarrow
\phantom{\s}           &      &
\\   0 & \fle & Z^1 (K ; \C )     & \stackrel{\iota^1}{\fle} & Z^1 (C \, ; \,
\C) & \stackrel{\pi^1}{\fle} &
\bigoplus _{j = 1}^v \, \C\cdot {s_js_0} & \fle & 0
\end{array}
\end{equation}
\noindent avec $\pi^0$ et $\pi^1$ les projections lin\'eaires canoniques,
$\iota^0$ et $\iota^1$ les injections lin\'eaires induites par les inclusions de
sous-complexes, et
${\d_C^1}$ l'application de cobord $\;\sum {\l_s \cdot s} \efle \sum (\l_s - \l_{s'}
)Ê\cdot ss'\,$  dont la restriction \`a $\bigoplus _{j = 0}^v \,
\C \cdot {s_j}$ est
$$\s (\, \sum_{j = 0}^v \l_j\cdot s_j\, ) \; = \; \sum_{j = 0}^v ( \l_j -
\l_0 ) \cdot s_js_0\, .$$
\noindent On obtient la suite exacte~:
\begin{equation}\label{sexsipli}
\begin{array}{cccccccc} ker({\d_C^1})  & \fle & \ker(\s)  & \fle & \cal  H^1 (K) &
\fle H^1 (C\, ; \, \C) & \fle & coker(\s)
\end{array}
\end{equation}
\noindent On voit facilement~:
$$ ker({\d_C^1}) = \C \cdot \sum_{s \in {\gS}(C)} s\,, \qquad \ker(\s) = \C
\cdot \sum_{s = 0}^v s_j
$$ o\`u ${\gS}(C)$ est l'ensemble des sommets de $C\,$, et la premi\`ere fl\`eche  de
(\ref{sexsipli}) est un isomorphisme. Comme
$\s$ est visiblement surjective, on obtient la conclusion.
\end{proof}

Ceci ach\`eve la d\'emonstration du th\'eor\`eme de codimension (\ref{critfr})

\subsection{Feuilletages de deuxi\`eme esp\`ece non-d\'eg\'en\'er\'es}
\addcontentsline{toc}{section}{\hspace{0,8em} {}\thesubsection .  Feuilletages de deuxi\`eme esp\`ece non-d\'eg\'en\'er\'es}
Mettons maintenant en \'evidence une classe  de feuilletages t.f.f. pour
lesquels $\htau (\F )$ peut se calculer seulement
 \`a partir de la donn\'ee de l'arbre dual $\aA (\F ) \,$.\\

\indent Pour cela munissons $\aA (\F)\,$ de la m\'etrique pour laquelle les ar\^etes sont
isom\'etriques \`a  l'intervalle
$[0 , 1]\,$. Rappelons que la {\it valence d'un sommet}
$s$ de
$\aA (\F)\,$, not\'ee
$v(s)\,$, est le nombre de fl\`eches et d'ar\^etes attach\'ees \`a $s\,$. On  appellera
aussi {\it valence d'une compo\-sante}
$D$ du diviseur exceptionnel $\D_\F$ et l'on notera $v(D)$ la valence du sommet  de $\aA (\F)$  corres\-pondant. Nous
appellerons {\it chaine de $\aA (\F)\,$} toute g\'eod\'esique de $\aA (\F)\,$ qui joint
deux sommets de valence
$\geq 3\,$. Une chaine de $\aA(\F)$ peut \^etre \'eventuellement r\'eduite \`a un
segment de longueur 1. A une chaine de $\aA (\F)\,$ correspond biunivoquement, soit une
union connexe maximale de diviseurs de $\D_\F$ de valence 2, soit un point d'intersection
de deux composantes de $\D_\F$ de valence $\geq 3$, appel\'ee dans les deux cas {\it
chaine de
$\D_\F$}.

\begin{defi}\label{nondeg} Un feuilletage formel de deuxi\`eme esp\`ece (\ref{especes}) $\F$ \`a l'origine de
$\C^2$  est dit {\it non-d\'eg\'en\'er\'e} s'il satisfait les conditions  suivantes~:
\begin{enumerate}
\item le groupe d'holonomie $H_D$ de toute composante $D$ de $\D_\F$ de  valence $\geq 3$
est non-commutatif,
\item pour toute chaine $\cal C$ de $\D_\F$ on a~: $\hO_\wF (\cal C) = \C\,$.
\end{enumerate}
\end{defi}
\noindent Remarquons que la condition de deuxi\`eme esp\`ece  implique que tout germe  d'int\'egrale premi\`ere
$f \in \hO_{\F , \, m}$ en un point singulier
$m$ situ\'e sur une chaine
$\cal C$ de $\D_\F\,$, se prolonge \`a toute la chaine~: $f \in \hO_\F (\cal  C ) \,$.
Ainsi on obtient une d\'efinition \'equivalente en rempla\c cant la condition 3. par~:
\begin{enumerate}
\item[2'.]\label{condsup}  pour toute chaine $\cal C$ de $\D_\F\,$ et tout un ouvert $U$ d'un recouvrement distingu\'e $\U$
de
$\D_\F$  intersectant $\cal C$ on a~:
$\hO_\F (U ) = \C\,$.
\end{enumerate}
Sur les chaines $\cal C$ de $\D_\F$ non-r\'eduites \`a un point, cette condition
est aussi \'equivalente \`a la non-finitude du groupe d'holonomie d'une composante
irr\'eductible quelconque de
$\cal C$.\\

On a le corollaire suivant au th\'eor\`eme de codimension  (\ref{critfr})~:
\begin{coro}\label{codimsshyp} Soit $\F$ un feuilletage de deuxi\`eme esp\`ece non-d\'eg\'en\'er\'e \`a l'origine de
$\C^2\,$. Alors
$\F$ est t.f.f., satisfait l'\'egalit\'e
$\,\hb (\F ) = \htau (\F) + \hd (\F)\,$ et $\,\htau (\F)$ est \'egal \`a la  somme du
nombre de chaines de $\aA (\F)\,$.
\end{coro}
\begin{preuvede}{du corollaire}  La condition 1. de la d\'efinition ci-dessus exclu, \-d'apr\`es
\cite{C-M}, l'existence d'un facteur int\'egrant formel de $\F$. La proposition (\ref{splitcod}) donne alors l'\'egalite
$\,\hb (\F ) = \htau (\F) + \hd (\F)\,$.
Pour achever la
d\'emonstration, il suffit de montrer  que
$\hgR^\ast (\F)$ est connexe r\'epulsif et de calculer
$\htau (\F)\,$.\\

\indent Examinons maintenant les propri\'et\'es du nerf $\hgN^\ast (\F)$  induites par les
conditions de non-d\'eg\'en\'erescence. Nous reprenons les notions de  branches s\'ecable et
d'\'elagage introduites dans la d\'emonstration du th\'eor\`eme (\ref{critfr}).
D'apr\`es
ce qui pr\'ec\`ede, en  examinant les listes (\ref{casposs}) et (\ref{posco}) on voit
que~:
\begin{enumerate}
\item[(a)]\label{aabbcc} tout sommet de valence $\geq 3$ est rouge et de poids 0, tout sommet rouge de
valence $\leq 2$  est de poids 1.
\item[(b)] toute ar\^ete verte est~: soit une simple fl\`eche orient\'ee vers son  sommet de type
0, qui est aussi vert, soit une double fl\`eche,
\item[(c)] toute ar\^ete verte de $\hgN^\ast (\F)$ qui joint deux sommets de  valence $\leq 2$
est une double fl\`eche joignant deux sommets verts,
\end{enumerate}
\noindent On en d\'eduit que les sommets $s_j\,$, $j = 1 , \ldots , r\,$,  extr\'emit\'es
de $\hgN^\ast (\F)\,$ sont aussi chacun l'extr\'emit\'e d'une branche s\'ecable rouge ou bien
verte $B_j\,$. Les sommets  d'extr\'emit\'es du sous-graphe $\hgN_1$ de $
\hgN^\ast (\F)$ obtenu en \'elagant $
\hgN^\ast (\F)$ des branches $B_1, \ldots , B_r$, \'etaient visiblement, dans
$\hgN^\ast (\F)\,$, des sommets de valence
$\geq 3\,$ et donc des sommets rouges de poids 0. Ainsi l'\'elagage total est effectu\'e en une seule
\'etape~:
$$ El\left( \hgN^\ast (\F) \right) = \hgN_1\,.
$$
\noindent  Chaque ar\^ete de $\hgN_1$ fait visiblement partie d'une
g\'eod\'esique $\cal C$ dont les sommets sont de valence 2 dans $
\hgN^\ast (\F)\,$, sauf les sommets d'extr\'emit\'es qui sont de valence
$\geq 3$ dans $ \hgN^\ast (\F)$. Une telle g\'eod\'esique correspond \`a  une
chaine de $\aA(\F)$.  D'apr\`es la condition 2'. de (\ref{nondeg}) les sommets et
les ar\^etes de $\cal C$ sont rouges  et d'apr\`es (a) sont tous de poids 1,
sauf les sommets  d'extr\'emit\'es qui sont de poids 0. Ainsi tous les sommets et ar\^etes de $\hgN_1$ sont rouges. Les
sommets de valence $\geq 3$ ainsi que les sommets d'extr\'emit\'es sont de poids 0. Tous les autres sommets ainsi
que les ar\`etes sont de poids 1. Les g\'eod\'esiques de $\hgN_1$ dont les sommets sont tous de poids 1 sauf les
extr\'emit\'es qui sont de poids 0, correspondent biunivoquement aux chaines de $\aA(\F)$. La formule
(\ref{nbrbousph}) et le th\'eor\`eme de codimension  (\ref{critfr}) montrent que le nombre
de ces g\'eod\'esiques est
$\htau(\F)$. Ceci ach\`eve la d\'emonstration.
\end{preuvede}

\subsection{D\'emonstration du th\'eor\`eme de pr\'eparation ($\ref{enimdir}$)}\label{primdir}
\addcontentsline{toc}{section}{\hspace{0,8em} {}\thesubsection .  D\'emonstration du th\'eor\`eme de pr\'eparation}
Conservons les notations fix\'ees en t\^ete de ce chapitre 4. En particulier le compos\'e de l'application
$E_{\sF_P} : \M_{\sF_P} \fle \C^2\times P$ d'\'equir\'eduction  de $\sF_P$ et de la
projection de $\C^2 \times P$ sur $P$ est not\'e
$$
\pi_{\sF_P} : \M_{\sF_P} \fle P \simeq \left( \C^p, 0 \right),\quad %\hbox{ou encore}\quad
\pi_{\sF_P} =(t_1, \ldots ,t_p)\,, \quad t_j := pr_j\circ\pi_{\sF_P}\,,
$$ o\`u $pr_j : \C^p \fle \C$ d\'esigne la $j$-i\`eme projection. Nous identifions
toujours la cime $\M_\F$ de la r\'eduction de $\F$ \`a la fibre $\pi_{\sF_P}^{-1}(0)\,$, et
le diviseur exceptionnel $\D_\F
\subset \M_{\F}$ \`a l'intersection de
$\pi_{\sF_P}^{-1}(0)$ et du diviseur exceptionnel $\D_{\sF_P}
\subset \M_{\sF_P}$. On peut donc \'ecrire~:
$\M_\F \subset \M_{\sF_P}$ et $\D_\F \subset \D_{\sF_P} \subset \M_{\F_P}$. Par
hypoth\`ese, pour chaque ouvert $U$ du recouvrement distingu\'e
$\U$ de $\D_\F$, le germe de $\wF_P$ le long de $U$ est conjugu\'e au germe de la
d\'eformation constante. On dispose ainsi de germes de r\'etractions
$$ R_U := \left( \M_{\sF_P} ,  U \right) \fle \left( \M_{\F} ,  U \right) \, , \qquad U
\in \U\,,
$$ tels que l'image r\'eciproque
\footnote{
Il s'agit du faisceau engendr\'e par les images r\'eciproques par $R_U$ des germes
formes diff\'erentielles de $\hL_\wF$, cf. (\ref{trans.strict})
}
$R_U^\ast\left( \wF \right)$ de $\wF$ et l'application
tangente de $\pi_{\sF_P}$  d\'efinissent
$\wF_P$ au voisinage de
$U$, c'est \`a dire~: $
\hL_{\wF_P} = R_U^\ast\left(\hL_\wF \right) + \left( dt_1 , \ldots ,dt_p\right)$.
Consid\'erons le faisceau d'anneaux
$\underline{\hO}_{\wF_P}$ de base $\D_{\sF_P}$ constitu\'e des germes $f \in
\hO_{\M_{\wF_P}, m}$ de fonctions transversalement formelles le long de
$\D_{\sF_P}
$ qui sont des {\it int\'egrales premi\`eres} de $\wF_P$, c'est \`a dire qui v\'erifient
$df \wedge \n =0$, pour tout germe
$\n \in\hL_{\wF_P, m}$. Restreignons ce faisceau au dessus de $0$ et consid\'erons le
faisceau suivant
$$
\hO_{\wF_P} := i^{-1}\left(\underline{\hO}_{\wF_P} \right)\, , \quad \hbox{o\`u}\quad i:
\D_\F \ifle
\D_{\sF_P}\,,
\quad i(m) := m\,
$$ qui a m\^emes fibres que $\underline{\hO}_{\wF_P}$ mais est de base $\D_\F$. Les
co-morphismes
$$ R_U^\ast : \hO_\wF (U) \ifle \hO_{\wF_P}(U) \quad \hbox{et} \quad \pi_{\sF_P}^\ast :
\O_P \ifle
\hO_{\wF_P}(U)\, , \quad U \in \U \,,
$$ permettent de  consid\'erer $\O_P$ et $\hO_\wF(U)$ comme des sous-anneaux de
$\hO_{\wF_P}(U)$.  On voit que les seules \'eventualit\'es sont~:
\begin{itemize}
\item $\hO_\wF(U) = \C \quad \hbox{et} \quad \hO_{\wF_P}(U) = \pi_{\sF_P}^\ast\left(
\O_P\right) = \C\left\{ t  \right\} $,
\item $\hO_\wF(U) = \CÊ[[f_U]]\quad \hbox{et} \quad
\hO_{\wF_P}(U) = \pi_{\sF_P}^\ast\left(
\O_P\right)[[ f_U\circ R_U ]] = \C\left\{ t  \right\} [[ f\circ R_U ]] $,
\end{itemize} o\`u  $\pi_{\sF_P}^\ast\left(
\O_P\right)[[f_U\circ R_U]]$ d\'esigne le compl\'et\'e de
$\pi_{\sF_P}^\ast\left(\O_P\right)[f_U\circ R_U]$  dans $\hO_{\wF_P}(U)$ pour la topologie
$\left(f_U\circ R_U
\right)$-adique et o\`u $\C\left\{ t \right\}$ est l'anneau des s\'eries convergentes en
les variables $t_1,  \ldots , t_p$. \\

Tout champ de vecteurs basique $X\in \hB_\wF(U)$, $U \in \U$, se rel\`eve par
$R_U$ en un unique champ basique et vertical not\'e $R_U^\ast(X)$, gr\^ace encore \`a la
trivialit\'e de $\wF_P$ le long de U. Pour le $\C\{ t \}$-module $\hT_{\wF_P}(U)$, les
seules
\'eventualit\'es sont~:

\begin{enumerate}
\item
$\hO_\F(U) = \C[[f_U]]$ et $\hT_\wF(U) = \C[[f_U]]\cdot \left\{X_U \right\}$, o\`u $f_U $
est une int\'egrale premi\`ere non-constante et o\`u $X_U \in \hB_\F(U)$ d\'efinit un
\'el\'ement non-nul apropri\'e $\left\{X_U\right\}$ de $\hT_\F(U)$ . Alors
$\hO_{\wF_P}(U) =\C\{t\}[[F_U]]$ et  $\hT_{\wF_P}(U) = \C\{t\}[[F_U]]   \cdot \left\{
Z_U\right\} $, avec $F_U := f_U
\circ R_U \notin \C\{ t\}$ et $Z_U := R_U^\ast (X_U) $ d\'efinissant un \'el\'ement
$\left\{ Z_U\right\}$ non-nul de
$\hT_{\wF_P}(U)$.
\item
$\hO_\F(U) = \C$ et $\hT_\wF(U) = \C \cdot \left\{X_U \right\}$ avec $\left\{ X_U\right\}
\not=\{ 0\}$. Alors
$\hO_{\wF_P}(U) =\C\{t\}$ et  $\hT_{\wF_P}(U) = \C\{t\} \cdot \left\{ Z_U\right\} $,
avec  $Z_U :=  R_U^\ast (X_U) $ et $\left\{ Z_U \right\} \not= \{ 0 \}$.
\item
$\hO_\F(U) = \C$ et $\hT_\wF(U) = \left\{ 0 \right\}$. Alors
$\hO_{\wF_P}(U) =\C\{t\}$ et $\hT_{\wF_P}(U) = \left\{0 \right\}$.
\end{enumerate}

En d\'esignant toujours par $\ast$ un sommet $s$ ou une ar\^ete $ss'$ du nerf
$\hgN^\ast(\F)$ de $\F$, posons maintenant~:
$$
\cal A_\ast := \hO_{\wF_P}(U_\ast) \qquad \hbox{et} \qquad \cal E_\ast := \hT_{\wF_P}
(U_\ast )\,,
$$
o\`u $U_\ast$ est d\'efini comme en (\ref{subsenerf}). Nous avons seulement, comme en
(\ref{subsenerf}), les trois possibilit\'es~:
\begin{enumerate}
\item $\cal A_\ast = \C\{ t\}\,$ et $\cal E_\ast = \{ 0 \}\,$,
\item $\cal A_\ast = \C\{ t\}\,$ et $\cal E_\ast = \C\{ t\} \cdot \cal Z_\ast\,$, avec $
\cal Z_\ast \not=  0\,$
\item $\cal A_\ast = \C\{ t\}[[F_\ast]]\,$ et $\cal E_\ast = \C\{ t\}[[F_\ast]] \cdot \cal
Z_\ast\,$,  avec
$\cal A_\ast \not= \C\{ t\} \,$, $ \cal Z_\ast \not= 0\,$.
\end{enumerate}
Pour chaque ar\^ete $ss'$ nous disposons, comme en (\ref{aplrstr}),
d'op\'erations de  restrictions
$$
\goth r_{ss'}^s : \cal E_s \ifle \cal E_{ss'}\, \qquad \hbox{et}\qquad \goth s_{ss'}^s :
\cal A_s
\ifle \cal A_{ss'}\,,
$$ qui sont maintenant des morphismes injectifs de $\C\{ t\}$-modules. De plus, toujours
\`a cause de la trivialit\'e de $\wF_P$  le long de
$U_\ast$, le morphisme $\goth r_{ss'}^s$, resp. $\goth s_{ss'}^s$, est bijectif si et
seulement si $\r_{ss'}^s$, resp.
$\s_{ss'}^s$, l'est. La {\it $\wF_P$-application de cohomologie} associ\'ee \`a un un
sous-graphe $K$ de
$\hgN^\ast(\F)$, d\'efinie comme en (\ref{eqchomo}) par~:
$$
\begin{array}{rcl}
\goth D_K : \prod_{s \in \gS(K)} \cal E_s  =:   \goth Z^0(K) &\fle &\goth Z^1(K)  :=
\prod_{(s, s') \in  {\gS}(K)^{\check{2}}}
\cal E_{ss'}\, ,\\
\left( \cal Z_s\right)_s & \efle &\left( \cal Z_{s'} - \cal Z_s\right)_{(s, s')}\,,
\end{array}
$$
est maintenant un morphisme de $\C\{ t\}$-modules. En notant $\cal H_P^1(K) :=
coker(\goth D_K)$, il vient~:
$$ H^1 \left( \U ; \, \hT_{\wF_P}\right) = \cal H_P^1\left( \hgN_P^\ast \left( \F
\right)\right)\,.
$$ Le proc\'ed\'e d'\'elagage des branches vertes et rouges s\'ecables utilis\'e dans la
d\'emonstration du th\'eor\`eme (\ref{critfr}) se retranscrit ici litt\'eralement~: le lemme
(\ref{secvert}) reste vrai lorqu'on remplace $\cal H^1( K )$, $\H^1(K')$ par
$\cal H_P^1( K )$, $\H_P^1(K')$ et, sous l'hypoth\`ese que $\F$ est t.f.f., on a encore
$$
\H_P^1 \left( \hgN^\ast \left( \F \right) \right) \simeq
\H_P^1 \left(\hgR^\ast \left( \F \right)\right) \simeq
\H_P^1 \left( El \left( \hgR^\ast \left( \F \right) \right) \right) \,.
$$ Les calculs de la d\'emonstration du lemme (\ref{cohorouge}) se retranscrivent aussi. Les $\C\{
t\}$-modules
$\goth Z^0 \left( K' \right)$, $\goth Z^1 \left( K' \right)$ et  $\goth Z^1  \left( K
\right)$ sont respectivement isomorphes  aux espaces de cochaines simpliciales
$$
 Z^0 \left( K' ; \, \C \{ t \} \right) := \bigoplus_{s \in \gS(K')}\, \C\{ t \}\cdot
{s}\,,\quad
Z^1 \left( K' ; \, \C \{ t \} \right) := \bigoplus _{ss' \in {\gS}(K')^{\check{2}} } \,
\C\{ t \}\cdot {ss'} \,,$$
$$
Z^1\left( K ;\, \C \{ t \} \right) := \bigoplus _{ss' \in {\gS}(K)^{\check{2}} } \, \C\{
t \} \cdot {ss'} \,.
$$
On obtient un diagramme similaire au diagramme (\ref{ssuuiiex}), dont la suite exacte
longue associ\'ee produit un isomorphisme
$\C\{ t \}$-lin\'eaire entre
$\cal H_P^1 \left( K \right)$ et $H^1 \left( C ; \, \C\{ t\} \right)$. Comme $H^1 \left( C
; \, \C\{ t\} \right) \simeq H^1 \left( C ; \, \C \right) \otimes_\C \C\{ t\}$, on peut
finalement conclure que $H^1 \left( \U ; \, \hT_{\wF_P} \right)$ est un $\C\{ t\}$-module
libre de rang \'egal \`a $dim_\Z H^1\left( \goth C(\F) ;\, \Z\right)$, c'est \`a dire de
rang $\widehat{\t} (\F)$, o\`u $\goth C(\F)$ est le bouquet de cercles d\'efini en (\ref{bouquet}).\\

Pour voir que $H^1\left( \U ; \, \hB_{\wF_P}^v \right)$ est de type fini sur $\C\{t \}$,
il suffit maintenant de rappeler \cite{Minv} que
$H^1\left( \U ; \, \hX_{\wF_P} \right)$ est libre de type fini et de consid\'erer les
derniers termes de la suite exacte longue  (\ref{dlousski}). Lorsque $\F$ n'admet pas de
facteur int\'egrant formel, on voit facilement \`a l'aide de (\ref{ffint}) que
$H^0\left( \U ;
\hT_{\wF_P}
\right)$ est nul. Ainsi (\ref{dlousski}) montre que $H^1\left( \U ; \, \hB^v_{\wF_P}
\right)$  est la somme directe de
$H^1\left(
\U ;
\,
\hX_{\wF_P}
\right)$ et de $H^1\left( \U ; \, \hT_{\wF_P} \right)$.
%%%%%%%%%%%%%%%%%%%%%%%%%%%%%%%%%%%%%%%%%%%%%%%%%%%%%%%%%%%%%%%%%%%%%

\section{G\'en\'ericit\'e des singularit\'es t.f.f.}\label{tffgen}
\subsection{Jets d\'eterminant les singularit\'es de deuxi\`eme esp\`ece}
\addcontentsline{toc}{section}{\hspace{0,8em} {}\thesubsection .  Jets d\'eterminant la deuxi\`eme esp\`ece}
Consid\'erons l'espace $\hL_{\C^2,\,0}$ comme la limite projective des espaces de jets, via les applications canoniques
$ \,j^k :  \hL_{\,\C^2,\,0}\fl  J_0^k\L := \left. \hL_{\,\C^2,\,0}
\,\right/\widehat{\goth m}_2\!\cdot\!\hL_{\C^2,\,0}\,,
$
o\`u $\widehat{\goth m}_2$ d\'esigne l'id\'eal maximal de $\hO_{\C^2,\,0}$. A l'aide des coordonn\'ees canoniques identifions
$J_0^k\L$ \`a l'espace vectoriel
$$
\cal P_k := \left\{ \left. Ê(P, Q) \in \C [X, Y] \;\right/ \; deg(P), deg(Q) \leq k \right\}\,.
$$
 \noindent Classiquement pour chaque entier
$r\,$, le sous-ensemble $L_{k, \, r} \subset \cal P_k\,$ des couples $(P, Q)\,$ qui engendrent un id\'eal de $\hO_{\C^2 , \, 0}$
de {\it colongueur}
$$
\mu (P, Q) := dim_\C \left( \hO_{\C^2 , 0} / (P, Q) \right)
$$
sup\'erieure ou \'egale \`a $r$, est alg\'ebrique ferm\'e. En raffinant la partition de $\cal P_k$ par les composantes connexes
des diff\'erences $L_{k, \,r + 1} - L_{k\, , r}$, on obtient une  une stratification $\underline{\cal
E}_\mu^k$ de
$\cal P_k\,$, localement finie pour la topologie de zariski, telle que chaque strate est un ensemble constructible lisse
le long duquel
$\mu (P, Q)$ est constant.\\

\indent Consid\'erons maintenant la famille "lin\'eaire" de feuilletages $\sF_{\cal P_k}\,$, d'espace (global) de param\`etres
$\cal P_k\,$, d\'efinie par~:
$$
\Omega_{P,\, Q} :=  \widetilde{E}_P(x,y)\,dx + \widetilde{E}_Q(x,y)\,dy \quad\mbox{o\`u}\quad  \widetilde{E}_F(x, y)  :=
F(x, y)\,,Ê\quad F \in \cal P_k\,.
$$
Il d\'ecoule\footnote{
La suite de sous-espaces vectoriels  de codimension finie  $F_j:= (P, Q) + \widehat{\goth m}_2^jÊ\subset \hL_{\,\C^2,
0}$ est d\'ecroissante et
$F_r = F_{r+1}$ pour un entier $r \leq \mu (P,Q)$.  On conclut par le lemme de Nakayama.
}
de l'inclusion $\,\widehat{\goth m}_2^{\mu (P,Q)} \subset
(P, Q)\,$ que l'ouvert
$$
W_{\mu.det}^k \; := \; J_0^k\hL \; - \; L_{k , \, k + 1} \, ,
$$
est {\it $k$-d\'eterminant pour le nombre de Milnor}, c'est \`a dire~: $\mu (\w)   =   \mu
(\n)   <   \infty$ d\`es que $
j^k\w =  j^k\n \in   W_{\mu.det}^k\,$.
Visiblement $W_{\mu.det}^k $ est une union de strates de $\underline{\cal E}_\mu^k$ et la famille $\cal W^k_{\mu.det} :=
(j^k)^{- 1}(W_{\mu.det}^k)$ v\'erifie~:
$$
\cal W_{\mu.det}^k \, \subset \cal W_{\mu.det}^{k + 1}\,,\quad k \, \in \N
\quad \hbox{et} \quad \bigcup_{j\in\N} \cal W_{\mu.det}^j\; = \;
\hL_{sat}\;,
$$
\noindent o\`u $\hL_{sat}$ d\'esigne l'ensemble des 1-formes formelles \`a singularit\'e isol\'ee.\\

D\'esignons par $\nu_\w$ la multiplicit\'e alg\'ebrique \`a l'origine de $\w \in \hL_{sat}$ et par $h'_\w$ la hauteur de son arbre de
pr\'e-r\'eduction (\ref{preredu}). Rappelons qu'apr\`es un  \'eclatement la multiplicit\'e du transform\'e strict de $\w$ en un point
quelconque du diviseur exceptionnel est au plus $\nu_\w + 1$.
 D\'efinissons alors la fonction
$\s : {\N^{\ast}}^2 \fle \N$  par la relation de r\'ecurence~:
$$
\left\{
\begin{array}{l}
\s (v, 1) := v + 1\,, \qquad \s (v , 2 ) := 2 v +1 \, ,\\
\s (v , h+2 ) = \s ( v , h + 1 ) + \s ( v , h ) + v + h\,.
\end{array}
\right.
$$
\noindent On voit facilement
qu'en chaque point $m$ du diviseur de c\^{\i}me $\wD'_{\w}$ de la
pr\'e-r\'eduction de $\w\,$, le transform\'e strict de $\F_{\w}$ s'obtient en divisant l'image r\'eciproque de
$ \w$ par un facteur $u^pv^{\e q}$, $\e = 0$ ou 1, avec~:
$$
p, \; q \; \leq \; \s \left(\nu_\w , h'_\w\right)\,.
$$
\noindent $uv^\e =0$ \'etant une \'equation r\'eduite locale de $\wD'_\w$ en $m$. Comme pour tout entier
$r\,$, le
$r$-jet\footnote{
Deux formes formelles le long d'un diviseur $\wD$ {\it ont m\^eme $r$-jet le long de $\wD$} si elles diff\`erent d'une forme \`a
coefficients dans $\cal I_{\wD}^{r + 1}$, o\`u $\cal I_{\wD}$ d\'esigne le faisceau des fonctions transversallement formelles le
long de
$\wD$ et nulles sur $\wD$.
}
de l'image r\'eciproque de
$\w$ le long de
$\cal D'_\w$ ne d\'epend que du
$r$-jet de
$\w$ \`a l'origine, on obtient~:

\begin{lemm}\label{jdetcri}
Soit $\w$ une 1-forme formelle \`a singularit\'e isol\'ee et $k\,$, $l\,$ des entiers qui v\'erifient~:
$
0 \; <\; l < \;  k\; -\; \s \left(\nu_\w ,  h'_\w \right)\,
$. Alors $j^k\w$ est $l$-d\'eterminant pour la pr\'e-r\'eduction des singularit\'es, i.e. $\n$ et $\w$ ont m\^eme pr\'e-r\'eduction des
singu\-larit\'es d\`es que $j^k\n = j^kÊ\w$.
\end{lemm}

A l'aide du le th\'eor\`eme (\ref{geredd}) appliqu\'e \`a la famille
$\left( \Omega_{P,
\,Q}\right)_{(P, Q)\in S}$ on peut construire une filtration
d\'ecroissante de chaque strate $S \in \underline{\cal E}_\mu^k$ par des sous-ensembles alg\'ebri\-ques ferm\'es
$L'_{k, r+1}(S)$ telle que les diff\'erences $\Delta_{k, \, r} (S) := L'_{k, r}(S) -
L'_{k, \, r - 1 } (S)$ sont lisses et telle ,que lorsque $(P, Q)$ varie dans $\Delta_{k, \, r} (S)\,$, le germe de
$\Omega_{P,
\,Q}\,$ en chaque point de $\Delta_{k, \, r} (S)$ d\'efinit une d\'eformation \'equi-pr\'e-r\'eductible . Notons $\underline{\cal
E}_{p.red}^k$ la stratification de $\cal P_k$ dont les strates sont les composantes connexes des
$\Delta_{k,
\, r} (S)\,$, $S\in \underline{\cal E}_\mu^k$. \\

\indent Pour chaque strate $S\in\underline{\cal E}_{p.red}^k\,$, notons respectivement $h'_S$ et $\nu_S$ la hauteur de
l'arbre de pr\'e-r\'eduction et la multiplicit\'e \`a l'origine d'une
1-forme quelconque $\Omega_{P, \, Q}\,$, $(P, Q)\in S\,$.
Pour chaque entier
$l\,$, $0 \leq l < k\,$, consid\'erons l'ensemble constructible $L"_{k, \, l}\,$ form\'e de l'union des strates $S$ de
$\underline{\cal E}_{p.red}^k$ telles que
$
k\; -\; \s \left(\nu_S ,  h'_S\right)\; \leq \; l\,$.
Notons~:
$$
W_{p.red}^{k,\,l} \; := \; J_0^{k}\hL\, - \, (L_{k ,\, k + 1} \cup L"_{k, \, l}) \quad \mbox{et} \quad
\cal W_{p.red}^{k,\,l} := \left( j^k\right )^{- 1} \left( W_{p.red}^{k,\,l} \right)\, .
$$
\noindent Fixons l'entier $l \geq 1\,$ et notons $j^k_{k'} : J_0^{k'}\hL \fle J_0^k\hL$, $k \leq
k'$, les applications canoniques.  Il
d\'ecoule du lemme (\ref{jdetcri})~:
\begin{enumerate}
\item Si $j^k\w \in W_{p.red}^{k,\,l}\,$ alors $\w$ est \`a singularit\'e isol\'ee et son $k$-jet est
$l$-d\'eterminant pour la pr\'e-r\'eduction des singularit\'es;
\item Si $S$ est une strate de $\underline{\cal E}_{p.red}^k$ contenue dans $ W_{p.red}^{k,\,l}$ et  $k\leq k'\,$, alors
$(j^k_{k'})^{-1}(S)\,$ est une strate de $\underline{\cal E}_{p.red}^{k'}\,$ contenue dans $ W_{p.red}^{k',\,l}$; En
parti\-culier
$(j^k_{k'})^{-1}(W_{p.red}^{k,\,l})
\subset  W_{p.red}^{k',\,l}\,$;
\item On a ~: $\;\cal W_{p.red}^{k,\,l} \; \subset \cal W_{p.red}^{k' ,\,l}\, $, $\,k \leq k'\;\;$ et
$\;\;\bigcup_{\stackrel{\scriptstyle r \in \N}{\scriptstyle r > l}}
\cal W_{p.red}^{r,\,l}\; = \; \hL_{sat}\,$.
\end{enumerate}
\noindent Ainsi deux \'el\'ements distincts de la famille
$$
{\cal E}_{p.red} \, := \, \left\{ \left. (j^k)^{-1} (S) \subset \hL_{sat} \right/ S \in \underline{\cal E}^k_{p.red} \, , \,
k - \s(\nu_S , h'_S ) > 1\, , k \in \N
\right\}
$$
\noindent sont d'intersection vide. Nous pouvons maintenant appliquer le th\'eor\`eme (\ref{geredd}) \`a chaque
$S \in \underline{\cal E}^k_{p.red}$ avec $k - \s(\nu_S , h'_S ) > 1\,$. Finalement nous obtenons~:
\begin{theo}\label{prstrpred}
Pour chaque
$\cal S \in {\cal E}_{p.red}$ il existe une fonction croissante
$k_S : \N \fle \N$ telle que pour  $k \geq k_S(l)$  on a~:
\begin{enumerate}
\item l'ensemble $\underline S_k := j^k (\cal S)$ est constructible et $\cal S = (j^k)^{-1}(\underline S_k)\,$,
\item les jets d'ordre $k$ des \'el\'ements de $\cal S$ sont $l$-d\'eterminants pour la pr\'e-r\'eduction des singularit\'es,
\item une famille $\{ \w_t\}_t$ de 1-formes formelle d\'ependant analytiquement d'un param\`etre $t$ d\'efinit une
d\'eformation \'equi-pr\'e-r\'eductible de $\w_0\,$ si et seulement si $\w_t \in \cal S\,$.
\item le sous-ensemble $\cal S^{2}\subset\cal S$ des \'el\'ements de $\cal S$ qui sont des 1-formes de deuxi\`eme esp\`ece v\'erifie~:
\begin{enumerate}
\item $\cal S^{2} = (j^k)^{-1}(\underline{S}_k^{2})$ avec $\underline{S}_k^{2} := j^k(\cal S^{2})$
\item la diff\'erence $ \underline{S}_k - \underline{S}_k^{2}$ est une union d\'enombrable localement finie (dans
$\underline{S}_k$, pour la topologie usuelle) de sous ensembles alg\'ebriques et l'adh\'erence $\overline{\underline{S}_k -
\underline{S}_k^{2}}$ est semi-alg\'ebrique r\'eelle.
\end{enumerate}
\end{enumerate}
\end{theo}

%%%%%%%%%%%%%%%%%%%%%%%%%%%%%%%%%%%%%%%%%%%%%%%%%%%%%%%%%%%%%%%%%%%%%
\subsection{Krull-densit\'e des singularit\'es t.f.f.}
\addcontentsline{toc}{section}{\hspace{0,8em} {}\thesubsection .  Krull-densit\'e des singularit\'es t.f.f.}
Nous \'etudions ici la Krull densit\'e des feuilletages
t.f.f.. Nous nous restreignons ici aux feuilletages de deuxi\`eme esp\`ece, qui donnent des en\-nonc\'es
plus concis. Mais la Krull densit\'e peut s'obtenir pour des classes plus larges de feuilletages, en
permettant par exemple des singularit\'es de type selle-n\oe ud tangent r\'esonant (\ref{sntagent}). \\

D\'esignons par
$\widehat{\goth E}^{2}$,  resp. $\widehat{\goth E}_{tff}^{2}$, resp. $\widehat{\goth{ND}}$, les sous-ensembles de
$\hL_{\C^2,\,0}$ des 1-formes formelles de deuxi\`eme esp\`ece, resp. de deuxi\`eme esp\`ece t.f.f., resp.
non-d\'eg\'en\'er\'ees (\ref{nondeg}). D'apr\`es (\ref{codimsshyp}) on a $\widehat{\goth{ND}} \subset \widehat{\goth E}_{tff}^2$. Une
cons\'equence im\-m\'ediate du th\'eor\`eme
 pr\'ec\'edent (\ref{prstrpred}) est que $\widehat{\goth E}^2$ et
$\widehat{\goth{ND}}$ sont des ouverts de
$\hL_{\C^2,\, 0}$ pour la topologie de Krull.

\begin{enonce}{Th\'eor\`eme de g\'en\'ericit\'e}\label{mthm.gen.}
L'ouvert $\widehat{\goth{ND}}$ - et \`a fortiori $\widehat{\goth E}_{tff}^2$, est dense dans $\widehat{\goth{E}}^2$ pour
la topologie de Krull.
\end{enonce}
\noindent De mani\`ere plus explicite~:
\begin{theo}\label{genexpl}
Soit $\F$  un feuilletage formel de deuxi\`eme esp\`ece (\ref{especes})\`a l'origine de $\C^2$ donn\'e par une 1- forme formelle $
\w \in
\hL_{\C^2,0}
$ et soit $k$ un entier positif. Alors il existe une 1-forme $\w'\in \hL_{\C^2,\,0}$ et un entier $k' \geq k$ tels que~:
\begin{enumerate}
\item $\w'$ est non-d\'eg\'en\'er\'ee au sens de (\ref{nondeg}), \`a fortiori t.f.f., et a m\^eme $k$-jet que $\w'$~: $j^k\w' =
j^k\w$,
\item toute 1- forme $\n \in \hL_{\C^2,0}$ v\'erifiant l'\'egalit\'e $j^{k'}\n = j^{k'}\w'$ est aussi de deuxi\`eme esp\`ece
et non-d\'eg\'en\'er\'ee.
\end{enumerate}
\end{theo}
\begin{proof} L'assertion 2. est une reformulation du fait que $\widehat{\goth{ND}}$ est un ouvert pour la topologie de
Krull, ce qui se d\'eduit sans peine de th\'eor\`eme (\ref{prstrpred}). D\'emontrons l'assertion 1.\\

Notons ici
$\bC(\wF)$ l'ensemble
$\bC\left(\AÊ[\F]\right)$ des \'el\'ements critiques de l'arbre de r\'eduction de $\F$, d\'efini en (\ref{doncrit}). Les propri\'et\'es
de non-d\'eg\'en\'erescence sont "semi-locales", dans le sens o\`u elles se v\'erifient \`a partir seulement de la collection
 des germes $\wF_K$ le long de chaque ensemble $K \in \bC(\wF)$, du transform\'e strict $\wF$ de la r\'eduction de $\F$. La
d\'emonstration va consister \`a construire un syst\`eme semi-local coh\'erent (\ref{coherent}) de d\'eformations de~$\wF$
\begin{equation}\label{sslgen}
\bS := \left(\wF_{\C, \,K}\right)_{K\, \in\, \bC(\wF)}\,,
\end{equation}
d'espace de param\`etres $(\C, 0)$ tel que, lorsque l'on fixe le param\`etre $t$ assez petit  et non-nul, les feuilletages
t.f.
$\wF_{\C, \,K}(t)$ satisfont les proprit\'e\'es de non-d\'eg\'en\'erescence. De plus ces d\'eformations seront tangentes
\`a la d\'eformation constante,
\`a un ordre $l$ suffisamment grand pour que le th\'eor\`eme de r\'ealisation tangente (\ref{de.real.tg.}) donne le
r\'esultat.\\

Avant de commencer cette construction, fixons le vocabulaire que nous allons employer ici. Soit $L \subset  \D_\w$ une union
de composantes ir\'eductible du diviseur $\D_\w$ de la r\'eduction de $\F$. Notons $Comp(L)$ l'ensemble des composantes
irr\'eductibles de $L$. Pour $D \in Comp (L)$ nous appelons respectivement { \it $L$-valence de} $D$ et {\it $\wF$-valence de
$D$} les entiers
$$
v_L(D) := Card \left( D \cap Sing(L)\right) \,, \quad \mbox{resp.} \quad v_\wF(D) := Card \left( D \cap Sing(\wF)\right)\,.
$$
\noindent Les \'el\'ements de $Comp(L)$ de $L$-valence 1 s'appellent {\it extr\'emit\'es de} $L$. Une extr\'emit\'e $D$ de $L$
 sera appel\'ee {\it diviseur de t\^ete de $L$} si $v_\wF(D) \geq 3$. Soit $C$ une union connexe
maximale de composantes irr\'eductibles de $L$ de $\wF$-valences \'egales \`a 1 ou \`a 2. Si l'une des deux extr\'emit\'es de
$C$ est aussi une extr\'emit\'e de $L$ nous dirons que $C$ est une {\it $\wF$-branche de} $L$. Dans ce cas l'adh\'erence de
$L-C$ est encore connexe et $(\overline{L-C})Ê\cap L$ est r\'eduit \`a un seul point, appel\'e {\it point d'attache de $C$ sur
$L$}. Par contre, si $C$ ne contient pas d'extr\'emit\'e de $L$, on dit que $C$ est une {\it $\wF$-chaine de $L$}. Alors
$C$ est aussi une {\it chaine de
$L$}
\footnote
{i.e. une union connexe maxi\-male de composantes de $L$ de $L$-valences 2.
}, l'adh\'erence de $L-C$
n'est plus connexe et $(\overline{L-C})Ê\cap L$ est compos\'e de deux points appel\'es encore {\it points d'attache de $C$
sur $L$}.\\

 Pour toute la suite de la d\'emonstration fixons un \'el\'ement $D'_0Ê\in Comp(\D_\w)$ tel que~:
$ v_\wF(D'_0) > v_{\D_\w} (D'_0) $ et $ v_\wF (D'_0) \geq 3 $. Consid\'erons la filtration (finie) de $\D_\w$
$$
\D_\w =: L^0 \,\nsupset \, L^1 \, \nsupset\, \cdots \,\nsupset\, L^q = D'_0
$$
\noindent d\'efinie par~:
\begin{enumerate}
\item[a.] si l'union $Br(L^k)$ des $\wF$-branches de $L^k$ est non-vide, nous posons $L^{k+1} := L^k - Br( L^k)$,
\item[b.]si $Br(L^k) = \emptyset$, alors l'ensemble des diviseurs de t\^ete $\not= D'_0$
$$
{Tet}(L^k) := \left\{ D \in Comp(L^k)\; \left/ \;v_{L^k}(D) = 1 ,
\,v_\wF (D) \geq 3, \,D \not= D'_0  \right. \right\}
$$
\noindent est non-vide et nous posons $L^{k+1} := \overline{L^k - \cup_{D \in Tet(L^k)} \,D}$.
\end{enumerate}
\noindent L'op\'eration $L^k \efle  L^{k+1}$ sera appel\'ee un {\it \'elagage de $L^k$} dans le cas a. et un {\it \'et\^etage
de
$L^k$} dans le cas b.. La filtration $(L^j)_j$ de $\D_\w$ induit le filtration $\,\cal C^0 \nsubset \cal C^1 \nsubset
\cdots
\nsubset
\cal C^q\,$ de
$\bC(\wF)$ d\'efinie par ~:
$$
\cal C^k := \left\{  \left. K \in \bC(\wF) \,\right/\, K \subset \overline{\D_\w - L^{k+1} } \,
\right\}\,.
$$
\noindent C'est par induction suivant cette derni\`ere cette filtration que nous allons construire la col\-lec\-tion
(\ref{sslgen}) .  Nous nous appuierons sur la proposition (\ref{constrcol}) et sur le lemme suivant.
\begin{lemm}\label{lemmeun}
Soit $C \subset \D_\w$ une $\wF$-chaine de $\D_\w$ et $K_0$ un \'el\'ement critique de $\wF$ contenu dans $C$. Soit $l \in
\N$ et $\wF'_{\C, \,K_0}$ une  d\'eformation t.f. du germe $\wF_{K_0}$ de $\wF$ le long de $K_0$, r\'ealis\'ee sur
le germe d'espace produit $(\M_\F \times \C,\, K_0 \times \{0\}) $. Supposons de plus que $\wF'_{\C, \, K_0}$ est $l$-tangente
le long de
$(\D_\w\times \C\,, K_0 \times \{0\} )$
\`a la d\'eformation constante. Alors il existe une collection
\begin{equation}\label{deflem}
\wF'_{\C, \,K} \,, \quad \hbox{o\`u~:}\, \quad \quad K \subset C\,, \quad K \in \bC( \wF ) \, ,
\end{equation}
de germes, le long des ensembles critiques $K\subset C$, de d\'eformations t.f. des germes $\wF_K$ de $\wF$ le long de $K$,
qui sont $l$-tangentes \`a l'identit\'e le long de $\D_\w$ et qui v\'erifient les relations de coh\'erence (\ref{coherent}).
De plus les \'el\'ements de cette collection sont uniques  \`a conjugaison t.f. pr\`es.
\end{lemm}
\begin{proof}
La construction se fait de proche en proche gr\^ace  \`a la proposition (\ref{constrcol}) et au sous-lemme suivant dont la
d\'emonstration r\'esulte imm\'ediatement de la classification t.f. des singularit\'es r\'eduites, par l'holonomie de l'une des
vari\'et\'es invariantes (la vari\'et\'e forte dans le cas d'un selle n\oe ud), cf. \cite{Ra-Ma} et \cite{Ra-Ma.s.n.}.
\end{proof}

\begin{enonce}{Sous-lemme}\label{slemme2}
Soit $\F'$ un germe \`a l'origine d'un feuilletage t.f. le long de $S := \C \times 0 \subset \C^2$. Supposons que $S$ est un
ensemble  invariant de $\F'$ et que $\F'$ est r\'eduit mais n'est pas un selle-n\oe ud tangent \`a $S$. Alors, pour toute
d\'eformation t.f.
\footnote{i.e. toute s\'erie $h(z, t) = \underline h(z) + \sum_{j=1}^\infty  h_j(t)z^j \in \C\{ t\}[[z]]$ telle que les rayons de
convergences des coefficients $h_j(t)$  sont  uniform\'ement minor\'es par une constante $>0$.}
de l'holonomie $\uh(z) \in \C[[z]]$ de $\F'$ le long d'un petit lacet $\gamma (s) := (\e e^{2i\pi s} , 0)$, $s \in
[0, 1 ]$,  il existe un germe $\F'_\C$ de d\'eformation de $\F'$, t.f. le long de $S$, d'espace de param\`etres $(\C, 0)$,
unique
\`a conjugaison t.f.  pr\`es, telle que $h(z,t)$ soit l'holonomie de $\F'_\C$ le long de $\gamma $. De plus, si pour un entier
 $ r $ donn\'e on a $h(z, t) = \uh(z)$ $mod (z^{r+1})$, alors la d\'eformation $\F'_\C$ peut \^etre choisie $r$-tangente le
long de $S$ \`a la d\'eformation constante.
\end{enonce}
\noindent Nous sommes maintenant en mesure de commencer la construction de la collection (\ref{sslgen}).\\

\indent {$\bullet $ \it \underline{Premi\`ere \'etape ~:} la construction de $\wF_{\C , K}$ pour $K \in \cal C^0$ lorsque
$Br(\D_\w)
\not=
\emptyset$}. Pour $B \in Br(\D_\w)$, deux cas se pr\'esentent~:

\indent {\it a) les deux extr\'emit\'es de $B$ sont de valence \'egale  \`a 2.} Visiblement l'une d'elle (celle qui est
aussi extr\'emit\'e de $\D_\w$) porte un point singulier $c_0$ de $\wF$ qui n'est pas  une singularit\'e de $\D_\w$. Nous
choisissons pour $\wF_{\C, \,c_0}$ une d\'eformation t.f. telle que pour $t$ asez petit, le
feuilletage $\wF_{\C, \,c_0}(t)$ ne poss\`ede au voisinage de $c_0$ qu'un seul point singulier et qu'en ce point  il soit
r\'eduit, sans int\'egrale premi\`ere t.f. non-constante et tangent  \`a $\D_\w$. Nous exigeons de plus que $\wF_{\C,
\,c_0}$ soit tangente \`a la d\'eformation constante le long de $\D_\w$ \`a l'ordre $l$ voulu. Pour les autres ensembles
critiques $K \in
\cal C^0$, $K \subset B$, nous prenons pour $\wF_{\C, \,K}$ les d\'eformations donn\'ees par le lemme (\ref{lemmeun}).\\
\indent{\it b) $B$ poss\`ede une extr\'emit\'e de valence \'egale \`a 1.}  Il est facile de voir que $\wF$ poss\`ede alors le long de
tout
$B$ un germe d'int\'egrale premi\`ere t.f. non-constante $f_B \in \hO_{\wF}(B)$ qui s'annule sur le diviseur. Nous prenons
alors pour
$\wF_{\C, \,K}$,
$K
\in \cal C^0$, $K \subset B$, les d\'eformations constantes.\\
\indent Dans les deux cas a) et b), tra\c cons sur la composante $D_B$ de $L^1$ qui porte le point d'attache $c_B$ de $B$ \`a
$L^1$, un lacet $\wg$ bordant un disque conforme $W'$ tel que $\overline{W'} \cap Sing(\wF) = \{c_B\} $

\begin{lemm}\label{lemmetrois}
l'holonomie $\widetilde{h}(z,t) $ de $\wF_{\C,\, c_B}$ le long de $\wg$ n'est pas \'egale \`a l'identit\'e.
\end{lemm}
\begin{proof}
L'assertion est triviale dans le cas a). Dans le cas b) ordonnons les composantes de $B$ , $B= D_0 \cup \cdots \cup D_q$ de
sorte que $v_{\wF}(D_0) = 1$, $D_j \cap D_{j-1}$ est r\'eduit \`a un point $c_j$, $j = 1, \ldots , q$ et $D_q \cap D_B = \{
c_B \}$. \\

La formule d'indice de Camacho-Sad \cite{C-S} donne les relations suivantes entre les multiplicit\'es $m_j$ des $D_j$ comme
composantes irr\'eductibles des z\'eros de l'int\'egrale premi\`ere $f_B$~:
$$
\left\{
\begin{array}{r}
m_0 e_0 + m_1 \phantom{{}_{j+}}\,= 0 \\
m_0 + m_1 e_1 + m_2\phantom{{}_{j+}} = 0\\
\cdots\cdots\cdots\cdots\cdots\cdots\cdots\\
m_{j-1} + m_je_j + m_{j+1} = 0\\
\cdots\cdots\cdots\cdots\cdots\cdots\cdots\\
m_{q - 1} +  m_q  e_q +\, m\phantom{{}_{j+}} = 0
\end{array}
\right.
$$
\noindent o\`u $e_j$ d\'esigne l'auto-intersection $(D_j , D_j)$ et $m$ la multiplicit\'e du germe de courbe $(D_B, c_B)$
comme z\'ero de $f_B$. Il vient
$$
m_0 (e_0 + 2) + \cdots + m_q (e_q +2) =  m_0 + m_q - m \, .
$$
\noindent la minimalit\'e de la r\'eduction donne~: $e_j \leq - 2\,,\quad j = 0, \ldots ,q$. D'o\`u~:
$$
\frac{m_q}{m} \leq 1 - \frac{m_0}{m} < 1 \,.
$$
\noindent Mais l'holonomie $\widetilde{h}(z,t)$ est conjugu\'ee \`a $(z, t) \efle exp(2i\pi \frac{m_q}{m} \cdot z)$,  car la
d\'eformation $\wF_{\C, \,c_B}$ est constante. Ainsi $\widetilde{h}(z,t)$ est bien distinct de l'identit\'e.
\end{proof}

\indent {$\bullet $ \it \underline{Deuxi\`eme \'etape}~: la construction de $\wF_{\C , K}$ pour $K \in \cal C^0$ lorsque
$Br(\D_\w) =
\emptyset$}. Alors $\overline{\D_\w - L^1}$ est une union de composantes de t\^ete de $\D_\w$. Soit $D$ l'une d'elles. Notons
$c_1,
\ldots ,c_{r +1}$,
$r \geq 2$, les points singuliers de $\wF$ situ\'es sur $D$,
$c_{r+1}$ d\'esignant le point d'attache de $D$ \`a $L^1$. Donnons nous des lacets simples $\gamma_j$ d'origine $m_0$
commune, bordant des disques conformes $W_j \subset D$ tels que $\overline{W_j} \cap \overline{W_k} = \{m_0\}$, $j\not=k$ et
$\overline{W_j} \cap Sing(\wF) = \{m_j\}$, $m_j \in \inte W_j$, $j, k = 1, \ldots , r+1$. Nous choisissons ces
lacets pour que leurs classes $\ptsur\gamma_j$ dans
$\pi_1( D^\ast, m_0)$, $D^\ast := D -\{c_1, \ldots c_{r+1}\}$, v\'erifient : $\ptsur \gamma_1 \vee\cdots \vee \ptsur
\gamma_r = \ptsur \gamma_{r+1}^{\,\,\,-1}\,$. Ainsi les holonomies $\underline h_j(z)$ de $\wF$ le long des
$\gamma_j$v\'erifient~: $\underline h_1(z) \circ\cdots \circ \underline h_r(z) = \underline h_{r+1}^{\,-1}(z)\,$. Donn\'e un
entier $l \geq 1$, fixons des d\'eformations t.f. $h_j(z,t)$ des $\underline h_j(z)$, $j = 1, \ldots ,
r+1$, de sorte que $h_j(z, t) = \underline h_j(z) \quad mod (z)^{l+1}$ et que, pour $t\not= 0$ petit,
$h_{r+1,\, t} := \left( h_{1,\, t} \circ \cdots \circ h_{r,\, t}\right)^{-1}\,$  ne soit pas p\'eriodique  et que le groupe
engendr\'e par
$h_{1,
\,t}\,, \ldots , h_{r,\,  t}$ ne soit pas commutatif. Le lemme (\ref{slemme2}) et la proposition (\ref{constrcol})
d\'efinissent alors des d\'eformations $\wF_{\C, \,c_j}$, $j = 1, \ldots , r+1$ et $\wF_{\C, \,D^\ast}$ respectivements. Ces
d\'eformations,
$l$-tangentes \`a la d\'eformation constante, sont uniques \`a conjugaison t.f.  pr\`es. Nous d\'efinissons la collection $\wF_{K,
\C}$ $K \in \cal C^0$, en effectuant cette  construction pour chaque composante de $\overline{\D_\w - \cal L^1}$. \\

\begin{rema}\label{remetoi}  En chaque point d'attache $c'$ d'une composante de $\overline{\D_\w -
L^1}$ sur $ L^1$, pour $t\not=0$ petit, le feuilletage $\wF_{\C,\, c'}(t)$ ne poss\`ede pas d'int\'egrale premi\`ere t.f.
non-constante.
\end{rema}

\indent{\it  $\bullet $ \underline{Troisi\`eme \'etape }: L'induction.}  pour un entier $k$, $0 \leq k \leq q-1$, supposons
construite la collection
$\left( \wF_{\C,
\,K}\right)_{K \in \cal C^k}$de mani\`ere que les conditions de coh\'erence, de non-d\'eg\'enerescence et de tangence \`a la d\'eformation
constante soient satisfaites. Nous allons construire les d\'eformations
\begin{equation}\label{colldef}
\wF_{\C, \,K}\,,Ê\qquad K\; \in \; \cal C^{k+1} \,- \,\cal C^k\,,
\end{equation}
pour que la collection  $\left( \wF_{\C, \,K}\right)_{K \in \cal C^{k+1}}$ satisfasse aussi les m\^emes conditions.\\

Lorsque l'op\'eration $L^{k+1} \efle L^{k+2}$ est un \'elagage, chaque branche $B$ de $L^{k+1}$ porte sur l'une de ses extr\'emit\'es
un point singulier $c_B$ qui est un point d'attache d'une branche ou d'un diviseur de t\^ete de $L^k$. Ainsi $c_B$ est un
\'el\'ement de $\cal C^k$ et la d\'eformation $\wF_{ \C, \, c_B }$ est d\'ej\`a donn\'e par l'induction. Nous d\'efinissons alors  les
autres d\'eformations de la collection (\ref{colldef}) \`a l'aide du lemme (\ref{lemmeun}).\\

Supposons maintenant que l'op\'eration $L^{k+1} \efle L^{k+2}$ est un \'et\^etage. Consid\'erons un diviseur de t\^ete $D \subset
L^{k+1}$. Adoptons les m\^emes notations qu'\`a la deuxi\`eme \'etape. Ordonnons aussi les points singuliers port\'es par $D$ pour que~:
$$
c_1 ,\ldots , c_s \in \cal C^k \qquad \mbox{et}   \qquad  c_j \notin \cal C ^k \quad \mbox{pour}\quad j\geq s+1\,,
$$
\noindent et pour que $c_{r+1}$ soit le point d'attache de $D$ sur $L^{k+1}$. On a~: $s \geq 1$ et $r \geq 2$. Pour $j =
1, \ldots , s$, les d\'eformations $\wF_{\C,\, c_j}$ sont d\'eja donn\'ees par l'induction.  En plus du syst\`eme de lacets
$\gamma_1 , \ldots , \gamma_{r+1}$, donnons-nous aussi pour chaque $j = 1, \ldots , s$ un disque conforme
$\Delta_j\subset W_j\,$, $c_j \in \, \inte{\Delta}_j\,$ suffisament petit pour que
les germes
$\wF_{\C , c_j}(t)$, soient r\'ealis\'es et \`a singularit\'es isol\'ees $c_j(t)$ dans $\Delta_j$.
\\

Remarquons d'abord que les holonomies $\underline \wh_{j, t}$  de  $\wF_{\C , c_j}(t)$ le long des lacets $\b_j:= \partial
\Delta_j$,
$j=1,
\cdots , s$, sont toutes diff\'erentes de l'identit\'e. En effet, si $k>1$, les germes
$\wF_{\C , c_j}(t)$ en
$c_j(t)$ ne poss\`edent pas d'int\'egrale premi\`ere t.f. non-constante d'apr\`es les conditions de non-d\'eg\'en\'erescence de
l'induction~: les points $c_1, \ldots , c_s$ sont des points d'attache de chaines de $\D_\w$, ou bien des points
d'intersections de deux composantes de $\D\w$ de valence
$\geq 3$. Si $k=1$, les points $c_1, \ldots , c_s$ sont des points d'attache de branche ou de diviseurs de t\^ete. Le lemme
(\ref{lemmetrois}) ou la remarque $(\ref{remetoi})$ ci dessus permettent de conclure.\\

\underline{\it Premi\`ere \'eventualit\'e~$s=r$}. Pour construire $\wF_{\C ,\,K}$ pour $K = D - \{c_1, \ldots , c_{r+1}\}$, il
suffit de d\'eterminer des diff\'eomorphismes formels $h_{1, t},\ldots , h_{r+1, t}$ d\'e\-pen\-dant holomorphiquement d'un
petit param\`etre $t \in (\C, 0)$, qui satisfont les condition (a)--(c) ci dessous, puis d'appliquer la proposition
(\ref{constrcol})\`a ces diff\'eomorphismes~:
\begin{enumerate}
\item[(a)] pour $j = 1, \ldots , s$, chaque $h_{j, t}$ est t.f. conjugu\'e \`a $\underline{\wh}_{j, t}$ par un diff\'eomorphisme t.f.
$l$-tangent \`a l'identit\'e,
\item[(b)] $h_{r+1, t} = \left( h_{1, t} \circ \cdots Ê\circ h_{r, t} \right)^{-1}$,
\item[(c)] pour $t \not= 0$, $h_{r+1, t}$ n'est pas p\'eriodique et le groupe engendr\'e par $\,h_{1, t} , \ldots ,h_{r, t}$
n'est pas ab\'elien.
\end{enumerate}
\noindent Les diff\'eomorphismes t.f. $\underline{\wh}_{j, t}$ \'etant diff\'erents
de l'identit\'e pour $t\not=0$ nous d\'eterminons ais\'ement $h_{2, t},\ldots,h_{r, t}$ v\'erifiant (a) pour que
$g_t :=h_{2, t} \circ \cdots
h_{r, t}$ ne soit pas l'identit\'e. On utilise alors lemme suivant dont la d\'emonstration est laiss\'ee au lecteur.
\begin{lemm}
Pour tout entier $l \geq 1$ il existe une famille holomorphe de s\'eries formelles $\Phi_t(z)$ telle que~:
\begin{enumerate}
\item $\phi_t(z) = z \quad mod\, ( z^{l+1})$,
\item $\phi_t \circ \underline{\wh}_{1, t}\circ \Phi_t^{-1}$ et $g_t$ ne commuttent pas,
\item $\phi_t \circ \underline{\wh}_{1, t}\circ \Phi_t^{-1}\circ g_t$ n'est pas p\'eriodique.
\end{enumerate}
\end{lemm}
\noindent On pose $h_{1, t} := \phi_t \circ \underline{\wh}_{1, t}\circ \Phi_t^{-1}$ et  $h_{r+1, t}$ est donn\'e par la
condition (b). Ceci ach\`eve la construction de $\wF_{\C, \, K}$. \\

Pour construire $\wF_{\C , \, c_{r+1}}$, on applique le
lemme (\ref{slemme2}) \`a $\Psi \circ h_{r+1, t} \circ \Psi^{-1}$, o\`u $\Psi(z)$ d\'esigne le transport holonome, via le
feuilletage $\wF$, le long d'un chemin simple $\a$ trac\'e sur $W_{r+1} -\{c_{r+1}\}$ qui relie $m_0$ \`a l'origine de
$\b_{r+1}$.\\

\underline{\it Seconde \'eventualit\'e~: $s \leq r-1$}. Visiblement les points $c_{s+1},\ldots , c_r$ ne sont pas des singularit\'es
de
$\D_\w$. Il suffit de construire $h_{1, t} , \ldots, h_{r+1, t}$ v\'erifiant les conditions (a)--(c) pr\'ec\'edentes, ainsi que
la condition suivante~:
\begin{enumerate}
\item[(d)]pour $j = s+1, \ldots ,r$ chaque  $h_{j, t}$ est une d\'eformation t.f.
de l'holonomie $\underline{h}_j$ de $\wF$ le
long du lacet $\gamma_j$ et de plus $h_{j, t} =
\underline{h}_j\quad mod (z^{l+1})$.
\end{enumerate}
\noindent Les degr\'es de libert\'e suppl\'e
mentaires que nous avons ici rendent plus facile cette construction et nous laissons au
lecteur le soin de lachever.
\end{proof}
%%%%%%%%%%%%%%%%%%%%%%%%%%%%%%%%%%%%%%%%
%%%%%%%%%%%%%%%%%%%%%%%%%%%%%%%%%%%%%%%%%%%%%%%%%%
%%%%%%%%%%%%%%%%%%%%%%%%%%%%%%%%%%%%%%%%

\vspace{3em}

%%%%%%%%%%%%%%%%%%%%%%%%%%%%%%%%%%%%%%%%%%%%

%%%%%%%%%%%%%%%%%%%%%%%%%%%%%%%%%%%%%%%%%%%%
%%%%%%%%%%%%%%%%%%%%%%%%%%%%%%%%%%%%%%%%%%%%
\end{document}